\documentclass[12pt,reqno]{amsart}

\usepackage{amscd}
\usepackage{amssymb}
\usepackage{float}
\usepackage{url}
\usepackage[all]{xy}
\usepackage{color}
\usepackage{hyperref}
\usepackage{pgf,tikz}
\usepackage{mathrsfs}
\usetikzlibrary{arrows}
\usetikzlibrary[patterns]
\definecolor{ttqqqq}{rgb}{0.,0.,0.}
\newcommand{\Z}{{\mathbb Z}}
\newcommand{\Q}{{\mathbb Q}}
\newcommand{\C}{{\mathbb C}}
\newcommand{\R}{{\mathbb R}}
\renewcommand{\P}{{\mathbb P}}
\renewcommand{\H}{{\mathbb H}}

\newcommand{\BB}{{\mathcal B}}
\newcommand{\DD}{{\mathcal D}}

\newcommand{\HH}{{\mathcal H}}

\newcommand{\LL}{{\mathcal L}}

\newcommand{\OO}{{\mathcal O}}

\newcommand{\RR}{{\mathcal R}}

\newcommand{\www}{\widetilde}

\newcommand{\oooo}{\overline}
\newcommand{\uuuu}{\underline}

\newcommand{\ppp}{\partial}

\DeclareMathOperator{\Aut}{Aut}

\DeclareMathOperator{\Eiw}{Eiw}

\DeclareMathOperator{\Gal}{Gal}

\DeclareMathOperator{\Gr}{Gr}

\DeclareMathOperator{\id}{id}
\DeclareMathOperator{\Imm}{Im}
\DeclareMathOperator{\Lie}{Lie}

\DeclareMathOperator{\mult}{mult}
\DeclareMathOperator{\Norm}{Norm}
\DeclareMathOperator{\ord}{ord}
\DeclareMathOperator{\Ord}{Ord}

\DeclareMathOperator{\rank}{rank}

\DeclareMathOperator{\Sp}{Sp}
\DeclareMathOperator{\Stab}{Stab}
\DeclareMathOperator{\supp}{supp}

\DeclareMathOperator{\tr}{tr}

\begin{document}

\theoremstyle{plain}
\newtheorem{lemma}{Lemma}[section]
\newtheorem{definition/lemma}[lemma]{Definition/Lemma}
\newtheorem{theorem}[lemma]{Theorem}
\newtheorem{proposition}[lemma]{Proposition}
\newtheorem{corollary}[lemma]{Corollary}
\newtheorem{conjecture}[lemma]{Conjecture}
\newtheorem{conjectures}[lemma]{Conjectures}

\theoremstyle{definition}
\newtheorem{definition}[lemma]{Definition}
\newtheorem{withouttitle}[lemma]{}
\newtheorem{remark}[lemma]{Remark}
\newtheorem{remarks}[lemma]{Remarks}
\newtheorem{example}[lemma]{Example}
\newtheorem{examples}[lemma]{Examples}
\newtheorem{notations}[lemma]{Notations}

\title
[Torelli results for marked bimodal singularities]
{$\mu$-constant monodromy groups and Torelli
results for the quadrangle singularities and the 
bimodal series}

\author{Falko Gauss\and Claus Hertling}

\address{Falko Gauss\\Universit\"at Mannheim\\ 
Lehrstuhl f\"ur Mathematik VI\\
Seminargeb\"aude A 5, 6\\
68131 Mannheim, Germany}

\email{gauss@math.uni-mannheim.de}

\address{Claus Hertling\\Universit\"at Mannheim\\ 
Lehrstuhl f\"ur Mathematik VI\\
Seminargeb\"aude A 5, 6\\
68131 Mannheim, Germany}

\email{hertling@math.uni-mannheim.de}

\thanks{This work was supported by the DFG grant He2287/4-1 (SISYPH). Moreover the first author gratefully acknowledges support from the JSPS Postdoctoral Fellowship for Research in Japan ID No. PE17702.}


%% Mathematical classification (2000)

\date{October 10, 2017}

\begin{abstract}
This paper is a sequel to \cite{He11} and \cite{GH17}. 
In \cite{He11} a notion of
marking of isolated hypersurface singularities was defined, 
and a moduli space $M_\mu^{mar}$ 
for marked singularities in one $\mu$-homotopy class 
of isolated hypersurface singularities was established.
It is an analogue of a Teichm\"uller space.
It comes together with a $\mu$-constant monodromy group 
$G^{mar}\subset G_\Z$. Here $G_\Z$ is the group of automorphisms
of a Milnor lattice which respect the Seifert form.
It was conjectured that $M_\mu^{mar}$ is connected.
This is equivalent to $G^{mar}= G_\Z$.
Also Torelli type conjectures were formulated.
In \cite{He11} and \cite{GH17} $M_\mu^{mar}, G_\Z$ and $G^{mar}$
were determined and all conjectures were proved
for the simple, the unimodal and the exceptional bimodal
singularities.
In this paper the quadrangle singularities and the bimodal series
are treated. The Torelli type conjectures are true.
But the conjecture $G^{mar}= G_\Z$ and $M_\mu^{mar}$ connected
does not hold for certain subseries of the bimodal series.
\end{abstract}

\dedicatory{To the memory of Egbert Brieskorn}

\maketitle

\tableofcontents

\setcounter{section}{0}

\section{Introduction}\label{c1}

\noindent
We dedicate this paper to the memory of Egbert Brieskorn.
It has its roots in work which the second author, Claus Hertling,
had done as a student of Brieskorn in Bonn in the early 90'ies.

\subsection{Reminiscences of the second author}\label{c1.1}
Prof. Dr. Egbert Brieskorn accepted me as a diploma student
in the spring of 1989. On March 10 and 13, 1989, he gave two full days
(Friday + Monday) of lectures for his new diploma students (including me)
and doctoral students. I still have his handwritten manuscript of
52 pages. There he introduced us to isolated hypersurface singularities.
He talked about the Jacobi algebra, the universal unfolding with
its discriminant, the Milnor fibration, its monodromy,
local systems and integrable connections and systems of regular singular
linear differential equations in general, his own work
on the Gauss-Manin connection and especially the Brieskorn lattice,
and the mixed Hodge structure which it induces.
He strongly recommended to read \cite{AGV88}, \cite{SaM89}
and \cite{SS85}. He proposed to me to work on the moduli of
singularities using the Gauss-Manin connection. 

I followed his advice in my diploma thesis and my doctoral thesis
and beyond the doctoral thesis. 
The subject developed into a long-going project of mine,
which I took up again and again. 
The present paper is in some sense a final step of it. 

In the doctoral thesis \cite{He93}, I formulated the global Torelli type
conjecture that an isolated hypersurface singularity is determined
up to right equivalence by its Brieskorn lattice together with
the Milnor lattice and the Seifert form
(conjecture \ref{t1.1} (b) reformulates this conjecture).
I proved it in the doctoral thesis for all unimodal singularities,
the exceptional bimodal singularities, the bimodal quadrangle singularities,
and the bimodal series $E_{3,p}$.

For the other seven bimodal series, I made in the spring 1993,
some months after finishing the doctoral thesis,
long calculations (120 pages) which led to a proof 
of this Torelli type conjecture for all series 
except the three bimodal subseries
$S_{10,r}^\sharp, S_{1,10r}, Z_{1,14r}$. 
At that time I thought that I would never review 
and publish these results. The paper \cite{He95} recapitulated
the main results of the doctoral thesis and of these calculations
for the eight bimodal series, but it did not at all give all details
(only 2.5 pages are devoted to the bimodal series).

Later I constructed a classifying space $D_{BL}$ for
Brieskorn lattices \cite{He99} and a moduli space $M_\mu(f_0)$
of the right equivalence classes of all singularities in the
$\mu$-homotopy class of a reference singularity $f_0$ \cite{He02}.
More recently, in \cite{He11}, I defined the notion of a 
{\it marked singularity}, I constructed a classifying space
$M_\mu^{mar}(f_0)$ for marked singularities, and I formulated
a Torelli type conjecture for marked singualarities,
which is stronger than the Torelli type conjecture in the doctoral
thesis for unmarked singularities.

The three papers \cite{He11}, \cite{GH17} and the present paper
prove the Torelli conjecture for marked singularities
for all singularities with modality 0, 1 and 2. 
The present paper deals with the bimodal quadrangle singularities
and the eight bimodal series. It comprises the calculations
from the spring 1993 and adds a lot more arguments and calculations,
which are necessary for the marked version.

It is satisfying, that the Torelli type conjectures hold for
all singularities with modality 0, 1 and 2. For each family, 
the interplay between the variations of the Brieskorn lattices 
and the automorphism group of the Milnor lattice with Seifert form 
is fascinating and takes the best possible shape. 
I believe that Brieskorn would have liked these positive results
and the many techniques used for their proofs. 
I thank him for proposing to me in March 1989 to work on the
moduli of singularities using the Gauss-Manin connection.
It was a good advice.

\subsection{Notions, conjectures and results}\label{c1.2}
In this paper, a {\it singularity} is
a holomorphic function germ $f:(\C^{n+1}\to (\C,0)$ with
an isolated singularity at $0$.
Then its {\it Milnor lattice} $Ml(f)\cong\Z^\mu$ is the $\Z$-lattice
$H_n(f^{-1}(\tau),\Z)$ for some small $\tau\in\R_{>0}$
for a suitable representative of $f$. Its {\it Seifert form} 
is called $L:Ml(f)\times Ml(f)\to\Z$. 
Its {\it monodromy} is called $M_h:Ml(f)\to Ml(f)$.
The automorphism group of the Milnor lattice with the Seifert 
form is $G_\Z(f):=\Aut(Ml(f),L)$.
It will play a predominant role in this paper.

This paper is a sequel to \cite{He11} and \cite{GH17}.
In \cite{He11}, a {\it strongly marked singularity}
$(f,\rho)$ and a {\it marked singularity} $(f,\pm\rho)$
are defined. Here one has to fix first a reference singularity
$f_0$. Then $f$ is in the $\mu$-homotopy class of $f_0$,
i.e. a $\mu$-constant family of singularities exists
which contains $f_0$ and $f$. And 
$\rho:(Ml(f),L(f))\to (Ml(f_0),L(f_0))$ is a chosen isomorphism.
Two singularities $f_1$ and $f_2$ are right equivalent
if a coordinate change $\varphi$ with $f_1=f_2\circ\varphi$ exists.
Two strongly marked singularities $(f_1,\rho_1)$ and $(f_2,\rho_2)$
are right equivalent if a coordinate change $\varphi$ with
$f_1=f_2\circ\varphi$ and $\rho_1=\rho_2\circ (\varphi)_{hom}$
exists, where $(\varphi)_{hom}:Ml(f_1)\to Ml(f_2)$ is the 
induced isomorphism.

In \cite{He02} a moduli space $M_\mu(f_0)$ for the right
equivalence classes of all singularities in the $\mu$-homotopy
class of a reference singularity $f_0$ was constructed as an 
analytic geometric quotient. In \cite{He11}, this construction
was enhanced to the construction of moduli spaces
$M_\mu^{mar}(f_0)$ and $M_\mu^{smar}(f_0)$ of
marked and strongly marked singularities. Here $M_\mu^{smar}(f_0)$
is Hausdorff and an analytic space only if assumption 
\eqref{8.1} or assumption \eqref{8.2} holds. 
\begin{eqnarray*}
\textup{Assumption (8.1):}&&\textup{Any singularity in the }
\mu\textup{-homotopy}\\
&&\textup{class of }f_0\textup{ has multiplicity }\geq 3.\nonumber\\
\textup{Assumption (8.2):}&&\textup{Any singularity in the }
\mu\textup{-homotopy}\\
&&\textup{class of }f_0\textup{ has multiplicity }2.\nonumber
\end{eqnarray*}
We expect that one of them holds for any $\mu$-homotopy class
of singularities. This would be an implication of the
Zariski multiplicity conjecture. But that is not proved in general.

But $M_\mu^{mar}(f_0)$ is fine, independently of these assumptions.
Locally it is isomorphic to the $\mu$-constant stratum 
$S_\mu(f)$ of a singularity in the base space of a universal
unfolding of that singularity. The group $G_\Z(f_0)$
acts properly discontinuously on $M_\mu^{mar}(f_0)$.
The quotient is $M_\mu^{mar}(f_0)/G_\Z\cong M_\mu(f_0)$.
Therefore that space is locally isomorphic to the quotient
of $S_\mu(f)$ by a finite group. 
$M_\mu^{mar}(f_0)$ can be considered as a  
{\it Teichm\"uller space} for singularities, in analogy to
the Teichm\"uller spaces for closed complex curves.
It can also be considered as a global $\mu$-constant stratum,
simultaneously for all singularities in one 
$\mu$-homotopy class.

The papers \cite{He11}, \cite{GH17} and this paper determine
$M_\mu^{mar}(f_0)$ for all singularities with modality
0, 1 and 2. The second column of the following table \eqref{1.1}
gives their isomorphism classes.

\begin{eqnarray}\label{1.1}
\begin{array}{lll}
\textup{Singularity family} & M_\mu^{mar}(f_0) & D_{BL}(f_0) \\
 \hline 
\textup{ADE-singularities} & \textup{point} & \textup{point} \\ 
\textup{simple elliptic sing.} & \H & \H \\ 
\textup{hyperbolic sing.} & \C & \C \\
\textup{exc. unimodal sing.} & \C & \C \\
\textup{exc. bimodal sing.} & \C^2 & \C^2 \\
\textup{quadrangle sing.} & (\H-\textup{(a discrete set)}\times\C 
 & \H\times\C \\
\textup{the 8 series, for }m\not|p & \C^*\times\C & \C^{N_{BL}} \\
\textup{the 8 subseries with }m|p & \infty\textup{ many copies of }
\C^*\times\C & \H\times \C^{N_{BL}}
\end{array}
\end{eqnarray}

Here the eight series and the respective numbers $m$ are
given in the following table \eqref{1.2}. Here $p\in\Z_{\geq 1}$.
\begin{eqnarray}\label{1.2}
\begin{array}{lllllllll}
\textup{series} & 
W_{1,p}^\sharp & S_{1,p}^\sharp & U_{1,p} & E_{3,p} & Z_{1,p} & Q_{2,p}
& W_{1,p} & S_{1,p} \\
m & 12 & 10 & 9 & 18 & 14 & 12 & 12 & 10 
\end{array}
\end{eqnarray}

One sees that $M_\mu^{mar}(f_0)$ is simply connected for all singularities
with modality 0 and 1 and for the exceptional bimodal singularities.
For the quadrangle singularities and the series with $m\not|p$, 
it is connected, but not simply connected. And for the subseries with
$m|p$, it is not even connected, but has infinitely 
many components.
This last result is a counterexample to conjecture 3.2 (a) 
in \cite{He11}, which said that $M_\mu^{mar}(f_0)$ should be connected. 

In \cite{He11}, also two subgroups $G^{smar}(f_0)$ and 
$G^{mar}(f_0)$ of $G_\Z(f_0)$ were defined.
$G^{smar}(f_0)$ was defined as the subgroup which is 
generated by the
transversal monodromies of all $\mu$-constant families 
which contain $f_0$. Then $G^{mar}(f_0)$ is the group generated by
$G^{smar}(f_0)$ and $-\id$.
A rough way to talk about this description is to say 
that the elements of $G^{smar}(f_0)$ are of geometric origin.
$G^{mar}(f_0)$ can also be characterized as the subgroup 
of $G_\Z$ which maps the component $(M_\mu^{mar})^0$ of 
$M_\mu^{mar}(f_0)$, which contains $[(f_0,\pm\id)]$, to itself. 
This last characterization gives
\begin{eqnarray}\label{1.3}
G_\Z(f_0)/G^{mar}(f_0)\stackrel{1:1}{\longleftrightarrow}\{\textup{components of }
M_\mu^{mar}(f_0)\}.
\end{eqnarray}
In view of this, $M_\mu^{mar}(f_0)$ is connected if and only if 
$G_\Z(f_0)=G^{mar}(f_0)$. By table \eqref{1.1}, this holds for
all singularities with modality 0, 1 or 2 except the eight 
subseries with $m|p$. 
Obviously, it is important to control $G_\Z(f_0)$.
This was the major task in \cite{He11} and \cite{GH17}
for the singularities considered there, and it takes approximately
half of this paper for the singularities considered here,
the bimodal series and the quadrangle singularities. 
The rough outcome in all cases is that the pair $(Ml(f_0),L)$ is 
surprisingly rigid and that $G_\Z(f_0)$ is surprisingly small.
The next table \eqref{1.4} gives more information on $G_\Z(f_0)$
for all singularities with modality 0, 1 and 2.
Here $M_h\in G_\Z$ is the classical monodromy. It commutes
with all elements of $G_\Z$.
The only families in table \eqref{1.4} where 
$\{\pm M_h^k\, |\, k\in\Z\}$
is not finite, are the hyperbolic singularities $T_{pqr}$.

\begin{eqnarray}\label{1.4}
\begin{array}{ll}
\textup{Singularity family} & 
G_\Z(f_0)/\{\pm M_h^k\, |\, k\in\Z\} \\ \hline 
\textup{ADE-singularities} & 
\{\id\}\textup{ or }S_2\textup{ or }S_3  \\  
\textup{simple elliptic sing.} & 
\textup{a finite extension of }SL(2,\Z) \\
\textup{hyperbolic sing.} & \textup{a finite group} \\
\textup{exc. unimodal sing.} & \{\id\}\textup{ or }S_2\textup{ or }S_3  \\
\textup{exc. bimodal sing.} & \{\id\}\textup{ or }S_2\textup{ or }S_3  \\
\textup{quadrangle sing.} &  \textup{a triangle group}  \\
\textup{the 8 series, for }m\not|p & \textup{a cyclic finite group}  \\
\textup{the 8 subseries with }m|p & \textup{an infinite Fuchsian group}  
\end{array}
\end{eqnarray}

In the case of the eight subseries with $m|p$, $G^{mar}(f_0)$
is the finite subgroup of the infinite group $G_\Z(f_0)$
such that $G^{mar}(f_0)/\{\pm M_h^k\, |\, k\in\Z\}$ is
the finite cyclic group which is generated by one elliptic element.

If the $\mu$-homotopy class of $f_0$ contains at least 
one singularity with multiplicity two, then $-\id\in G^{smar}(f_0)$
and $G^{smar}(f_0)=G^{mar}(f_0)$. Conjecture 3.2 (b) in 
\cite{He11} complements this. It claims that $-\id\notin G^{smar}(f_0)$
if assumption \eqref{8.1} holds. This is true for all singularities
with modality 0, 1 and 2. For the bimodal series and the quadrangle
singularities, it is proved in this paper.

In \cite{He99} the second author defined a classifying space
$D_{BL}$ for Brieskorn like lattices (i.e. for objects which
are sufficiently similar to the Brieskorn lattice $H_0''(f_0)$,
see section \ref{c7} before theorem \ref{t7.11} for details).
The group $G_\Z(f_0)$ acts properly discontinuously on it.
The elements of $D_{BL}(f_0)$ are {\it marked Brieskorn like lattices},
and the elements of $D_{BL}(f_0)/G_\Z(f_0)$ are isomorphism classes of 
Brieskorn like lattices. One obtains a holomorphic period map
\begin{eqnarray}\label{1.5}
BL: M_\mu^{mar}(f_0)\to D_{BL}(f_0).
\end{eqnarray}
By \cite[Theorem 12.8]{He02} it is $G_\Z(f_0)$-equivariant,
and it is an immersion (this fact is an infinitesimal Torelli 
type result). Now the following Torelli type conjectures
are natural. Part (a) is for marked singularities.
Part (b) recasts the Torelli type conjecture in \cite{He93}.
Part (a) implies part (b).

\begin{conjecture}\label{t1.1}
(a) \cite[Conjecture 5.3]{He11}
The map $BL$ is injective.

(b) \cite[Kap. 2 d)]{He93} 
The map $BL/G_\Z(f_0): M_\mu(f_0)\to D_{BL}(f_0)/G_\Z(f_0)$
is injective.
\end{conjecture}

\begin{theorem}\label{t1.2}
(\cite{He93}\cite{He11}\cite{GH17} and the theorems
\ref{t9.1} and \ref{t10.1} in this paper)
Both Torelli type conjectures are true for all
singularities with modality 0, 1 and 2.
\end{theorem}

The proofs have in almost all cases two parts:

\begin{list}{}{}
\item[(1)] A good control of an (often multivalued)
period map $T\to D_{BL}(f_0)$, where $T$ is the 
parameter space of a well chosen family of 
normal forms.
\item[(2)] A good control of $G_\Z(f)$ and its action
on $M_\mu^{mar}(f_0)$ and $D_{BL}(f_0)$.
\end{list}

In all cases, (1) is less work than (2). 
For the ADE-singularities, (1) is empty as there $T$
is a point, but (2) is not. 

Part (b) of conjecture \ref{t1.1} 
was proved in \cite{He93} for the unimodal and
bimodal singularities except seven of the eight series.
For the seven series, the second author had unpublished
calculations shortly after \cite{He93}.
But for technical reasons, part (b) stayed open for the
subseries $S_{1,10r}^\sharp, S_{1,10r}, Z_{1,14r}$.
\cite{He93} and these unpublished calculations give (1)
and a part of (2).

In view of these old results,
the major point in \cite{He11}, \cite{GH17}
and in this paper is (2). But also some refinement
of (1) is needed in the case of the singularities
in this paper. The refinement is used for a 
better control of the transversal monodromy 
of the family of normal forms.

Finally, the conjecture $G_\Z(f_0)=G^{mar}(f_0)$ is 
probably wrong in general as it is wrong for the
subseries with $m|p$. But for all singularities
with modality 0, 1 and 2 except the eight series,
the Torelli result for marked singularities and \eqref{1.3}
require $G_\Z(f_0)=G^{mar}(f_0)$
to be true, as $BL$ is an immersion and there 
$\dim M_\mu^{mar}(f_0)=\textup{modality}(f_0)=\dim D_{BL}(f_0)$.
And there $G_\Z(f_0)=G^{mar}(f_0)$ holds indeed.
For the eight series, $\dim D_{BL}(f_0)>\dim M_\mu^{mar}(f_0)$,
so there is enough space in $D_{BL}$ for 
infinitely many copies of $(M_\mu^{mar}(f_0))^0$.

Open questions are now how to control
the subgroup $G^{mar}(f_0)\subset G_\Z(f_0)$ in general,
and how to attack the Torelli conjectures in greater
generality. For the second question, we plan to
thicken $M_\mu^{mar}(f_0)$ to a $\mu$-dimensional
$F$-manifold $M^{mar}(f_0)$ which is locally at each point of 
$M_\mu^{mar}(f_0)$ the base space of a universal unfolding.
Then we will try to embed the Torelli type conjecture
for $M_\mu^{mar}$ into a family of Torelli type conjectures
for all the $\mu$-homotopy strata of multigerms of
singularities in $M^{mar}(f_0)$.
We hope that this global point of view and the different
geometry there with Stokes structures will give us
new techniques. But this is a hope for the future.

\subsection{Structure of the paper}\label{c1.3}
{\bf Section \ref{c2}} is a collection of techniques
which are useful to control the automorphisms of a
pair $(\Lambda,L)$ or a pair $(\Lambda,M_h)$ where
$\Lambda$ is a $\Z$-lattice, $L$ is a unimodular 
bilinear form and $M_h$ is an automorphism of finite order.
We define {\it Orlik blocks} and study their
automorphisms (lemma \ref{t2.8} will be very useful),
and we cite classical algebraic facts on unit roots $\zeta$
and the rings $\Z[\zeta]$. All this is needed for
the control of $G_\Z(f_0)$ in the sections \ref{c5}
and \ref{c6}.

{\bf Section \ref{c3}} discusses infinite Fuchsian groups
which arise as subgroups of groups
$GL(2,\Z[\zeta])$ with $\zeta$ a unit root.
They are in fact arithmetic Fuchsian groups.
But our treatment is essentially self-contained.
Solutions of Pell equations with coefficients
in $\Z[\zeta]$ play a role.
For the quadrangle singularities, we need a precise
analysis of some of these groups. They are certain
triangle groups.

{\bf Section \ref{c4}} recalls some classical notions
and facts around singularities: Milnor fibration,
Milnor lattice $Ml(f)$, monodromy $M_h$, 
Seifert form $L$, Coxeter-Dynkin diagram,
Stokes matrix, Thom-Sebastiani type results, 
suspension, polarized mixed Hodge structure on
$H^\infty_\C$, its polarizing form.

{\bf Section \ref{c5}} is long. It studies $G_\Z(f_0)$
for the eight bimodal series. 
Theorem \ref{t5.1} states the results.
We start with a distinguished basis of the Milnor lattice
with Coxeter-Dynkin diagram in \cite{Eb81}.
We calculate the monodromy $M_h$ and find 2 or 3 (3 only for 
$Z_{1,p}$) Orlik blocks whose direct sum is of
index 1 or 2 in $Ml(f_0)$.
Then $G_\Z(f_0)$ is studied using these Orlik blocks
and their rigidity 
and the results from the sections \ref{c2} and \ref{c3}.
A lot of calculations are needed, the different
series behave differently. The singularities in the families
$Q_{2,p}, W_{1,6s-3}, S_{1,10}$ need special care.

{\bf Section \ref{c6}} gives similar results for
$G_\Z(f_0)$ for the quadrangle singularities.
Theorem \ref{t6.1} states the results.
Many, but not all, calculations and arguments in section
\ref{c5} are also valid in section \ref{c6}.
Therefore this section is much shorter.

{\bf Section \ref{c7}} gives a rather complete account
on the Gauss-Manin connection and the Brieskorn lattice
$H_0''(f)$ of a singularity $f$. It does not rewrite
the proofs in \cite{Br70} and other papers, but it cites
almost all known results. A highlight is the treatment
of the bilinear forms. The polarizing form of the
polarized mixed Hodge structure is connected with
the restriction of K. Saito's higher residue pairings
to $H_0''(f)$ and with Pham's intersection form for
Lefschetz thimbles. 
We need the Fourier-Laplace transform $FL(H_0''(f))$
for a Thom-Sebastiani formula for Brieskorn lattices.
We need this in the special case of a suspension
$f(z_0,...,z_n)+z_{n+1}^2$ because we want to treat 
the suspensions in a more conceptual way than in 
\cite{He93}\cite{He11}\cite{GH17}.

{\bf Section \ref{c8}} reviews the notions and results
from \cite{He11}, the (strongly) marked singularities
and their moduli spaces $M_\mu^{smar}(f_0)$ and 
$M_\mu^{mar}(f_0)$, the $\mu$-constant monodromy
groups $G^{smar}(f_0)$ and $G^{mar}(f_0)$,
and the Torelli conjectures.
Corollary \ref{t8.14} is an application of the Thom-Sebastiani
result for $FL(H_0''(f))$ in section \ref{c7}
and states that the marked Torelli conjecture for
$f_0$ is equivalent to the marked
Torelli conjecture for $f_0(z_0,...,z_n)+\sum_{j=n+1}^m z_j^ 2$
for any fixed $m\geq n+1$. This allows us to consider in 
the sections \ref{c9} and \ref{c10} only the surface
singularities.

{\bf Section \ref{c9}} proves the marked Torelli 
conjecture for the bimodal series (theorem \ref{t9.1}).
It establishes the good control (1) of the
multivalued period map $T\to D_{BL}(f_0)$
where $T=\C^*\times \C$ is the parameter space
of normal forms in \cite{AGV85}.
Theorem \ref{t5.1} provides crucial information 
on $G_\Z(f_0)$.

{\bf Section \ref{c10}} proves the marked Torelli conjecture
for the quadrangle singularities (theorem \ref{t10.1}).
It starts with a careful own choice of normal forms
with parameter space $T=(\C-\{0,1\})\times\C$.
It establishes the good control (1) of the
multivalued period map $T\to D_{BL}(f_0)$.
Theorem \ref{t6.1} provides crucial information
on $G_\Z(f_0)$.

\section{\texorpdfstring{$\Z$}{Z}-lattices with unimodal bilinear form and monodromy}\label{c2}

\noindent
This section provides tools for the study of the Milnor lattices with
Seifert form and monodromy for the bimodal series and the
quadrangle singularities, in the sections \ref{c5} and \ref{c6}. 
These lattices turn out to be quite rigid and to have
rather few automorphisms.
This is important for the global Torelli results in the sections
\ref{c9} and \ref{c10}. 
This section puts together elementary, but nontrivial observations
about $\Z$-lattices with a unimodal bilinear form and an (induced)
monodromy.

Let $\Lambda$ be a $\Z$-lattice of rank $\mu\in\Z_{\geq 1}$,
i.e. a free $\Z$-module of rank $\mu$.
Let $L:\Lambda\times \Lambda\to\Z$ be a unimodal bilinear form,
i.e. for any basis $\delta_1,\ldots,\delta_\mu$ we have 
$\det(L(\delta_i,\delta_j)_{i,j=1,\ldots,\mu})=\pm 1$.
We do not suppose that $L$ is symmetric or antisymmetric.
Let $M_h:\Lambda\to\Lambda$ be the automorphism which is 
uniquely determined by
\begin{eqnarray}\label{2.1}
L(M_h(a),b)=-L(b,a)\quad\textup{for }a,b\in\Lambda.
\end{eqnarray}
We call $L$ the {\it Seifert form} and $M_h$ the {\it monodromy}.
\eqref{2.1} implies
\begin{eqnarray}\label{2.2}
L(M_h(a),M_h(b))=L(a,b),
\end{eqnarray}
i.e. $L$ is $M_h$-invariant. We make the assumption that
\begin{eqnarray}\label{2.3}
M_h\textup{ is finite,}
\end{eqnarray}
i.e. $M_h$ is semisimple and its eigenvalues are unit roots.
Then the characteristic polynomial $p_\Lambda$ of $M_h$ is
a product of cyclotomic polynomials.

\begin{notations}\label{t2.1}
(a) For any subring $R\subset\C$ denote
$\Lambda_R:=\Lambda\otimes_\Z R$.
For any monodromy invariant subspace $V\subset \Lambda_\C$
denote by $E(V)\subset S^1$ the set of eigenvalues of $M_h$ on
$V$ and by $p_V$ its characteristic polynomial.
For $\lambda\in E(V)$ denote 
$V_\lambda:=\ker(M_h-\lambda\id:V\to V)\subset V$.
For any monodromy invariant sublattice $\Lambda^{(1)}\subset \Lambda$
write $E(\Lambda^{(1)}):=E(\Lambda^{(1)}_\C)$ and
$p_{\Lambda^{(1)}}:=p_{\Lambda^{(1)}_\C}$ and 
$\Lambda^{(1)}_\lambda:=(\Lambda^{(1)}_\C)_\lambda$.
For any product $p\in\Z[t]$ of cyclotomic polynomials
with $p|p_{\Lambda^{(1)}}$ denote
\begin{eqnarray}\label{2.4}
\Lambda^{(1)}_{\C,p}:=\bigoplus_{\lambda:\, p(\lambda)=0}
\Lambda^{(1)}_\lambda
\quad\textup{and}\quad \Lambda^{(1)}_p
:=\Lambda^{(1)}_{\C,p}\cap\Lambda^{(1)}.
\end{eqnarray}
Then $\Lambda^{(1)}_p$ is a primitive and monodromy invariant sublattice
of $\Lambda^{(1)}$.

\medskip
(b) Recall that a sublattice $\Lambda^{(1)}$ of $\Lambda$ is 
{\it primitive} (in $\Lambda$) 
if and only if $\Lambda/\Lambda^{(1)}$ has no torsion
and that for any sublattice $\Lambda^{(2)}\subset\Lambda$
there is a unique primitive sublattice $\Lambda^{(3)}$ with
$\Lambda^{(3)}_\Q=\Lambda^{(2)}_\Q$, that it is
$\Lambda^{(3)}=\Lambda^{(2)}_\Q\cap\Lambda$ and that
$[\Lambda^{(3)}:\Lambda^{(2)}]<\infty$.

\medskip
(c) For $n\in\Z_{\geq 1}$, the cyclotomic polynomial
$\Phi_n$ is $$\Phi_n=
\prod_{\lambda:\,\textup{ord}(\lambda)=n}(t-\lambda).$$
It is unitary and in $\Z[t]$ and irreducible in $\Z[t]$ and $\Q[t]$.

\medskip
(d) We define the square root on $S^1-\{-1\}$ by
$\sqrt{e^{2\pi i\alpha}}:=e^{\pi i \alpha}$ for 
$\alpha\in ]-\frac{1}{2},\frac{1}{2}[$.
\end{notations}

\begin{lemma}\label{t2.2}
(a) Let $\lambda\in E(\Lambda)-\{1\}$. Then the sesquilinear
(=linear$\times$semilinear) form
$h_\lambda:\Lambda_\lambda\times\Lambda_\lambda\to\C$ with 
\begin{eqnarray}\label{2.5}
h_\lambda(a,b):=\sqrt{-\lambda}\cdot L(a,\oooo{b})
\end{eqnarray}
is hermitian, i.e. $h_\lambda(b,a)=\oooo{h_\lambda(a,b)}$.
Especially, $\sqrt{-\lambda}\cdot L(a,\oooo{a})\in\R$.
Together, these forms define a hermitian form 
$h:=\bigoplus_{\lambda\in E(\Lambda)-\{1\}}h_\lambda$.

\medskip
(b) Let $V\subset \Lambda_\C$ be a monodromy invariant
subspace with $1\notin E(V)$. The following two properties
are equivalent.

\begin{list}{}{}
\item[$(\alpha)$]
$h|_V$ is positive definite.
\item[$(\beta)$]
The hermitian form on $V$ defined by 
$(a,b)\mapsto L(a,\oooo{b})+L(\oooo{b},a)$ is positive definite.
\end{list}
\end{lemma}

{\bf Proof:}
(a) For $a,b\in V_\lambda$
\begin{eqnarray*}
\sqrt{-\lambda}\cdot L(b,\oooo{a}) 
&=& -\sqrt{-\lambda}\cdot L(M_h(\oooo{a}),b)
= -\sqrt{-\lambda}\cdot\oooo{\lambda}\cdot L(\oooo{a},b)\\
&=& \sqrt{-\oooo{\lambda}}\cdot \oooo{L(a,\oooo{b})}
= \oooo{\sqrt{-\lambda}\cdot L(a,\oooo{b})}.
\end{eqnarray*}

(b) Consider some $\lambda\in E(V)$. Observe
$\sqrt{-\lambda}+\sqrt{-\oooo{\lambda}}>0$ and for $a,b\in V$
\begin{eqnarray*}
L(a,\oooo{b})+L(\oooo{b},a)
&=& L(a,\oooo{b})+\oooo{L(b,\oooo{a})} \\
&=& \sqrt{-\oooo{\lambda}}\cdot h_\lambda(a,b) 
+ \oooo{\sqrt{-\oooo{\lambda}}\cdot h_\lambda(b,a)}\\
&=& (\sqrt{-\oooo{\lambda}}+\sqrt{-\lambda})\cdot h_\lambda(a,b).
\hspace*{2cm}\Box
\end{eqnarray*}

\begin{remarks}\label{t2.3}
(i) The surface singularities considered in this paper do not have 1 
as an eigenvalue of their monodromy. Therefore we do not treat it here.

\medskip
(ii) Part (b) of lemma \ref{t2.2} connects to the polarization of
the polarized Hodge structure of these surface singularities
and rewrites it in different ways.
$(\beta)$ is the classical way, with $-L-L^t$ on $\Lambda_\R$
as intersection form and $L+L^t$ as polarizing form.
And $(\alpha)$ is the way used in the sections \ref{c3}, \ref{c5}
and \ref{c6}.
\end{remarks}

In 1972 Orlik formulated the beautiful conjecture \ref{t2.5} below
on the integral monodromy of quasihomogeneous singularities
\cite{Or72}. It is known to be true for the quasihomogeneous
curve singularities \cite{MW86} and for the quasihomogeneous
singularities with modality $\leq 2$ \cite{He95}.
But it is open for most other quasihomogeneous singularities.

A key observation for the treatment of the Milnor lattices of the
bimodal series singularities and the quadrangle singularities
is that they all have a structure close to Orlik's conjecture.
The following definition gives the ingredients.

\begin{definition}\label{t2.4}
Let $(\Lambda,L,M_h)$ be as above. An {\it Orlik block} is a 
primitive and monodromy invariant sublattice 
$\Lambda^{(1)}\subset\Lambda$ with $\Lambda^{(1)}\supsetneqq \{0\}$
and with a cyclic generator, i.e.
a lattice vector $e^{(1)}\in\Lambda^{(1)}$ with
\begin{eqnarray}\label{2.6}
\Lambda^{(1)} = \bigoplus_{j=0}^{\deg p_{\Lambda^{(1)}}-1}\Z\cdot 
M_h^j(e^{(1)}).
\end{eqnarray}
\end{definition}

\begin{conjecture}\label{t2.5}
\cite[conjecture 3.1]{Or72}
Let $(\Lambda,M_h)$ be the Milnor lattice with
monodromy of a quasihomogeneous singularity.
Let $k:=\max(\dim \Lambda_\lambda\, |\, \lambda\in E(\Lambda))$.
Then a decomposition 
$\Lambda=\bigoplus_{j=1}^k\Lambda^{(k)}$ into Orlik blocks
$\Lambda^{(1)},\ldots,\Lambda^{(k)}$ with
$p_{\Lambda^{(j+1)}}|p_{\Lambda^{(j)}}$ for $0\leq j<k$ exists.
\end{conjecture}

\begin{remarks}\label{t2.6}
(i) A cyclic monodromy module has only one Jordan block for each
eigenvalue. In this paper $M_h$ is semisimple. Therefore
in an Orlik block, each eigenvalue has multiplicity one.

\medskip
(ii) In Orlik's conjecture \ref{t2.5}, the polynomials 
$p_{\Lambda^{(1)}},\ldots,p_{\Lambda^{(k)}}$ are unique. They are
\begin{eqnarray}\label{2.7}
p_{\Lambda^{(j)}} =\prod_{\lambda\in E(\Lambda):\, \dim \Lambda_\lambda\geq j}
(t-\lambda)\quad\textup{ for }j=1,\ldots,k.
\end{eqnarray}

\medskip
(iii) In the sections \ref{c5} and \ref{c6}, we will work most
often with two Orlik blocks $\Lambda^{(1)}$ and $\Lambda^{(2)}$
such that $\Lambda^{(1)}+\Lambda^{(2)}=\Lambda^{(1)}\oplus
\Lambda^{(2)}$ and that it is 
either equal to $\Lambda$ or has index 2 in $\Lambda$
and such that $L(\Lambda^{(1)},\Lambda^{(2)})=
L(\Lambda^{(2)},\Lambda^{(1)})=0$.

\medskip
(iv) In all cases in section \ref{c5} with
$[\Lambda:\Lambda^{(1)}\oplus \Lambda^{(2)}]=2$
except $S_{1,10}$, we will show
\begin{eqnarray}\label{2.8}
\Aut(\Lambda,L)=\Aut(\Lambda^{(1)}\oplus\Lambda^{(2)},L).
\end{eqnarray}

In many of these cases, there is an element $\gamma_5
\in\Lambda^{(1)}_{\Phi_2}-\{0\}$ which is mapped by
any element $g$ of 
$\Aut(\Lambda,L)\cup \Aut(\Lambda^{(1)}\oplus\Lambda^{(2)},L)$ 
to $\pm \gamma_5$ and such that
\begin{eqnarray}\label{2.9}
\Lambda^{(1)}\oplus \Lambda^{(2)}=
\{a\in\Lambda\, |\, L(a,\gamma_5)\in 2\Z\}.
\end{eqnarray}
Then any $g\in \Aut(\Lambda,L)$ maps $\Lambda^{(1)}\oplus \Lambda^{(2)}$
to itself, so
$\Aut(\Lambda,L)\subset\Aut(\Lambda^{(1)}\oplus\Lambda^{(2)},L).$

If this inclusion $\subset$ holds, the following argument
shows that $\Aut(\Lambda,L)$ is either equal to or a subgroup
of index 2 in $\Aut(\Lambda^{(1)}\oplus\Lambda^{(2)},L)$.
Unfortunately it looks hard to exclude the second case.
Therefore in section \ref{c5} we show the equality \eqref{2.8}
in a different (and more laborious) way.

Let $\Lambda^{(0)}\subset\Lambda_\Q$ be the unique lattice
such that 
$$L:\Lambda^{(0)}\times (\Lambda^{(1)}\oplus \Lambda^{(2)})\to\Z$$
is unimodal.
Then $\Lambda^{(0)}\supset\Lambda\supset \Lambda^{(1)}\oplus\Lambda^{(2)}$
and $[\Lambda^{(0)}:\Lambda]=2$
and $\Aut(\Lambda^{(1)}\oplus\Lambda^{(2)},L)
=\Aut(\Lambda^{(0)},L)$. 

{\bf 1st case,} $\Lambda^{(0)}/(\Lambda^{(1)}\oplus\Lambda^{(2)})
\cong\Z/4\Z$. Then $\Lambda$ is the unique lattice between
$\Lambda^{(0)}$ and $\Lambda^{(1)}\oplus\Lambda^{(2)}$ with
$[\Lambda^{(0)}:\Lambda]=2$.
Then any $g\in \Aut(\Lambda^{(1)}\oplus\Lambda^{(2)},L)$
respects $\Lambda$, so \eqref{2.8} holds.

{\bf 2nd case,} $\Lambda^{(0)}/(\Lambda^{(1)}\oplus\Lambda^{(2)})
\cong\Z/2\Z\times\Z/2\Z$. Then there are three lattices
between $\Lambda^{(0)}$ and $\Lambda^{(1)}\oplus\Lambda^{(2)}$
with index 2 in $\Lambda^{(0)}$, one for each subgroup
of index 2 in $\Z/2\Z\times\Z/2\Z$.
One of them is $\Lambda$. Another one is 
$\{a\in\Lambda^{(0)}\, |\, L(a,\gamma_5)\in 2\Z\}$.
No element of $\Aut(\Lambda^{(0)},L)$ maps $\Lambda$ to 
this lattice. But it looks hard to exclude the possibility that
half of the elements of $\Aut(\Lambda^{(0)},L)$
map $\Lambda$ to the third lattice between $\Lambda^{(0)}$
and $\Lambda^{(1)}\oplus\Lambda^{(2)}$.

\medskip
(v) If $\Lambda^{(1)}\subset\Lambda$ is an Orlik block
with cyclic generator $e^{(1)}$ and if 
$p_{\Lambda^{(1)}}=p_1\cdot p_2$ with
$\deg p_1\geq 1$ and $\deg p_2\geq 1$, then
the sublattice $\Lambda^{(2)}:=\Lambda^{(1)}_{p_1}$ is also an Orlik block, and
a cyclic generator is
\begin{eqnarray}\label{2.10}
e^{(2)}:= p_2(M_h)(e^{(1)}).
\end{eqnarray}

\medskip
(vi) If $\Lambda^{(1)}\subset\Lambda$ is an Orlik block 
with generator $e^{(1)}$ and $\lambda\in E(\Lambda^{(1)})$ is an eigenvalue
of the monodromy on $\Lambda^{(1)}$, then an eigenvector is
\begin{eqnarray}\label{2.11}
v(e^{(1)},\lambda) := \frac{p_{\Lambda^{(1)}}}{t-\lambda}(M_h)
(e^{(1)}).
\end{eqnarray}
And then
\begin{eqnarray}
&&L(v(e^{(1)},\lambda),v(e^{(1)},\oooo{\lambda}))\nonumber\\
&=& 
L(v(e^{(1)},\lambda),\frac{p_{\Lambda^{(1)}}}{t-\oooo{\lambda}}(M_h)
(e^{(1)})) \nonumber \\
&=& L(\frac{p_{\Lambda^{(1)}}}{t-\oooo{\lambda}}(M_h^{-1})
v(e^{(1)},\lambda),e^{(1)})\nonumber \\
&=& \frac{p_{\Lambda^{(1)}}}{t-\oooo{\lambda}}(\oooo{\lambda})\cdot
L(v(e^{(1)},\lambda),e^{(1)})\nonumber\\
&=& \frac{p_{\Lambda^{(1)}}}{t-\oooo{\lambda}}(\oooo{\lambda})\cdot
L(\frac{p_{\Lambda^{(1)}}}{t-\lambda}(M_h)(e^{(1)}),e^{(1)}).
\label{2.12}
\end{eqnarray}
This calculation will be useful in section \ref{c5}.
\end{remarks}

The following two lemmata concern automorphisms of sums of Orlik blocks
(lemma \ref{t2.7}) or of a single Orlik block (lemma \ref{t2.8}).
They will be useful tools in order to show the rigidity
of the Milnor lattices in the sections \ref{c5} and \ref{c6}.

\begin{lemma}\label{t2.7}
Let $(\Lambda,M_h)$ be as above (we will not need $L$ here, only $M_h$).
Let $\Lambda^{(1)},\ldots,\Lambda^{(k)}\subset\Lambda$ be Orlik blocks with 
cyclic generators $e^{(1)},\ldots,e^{(k)}$ and with
$$\Lambda^{(1)}+\ldots+\Lambda^{(k)}
=\Lambda^{(1)}\oplus \ldots \oplus \Lambda^{(k)}.$$
Consider an element $g\in\Aut(\Lambda^{(1)}\oplus\ldots\oplus
\Lambda^{(k)},M_h)$. Then there are unique polynomials
$p_{ij}\in\Z[t]_{<\rank \Lambda^{(j)}}$ for $i,j=1,\ldots,k$
with
\begin{eqnarray}\label{2.13}
g(e^{(j)})=\sum_{i=1}^k p_{ij}(M_h)(e^{(i)}).
\end{eqnarray}

Suppose now that $p_0\in\Z[t]$ divides
$\gcd(p_{\Lambda^{(1)}},\ldots,p_{\Lambda^{(k)}})$ and that
\begin{eqnarray}\label{2.14}
g=\id\textup{ on }\quad 
\Lambda^{(j)}_{p_{\Lambda^{(j)}}/p_0}
\quad\textup{ for any }j,
\end{eqnarray}
so that $g$ acts nontrivial only on
$(\Lambda^{(1)}\oplus\ldots\oplus\Lambda^{(k)})_{p_0}$. Then
\begin{eqnarray}\label{2.15}
p_{ij} = \delta_{ij}+\frac{p_{\Lambda^{(i)}}}{p_0}\cdot q_{ij}
\end{eqnarray}
for suitable polynomials $q_{ij}\in\Z[t]_{<\deg p_0}$.

Suppose furthermore that a unit root $\xi$ satisfies $p_0(\xi)=0$.
Then $g$ with respect to the eigenvectors 
$v(e^{(1)},\xi)\in\Lambda^{(1)}_\xi,\ldots,
v(e^{(k)},\xi)\in\Lambda^{(k)}_\xi$ (defined in \eqref{2.11}) 
is given by 
\begin{eqnarray}\label{2.16}
g(v(e^{(j)},\xi))=\sum_{i=1}^k (\delta_{ij}+\frac{p_{\Lambda^{(j)}}}{p_0}
\cdot q_{ij})(\xi)\cdot v(e^{(i)},\xi)
\end{eqnarray}
\end{lemma}

{\bf Proof:}
Only the part after \eqref{2.13} is nontrivial.
Suppose that $p_0$ and $g$ are as stated above. By assumption
\begin{eqnarray*}
g(e^{(j)})-e^{(j)} &\in& 
(\Lambda^{(1)}\oplus\ldots \oplus\Lambda^{(k)})_{p_0}\\
&\subset& \bigoplus_{i=1}^k \Lambda^{(i)}_{\C,p_0}
=\bigoplus_{i=1}^k \frac{p_{\Lambda^{(i)}}}{p_0}(M_h)(\Lambda^{(i)}_\C).
\end{eqnarray*}
Thus 
$p_{ij}-\delta_{ij}\in \frac{p_{\Lambda^{(i)}}}{p_0}\cdot
\C[t]$, thus
$p_{ij}-\delta_{ij}\in \frac{p_{\Lambda^{(i)}}}{p_0}\cdot
\Z[t]_{<\deg p_0}$. 

The following calculation proves \eqref{2.16}.
\begin{eqnarray*}
g(v(e^{(j)},\xi)) &=& g\left(\frac{p_{\Lambda^{(j)}}}{t-\xi}(M_h)(e^{(j)})\right)
= \frac{p_{\Lambda^{(j)}}}{t-\xi}(M_h)\left(g(e^{(j)})\right)\\
&=& \frac{p_{\Lambda^{(j)}}}{t-\xi}(M_h)\left(
\sum_{i=1}^k \left(\delta_{ij}+
\frac{p_{\Lambda^{(i)}}}{p_0}\cdot q_{ij}\right)(M_h)(e^{(i)})\right)\\
&=& \sum_{i=1}^k \left(\left(\delta_{ij}+
\frac{p_{\Lambda^{(i)}}}{p_0}\cdot q_{ij}\right)\cdot 
\frac{p_{\Lambda^{(j)}}}{t-\xi}\right) (M_h)
(e^{(i)})\\
&=& \sum_{i=1}^k \left(\delta_{ij}+
\frac{p_{\Lambda^{(j)}}}{p_0}\cdot q_{ij}\right)(M_h)
(v(e^{(i)},\xi))\\
&=& \sum_{i=1}^k \left(\delta_{ij}+
\frac{p_{\Lambda^{(j)}}}{p_0}\cdot q_{ij}\right)(\xi)
\cdot v(e^{(i)},\xi).
\end{eqnarray*}
\hfill $\Box$ 

\bigskip

The following lemma is cited from \cite[lemma 8.2]{He11},
but it goes back to arguments in \cite[ch. 6]{He98}.

\begin{lemma}\label{t2.8}
Let $(\Lambda,L,M_h)$ be as above.
Suppose that $\Lambda$ is a single Orlik block.

We make the following nontrivial assumption on the set 
\begin{eqnarray}\label{2.17}
\Ord:=\{\ord\lambda\, |\, \lambda\textup{ eigenvalue of }M_h\}
\subset \Z_{\geq_1}
\end{eqnarray}
of orders of the eigenvalues of the monodromy $M_h$:
There exist four sequences $(m_i)_{i=1,\ldots,|\Ord|}$, $(j(i))_{i=2,\ldots,|\Ord|}$, 
$(p_i)_{i=2,\ldots,|\Ord|}$, $(k_i)_{i=2,\ldots,|\Ord|}$ of numbers in $\Z_{\geq 1}$
and two numbers $i_1,i_2\in\Z_{\geq 1}$ with $i_1\leq i_2\leq |\Ord|$ and 
with the properties: 
\begin{list}{}{}
\item $\Ord=\{m_1,\ldots,m_{|\Ord|}\}$,\
\item $p_i$ is a prime number, $p_i=2$ for $i_1+1\leq i\leq i_2$, 
$p_i\geq 3$ else, 
\item  $j(i)=i-1$ for $i_1+1\leq i\leq i_2$, $j(i)<i$ else,
\item $m_i=m_{j(i)}/p_i^{k_i}.$
\end{list}

Then
\begin{eqnarray}\label{2.18}
\Aut(\Lambda,L,M_h)=\{\pm M_h^k\, |\, k\in \Z\}.
\end{eqnarray}
\end{lemma}

We will need some basic facts for the unit roots
$\zeta=e^{2\pi i /m}$ with $m\in\{10,12,14,18\}$.
The following theorem \ref{t2.9} collects some facts 
for general unit roots.
Theorem \ref{t2.10} cites two classical results
on orders in algebraic number fields.
Lemma \ref{t2.11} puts together some specific properties
for the unit roots of the orders $m\in\{10,12,14,18\}$.

\begin{theorem}\label{t2.9}
Fix $m\in\Z_{\geq 3}$ and define $\zeta:=e^{2\pi i /m}$,
$p_1:=\zeta+\oooo{\zeta}$.

(a) 
\begin{eqnarray*}
\Eiw(\zeta)&:=& \{\pm \zeta^k\, |\, k\in\Z\}\\
&=& \{\textup{unit roots in }\Q(\zeta)\}
=\{\textup{unit roots in }\Z[\zeta]\}\\
&=& \{a\in\Z[\zeta]\, |\, |a|=1\}.
\end{eqnarray*}

(b) $\Z[\zeta]$ is the ring of algebraic integers of $\Q(\zeta)$.

(c) $\Z[p_1]$ is the ring of algebraic integers of $\Q(p_1)$.
And $\Q(p_1)$ is the maximal real subfield of $\Q(\zeta)$.

(d) $\Q(\zeta)$ has class field number 1 and thus $\Z[\zeta]$
is a principal ideal domain if and only if 
$m\in A_1\cup A_2\cup A_3$ where
\begin{eqnarray*}
A_1&=& \{1,3,5,\ldots,21\}\cup\{25,27,33,35,45\},\\
A_2&=& \{2n\, |\, n\in A_1\},\\
A_3&=& \{4n\ |\, n\in A_4\},\quad
A_4=\{1,2,3,\ldots,12\}\cup\{15,21\}.
\end{eqnarray*}

(e) If $\Q(\zeta)$ has class field number 1, then $\Q(p_1)$
has class field number $1$ and thus $\Z[p_1]$ is a 
principal ideal domain.

(f) $\zeta-1\in(\Z[\zeta])^*$ if $m\notin\{p^k\, |\, 
p\textup{ a prime number}, \ k\in\Z_{\geq 1}\}$.

$\zeta+1\in(\Z[\zeta])^*$ if $m\notin\{2\cdot p^k\, |\, 
p\textup{ a prime number}, \ k\in\Z_{\geq 1}\}.$
\end{theorem}

{\bf Proof:}
(a) \cite{Wa97} lemma 1.6 and exercise 2.3.
(b) \cite{Wa97} theorem 2.6.
(c) \cite{Wa97} proposition 2.16.
(d) \cite{Wa97} theorem 11.1.
(e) \cite{Wa97} theorem 4.10.
(f) \cite{Wa97} proposition 2.8. \hfill $\Box$

\begin{theorem}\label{t2.10}
Let $K$ be an algebraic number field of degree $n=s+2t$ over $\Q$
with $s$ real embeddings $\sigma_j:K\to \R$, $j=1,\ldots,s$, and
$2t$ complex embeddings $\sigma_j:K\to\C$, $j=s+1,\ldots,n$,
with $\sigma_{s+t+j}=\oooo{\sigma_{s+j}}$ for $j=1,\ldots,t$.

(a) \cite[Ch. 2, 3.1 Theorem 1]{BS66}
Define $\uuuu{\sigma}:=(\sigma_1,\ldots,\sigma_{s+t}):
K\to \R^s\times\C^t\cong \R^n$. 
Any $\Q$-basis of $K$ maps to an $\R$-basis of $\R^n$.
Thus the image under $\uuuu{\sigma}$ of any order 
$\OO\subset K$ is a lattice of rank $n$ in $\R^n$.

(b) (Dirichlet's unit theorem, \cite[Ch. 2, 4.3 Theorem 5]{BS66})
Let $\OO\subset K$ be an order. 
One can choose $r=s+t-1$ units $a_1,\ldots,a_r\in\OO^*$
such that any unit has a unique representation
$\xi\cdot a_1^{k_1}\cdot \ldots\cdot a_r^{k_r}$ with 
$k_1,\ldots,k_r\in\Z$ and $\xi$ a root of $1$ in $\OO$.
\end{theorem}

Of course, $n=\varphi(m)=2t$ in the case 
$\OO=\Z[\zeta]\subset K=\Q(\zeta)$, 
and $n=\frac{\varphi(m)}{2}$ in the case $\OO=\Z[p_1]\subset
K=\Q[p_1]$, where $\zeta=e^{2\pi i/m}$ and $p_1=\zeta+\oooo\zeta$.

The unit roots of orders $m\in\{10,12,14,18\}$ are most 
important in this paper. The next lemma collects specific
properties of $\Z[\zeta]$ for these orders.

\begin{lemma}\label{t2.11}
Fix $m\in\{10,12,14,18\}$ and define $\zeta=e^{2\pi i /m}$
and $p_1=\zeta+\oooo\zeta$.

$\Z[\zeta]$ and $\Z[p_1]$ are principal ideal domains
(by theorem \ref{t2.9} (d)+(e)).

(a) $m=10$: $\Phi_{10}(t)=t^4-t^3+t^2-t+1$, 
\begin{eqnarray*}
\Z[\zeta]^*&=&\Eiw(\zeta)\cdot \Z[p_1]^*\supset \{\zeta-1\},\\
\Z[p_1]^*&=& \{\pm 1\}\times \{p_1^k\, |\, k\in\Z\}
\supset\{p_1-2,p_1-1,p_1,p_1+1\},\\
p_1&=& \frac{\sqrt{5}+1}{2}>0,\quad 
p_3:=\zeta^3+\oooo{\zeta}^3=\frac{-\sqrt{5}+1}{2}<0,\\
\Gal(\Q(p_1):\Q)&=& \{\id,\varphi\},\quad
\varphi:p_1\mapsto p_3\mapsto p_1,\\
(x-p_1)(x-p_3)&=& x^2-x-1,\quad p_1+p_3=1,\ p_1p_3=-1,\ 
p_1^2=p_1+1.
\end{eqnarray*}

(b) $m=12$: $\Phi_{12}(t)=t^4-t^2+1$, 
\begin{eqnarray*}
\Z[\zeta]^*&=&\Eiw(\zeta)\cdot \Z[p_1]^*
\cup (\zeta+1)\cdot\Eiw(\zeta)\cdot\Z[p_1]^*\\
&=& \Eiw(\zeta)\cdot \{(\zeta+1)^k\, |\, k\in\Z\}
\supset \{\zeta-1,\zeta+1\},\\
\Z[p_1]^*&=& \{\pm 1\}\times \{p_1^k\, |\, k\in\Z\}
\supset\{p_1-2,p_1+2\},\\
p_1&=& \sqrt{3}>0,\quad 
p_5:=\zeta^5+\oooo{\zeta}^5=-\sqrt{3}<0,\\
\Gal(\Q(p_1):\Q)&=& \{\id,\varphi\},\quad
\varphi:p_1\mapsto p_5\mapsto p_1,\\
(x-p_1)(x-p_5)&=& x^2-3,\quad p_1+p_5=0,\ p_1p_5=-3,\ 
p_1^2=3.
\end{eqnarray*}

(c) $m=14$: $\Phi_{14}(t)=t^6-t^5+t^4-t^3+t^2-t+1$, 
\begin{eqnarray*}
\Z[\zeta]^*&=&\Eiw(\zeta)\cdot \Z[p_1]^*\supset \{\zeta-1\},\\
\Z[p_1]^*&=& \{\pm 1\}\times \{p_1^{k_1}p_3^{k_3}\, |\, 
k_1,k_3\in\Z\}\\
&\supset&\{p_1-2,p_1-1,p_1,p_1+1\},\\
p_1&>&0,\ p_3:=\zeta^3+\oooo{\zeta}^3>0,
\ p_5:= \zeta^5+\oooo{\zeta}^5<0,\\
%\end{eqnarray*}
%\begin{eqnarray*}
\Gal(\Q(p_1):\Q)&=& \{\id,\varphi,\varphi^2\},\quad
\varphi:p_1\mapsto p_3\mapsto p_5\mapsto p_1,\\
(x-p_1)(x-p_3)(x-p_5)&=& x^3-x^2-2x+1,\quad 
p_1+p_3+p_5=1,\\ 
p_1p_3p_5&=&-1,\quad p_1p_3=p_1-1,\quad p_1^2=-p_5+2.
\end{eqnarray*}

(d) $m=18$: $\Phi_{18}(t)=t^6-t^3+1$, 
\begin{eqnarray*}
\Z[\zeta]^*&=&\Eiw(\zeta)\cdot \Z[p_1]^*\supset \{\zeta-1\},\\
\Z[p_1]^*&=& \{\pm 1\}\times \{p_1^{k_1}p_5^{k_5}\, |\, 
k_1,k_5\in\Z\}\\
&\supset&\{p_1-2,p_1,p_1+1\},\\
p_1&>&0,\ p_5:=\zeta^5+\oooo{\zeta}^5<0,
\ p_7:= \zeta^7+\oooo{\zeta}^7<0,\\
\Gal(\Q(p_1):\Q)&=& \{\id,\varphi,\varphi^2\},\quad
\varphi:p_1\mapsto p_5\mapsto p_7\mapsto p_1,\\
(x-p_1)(x-p_5)(x-p_7)&=& x^3-3x-1,\quad 
p_1+p_5+p_7=0,\\ 
p_1p_5p_7&=&1,\quad p_1p_5=-p_5-1,\quad p_1^2=-p_7+2.
\end{eqnarray*}
\end{lemma}

{\bf Proof:} 
That the index $[\Z[\zeta]^*:\Eiw(\zeta)\cdot \Z[p_1]^*]$
is $1$ for $m\in\{10,14,18\}$ and $2$ for $m=12$,
follows from \cite[theorem 4.12 and corollary 4.13]{Wa97}.
That $\Z[p_1]^*$ is as stated, follows for $m\in\{10,14,18\}$
from \cite[theorem 8.2 and lemma 8.1 (a)]{Wa97}.
For $m=12$ \cite[\S 8.1]{Wa97} is not so useful,
but there the proof of $\Z[p_1]^*=\{\pm 1\}\cdot\{p_1^k\, |\,
k\in\Z\}$ is easy. Everything else is elementary.
\hfill$\Box$

\bigskip
Part (b) of the following lemma applies with 
$\Lambda=Ml(f)$ and $\Lambda^{(1)}=\www B_1\oplus B_2$
(see the theorems \ref{t5.1} and \ref{t6.1})
to most of the Milnor lattices in the sections \ref{c5} 
and \ref{c6}. We will need \eqref{2.19}.

\begin{lemma}\label{t2.12}
(a) Let $p=\prod_{i\in I}\Phi_{m_i}$ be a product of
cyclotomic polynomials. Then
$p(1)\equiv 1(2)$ if and only if all
$m_i\in\Z_{\geq 1}-\{2^k\, |\, k\in\Z_{\geq 0}\}$.

\medskip
(b) Let $(\Lambda,L,M_h)$ be as above
(we will not need $L$ here, only $M_h$).
Let $\Lambda^{(1)}\subset \Lambda$ be an $M_h$-invariant
sublattice with $[\Lambda:\Lambda^{(1)}]=2$. Write
\begin{eqnarray*}
p_\Lambda&=& p_1\cdot p_2\quad\textup{with }
p_j=\prod_{m\in J_j}\Phi_m\\
\textup{and} && 
J_1\subset \Z_{\geq 1}-\{2^k\, |\, k\in\Z_{\geq 0}\},\quad
J_1\subset \{2^k\, |\, k\in\Z_{\geq 0}\}.
\end{eqnarray*}
Then $J_2\neq\emptyset$, $p_2\neq 1$, and 
\begin{eqnarray}\label{2.19}
\Lambda_p&=& \Lambda_p^{(1)}\qquad \textup{ for any }
p\textup{ with }p|p_1,\\
{}[\Lambda_p:\Lambda_p^{(1)}]&=& 2\qquad \textup{ for any }
p\textup{ with }p_2|p. \label{2.20}
\end{eqnarray}
\end{lemma}

{\bf Proof:}
(a) Observe $\Phi_{2^k}(t)=t^{2^{k-1}}+1$ for $k\geq 1$ and 
\begin{eqnarray}\label{2.21}
t^{2^k\cdot q}-1&=& (t^{2^k}-1)(t^{2^k(q-1)}+t^{2^k(q-2)}+\ldots
+t^{2^k}+1).
\end{eqnarray}
For odd $q>1$, the second factor has at $t=1$ the odd value $q$.
Therefore $\Phi_m(1)\equiv 1(2)$ for any $m$ with 
$2^k|m|2^k\cdot q$ and $2^k\neq m$ with $q$ odd.

\medskip
(b) For an arbitrary element $\gamma\in\Lambda-\Lambda^{(1)}$,
\begin{eqnarray*}
\Lambda-\Lambda^{(1)}&=& \gamma+\Lambda^{(1)}.
\end{eqnarray*}
This set is $M_h$-invariant because $\Lambda^{(1)}$ is 
$M_h$-invariant.
Thus for any $k\in\Z_{\geq 1}$  
$M_h^k(\gamma)\in \Lambda-\Lambda^{(1)}$. 
By part (a) $p_1(1)\equiv 1(2)$.
Thus $p_1(M_h)(\gamma)\in \Lambda-\Lambda^{(1)}$ and
\begin{eqnarray*}
p_1(M_h)(\Lambda-\Lambda^{(1)})\subset \Lambda-\Lambda^{(1)}.
\end{eqnarray*}
On the other hand
\begin{eqnarray*}
p_1(M_h)(\Lambda_{p_1})&=&\{0\}\subset\Lambda^{(1)},
\quad\textup{thus }\Lambda_{p_1}\subset \Lambda^{(1)},
\quad\textup{thus }\eqref{2.19}.\\
p_1(M_h)(\Lambda)&\subset& \Lambda_{p_2},\quad
\textup{thus }\Lambda_{p_2}\cap(\Lambda-\Lambda^{(1)})\neq
\emptyset,\quad\textup{thus }\eqref{2.20}.
\end{eqnarray*}
\hfill  $\Box$

\section{Some Fuchsian groups}\label{c3}
\setcounter{equation}{0}

\begin{notations}\label{t3.1}
For any $m\in\Z_{\geq 3}$ define $\zeta:=e^{2\pi i/m}$
and $p_1:=\zeta+\oooo\zeta$. 
The letter $\xi$ will denote in this section 
a primitive $m$-th unit root.
An element of $\Q(\zeta)$ will be written 
as $a$ or $a(\zeta)$. Then $a(\xi)$ is the image 
$\varphi(a)$ for $\varphi\in\Gal(\Q(\zeta):\Q)$ with
$\varphi(\zeta)=\xi$.
\end{notations}

Any element $A=\begin{pmatrix}a&b\\c&d\end{pmatrix}\in GL(2,\C)$
acts on $\P^1\C$ by the linear transformation
$z\mapsto \frac{az+b}{cz+d}$, which is an automorphism
of $\P^1\C$. The {\it limit set} 
$L(\Gamma)\subset\P^1\C$ of a subgroup $\Gamma\subset GL(2,\C)$
is \cite[III 1B]{Le64}
\begin{eqnarray*}
L(\Gamma)&=&\{z\in\P^1\C\, |\, \exists\ z_0\in\P^1\C\textup{ and }
\exists\textup{ a sequence of different}\\
&&\textup{elements }\gamma_i\in\Gamma\textup{ with }
\gamma_i(z_0)\to z\}.
\end{eqnarray*}
A subgroup $\Gamma\subset GL(2,\C)$ and the induced subgroup
of $PGL(2,\C)$ are called {\it Fuchsian} if $\Gamma$
maps a certain circle $C\subset\P^1\C$ to itself and
$L(\Gamma)\subset C$. By a theorem of Poincar\'e
\cite[III 3I]{Le64}, a subgroup $\Gamma\subset GL(2,\C)$
is Fuchsian if it maps a certain circle $C\subset\P^1\C$
to itself and is discrete in $GL(2,\C)$.

In the sections \ref{c5} and \ref{c6} we will encounter
Fuchsian groups which arise in the following way.

\begin{theorem}\label{t3.2}
Let $m\in\Z_{\geq 3}$, $\zeta:=e^{2\pi i /m}$, 
$p_1:=\zeta+\oooo\zeta$, and $w=w(\zeta)\in\Q(\zeta)$ with
\begin{eqnarray}\label{3.1}
w(\zeta)>0&&(\textup{thus }
w(\zeta)=w(\oooo\zeta)\in\Q(p_1)),\\
w(\xi)<0&&\textup{for any primitive }m\textup{-th unit root }
\xi\notin\{\zeta,\oooo\zeta\}.\label{3.2}
\end{eqnarray}
Then the matrix group
\begin{eqnarray}\label{3.3}
\Gamma:=\{A\in GL(2,\Z[\zeta])\, |\, 
\begin{pmatrix}-1&0\\0&w\end{pmatrix}
=A^t\begin{pmatrix}-1&0\\0&w\end{pmatrix}
\oooo{A}\}
\end{eqnarray}
is an infinite Fuchsian group. It preserves the circle
\begin{eqnarray}\label{3.4}
C=\{z\in\C\, |\, |z|^2=w\}.
\end{eqnarray}
The map
\begin{eqnarray}
\{(a,c,\delta)\in\Z[\zeta]^2\times\Eiw(\zeta)
\, |\, |a|^2-1=w\cdot |c|^2\}
\to \Gamma\nonumber\\
(a,c,\delta)\mapsto
A:=\begin{pmatrix}a & w\cdot \oooo c\cdot 
\delta\\ c & \oooo a\cdot\delta\end{pmatrix}
\label{3.5}
\end{eqnarray}
is a bijection 
(here $\Eiw(\zeta)=\{\pm\zeta^k\, |\, k\in\Z\}$, see 
theorem \ref{t2.9} (a)).
\end{theorem}

{\bf Proof:}
The matrix $\begin{pmatrix}-1&0\\0&w\end{pmatrix}$
defines an indefinite hermitian form on $\C^2$.
The isotropic lines are $\C\cdot \begin{pmatrix}z\\1\end{pmatrix}$
with $z\in C$. Therefore any matrix $A\in\Gamma$
maps $C$ to itself.

The matrix equation which defines $\Gamma$ can be spelled out
as follows,
\begin{eqnarray}
\begin{pmatrix}-1&0\\0&w\end{pmatrix}
&=&\begin{pmatrix}a&c\\b&d\end{pmatrix}
\begin{pmatrix}-1&0\\0&w\end{pmatrix}
\begin{pmatrix}\oooo a&\oooo b\\ \oooo c&\oooo d\end{pmatrix}
\nonumber\\
&=&\begin{pmatrix}-a\oooo a+wc\oooo c & 
-a\oooo b+wc\oooo d \\ -\oooo a b+w\oooo c d & 
-b\oooo b+wd\oooo d\end{pmatrix}.\label{3.6}
\end{eqnarray}
The determinant $\delta=\det A=ad-bc$ is in
$\Z[\zeta]$ and has absolute value 1, so it is in
$\Eiw(\zeta)$ by theorem \ref{t2.9} (a).
The equations above give
\begin{eqnarray}\label{3.7}
\oooo a \delta &=& \oooo a(ad-bc)=(wc\oooo c+1)d-(w\oooo cd)c=d,\\
w\oooo c\delta &=& w\oooo c(ad-bc)=(\oooo ab)a-(a\oooo a-1)b=b.
\nonumber
\end{eqnarray}
This yields the bijection \eqref{3.5}.

The defining equation 
\begin{eqnarray}\label{3.8}
|a(\zeta)|^2-1=w(\zeta)\cdot |c(\zeta)|^2
\end{eqnarray}
for the pairs $(a(\zeta),c(\zeta))\in\Z[\zeta]^2$ on the left
hand side of \eqref{3.5} is in the case 
$(a,c)\in\Z[p_1]^2$ and $w(\zeta)\in\Z[p_1]$ a {\it Pell equation}.
We obtain the inequalities
\begin{eqnarray}
0&\leq& |c(\zeta)|^2=w(\zeta)^{-1}(|a(\zeta)|^2-1),\nonumber\\
|a(\zeta)|&\geq& 1\label{3.9}
\end{eqnarray}
and
\begin{eqnarray}
0&\leq& |c(\xi)|^2=(-w(\xi))^{-1}(1-|a(\xi)|^2)<(-w(\xi))^{-1},
\nonumber\\
|a(\xi)|&\leq& 1 \quad\textup{for any primitive }m\textup{-th
unit root }\xi\notin\{\zeta,\oooo\zeta\}.\label{3.10}
\end{eqnarray}

$\Gamma$ maps $C$ to itself. Therefore by Poincar\'e's
theorem, it is a Fuchsian group if it is a discrete matrix
group. This holds if the set
\begin{eqnarray*}
P_1:=\{a\in\Z[\zeta]\, |\, \exists\ c\in\Z[\zeta]
\textup{ with }|a|^2-1=w\cdot |c|^2\}
\end{eqnarray*}
intersects each compact set $K\subset \C$ in a finite set.

The embedding $\uuuu{\sigma}:\Q(\zeta)\to\R^{\varphi(n)}$
from theorem \ref{t2.10} (a) maps $\Z[\zeta]$ to a lattice
in $\R^{\varphi(n)}$. 
Because of \eqref{3.10}, it maps $P_1\cap K$ to a subset of
\begin{eqnarray*}
\uuuu\sigma(\Z[\zeta])\cap \left( K\times
\{z\in\C\, |\, |z|\leq 1\}^{\varphi(n)/2-1}\right) .
\end{eqnarray*}
This is a finite set. Therefore $\Gamma$ is a Fuchsian group.

The next lemma shows that the set $P_1$ and the group
$\Gamma$ contain infinitely many elements.
\hfill $\Box$

\begin{lemma}\label{t3.3}
Let $m\in\Z_{\geq 3}$, $\zeta,p_1$ and $w\in\Q(p_1)$ be
as in theorem \ref{t3.2}. Then the set
\begin{eqnarray}\label{3.11}
P_2:=\{(a,c)\in\Z[p_1]\, |\, a^2-1=w\cdot c^2\}
\end{eqnarray}
contains infinitely many pairs. 
If $w\in\Z[p_1]$, then $P_2$ contains pairs $(a,c)$ with 
$w|(a-1)$.
\end{lemma}

{\bf Proof:}
If $\www w=w\cdot u^2$ for some $u\in\Z[p_1]-\{0\}$ then
a pair $(a,\www c)\in\Z[p_1]^2$ with $a^2-1=\www w\cdot\www c^2$
induces a pair $(a,c)=(a,\www c\cdot u)$ in $P_2$. 
Therefore we can suppose $w\in\Z[p_1]$.

We will now construct infinitely many units in 
$\Z[\sqrt{w},p_1]^*-\Z[p_1]^*$
and from them infinitely many pairs $(a,c)$ in $P_2$. 

The algebraic number field $\Q(\sqrt{w},p_1)$ has degree
$\varphi(m)$ over $\Q$ and two real embeddings and 
$\varphi(m)-2$ complex embeddings, because of 
\eqref{3.1} and \eqref{3.2}.
By Dirichlet's unit theorem (theorem \ref{t2.10} (b)),
the unit group $\Z[\sqrt{w},p_1]^*$ of the 
order $\Z[\sqrt{w},p_1]$ in $\Q(\sqrt{w},p_1)$
contains a free abelian group of rank 
$2+\frac{\varphi(m)-2}{2}-1=\frac{\varphi(m)}{2}$.

The unit group $\Z[p_1]^*$ contains only a free abelian group
of rank $\frac{\varphi(m)}{2}-1$. Therefore infinitely many
units $a_1+\sqrt{w}c_1\in\Z[\sqrt{w},p_1]^*$
with $a_1\neq 0$ and $c_1\neq 0$ exist.
Then also $a_1-\sqrt{w}c_1$,
\begin{eqnarray*}
(a_1+\sqrt{w}c_1)^2&=&(a_1^2+wc_1^2)+\sqrt{w}(2a_1c_1)=:
a_2+\sqrt{w}c_2,\\
\textup{and}\quad 
h&:=&(a_1+\sqrt{w}c_1)(a_1-\sqrt{w}c_1)=a_1^2-wc_1^2
\end{eqnarray*}
are units, $h$ being in $\Z[p_1]^*$. Then
\begin{eqnarray}\label{3.12}
(a_3,c_3):=(\frac{a_2}{h},\frac{c_2}{h})\in P_2
\end{eqnarray}
because
\begin{eqnarray*}
a_3^2-wc_3^2&=& h^{-2}(a_2^2-wc_2^2)
=h^{-2}(a_2+\sqrt{w}c_2)(a_2-\sqrt{w}c_2)\\
&=& h^{-2}(a_1+\sqrt{w}c_1)^2(a_1-\sqrt{w}c_1)^2=1.
\end{eqnarray*}
Only finitely many units $a_1+\sqrt{w}c_1$ can give the 
same pair $(a_3,c_3)$. Therefore there are infinitely many
pairs $(a_3,c_3)$ in $P_2$.

For the last statement, suppose that $(a_4,c_4)\in P_2$
with $c_4\neq 0$. Then the pair 
$(a_5,c_5):=(a_4^2+wc_4^2,2a_4c_4)$ is also in $P_2$,
\begin{eqnarray*}
a_5^2-wc_5^2 &=& (a_5+\sqrt{w}c_5)(a_5-\sqrt{w}c_5)\\
&=& (a_4+\sqrt{w}c_4)^2(a_4-\sqrt{w}c_4)^2 
= (a_4-\sqrt{w}c_4)^2=1.
\end{eqnarray*}
And it satisfies $w|(a_5-1)$ because of
\begin{eqnarray*}
a_5-1 =a_4^2+wc_4^2-1=2wc_4^2.
\end{eqnarray*}
\hfill $\Box$

\begin{remarks}\label{t3.4}
(i) The equation $a^2-1=wc^2$ is for $w\in\Z[p_1]$ a 
{\it Pell equation}. A generalization of lemma \ref{t3.3}
is theorem 3 in \cite{Sch06}.

\medskip
(ii) The notion of an {\it arithmetic Fuchsian group}
is defined in \cite[ch 9.2]{Sh71}.
The group $\Gamma$ in theorem \ref{t3.2} is in fact
an arithmetic Fuchsian group. This would follow immediately
from \cite[theorem 2]{Ta75}, if it were clear a priori
that $\Gamma$ is a Fuchsian group of the first kind,
i.e. a Fuchsian group with limit set $L(\Gamma)=C$.
It follows with some work from a comparison of the data
in theorem \ref{t3.2} with the data in 
\cite[ch. 9.2]{Sh71}.

\medskip
(iii) The five triangle groups below in theorem \ref{t3.6}
are arithmetic triangle groups. They are in the list in
\cite[theorem 3]{Ta77} of all 85 arithmetic triangle groups.

\medskip
(iv) Theorem \ref{t3.2} and lemma \ref{t3.3} will be used
in the steps 2 and 4 in the proof of theorem \ref{t5.1}
on the groups $G_\Z$ for the bimodal series.
\end{remarks}

\begin{remarks}\label{t3.5}
(i) The triangle groups  below in theorem \ref{t3.6} will
arise in theorem \ref{t6.1} as quotients of the groups $G_\Z$ 
for the quadrangle singularities.

\medskip
(ii) There the first six of the eight elements $w(\zeta)$
in table \eqref{5.71} in the case $r=0$ will be used.
So here $W_{1,0}$ and $S_{1,0}$ are seen as 0-th members
of the series $W_{1,p}^\sharp$ and $S_{1,p}^\sharp$,
not the series $W_{1,p}$ and $S_{1,p}$.

\medskip
(iii) Using the notations and formulas from lemma \ref{t2.11},
the first six of the eight elements $w(\zeta)$ in table 
\eqref{5.71}  in the case $r=0$ can be written as follows.
In the case $U_{1,0}$ we change from $m=9$ to $m=18$,
so below $\zeta=e^{2\pi i/18}$ for $E_{3,0}$ and $U_{1,0}$.
\begin{eqnarray}
W_{1,0}&:& w(\zeta)=\frac{6}{(2-p_1)p_1}
=\frac{1}{(2-p_1)(2+p_1)}
\cdot 2p_1(p_1+2).\nonumber\\
S_{1,0}&:& w(\zeta)=\frac{-2}{(-p_3)(-p_3-1)}
=1\cdot 2p_1^3.\nonumber\\
U_{1,0}&:& w(\zeta)=\frac{-3}{(2+p_7)(1-p_1)}
=1 \cdot p_1(p_1+2).\nonumber\\
E_{3,0}&:& w(\zeta)=\frac{3(2-p_1)}{(p_1+2)(p_1-1)}
=(2-p_1)^2 \cdot p_1(p_1+2).\nonumber\\
Z_{1,0}&:& w(\zeta)=\frac{1}{-p_5}= 1\cdot (-p_5)^{-1}
=1\cdot (p_1-1).\nonumber\\
Q_{2,0}&:& w(\zeta)=\frac{2-p_1}{p_1+1}=(2-p_1)
\cdot \frac{1}{p_1+1}.\label{3.13}
\end{eqnarray}

\medskip
(iv) In theorem \ref{t3.2} one can replace $w$ by 
$\www w:=w\cdot u\oooo u$ for any 
$u\in\Z[\zeta]^*$. The group $\Gamma$ for $w$ and
the group $\www\Gamma$ for $\www w$ are isomorphic, 
and the triples in \eqref{3.5}
are related by 
\begin{eqnarray*}
(\www a,\www c,\www \delta)=(a,c\cdot u^{-1},\delta).
\end{eqnarray*}
We can choose $u$ such that $\www w$ is simpler to work
with than $w$.  In the products for $w$ in (iii), the left
terms are of the form $u\oooo u$ for a suitable
unit $u\in\Z[\zeta]^*$. The right terms are
$\www w$. We will work with the terms $\www w$ 
in theorem \ref{t3.6}.
\end{remarks}

\begin{theorem}\label{t3.6}
The image in $PGL(2,\C)$ of the group $\Gamma$ in theorem 
\ref{t3.2} for the following values of $m$ and $w$
\begin{eqnarray}\label{3.14}
\begin{array}{l|l|l|l|l|l}
&W_{1,0} & S_{1,0} & E_{3,0}\ \&\ U_{1,0} & Z_{1,0} & Q_{2,0} \\
m& 12 & 10 & 18  & 14 & 12 \\
w & 2p_1(p_1+2) & 2p_1^3 & p_1(p_1+2) & (-p_5)^{-1} & (p_1+1)^{-1}
\end{array}
\end{eqnarray}
is a Schwarzian triangle group of the following type:
\begin{eqnarray}\label{3.15}
\begin{array}{l|l|l|l|l}
W_{1,0} & S_{1,0} & E_{3,0}\ \&\ U_{1,0} & Z_{1,0} & Q_{2,0} \\
(2,12,12) & (2,10,10) & (2,3,18) & (2,3,14) & (2,3,12) 
\end{array}
\end{eqnarray}
\end{theorem}

{\bf Proof:}
The proof has three steps. In step 1, we will present
two matrices $A_1$ and $A_2$ in $\Gamma$ whose images in
$PGL(2,\C)$ are elliptic and generate in each case
a Schwarzian triangle group of the claimed type.
We will prove this.
In step 2, we will show that no matrix in $\Gamma$ is closer
to $A_1$ than $A_2$. This will be used in step 3 to prove
that the images in $PGL(2,\C)$ of $A_1$ and $A_2$ 
generate the image of $\Gamma$ in $PGL(2,\C)$.
The steps 1 and 3 together give theorem \ref{t3.6}.

\medskip
{\bf Step 1:}
One checks easily with \eqref{3.5} that the following 
matrices $A_1$ and $A_2$ are in $\Gamma$.
\begin{eqnarray}\label{3.16}
A_1=\begin{pmatrix}\zeta & 0 \\ 0 & 1\end{pmatrix}
\qquad\textup{for all 5 cases.}
\end{eqnarray}
\begin{eqnarray}
\begin{array}{ll}
W_{1,0}: & A_2=
\begin{pmatrix} p_1+2 & -2p_1(p_1+2)\\ 1 & -(p_1+2)\end{pmatrix},
\quad \det A_2=-1,\\
S_{1,0}: & A_2=
\begin{pmatrix} (\zeta+1)p_1 & -2p_1^3\zeta\\ 
1 & -(\zeta+1)p_1\end{pmatrix},
\quad \det A_2=-\zeta,\\
E_{3,0}\ \&\ U_{1,0}: & A_2=
\begin{pmatrix} p_1+1 & -p_1(p_1+2)\\ 1 & -(p_1+1)\end{pmatrix},
\quad\det A_2=-1,\\
Z_{1,0}: & A_2=p_1(1-\zeta^3)\cdot
\begin{pmatrix} 1& -(-p_5)^{-1} \\
1 & -1 \end{pmatrix},
\quad\det A_2=\zeta^3,\\
Q_{2,0}: & A_2=
\begin{pmatrix} \zeta+1 & -\zeta\\ 
p_1+1 & -(\zeta+1)\end{pmatrix}, 
\quad\det A_2=-\zeta.
\end{array}\label{3.17}
\end{eqnarray}

A matrix $A\in GL(2,\C)$ is elliptic if its eigenvalues
$\lambda_1$ and $\lambda_2$ satisfy
$\frac{\lambda_2}{\lambda_1}\in S^1$. 
Let $\begin{pmatrix}z_j\\1\end{pmatrix}$ be an eigenvector
with eigenvalue $\lambda_j$ for $j=1,2$ 
(possibly $z_1=0$ and $z_2=\infty$).
Then the linear transformation of $A$ is a rotation around
the fixed point $z_1$ with angle 
$\alpha(A)=\arg\frac{\lambda_2}{\lambda_1}$.
For $A\in\Gamma$ elliptic we number the eigenvalues
$\lambda_1,\lambda_2$ such that $|z_1|<|z_2|$, so then
$|z_1|^2<w$ and $z_1$ is in the interior of the circle $C$.
One sees in all 5 cases
\begin{eqnarray}\label{3.18}
\lambda_1(A_1)=1,\quad \lambda_2(A_1)=\zeta,\quad
\alpha(A_1)=\frac{2\pi}{m},\\
\tr(A_2)=0,\quad \alpha(A_2)=\pi.\label{3.19}
\end{eqnarray}
The following table lists for the product $A_1A_2$
the eigenvalues $\lambda_1,\lambda_2$ and the
angle $\alpha_1(A_3)$.
\begin{eqnarray}\label{3.20}
\begin{array}{llll}
 & \lambda_1&\lambda_2&\alpha \\
W_{1,0} & \zeta^4&\zeta^3&\frac{-2\pi}{12} \\
S_{1,0} & \zeta^4&\zeta^3&\frac{-2\pi}{10} \\
E_{3,0}\ \&\ U_{1,0}& \zeta^8&\zeta^2&\frac{-2\pi}{3} \\
Z_{1,0} & e^{2\pi i/6}\zeta^2 & e^{-2\pi i/6}\zeta^2 & 
\frac{-2\pi}{3} \\
Q_{2,0} & \zeta^6&\zeta^2&\frac{-2\pi}{3}
\end{array}
\end{eqnarray}
Therefore the images of $A_1$ and $A_2$ in $PGL(2,\C)$
generate a Schwarzian triangle group of the type in
table \eqref{3.15} \cite[VII 1G]{Le64}.

\medskip
{\bf Step 2:}
Write $A_2=\begin{pmatrix}a_2&b_2\\c_2&d_2\end{pmatrix}$
and write $A=\begin{pmatrix}a&b\\c&d\end{pmatrix}$
for any $A\in\Gamma$.

\medskip
{\bf Claim 1:} Any $A\in\Gamma$ with $c\neq 0$ satisfies
$|a|\geq |a_2|$.

\medskip
The proof consists in making the proof of theorem \ref{t3.2}
more constructive. 

First we look for candidates $f\in\Z[p_1]$
of $|a|^2$ which are compatible with the inequalities
\eqref{3.9} and \eqref{3.10} and which satisfy 
$f<|a_2|$. 
Then we will show that these candidates are not compatible
with the equality $|a|^2=1+w\cdot|c|^2$.

Denote by $\uuuu\sigma^\R=
(\sigma^\R_1,\ldots,\sigma^\R_{\varphi(m)/2}): 
\Q(p_1)\to\R^{\varphi(m)/2}$ the tuple of the embeddings
$\sigma_j^\R:\Q(p_1)\to\R$. 
Then $\uuuu\sigma^\R(\Z[p_1])$ is a $\Z$-lattice in
$\R^{\varphi(m)/2}$. The candidates are the numbers
$f=f(p_1)$ in $\Z[p_1]$ with 
\begin{eqnarray}\label{3.21}
\uuuu\sigma^\R(f)&\in&  ]1,|a_2|^2[\ \times\ 
\ ]0,1[^{\varphi(m)/2-1}.
\end{eqnarray}
This follows from the inequalities \eqref{3.9} and \eqref{3.10}.
With sufficient numerical precision of the numbers  $p_j$
in lemma \ref{t2.11}, it is easy to find these candidates.
They are as follows.
\begin{eqnarray*}
W_{1,0}&:& f(p_1)=\alpha\cdot 1+\beta\cdot p_1,\quad
(\alpha,\beta)\in\{(2,1),(4,2),(6,3)\}.\\
S_{1,0}&:& f(p_1)=\alpha\cdot 1+\beta\cdot p_1,
\quad (\alpha,\beta)\in\{(2,2),(2,3)\}.\\
E_{3,0}& \& & U_{1,0}\ :\ \emptyset.\\
Z_{1,0}&:& \emptyset.\\
Q_{2,0}&:& \emptyset.
\end{eqnarray*}
All these candidates will be excluded with the help of the 
condition 
\begin{eqnarray*}
\Norm(|a|^2-1)=\Norm(w\cdot |c|^2)=\Norm(w)\cdot \Norm(|b|^2).
\end{eqnarray*}
Here the norm is the norm in $\Q(p_1)$ and $\Z[p_1]$ 
with values in $\Q$ respectively $\Z$. 

The case $W_{1,0}$: $\Norm(w)=-12,\quad \Norm(1+p_1)=-2,\quad
\Norm(3+2p_1)=-3,\quad \Norm(5+3p_1)=-2.$

The case $S_{1,0}$: $\Norm(w)=-4,\quad \Norm(1+2p_1)=-1,\quad 
\Norm(1+3p_1)=-5$.

\medskip
{\bf Step 3:} It is sufficient to show the following claim 2.

\medskip
{\bf Claim 2:} For any matrix $A_3\in \Gamma$ with 
$c_3\neq 0$, a number $k\in\Z$ exists such that the product
\begin{eqnarray}\label{3.22}
A_4:=A_3\cdot A_1^{-k}A_2A_1^k=
\begin{pmatrix}a_3&b_3\\c_3&d_3\end{pmatrix}
\begin{pmatrix}a_2&\zeta^{-k}b_2\\ \zeta^kc_2&d_2\end{pmatrix}
\end{eqnarray}
satisfies
\begin{eqnarray}\label{3.23}
|c_4|<|c_3|,\qquad\textup{here }
c_4=c_3a_2+\zeta^kd_3c_2.
\end{eqnarray}

\medskip
We can choose $k\in\Z$ such that
\begin{eqnarray}\label{3.24}
\beta:=|\arg(c_3a_2)-\arg(-\zeta^kd_3c_2)|\leq\frac{\pi}{m}.
\end{eqnarray}
Observe
\begin{eqnarray}\label{3.25}
\frac{|\zeta^kd_3c_2|^2}{|c_3a_2|^2}
=\frac{|a_3|^2\frac{|a_2|^2-1}{w(\zeta)}}
{\frac{|a_3|^2-1}{w(\zeta)}|a_2|^2}
=\frac{1-|a_2|^{-2}}{1-|a_3|^{-2}}.
\end{eqnarray}
The trivial inequality $1-|a_3|^{-2}<1$ and the
inequality $|a_3|\geq |a_2|$ from step 2 give the
inequalities
\begin{eqnarray}\label{3.26}
\left(1-|a_2|^{-2}\right) |c_3a_2|^2
<|\zeta^kd_3c_2|^2\leq |c_3a_2|^2.
\end{eqnarray}
Observe also
\begin{eqnarray}\label{3.27}
\sqrt{1-|a_2|^{-2}}<\cos\frac{\pi}{m}.
\end{eqnarray}
Therefore 
\begin{eqnarray}
|c_4|&=& 
|c_3a_2|^2(\sin\beta)^2+
\left(|c_3a_2|\cos\beta-|d_3c_2|\right)^2\nonumber\\
&<& |c_3a_2|^2(\sin\frac{\pi}{m})^2+
\left(1-\sqrt{1-|a_2|^{-2}}\right)^2\cdot |c_3a_2|^2
\nonumber\\
&=& |c_3|^2\cdot |a_2|^2\left((\sin\frac{\pi}{m})^2
+\left(1-\sqrt{1-|a_2|^{-2}})\right)^2\right)\nonumber\\
&\stackrel{(*)}{<}& |c_3|^2.\label{3.28}
\end{eqnarray}
$\stackrel{(*)}{<}$ follows in all 5 cases by an explicit
calculation.
\hfill $\Box$

\section{Review on the topology of singularities}\label{c4}
\setcounter{equation}{0}

\noindent
In this section, we recall some classical facts about the
topology of singularities, and we fix some notations.

An {\it isolated hypersurface singularity}
(short: {\it singularity})
is a holomorphic function germ $f:(\C^{n+1},0)\to (\C,0)$ 
with an isolated singularity at $0$. 
Its {\it Jacobi ideal} is 
$$J(f):=(\frac{\ppp f}{\ppp x_0},\ldots,
\frac{\ppp f}{\ppp x_n})\subset \OO_{\C^{n+1},0}.$$
Its {\it Jacobi algebra} is $\OO_{\C^{n+1},0}/J(f)$.
Its {\it Milnor number}  $\mu:=\dim\OO_{\C^{n+1},0}/J(f)$
is finite. 
For the following notions and facts compare \cite{AGV88} and 
\cite{Eb07}.
A {\it good representative} of $f$ has to be defined with some 
care \cite{Mi68}\cite{AGV88}\cite{Eb07}. It is $f:X\to \Delta$
with $\Delta=\{\tau\in\C\, |\, |\tau|<\delta\}$ 
a small disk around 0 and 
$X=\{x\in\C^{n+1}\, |\, |x|<\varepsilon\}\cap f^{-1}(\Delta)$ 
for some sufficiently small $\varepsilon>0$
(first choose $\varepsilon$, then $\delta$).
Then $f:X'\to \Delta'$ with $X'=X-f^{-1}(0)$ and 
$\Delta'=\Delta-\{0\}$ is a locally trivial $C^\infty$-fibration,
the  {\it Milnor fibration}. Each fiber has the
homotopy type of a bouquet of $\mu$ $n$-spheres \cite{Mi68}.

Therefore the (reduced for $n=0$) 
middle homology groups are {}\\{}
$H_n^{(red)}(f^{-1}(\tau),\Z) \cong \Z^\mu$ for $\tau\in \Delta'$.
Each comes equipped with an intersection form $I$, 
which is a datum of one fiber,
a monodromy $M_h$ and a Seifert form $L$, which come from the 
Milnor fibration,
see \cite[I.2.3]{AGV88} for their definitions.
$M_h$ is a quasiunipotent automorphism, $I$ and $L$ are 
bilinear forms with values in $\Z$,
$I$ is $(-1)^n$-symmetric, and $L$ is unimodular. $
L$ determines $M_h$ and $I$ because of the formulas
\cite[I.2.3]{AGV88}
\begin{eqnarray}\label{4.1}
L(M_ha,b)&=&(-1)^{n+1}L(b,a),\\ \label{4.2}
I(a,b)&=&-L(a,b)+(-1)^{n+1}L(b,a)=L((M-\id)a,b).
\end{eqnarray}
\eqref{4.2} tells especially that $\ker(M_h-\id)$ is the radical
of $I$ and that $L$ is $(-1)^{n+1}$-symmetric on this radical.
The semisimple part of $M_h$ is called $M_s$, the unipotent
part $M_u$, the nilpotent part $N=\log M_u$.

The Milnor lattices $H_n(f^{-1}(\tau),\Z)$ for all 
Milnor fibrations $f:X'\to \Delta'$ and then all 
$\tau\in\R_{>0}\cap T'$ are canonically isomorphic,
and the isomorphisms respect $M_h$, $I$ and $L$. 
This follows from Lemma 2.2 in \cite{LR73}. 
These lattices are identified and called 
{\it Milnor lattice} $Ml(f)$.

The group $G_\Z$ is 
\begin{eqnarray}\label{4.3}
G_\Z=G_\Z(f):= \Aut(Ml(f),L)=\Aut(Ml(f),M_h,I,L),
\end{eqnarray}
the second equality is true because $L$ determines $M_h$ and $I$.
A good control of this group for the bimodal series
and the quadrangle singularities will be crucial
in this paper. It is the task of the sections
\ref{c5} and \ref{c6}.

The Milnor lattice comes equipped with a 
set $\BB$ of {\it distinguished bases},
certain tuples $\uuuu{\delta}=(\delta_1,\ldots,\delta_\mu)$
of $\Z$-bases of the Milnor lattice.
Each one is defined with a generic deformation of $f$
which has $\mu$ $A_1$-singularities which have all
different critical values. One chooses a {\it distinguished
system of paths} in $\Delta$ from the critical values
to $\delta\in\ppp\Delta$ and pushes vanishing cycles
along these paths to $H_n(f^{-1}(\delta),\Z)=Ml(f)$.
See \cite{AGV88} or \cite{Eb07} for details.
In all cases except the simple singularities, the
set $\BB$ is infinite. Each distinguished basis
determines the monodromy by the formula
\begin{eqnarray}\label{4.4}
M_h=s_{\delta_1}\circ \ldots\circ s_{\delta_\mu}
\end{eqnarray}
where
\begin{eqnarray}\label{4.5}
s_\delta&:&Ml(f)\to Ml(f),\nonumber\\
s_\delta(b)&:=& b-(-1)^{n(n+1)/2}\cdot I(\delta,b)\cdot \delta,
\end{eqnarray}
is the {\it Picard-Lefschetz transformation} of a vanishing
cycle $\delta$, a reflection for even $n$ and a symplectic
transvection for odd $n$.

The matrix of the Seifert form with respect to a distinguished
basis is lower triangular with 
$(-1)^{(n+1)(n+2)/2}$ on the diagonal. This motivates
two definitions, the {\it normalized Seifert form}
\begin{eqnarray}\label{4.6}
L^{hnor}&:=& (-1)^{(n+1)(n+2)/2}\cdot L,
\end{eqnarray}
and the {\it Stokes matrix} $S$ of the distinguished basis with
\begin{eqnarray}\label{4.7}
S&:=& (-1)^{(n+1)(n+2)/2}\cdot L(\uuuu{\delta}^t,\uuuu{\delta})^t
=L^{hnor}(\uuuu{\delta}^t,\uuuu{\delta})^t.
\end{eqnarray}
$S$ is an upper triangular matrix in $GL(\mu,\Z)$
with $1$'s on the diagonal. 

The {\it Coxeter-Dynkin diagram} of a distinguished basis encodes
$S$ in a geometric way. It has $\mu$ vertices which are numbered
from 1 to $\mu$. Between two vertices $i$ and $j$ with $i<j$
one draws

\begin{tabular}{ll}
no edge & if $S_{ij}=0$, \\
$|S_{ij}|$ edges & if $S_{ij}<0$, \\
$S_{ij}$ dotted edges & if $S_{ij}>0$. 
\end{tabular}

Coxeter-Dynkin diagrams for the 8 bimodal series 
will be given in section \ref{c5}, following 
\cite{Eb81}.

A result of Thom and Sebastiani
compares the Milnor lattices and monodromies of 
the singularities $f=f(x_0,\ldots,x_n),g=g(y_0,\ldots,y_m)$ and
$f+g=f(x_0,\ldots,x_n)+g(x_{n+1},\ldots,x_{m+n+1})$.
There are extensions by Deligne for the Seifert form
and by Gabrielov for distinguished bases. 
All results are in \cite[I.2.7]{AGV88}. 
They are restated here.
There is a canonical isomorphism
\begin{eqnarray}\label{4.8}
\Phi:Ml(f+g)&\stackrel{\cong}{\longrightarrow} 
&Ml(f)\otimes Ml(g),\\
\textup{with } M_h(f+g)&\cong & M_h(f)\otimes M_h(g) 
\label{4.9}\\
\textup{and } 
L^{hnor}(f+g)&\cong& L^{hnor}(f)\otimes L^{hnor}(g).
\label{4.10}
\end{eqnarray}
If $\uuuu{\delta}=(\delta_1,\ldots,\delta_{\mu(f)})$
and $\uuuu{\gamma}=(\gamma_1,\ldots,\gamma_{\mu(g)})$ are
distinguished bases of $f$ and $g$  with Stokes matrices
$S(f)$ and $S(g)$, then 
$$\Phi^{-1}(\delta_1\otimes \gamma_1,\ldots,
\delta_1\otimes \gamma_{\mu(g)},
\delta_2\otimes \gamma_1,\ldots,
\delta_2\otimes \gamma_{\mu(g)},
\ldots,
\delta_{\mu(f)}\otimes \gamma_1,\ldots,
\delta_{\mu(f)}\otimes \gamma_{\mu(g)})$$
is a distinguished basis of $Ml(f+g)$,
that means, one takes the vanishing cycles 
$\Phi^{-1}(\delta_i\otimes \gamma_j)$ in the lexicographic order.
Then by \eqref{4.7} and \eqref{4.10}, the matrix 
\begin{eqnarray}\label{4.11}
S(f+g)=S(f)\otimes S(g)
\end{eqnarray}
(where the tensor product is defined
so that it fits to the lexicographic order) 
is the Stokes matrix of this distinguished basis.

In the special case $g=x_{n+1}^2$,
the function germ 
$f+g=f(x_0,\ldots,x_n)+x_{n+1}^2\in \OO_{\C^{n+2},0}$
is called {\it stabilization} or {\it suspension} of $f$. 
As there are only two isomorphisms $Ml(x_{n+1}^2)\to\Z$, 
and they differ by a sign, there are two equally canonical
isomorphisms $Ml(f)\to Ml(f+x_{n+1}^2)$,
and they differ just by a sign. 
Therefore automorphisms and bilinear forms on $Ml(f)$ 
can be identified with automorphisms and bilinear forms on 
$Ml(f+x_{n+1}^2)$. In this sense \cite[I.2.7]{AGV88}
\begin{eqnarray}\label{4.12}
L^{hnor}(f+x_{n+1}^2) &=& L^{hnor}(f),\\ 
M(f+x_{n+1}^2)&=& - M(f),\label{4.13}\\
G_\Z(f+x_{n+1}^2)&=& G_\Z(f).\label{4.14}
\end{eqnarray}
The image in $Ml(f+x_{n+1}^2)$ of a distinguished
basis in $Ml(f)$ under either of the both isomorphisms
$Ml(f)\to Ml(f+x_{n+1}^2)$ is again a distinguished basis,
and it has the same Stokes matrix.

Denote by $H^\infty_\C$ the $\mu$-dimensional vector space
of global flat multi-valued sections in the flat cohomology
bundle $\bigcup_{\tau\in\Delta'}H^n(f^{-1}(\tau),\C)$
(reduced cohomology for $n=0$). 
It comes equipped with a $\Z$-lattice $H^\infty_\Z$, 
a real subspace $H^\infty_\R$, a monodromy which
is also denoted by $M_h$, and the dual 
$L^{nor}$ of the normalized Seifert form $L^{hnor}$.
It is a unimodular form on $H^\infty_\Z$,
and the analogue of \eqref{4.1},
\begin{eqnarray}\label{4.15}
L^{nor}(M_ha,b)=(-1)^{n+1}L^{nor}(b,a) \quad\textup{for }
a,b\in H^\infty_\Z,
\end{eqnarray}
holds.

We apply the notations \ref{t2.1} (a) to
$Ml(f)$ and to $H^\infty_\Z$ and extend them slightly:
\begin{eqnarray}\label{4.16}
Ml(f)_\lambda&:=& \ker(M_h-\lambda\id)^\mu:Ml(f)_\C\to
Ml(f)_\C,\\
Ml(f)_{\neq 1}&:=& \bigoplus_{\lambda\neq 1}Ml(f)_\lambda,\quad
Ml(f)_{\neq -1}:= \bigoplus_{\lambda\neq -1}Ml(f)_\lambda,
\nonumber\\
Ml(f)_p&:=& \bigoplus_{\lambda:\, p(\lambda)=0}
Ml(f)_\lambda,\quad Ml(f)_{p,\Z}:=Ml(f)_p\cap Ml(f).\nonumber
\end{eqnarray}
$H^\infty_\lambda$, 
$H^\infty_{\neq 1}$, $H^\infty_{\neq -1}$,
$H^\infty_p$ and $H^\infty_{p,\Z}$ are defined
analogously.

There are a natural Hodge filtration $F^\bullet_{St}$
on $H^\infty_\C$ and a weight filtration $W_\bullet$ on 
$H^\infty_\Q$ such that 
$(H^\infty_{\neq 1},H^\infty_{\neq 1,\Z},F^\bullet_{St},
W_\bullet,-N,S)$ is a polarized mixed Hodge structure
of weight $n$ and 
$(H^\infty_1,H^\infty_{1,\Z},F^\bullet_{St},
W_\bullet,-N,S)$ is a polarized mixed Hodge structure 
of weight $n+1$ \cite[Theorem 10.30]{He02}.

In the case of a singularity with semisimple monodromy,
so $N=0$, the weight filtrations become trivial,
and the polarized mixed Hodge structures are
polarized pure Hodge structures. This holds for all
bimodal singularities. Therefore we do not care here
about the weight filtration. 
We will define the Hodge filtration using the
Brieskorn lattice in theorem \ref{t7.7}
(following Varchenko, Scherk{\&}Steenbrink and M. Saito).

The pure Hodge structure of weight $n$ on 
$H^\infty_{\neq 1}$ for any singularity with semisimple
monodromy has the following properties.
The Hodge filtration is $M_s$-invariant and satisfies
\begin{eqnarray}\label{4.17}
H^\infty_\lambda&=&
\bigoplus_{p\in\Z}H^{p,n-p}_\lambda\quad\textup{for }
\lambda\neq 1,\\ 
\textup{where \ }
H^{p,n-p}_\lambda&:=&(F^p\cap \oooo{F^{n-p}})H^\infty_\lambda,
\quad(\Rightarrow H^{n-p,p}_{\oooo\lambda} 
=\oooo{H^{p,n-p}_\lambda},)\nonumber\\
\textup{equivalently \ }
H^\infty_\lambda&=& F^pH^\infty_\lambda\oplus 
\oooo{F^{n+1-p}H^\infty_\lambda}.\nonumber
\end{eqnarray}
The polarizing form carries an isotropy and a positivity
condition,
\begin{eqnarray}\label{4.18}
S(H^{p,n-p}_\lambda,H^{q,n-q}_{\oooo\lambda})
&=& 0\quad\textup{if }p+q\neq 0,\\
i^{p-(n-p)}\cdot S(a,\oooo{a})&>0&\quad\textup{for }
a\in H^{p,n-p}_\lambda-\{0\}.\label{4.19}
\end{eqnarray}
The pure Hodge structure of weight $n+1$ on $H^\infty_1$
has analogous properties, with $n$ replaced by $n+1$.

The polarizing form 
$S:H^\infty_\Q\times H^\infty_\Q\to\Q$ 
is defined by \cite[10.6]{He02}. 

\begin{eqnarray}\label{4.20}
S(a,b):=-L^{nor}(a,\nu b)
\end{eqnarray}
where $\nu:H^\infty_\Q\to H^\infty_\Q$ is the $M_h$-invariant automorphism
\begin{eqnarray}\label{4.21}
\nu :=\left\{\begin{array}{ll}
\frac{1}{M_h-\id}&\textup{on }H^\infty_{\neq 1},\\
\frac{-N}{M_h-\id}&\textup{on }H^\infty_1,\end{array}\right.
\end{eqnarray}
$S$ is nondegenerate and $M_h$-invariant.
It  is $(-1)^n$-symmetric on $H^\infty_{\neq 1}$ and 
$(-1)^{n+1}$-symmetric on $H^\infty_1$.
The restriction to $H^\infty_{\neq 1}$ is
$(-1)^{n(n+1)/2}\cdot I^\vee$, where $I^\vee$ on 
$H^\infty_{\neq 1}$ is dual to $I$
(which is nondegenerate on $Ml(f)_{\neq 1}$).

\section{The group \texorpdfstring{$G_\Z$}{GZ} for the bimodal 
series singularities}\label{c5}
\setcounter{equation}{0}

\noindent
The normal forms from \cite[\S 13]{AGV85} for the eight bimodal series
will be listed below in section \ref{c9}. The following table
gives their names, the Milnor numbers, certain polynomials $b_1,b_2$ or, in the
case of the series $Z_{1,p}$, polynomials $b_1,b_2,b_3$ such that
$b_1b_2$ respectively $b_1b_2b_3$ is the characteristic polynomial
of the surface singularities, and two important numbers $m$  and $r_I$.
In the series $p\in\Z_{\geq 1}$.

\begin{eqnarray}\label{5.1}
\begin{array}{lllllll}
\textup{series} & \mu & b_1 & b_2 & b_3 & m & r_I\\   \hline 
W_{1,p}^\sharp & 15+p & \Phi_{12}& (t^{12+p}-1)/\Phi_1& - & 12 & 1\\
S_{1,p}^\sharp & 14+p & \Phi_{10}\Phi_2 & (t^{10+p}-1)/\Phi_1 & - & 10 & 1\\
U_{1,p} & 14+p & \Phi_9 & (t^{9+p}-1)/\Phi_1 & - & 9 & 1 \\
E_{3,p} & 16+p & \Phi_{18}\Phi_2 & t^{9+p}+1& - & 18 & 2\\
Z_{1,p} & 15+p & \Phi_{14}\Phi_2 & t^{7+p}+1 & \Phi_2 & 14 & 2\\
Q_{2,p} & 14+p & \Phi_{12}\Phi_4\Phi_3 & t^{6+p}+1 & - & 12 & 2\\
W_{1,p} & 15+p & \Phi_{12}\Phi_6\Phi_3\Phi_2 & t^{6+p}+1 & - & 12 & 2 \\
S_{1,p} & 14+p & \Phi_{10}\Phi_5\Phi_2 & t^{5+p}+1 & - & 10 & 2
\end{array}
\end{eqnarray}

The following theorem on the group $G_\Z$ will be proved
in two steps. Directly after the theorem, the arguments and properties
which hold for all eight series will be given. 
Then in eight subsections, one for each series, the corresponding
objects will be made explicit and some specific details will be given.
For each series, denote $\zeta:=e^{2\pi i/m}\in S^1\subset\C$.

\begin{theorem}\label{t5.1}
For any surface singularity $f$ in any of the eight bimodal
series, the following holds.

(a) (See definition \ref{t2.3} for the notion {\rm Orlik block})
For all series except $Z_{1,p}$, there are Orlik blocks
$B_1,B_2\subset Ml(f)$, and for the series $Z_{1,p}$, there are
Orlik blocks $B_1,B_2,B_3\subset Ml(f)$ 
with the following properties.
The characteristic polynomial $p_{B_j}$ of the monodromy on $B_j$
is $b_j$. The sum $\sum_{j\geq 1}B_j$ 
is a direct sum $\bigoplus_{j\geq 1}B_j$, 
and it is a sublattice of $Ml(f)$ of full rank $\mu$ and of index $r_I$. 
Define
\begin{eqnarray}\label{5.2}
\www B_1:= \left\{\begin{array}{ll}
B_1 & \textup{for all series except }Z_{1,p},\\
B_1\oplus B_3 & \textup{for the series }Z_{1,p}.
\end{array}\right. 
\end{eqnarray}
Then 
\begin{eqnarray}\label{5.3}
L(\www B_1,B_2)&=&0=L(B_2,\www B_1)\qquad\textup{for all series},\\
G_\Z &=& \Aut(\bigoplus_{j\geq 1}B_j,L)
\qquad\textup{for all series except }S_{1,10}.\label{5.4}
\end{eqnarray}
In the case $S_{1,10}$, a substitute for \eqref{5.4} is
\begin{eqnarray}\label{5.5}
g\in G_\Z\textup{ with }g((B_1)_{\Phi_{10}})=(B_1)_{\Phi_{10}}
\Rightarrow g(B_j)=B_j\textup{ for }j=1,2.
\end{eqnarray}

\medskip
(b) $\Phi_m\not|\, b_2\iff m\not| \, p$. In that case 
\begin{eqnarray}
G_\Z &=& \{(\pm M_h^{k_1}|_{\www B_1})\times 
(\pm M_h^{k_2}|_{B_2})\, |\, k_1,k_2\in\Z\}. \label{5.6}
\end{eqnarray}

\medskip
(c) In the case of the subseries with $m | p$, the
eigenspace $Ml(f)_{\zeta}\subset Ml(f)_\C$ is 2-dimensional.
The hermitian form $h_\zeta$ on it from lemma \ref{t2.2} (a) with 
$h_\zeta(a,b):=\sqrt{-\zeta}\cdot L(a,\oooo{b})$ for $a,b\in
Ml(f)_\zeta$ is nondegenerate and indefinite, so $\P (Ml(f)_{\zeta})\cong \P^1$ contains a half-plane 
\begin{eqnarray}\label{5.7}
\HH_\zeta:=\{\C\cdot a\, |\, a\in Ml(f)_{\zeta}\textup{ with }h_\zeta(a,a)<0\}
\subset\P(Ml(f)_\zeta).
\end{eqnarray}
Therefore the group $\Aut(Ml(f)_{\zeta},h_\zeta)/S^1\cdot\id$ is isomorphic
to $PSL_2(\R)$. The homomorphism
\begin{eqnarray}\label{5.8}
%G_{\Z,2}&:=& \{g\in G_\Z\, |\, g=\id \textup{ on }Ml(f)_{e(1/(9+p))}\},\\
\Psi:G_\Z\to \Aut(Ml(f)_{\zeta},h_\zeta)/S^1\cdot \id,
\quad g\mapsto g|_{Ml(f)_\zeta}\textup{mod}\, S^1\cdot \id,
\end{eqnarray}
is well-defined. $\Psi(G_\Z)$ 
is an infinite Fuchsian group acting on the half-plane 
$\HH_\zeta$. And 
\begin{eqnarray}\label{5.9}
\ker\Psi &=& \{\pm M_h^k\, |\, k\in\Z\} .
\end{eqnarray}
\end{theorem}

{\bf Proof:} Here we explain the common arguments of the proof,
which hold for all eight series. We will announce definitions and properties
of several objects. In the following eight subsections, one for each series,
the objects will be defined, and their properties will be shown.

\medskip
(a) For each of the eight series of surface singularities,
a distinguished basis $e_1,\ldots,e_\mu$ with the 
Coxeter-Dynkin diagram in the corresponding figure will be given
in the subsections \ref{c5.1} to \ref{c5.8}.
The distinguished basis is the one in 
\cite[Tabelle 6 \& Abb. 16]{Eb81}, with a small change in the 
cases $W_{1,1}$ and $S_{1,1}$. They are exceptional in 
\cite{Eb81}. With the actions of the braids
$\alpha_1,\ldots,\alpha_{\mu-1}$ (see \cite[5.7]{Eb07}
for these braids and their actions) and a sign change,
we arrive at a new numbering of the same unnumbered diagram,
such that $W_{1,1}$ and $S_{1,1}$ are no longer exceptional
(i.e. the top vertex has the number $p+q+r+3$ in the notation
of \cite[Abb. 16]{Eb81} even for $W_{1,1}$ and $S_{1,1})$.
We thank Wolfgang Ebeling for the explanation how to arrive
at this numbering.

Recall that for a surface singularity
(then $n=2$) the reflection along a vanishing cycle $\delta$ is
$$s_\delta(b)=b+I(\delta,b)\cdot\delta\qquad\textup{for any }b\in Ml(f).$$ 
The Coxeter-Dynkin diagram has between the vertices $i$ and $j$ with $i<j$
no edge if $S_{ij}=0$, $|S_{ij}|$ edges if $S_{ij}<0$ and
$S_{ij}$ dotted edges if $S_{ij}>0$. Here for $i<j$ 
\begin{eqnarray}
\begin{array}{lll}
I(e_i,e_j)= I(e_j,e_i)&=-S_{ij},&I(e_i,e_i)=-2, \\
L(e_i,e_j)=0,&L(e_j,e_i)= S_{ij}, & L(e_i,e_i)=1.
\end{array} \label{5.10}
\end{eqnarray}
The monodromy can be calculated fairly efficiently by hand
(one should write down some intermediate steps)
with the formula 
\begin{eqnarray}\label{5.11}
M_h=s_{e_1}\circ \ldots\circ s_{e_\mu}.
\end{eqnarray}
The cyclic sublattices $B_j\subset Ml(f)$ are chosen
by choosing the generating lattice vectors $\beta_j$ with
\begin{eqnarray}\label{5.12}
B_j:=\sum_{i\geq 0}\Z\cdot M_h^i(\beta_j).
\end{eqnarray}
The following table gives them.
\begin{eqnarray}\label{5.13}
\begin{array}{llll}
\textup{series} &  \beta_1 & \beta_2  & \beta_3 \\   \hline 
W_{1,p}^\sharp & e_3 & e_8 & -\\
S_{1,p}^\sharp & e_8 & e_9 & -\\
U_{1,p} &  e_8 & e_{10} & -\\
E_{3,p} & e_3 & e_{10} & -\\
Z_{1,p} & e_8 & e_{11} & e_3-e_4-e_9\\
Q_{2,p} & e_8 & e_{11} & -\\
W_{1,p} & e_3+e_9+e_{11} & e_{16} & - \\
S_{1,p} & -e_{8}+e_{13}& e_{15} & -
\end{array}
\end{eqnarray}
We will write down the action of the powers of the monodromy,
\begin{eqnarray}\label{5.14}
\beta_j\mapsto M_h(\beta_j)\mapsto M_h^2(\beta_j)\mapsto \ldots \mapsto M_h^{\deg b_j}(\beta_j),
\end{eqnarray}
in the subsections.
Verifying $b_j(M_h)(\beta_j)=0$ will show that the characteristic
polynomial of $M_h$ on $B_j$ is $b_j$.
We will also write down nice generators of $B_j$. This will show
that $B_j$ is a primitive sublattice of $Ml(f)$,
that $\sum_{j\geq 1}B_j=\bigoplus_{j\geq 1}B_j$ is a direct sum
and that it is a sublattice of full rank and of index $r_I$ in $Ml(f)$.
In all cases except $W_{1,p}$ and $S_{1,p}$, 
the index $r_I$ is obvious from the nice generators, in the two cases 
$W_{1,p}$ and $S_{1,p}$, it requires the calculation
of a determinant. 

The left and right $L$-orthogonality of $\www B_1$ and $B_2$ in \eqref{5.3} will be proved now. 
$e_\mu$ is a cyclic generator for $B_2$ in all
eight series. 
The nice generators for $\www B_1$ show 
$\www B_1\subset\bigoplus_{j=1}^{\mu-2}\Z\cdot e_j$
for all cases except $W_{1,1}$ and $S_{1,1}$.
This and $L(e_i,e_\mu)=0$ for $i<\mu$ show $L(\www B_1,e_\mu)=0$, thus $L(\www B_1,B_2)=0$. From the CDD one sees easily 
$L(e_\mu,e_i)=0$ for $i\leq \mu-2$ for all cases except $W_{1,1}$ and $S_{1,1}$, thus $L(e_\mu,\www B_1)=0$ and $L(B_2,\www B_1)=0$.
For the cases $W_{1,1}$ and $S_{1,1}$, 
$L(B_1,e_\mu)=0=L(e_\mu,B_1)$ and thus 
$L(B_1,B_2)=0=L(B_2,B_1)$ hold also. 

\eqref{5.5} for $S_{1,10}$ will be shown in subsection
\ref{c5.8}.
With respect to part (a), it rests to show \eqref{5.4}.
It is trivial for the 3 series with $r_{I}=1$. It will be shown
in subsection \ref{c5.6} for the series $Q_{2,p}$
and in subsection \ref{c5.7} for the subseries
$W_{1,6s-3}$ ($s\in\Z_{\geq 1})$ of the series $W_{1,p}$.
For all other series, it will be shown below.
It requires a study of smaller Orlik blocks.
$\Phi_2|b_1$ holds in the series 
$S_{1,p}^\sharp$, $E_{3,p}$, $Z_{1,p}$, $W_{1,p}$ and $S_{1,p}$.
In these cases define (see \eqref{2.11} for the notion
$v(\beta_1,-1)$)
\begin{eqnarray}\label{5.15}
\gamma_1 := v(\beta_1,-1):=\frac{b_1}{\Phi_2}(M_h)(\beta_1)
\end{eqnarray}
and calculate $L(\gamma_1,\gamma_1)$ using \eqref{2.12}:
$L(\gamma_1,\gamma_1)=\frac{b_1}{\Phi_2}(-1)\cdot L(\gamma_1,\beta_1).$

\begin{eqnarray}\label{5.16}
\begin{array}{lll}
\textup{series} &  \gamma_1 & L(\gamma_1,\gamma_1)\\   \hline 
S_{1,p}^\sharp & \Phi_{10}(M_h)(e_8)=2e_1+e_2-e_4-e_5-e_6+e_8 & 5\\
E_{3,p} & \Phi_{18}(M_h)(e_3)=-e_2+2e_3+e_6-e_7+e_9 & 6\\
Z_{1,p} &  \Phi_{14}(M_h)(e_8)& \\
&  =e_2+e_3-3e_4 -e_6+e_7-3e_9-e_{10} & 21 \\
W_{1,p} &  (\Phi_{12}\Phi_6\Phi_3)(M_h)(e_3+e_9+e_{11}) & \\
& = e_4-e_5+e_9+e_{11}-e_{13}-e_{15} & 6 \\
S_{1,p} & (\Phi_{10}\Phi_5)(M_h)(-e_{8}+e_{13})\\
& =-2e_1+e_7-e_8-e_9-e_{11}-e_{12}-e_{14}  & 10
\end{array}
\end{eqnarray}
In the case of the series $Z_{1,p}$, define $\gamma_3:=\beta_3$
and calculate 
\begin{eqnarray}\label{5.17}
L(\gamma_3,\gamma_3)=3,\quad 
L(\gamma_1,\gamma_3)=L(\gamma_3,\gamma_1)=7.
\end{eqnarray}

$\Phi_2|b_2$ holds in certain subseries of the series
$S_{1,p}^\sharp$, $E_{3,p}$, $Z_{1,p}$, $W_{1,p}$ and $S_{1,p}$.
In these cases define 
\begin{eqnarray}\label{5.18}
\gamma_2 := v(\beta_2,-1):=\frac{b_2}{\Phi_2}(M_h)(\beta_2)
\end{eqnarray}
and calculate $L(\gamma_2,\gamma_2)$ using \eqref{2.12}:
$L(\gamma_2,\gamma_2)=\frac{b_2}{\Phi_2}(-1)\cdot L(\gamma_2,\beta_2).$

\begin{eqnarray}\label{5.19}
\begin{array}{lll}
\textup{series} & \textup{Condition for }\Phi_2|b_2 & 
L(\gamma_2,\gamma_2) \\ \hline 
S_{1,p}^\sharp & p\equiv 0(2) & 5+\frac{p}{2}\\
E_{3,p} &  p\equiv 0(2) &   18+2p \\
Z_{1,p} & p\equiv 0(2) & 14+2p \\
W_{1,p} &  p\equiv 1(2)  & 12+2p \\
S_{1,p} &  p\equiv 0(2)& 10+2p  
\end{array}
\end{eqnarray}
In table \eqref{5.20}, the first line for $S_{1,p}^\sharp$ is the
case $p\equiv 0(4)$, the second line is the case $p\equiv 2(4)$.
\begin{eqnarray}\label{5.20}
\begin{array}{lll}
\textup{series}  & \gamma_2 \\ \hline 
S_{1,p}^\sharp 
 & -e_2+e_4+e_5+e_6-e_7+\sum_{j=1}^{2+p/4}(e_{7+2j}+e_{10+\frac{p}{2}+2j})\\
 & -e_4+e_5+\sum_{j=1}^{(6+p)/4} (-e_{8+2j}+
 e_{11+\frac{p}{2}+2j}) \\
E_{3,p} & -e_2+2e_5+e_6-e_7+e_9
+2\sum_{j=1}^{4+p/2}e_{8+2j}  \\
Z_{1,p} & -e_2+2e_5 + e_6-e_7+e_{10} + 
2\sum_{j=1}^{3+p/2} e_{9+2j} \\
W_{1,p}& -2e_3+e_4+e_5+e_9+e_{11}
+e_{13}+e_{15}+2\sum_{j=1}^{(1+p)/2} e_{14+2j}  \\
S_{1,p} &  2(-e_1-e_2+e_4+e_5+e_6)-e_7-e_8\\
 & \hspace*{1cm}+e_9+e_{11}+e_{12}+e_{14}-2\sum_{j=1}^{p/2} e_{14+2j} 
\end{array}
\end{eqnarray}

In the subseries of $E_{3,p}, W_{1,p}$ and $S_{1,p}$
with $\Phi_2|b_2$, one sees
\begin{eqnarray}\label{5.21}
\www\gamma_2 :=\frac{1}{2}(\gamma_1+\gamma_2)\stackrel{!}{\in}
Ml(f).
\end{eqnarray}
In the subseries of $Z_{1,p}$ with $\Phi_2|b_2$, one sees
\begin{eqnarray}\label{5.22}
\www\gamma_2 :=\frac{1}{2}(\gamma_1+\gamma_2-3\gamma_3)\stackrel{!}{\in}
Ml(f).
\end{eqnarray}
Together with $[Ml(f):B_1\oplus B_2]=2$ for these
subseries, this shows
\begin{eqnarray}\label{5.23}
Ml(f)_{\Phi_2} &=& \Z\gamma_1\oplus\Z\www\gamma_2
\qquad \textup{for }E_{3,2q},W_{1,2q-1},S_{1,2q},\\
Ml(f)_{\Phi_2} &=& \Z(\gamma_1-2\gamma_3)\oplus 
\Z\www\gamma_2\oplus\Z\gamma_3 \qquad\textup{for }Z_{1,2q}.
\label{5.24}
\end{eqnarray}

For $S_{1,2q}^\sharp$, $Ml(f)=B_1\oplus B_2$ gives
$Ml(f)_{\Phi_2}=\Z\gamma_1\oplus\Z\gamma_2$.
The matrices of $L$ for these bases of $Ml(f)_{\Phi_2}$ in these
cases are
\begin{eqnarray}\label{5.25}
\begin{array}{ccc}
S_{1,2q}^\sharp & E_{3,2q} & Z_{1,2q} \\
\begin{pmatrix} 5 & 0 \\ 0 & 5+q \end{pmatrix} & 
\begin{pmatrix} 6 & 3 \\ 3 & 6+q \end{pmatrix} & 
\begin{pmatrix} 5&2&1 \\2&5+q&-1 \\1&-1& 3 \end{pmatrix} \\ 
  W_{1,2q-1} & S_{1,2q} & \\
\begin{pmatrix} 6 & 3 \\ 3 & 4+q \end{pmatrix} & 
\begin{pmatrix} 10 & 5 \\ 5 & 5+q \end{pmatrix}  & 
\end{array}
\end{eqnarray}
These matrices are positive definite. The corresponding 
quadratic forms 
$(x_1\, x_2)(\textup{matrix})\begin{pmatrix}x_1\\x_2\end{pmatrix}$
respectively 
$(x_1\, x_2\, x_3)(\textup{matrix})\begin{pmatrix}x_1\\x_2\\x_3\end{pmatrix}$
are
\begin{eqnarray}\
5x_1^2+(5+q)x_2^2 & \textup{for }S_{1,2q}^\sharp \nonumber\\
3x_1^2+3(x_1+x_2)^2 + (3+q)x_2^2 & \textup{for }E_{3,2q} \nonumber \\
(2x_1+x_2)^2+(x_1+x_3)^2 \hspace*{1cm}& \nonumber\\ 
+(x_2-x_3)^2+(3+q)x_2^2+x_3^2 & 
\textup{for }Z_{1,2q} \label{5.26}\\
3x_1^2+3(x_1+x_2)^2+(1+q)x_2^2 & \textup{for }W_{1,2q-1}\nonumber\\
5x_1^2+5(x_1+x_2)^2+qx_2^2& \textup{for }S_{1,2q} \nonumber
\end{eqnarray}
This shows 
\begin{eqnarray}\label{5.27}
\{a\in Ml(f)_{\Phi_2}\, |\, L(a,a)=L(\gamma_1,\gamma_1)\}
=\{\pm\gamma_1\}.
\end{eqnarray}
for $W_{1,2q-1}$ with $q\neq 2$,
for $S_{1,2q}$ with $q\neq 5$, and for all
$S_{1,2q}^\sharp$ and $E_{3,2q}$.
It shows for $Z_{1,2q}$
\begin{eqnarray}\label{5.28}
\{a\in Ml(f)_{\Phi_2}\, |\, L(a,a)=3\}
&=&\{\pm\gamma_3\},\\
\{a\in Ml(f)_{\Phi_2}\, |\, L(a,a)=5\}
&=&\{\pm(\gamma_1-2\gamma_3)\}.\label{5.29}
\end{eqnarray}
All this implies
\begin{eqnarray}
\Aut(Ml(f)_{\Phi_2},L) &=&
\{\pm\id|_{\Z\gamma_1}\}\times \{\pm\id|_{\Z\gamma_2}\}
\qquad\textup{for }S_{1,2q}^\sharp,\nonumber\\
&& \textup{for }E_{3,2q}, \textup{ for }S_{1,2q}
\textup{ with }q\neq 5,\nonumber\\
&& \textup{and for }W_{1,2q-1}\textup{ with }q\neq 2,\label{5.30}\\
\Aut(Ml(f)_{\Phi_2},L) &=& 
\{\pm\id|_{\Z\gamma_1\oplus\Z\gamma_3}\}\times
 \{\pm\id|_{\Z\gamma_2}\}
\ \textup{for }Z_{1,2q}.\hspace*{0.5cm}\label{5.31}
\end{eqnarray}
In the cases $S_{1,2q-1}^\sharp$, $E_{3,2q-1}$, $Z_{1,2q-1}$, 
$W_{1,2q}$ and $S_{1,2q-1}$ with $\Phi_2\not| \, b_2$, 
\begin{eqnarray}\label{5.32}
Ml(f)_{\Phi_2}= (\www B_1)_{\Phi_2}
\quad \textup{and}\quad
\Aut(Ml(f)_{\Phi_2},L) =\{\pm\id\}.
\end{eqnarray}
Define 
\begin{eqnarray}\label{5.33}
\gamma_4 := \left\{\begin{array}{ll} \gamma_1 & 
\textup{for }E_{3,p},W_{1,p},S_{1,p}\\
\gamma_1-3\gamma_3 & \textup{for }Z_{1,p}.\end{array}\right. 
\end{eqnarray}
Then for $E_{3,p}$, $W_{1,p}$ with $p\neq 3$, $S_{1,p}$ with $p \neq 10$, $Z_{1,p}$
\begin{eqnarray}\label{5.34}
g(\gamma_4)=\pm\gamma_4\qquad \textup{for }g\in G_\Z,
\end{eqnarray}
and for $E_{3,p}$, $W_{1,p}$ (including $p=3$), $S_{1,p}$ (including $p=10$), $Z_{1,p}$
\begin{eqnarray}\label{5.35}
\www B_1\oplus B_2 =\{a\in Ml(f)\, |\, L(a,\gamma_4)\equiv 0(2)\}.
\end{eqnarray}
Here $\subset$ \eqref{5.35} follows from $L(B_2,\gamma_4)=0$ and
$L(\beta_1,\gamma_4)\equiv 0(2)$ and in the case
of $Z_{1,p}$ $L(\beta_3,\gamma_4)=4$.
Now $=$ in \eqref{5.35} follows from $L(Ml(f),\gamma_4)=\Z$ and
$[Ml(f):\www B_1\oplus B_2]=2$.
Together \eqref{5.34} and \eqref{5.35}
show that any $g\in G_\Z$ respects $\www B_1\oplus B_2$, so
\begin{eqnarray}\label{5.36}
G_\Z\subset\Aut(\www B_1\oplus B_2,L)
\end{eqnarray}
for $E_{3,p}$, $W_{1,p}$ with $p\neq 3$, $S_{1,p}$ with $p \neq 10$ and $Z_{1,p}$.
We claim that \eqref{5.34} and thus \eqref{5.36} hold also
for $W_{1,3}$. That will be proved in the subsection \ref{c5.7}.

It rests to show $\Aut(\www B_1\oplus B_2,L)\subset G_\Z$
for the series $E_{3,p}$, $Z_{1,p}$, $W_{1,p}$, $S_{1,p}$. 
We will extend the definition of $\www\gamma_2$ in such a way
to the cases with $\Phi_2\not| b_2$ 
that $(\www B_1\oplus B_2)+\Z\cdot\www\gamma_2=Ml(f)$.
And we will show $g(\www\gamma_2)\in Ml(f)$ for any
$g\in\Aut(\www B_1\oplus B_2,L)$.
This implies $\Aut(\www B_1\oplus B_2,L)\subset G_\Z$.
The proof of $g(\www\gamma_2)\in Ml(f)$ requires a better
control of $\Aut(\www B_1\oplus B_2,L)$.

Consider all eight series and define 
\begin{eqnarray}\label{5.37}
b_4 := \frac{\gcd(b_1,b_2)}{\gcd(b_1,b_2,\Phi_m)}
=\gcd(\frac{b_1}{\Phi_m},b_2)\in\Z[t].
\end{eqnarray}
Then
\begin{eqnarray}\label{5.38}
b_4=\left\{\begin{array}{ll} 1 & \textup{for }
W_{1,p}^\sharp,S_{1,2q-1}^\sharp,U_{1,p},E_{3,2q-1},Z_{1,2q-1},\\
& Q_{2,p}\textup{ with }p\not\equiv 0(4),W_{1,2q},S_{1,2q-1},\\
\Phi_2 & \textup{for }S_{1,2q}^\sharp, E_{3,2q},Z_{1,2q},
W_{1,2q-1}\textup{ with }q\not\equiv 2(3),S_{1,2q},\\
\Phi_4 & \textup{for }Q_{2,4s},\\
\Phi_6\Phi_2 &  \textup{for }W_{1,6s-3}.\end{array}\right. 
\end{eqnarray}
We claim that in all cases except $S_{1,10}$, 
any $g\in G_\Z\cup\Aut(\www B_1\oplus B_2,L)$ maps
$(\www B_1)_{b_4}$ to $(\www B_1)_{b_4}$ and $(B_2)_{b_4}$
to $(B_2)_{b_4}$.
In the cases with $b_4=1$ this is an empty statement as then
$(\www B_1)_{b_4}=\{0\}=(B_2)_{b_4}$.
In the cases $Q_{2,p}$ with $p\equiv 0(4)$ and $W_{1,6s-3}$, this
will be shown in the subsections \ref{c5.6} and \ref{c5.7}.
In all other cases $b_4=\Phi_2$ and $(B_2)_{b_4}=\Z\cdot \gamma_2$ and
\begin{eqnarray}\label{5.39}
(\www B_1)_ {b_4} = \left\{\begin{array}{ll}
\Z\cdot(\gamma_1-2\gamma_3)\oplus\Z\cdot \gamma_3 & \textup{for }Z_{1,2q},\\
\Z\cdot\gamma_1&\textup{else}. 
\end{array}\right. 
\end{eqnarray}
Because $(\www B_1\oplus B_2)_{\Phi_2}\subset Ml(f)_{\Phi_2}$,
\eqref{5.27}--\eqref{5.29} hold also with $(\www B_1\oplus B_2)_{\Phi_2}$
instead of $Ml(f)_{\Phi_2}$. They characterize 
$(\www B_1)_{\Phi_2}$ within $Ml(f)_{\Phi_2}$ and within 
$(\www B_1\oplus B_2)_{\Phi_2}$.
Thus any $g\in G_\Z\cup\Aut(\www B_1\oplus B_2,L)$ maps $(\www B_1)_{\Phi_2}$
to itself, and then it maps also the $L$-orthogonal sublattice
$(B_2)_{\Phi_2}$ to itself.

For all eight series except $S_{1,10}$, 
this implies the following.
For any $g\in G_\Z\cup \Aut(\www B_1\oplus B_2,L)$
\begin{eqnarray}\label{5.40}
g:\www B_1\to \www B_1 \textup{ and }
B_2\to B_2 & \textup{ if }m\not|\, p,\\
\left. \begin{array}{r}
g:(\www B_1)_{b_1/\Phi_m}\to (\www B_1)_{b_1/\Phi_m}\\ 
g:(B_2)_{b_2/\Phi_m}\to (B_2)_{b_2/\Phi_m}\end{array}\right\}
& \begin{array}{r}\textup{ if }m|p\textup{ and the}\\
\textup{type is not }S_{1,10}.\end{array}  \label{5.41}
\end{eqnarray}

Now we want to apply lemma \ref{t2.8} to these Orlik blocks.
One checks easily that all hypotheses are satisfied.
Therefore
\begin{eqnarray}\label{5.42}
&&\Aut(\www B_1\oplus B_2,L) \\
&=& \{\pm M_h^k|_{\www B_1}\, |\, k\in\Z\}\times
\{\pm M_h^k|_{B_2}\, |\, k\in\Z\}\qquad\textup{if }m\not|\, p,
\nonumber
\end{eqnarray}
and if $m|\, p$  and the type is not $S_{1,10}$,  
then $\Aut(\www B_1\oplus B_2,L)$ projects to a subgroup of
\begin{eqnarray}\label{5.43}
&&\Aut((\www B_1)_{b_1/\Phi_m},L)\times\Aut((B_2)_{b_2/\Phi_m},L) \\
&=& \{\pm M_h^k|_{(\www B_1)_{b_1/\Phi_m}}\, |\, k\in\Z\}\times
\{\pm M_h^k|_{(B_2)_{b_2/\Phi_m}}\, |\, k\in\Z\}.\nonumber
\end{eqnarray}
The group $\Aut(\www B_1\oplus B_2,L)$ for $m\not| p$ is 
generated by $M_h,-\id$, $M_h|_{\www B_1}\times \id|_{B_2}$
and $(-\id|_{\www B_1})\times \id|_{B_2}$, 
and analogously for the group in \eqref{5.43} if $m|\, p$.

Now we extend the definition of $\gamma_2$. 
For $E_{3,2q-1}$, $Z_{1,2q-1}$ and $S_{1,2q-1}$ define it as
follows:
\begin{eqnarray}\label{5.44}
\gamma_2&:=& e_2-e_6+e_7+e_9\quad\textup{for }E_{3,2q-1},\\
\gamma_2&:=& e_2-e_6+e_7+e_{10}\quad\textup{for }Z_{1,2q-1},
\nonumber \\
\gamma_2&:=& 2(-e_1-e_2+\sum_{j\in\{4,5,6\}}e_j)
-e_7-e_8+\sum_{j\in\{9,11,12,14\}}e_j  
\quad\textup{for }S_{1,2q-1}.\nonumber
\end{eqnarray}
\eqref{5.4.4}, \eqref{5.5.5} and \eqref{5.8.2} show 
$\gamma_2\in B_2$.
For $W_{1,2q}$ (so $p=2q$) define
\begin{eqnarray}\label{5.45}
\gamma_2 &:=& (t^p(t+1)\Phi_{12}+
\sum_{j=0}^{p-1}t^j)(M_h)(e_{16})\\
&=& (t^p(1+t-t^2-t^3+t^4+t^5)+\sum_{j=0}^{p-1}t^j)(M_h)(e_{16})
\nonumber \\
&=& -2e_2+2e_6-2e_7+e_4+e_5+e_9-e_{11}+e_{13}-e_{15}.\nonumber
\end{eqnarray}
Observe that in the case $12|p$, $\Phi_{12}$ divides 
$\sum_{j=0}^{p-1}t^j$ 
so that then 
$\gamma_2\in \Phi_{12}(M_h)(B_2)=(B_2)_{b_2/\Phi_{12}}$.
In all four cases $\frac{1}{2}(\gamma_4+\gamma_2)\in Ml(f)$.

Now for the series $E_{3,p},Z_{1,p},W_{1,p}$ and $S_{1,p}$
\begin{eqnarray}\label{5.46}
\gamma_4\in (B_1)_{\Phi_2},\qquad 
\left\{\begin{array}{ll}
\gamma_2\in B_2&\textup{ if }m\not|\, p,\\
\gamma_2\in (B_2)_{b_2/\Phi_m}&\textup{ if }m|\, p,\end{array}\right. \\
\www\gamma_2 :=\frac{1}{2}(\gamma_4+\gamma_2)
\stackrel{!}{\in} Ml(f),\label{5.47}\\
Ml(f)=(\www B_1\oplus B_2)+\Z\www\gamma_2,
\label{5.48}\\
(M_h|_{\www B_1}\times\id|_{B_2})(\www\gamma_2) 
=((-\id|_{\www B_1})\times \id|_{B_2})(\www\gamma_2) \nonumber\\
=\frac{1}{2}(-\gamma_4+\gamma_2)=-\gamma_4+\www\gamma_2\in Ml(f).\label{5.49}
\end{eqnarray}
Therefore any $g\in \Aut(\www B_1\oplus B_2,L)$ maps
$\www\gamma_2$ to an element of $Ml(f)$.
Thus it maps $Ml(f)$ to $Ml(f)$, thus $g\in G_\Z$.
This finishes the proof of \eqref{5.4} and of part (a)
for all series except $Q_{2,p}$ and $W_{1,6s-3}$ and $S_{1,10}$.
For $Q_{2,p}$ and $W_{1,6s-3}$ and $S_{1,10}$ 
see the subsections \ref{c5.6}, \ref{c5.7} and \ref{c5.8}.

\medskip
(b) This follows immediately from \eqref{5.4} and \eqref{5.42}.
The subsections \ref{c5.6} and \ref{c5.7}
establish \eqref{5.4} and \eqref{5.42} also for the series
$Q_{2,p}$ and $W_{1,6s-3}$.

\medskip
(c) Now we consider the eight subseries with $m|p$. 
Write $p=m\cdot r$ with $r\in\Z_{\geq 1}$. 
Recall $\zeta=e^{2\pi i /m}$, and recall that $\Z[\zeta]$
is a principal ideal domain (lemma \ref{t2.11}).
In the following, $\xi$ will be any primitive $m$-th unit root.

Formula \eqref{2.19} in lemma \ref{t2.12} (b) applies with 
$\Lambda=Ml(f), \Lambda^{(1)}=\www B_1\oplus B_2, p=\Phi_m,$ 
and gives
\begin{eqnarray}\label{5.50n}
Ml(f)_{\Phi_m}=(\www B_1\oplus B_2)_{\Phi_m}
=(B_1\oplus B_2)_{\Phi_m}=(B_1)_{\Phi_m}\oplus
(B_2)_{\Phi_m}.
\end{eqnarray}
Therefore the space
\begin{eqnarray}\label{5.50}
Ml(f)_{\xi,\Z[\zeta]}&:=& Ml(f)_\xi \cap Ml(f)_{\Z[\zeta]}
\end{eqnarray}
is a free $\Z[\zeta]$-module of rank 2 with basis
$v_{1,\xi},v_{2,\xi}$ 
with
\begin{eqnarray}\label{5.51}
v_{j,\xi}:= v(\beta_j,\xi)=\frac{b_j}{t-\xi}(M_h)(\beta_j)\qquad \textup{for }
j=1,2
\end{eqnarray}
(see \eqref{2.11} for the notion $v(\beta_j,\xi)$).
Observe $v_{j,\oooo{\xi}}=\oooo{v_{j,\xi}}$. 

The proof of part (c) will consist of four steps.
Step 1 calculates the values of the hermitian form $h_\xi$ from
lemma \ref{t2.2} on a suitable $\Z[\zeta]$-basis of 
$Ml(f)_{\xi,\Z[\zeta]}$.
Step 2 analyzes what this implies for automorphisms of the pair
$(Ml(f)_{\xi,\Z[\zeta]},L)$ and thus gives a first approximation
to $\Psi(G_\Z)$. Step 3 uses \eqref{5.5} for $S_{1,10}$
and \eqref{5.41} for all other singularities and
the Orlik block structure of the blocks $B_j$ to control
the action of $g\in G_\Z$ on all eigenspaces simultaneously.
It will prove \eqref{5.9}.
Step 4 combines the steps 2 and 3 with results from section \ref{c3}
and shows that $\Psi(G_\Z)$ is an infinite Fuchsian group.

\medskip
{\bf Step 1:}
The form 
\begin{eqnarray*}
h_\xi: Ml(f)_\xi\times Ml(f)_\xi\to\C,\quad (a,b)\mapsto \sqrt{-\xi}\cdot 
L(a,\oooo{b})
\end{eqnarray*}
from lemma \ref{t2.2} 
is hermitian. In this step it will be calculated with respect to the
$\Z[\zeta]$-basis $v_{1,\xi},v_{2,\xi}$ of 
$Ml(f)_{\xi,\Z[\zeta]}$. For $i\neq j$ 
\begin{eqnarray}\label{5.52}
h_\xi(v_{i,\xi},v_{j,\xi})=\sqrt{-\xi}\cdot L(v_{i,\xi},v_{j,\oooo\xi})=0
\end{eqnarray}
because of \eqref{5.3}.
$L(v_{j,\xi},v_{j,\oooo{\xi}})$ will be calculated with \eqref{2.12},
\begin{eqnarray}\label{5.53}
L(v_{j,\xi},v_{j,\oooo{\xi}}) = \frac{b_j}{t-\oooo{\xi}}({\oooo{\xi}})\cdot
L(\frac{b_j}{t-\xi}(M_h)(\beta_j),\beta_j),
\end{eqnarray}
first for $j=2$, then for $j=1$.

One calculates for all eight subseries:
\begin{eqnarray*}
\begin{array}{c|c|c|c|c|c|c|}
k & 0 & 1 & 2 & \cdots & \deg b_2-1 & \deg b_2 \\ \hline 
L(M_h^k(\beta_2),\beta_2) & 1 & -1 & 0 & \cdots & 0 & 
0\textup{ if }r_I=1,\, -1\textup{ if }r_I\geq 2
\end{array}
\end{eqnarray*}

For the three subseries with $r_I=1$ (so $W^\sharp_{1,12r},
S^\sharp_{1,10r},U_{1,9r}$)
\begin{eqnarray}\label{5.54}
\frac{b_2}{t-\xi} &=& \frac{t^{m+p}-1}{(t-\xi)\cdot \Phi_1}
= \Phi_1^{-1}\cdot \sum_{j=0}^{m+p-1}\xi^{m+p-1-j}\cdot t^j,\\
\frac{b_2}{t-\oooo{\xi}}(\oooo{\xi}) 
&=& (\oooo{\xi}-1)^{-1}\cdot (m+p)\cdot\xi 
= m(1+r)(\oooo{\xi}-1)^{-1}\cdot \xi,\label{5.55}
\end{eqnarray}
\begin{eqnarray}
L(\frac{b_2}{t-\xi}(M_h)(\beta_2),\beta_2) &=& 
(\xi-1)^{-1}\cdot \oooo{\xi}\cdot (1-\oooo{\xi})={\oooo{\xi}}^2,
\label{5.56}\\
h_\xi(v_{2,\xi},v_{2,\xi})&=& m(1+r)\cdot (1-\xi)^{-1}\cdot 
\sqrt{-\xi}>0.\label{5.57}
\end{eqnarray}
For the five subseries with $r_I= 2$
\begin{eqnarray}\label{5.58}
\frac{b_2}{t-\xi} &=& \frac{t^{m/2+p}+1}{t-\xi}
= \sum_{j=0}^{m/2+p-1}\xi^{m/2+p-1-j}\cdot t^j,\\
\frac{b_2}{t-\oooo{\xi}}(\oooo{\xi}) 
&=& (\frac{m}{2}+p)(-\xi) 
= \frac{m}{2}(1+2r)(-\xi),\label{5.59}
\end{eqnarray}
\begin{eqnarray}
L(\frac{b_2}{t-\xi}(M_h)(\beta_2),\beta_2) &=& 
-\oooo{\xi}(1-\oooo{\xi}),\label{5.60}\\
h_\xi(v_{2,\xi},v_{2,\xi})&=& \frac{m}{2}(1+2r)
\cdot (1-\oooo{\xi})\cdot \sqrt{-\xi}>0.\label{5.61}
\end{eqnarray}

Now we turn to $h_\xi(v_{1,\xi},v_{1,\xi})$. One calculates for all
eight series 
\begin{eqnarray*}
\begin{array}{c|c|c|c|c|c|c|c|c|c|c|c|c|}
k & 0 & 1 & 2 & 3 & 4 & 5 & 6 & 7 & 8 & 9 & 10 & 11  \\
L(M_h^k(\beta_1),\beta_1) & & & & & & & & & & & &  \\ \hline 
\textup{for }W_{1,p}^\sharp & 1 & -1 & 1 & 0 & 0 & 1 & & & & & &  \\
\textup{for }S_{1,p}^\sharp & 1 & -1 & 0 & 1 & 0 & & & & & & &  \\
\textup{for }U_{1,p} & 1 & -1 & 0 & 0 & 1 & 0 & -1 & 0 & 0 & & & \\
\textup{for }E_{3,p} & 1 & -1 & 1 & 0 & 1 & 0 & 1 & 0 & 1 & & &  \\
\textup{for }Z_{1,p} & 1 & -1 & 0 & 0 & 1 & 0 & 0 & & & & & \\
\textup{for }Q_{2,p} & 1 & -1 & 0 & 1 & 0 & 0 & 0 & 0 & 0 & 0 & -1 & 0  \\
\textup{for }W_{1,p} & 3 & -3 & 2 & -1 & 0 & 1 & -1 & 1 & -1 & 0 & 1 & -2  \\
\textup{for }S_{1,p} & 2 & -2 & 0 & 1 & 0 & -1 & 1 & 0 & -1 & 0 & & 
\end{array}
\end{eqnarray*}
and 
\begin{eqnarray*}
\textup{for }W_{1,p}^\sharp && 
\frac{b_1}{t-\xi} = \frac{\Phi_{12}}{t-\xi} = t^3+\xi t^2+(\xi^2-1)t
+(\xi^3-\xi),\\
\textup{for }S_{1,p}^\sharp && 
\frac{b_1}{t-\xi} = \frac{\Phi_{10}\Phi_2}{t-\xi} = 
\frac{t^5+1}{t-\xi} = t^4+\xi t^3+\xi^2t^2+\xi^3t+\xi^4,\\
\textup{for }U_{1,p} && 
\frac{b_1}{t-\xi} = \frac{\Phi_9}{t-\xi} = 
\frac{t^6+t^3+1}{t-\xi} \\
&&= t^5+\xi t^4+\xi^2t^3+(\xi^3+1)t^2 +(\xi^4+\xi)t +(\xi^5+\xi^2),
\end{eqnarray*}
\begin{eqnarray*}
\textup{for }E_{3,p} && 
\frac{b_1}{t-\xi} = \frac{\Phi_{18}\Phi_2}{t-\xi} 
=\frac{t^7+t^6-t^4-t^3+t+1}{t-\xi} 
= t^6+(\xi+1)t^5\\
&& +(\xi^2+\xi)t^4 +(\xi^6+\xi^2)t^3+(\xi^7+\xi^6)t^2 
+(\xi^8+\xi^7)t+\xi^8,\\
\textup{for }Z_{1,p} && 
\frac{b_1}{t-\xi} = \frac{t^7+1}{t-\xi} 
 = t^6+\xi t^5+\xi^2 t^4+\xi^3t^3+\xi^4t^2+\xi^5t+\xi^6,\\
\textup{for }Q_{2,p} && 
\frac{b_1}{t-\xi} = \frac{\Phi_{12}\Phi_4\Phi_3}{t-\xi} = 
\frac{t^8+t^7+t^6+t^2+t+1}{t-\xi} \\
&& = t^7+(\xi+1)t^6 +(\xi^2+\xi+1)t^5 +(\xi^3+\xi^2+\xi)t^4\\
&& +(\xi^4+\xi^3+\xi^2)t^3+(\xi^5+\xi^4+\xi^3)t^2 
+(\xi^5+\xi^4)t+\xi^5,
\end{eqnarray*}
\begin{eqnarray*}
\textup{for }W_{1,p} && 
\frac{b_1}{t-\xi} = \frac{\Phi_{12}\Phi_6\Phi_3\Phi_2}{t-\xi} = 
\frac{t^9+t^8+t^5+t^4+t+1}{t-\xi} \\
&& = t^8+(\xi+1)t^7+(\xi^2+\xi)t^6 + (\xi^3+\xi^2)t^5 + (\xi^3+\xi^2)t^4\\
&& + (\xi^3+\xi^2)t^3+(\xi^4+\xi^3)t^2 +(\xi^5+\xi^4)t+\xi^5,\\
\textup{for }S_{1,p} && 
\frac{b_1}{t-\xi} = \frac{\Phi_{10}\Phi_5\Phi_2}{t-\xi} = 
%\frac{t^9+t^8+t^7+t^6+t^5+t^4+t^3+t^2+t+1}{t-\xi} \\
\frac{\sum_{j=0}^9 t^j}{t-\xi}\\
&& = t^8+(\xi+1)t^7+(\xi^2+\xi+1)t^6+(\xi^3+\xi^2+\xi+1)t^5\\
&&+(\xi^4+\xi^3+\xi^2+\xi+1)t^4+(\xi^4+\xi^3+\xi^2+\xi)t^3\\
&&+(\xi^4+\xi^3+\xi^2)t^2+(\xi^4+\xi^3)t+\xi^4 .
\end{eqnarray*}
This table and this list give the following values.
\begin{eqnarray*}
\begin{array}{c|c|c}
 & \frac{b_1}{t-\oooo{\xi}}(\oooo\xi) & 
 L(\frac{b_1}{t-\xi}(M_h)(\beta_1),\beta_1) \\
 \hline 
W_{1,p}^\sharp & 4\oooo\xi^3-2\oooo\xi=-2(\xi+\oooo\xi)\xi^2 & 
\xi^3(1-\xi)\\
S_{1,p}^\sharp & 5\oooo\xi^4=-5\xi &
-\xi(\xi^2+\oooo\xi^2-1)\\
U_{1,p} & 6\oooo\xi^5+3\oooo\xi^2=3\xi(\xi^3-1) & 
-\xi^6(\xi^2+\oooo\xi^2) \\
E_{3,p} & 3(\oooo\xi^6+\oooo\xi^5+\oooo\xi^9+\oooo\xi^8)
=-3(\xi+1)(\xi^3+1) & \xi^2(\xi+\oooo\xi)(\xi^2+\oooo\xi^2)\\
Z_{1,p} & 7\oooo\xi^6 =-7\xi & 
\xi^2(\xi^4+\oooo\xi^4+1)\\
Q_{2,p} & 6(\oooo\xi^7+\oooo\xi^6+\oooo\xi^5)=-6(\xi+\oooo\xi+1) & 
\xi^2(\xi+1)=(1-\xi)^{-1}\\
W_{1,p} & 4(\oooo\xi^8+\oooo\xi^7+\oooo\xi^6+\oooo\xi^5)
=4\oooo\xi^7(1+\xi)(\xi+\oooo\xi) & \xi^3(\xi-1)(\oooo\xi-1) \\
S_{1,p} & 5(\oooo\xi^8+\oooo\xi^7+\oooo\xi^6+\oooo\xi^5+\oooo \xi^4) & -1+\xi+\xi^2-2\xi^3+\xi^4  
\end{array}
\end{eqnarray*}
With $h_\xi(v_{1,\xi},v_{1,\xi})=\sqrt{-\xi}\cdot
L(v_{1,\xi},v_{1,\oooo\xi})$ and \eqref{5.53} 
and the information on the rings $\Z[\zeta]$ in lemma \ref{t2.11}, 
we obtain the following values.
\begin{eqnarray}\label{5.62}
\begin{array}{c|r}
 & h_\xi(v_{1,\xi},v_{1,\xi}) \\ \hline 
W_{1,p}^\sharp & (-2)(\xi+\oooo\xi)\cdot(1-\oooo\xi)\sqrt{-\xi} \\
S_{1,p}^\sharp & 5(\xi^2+\oooo\xi^2)(\xi^2+\oooo\xi^2-1)\cdot(1-\xi)^{-1}
\sqrt{-\xi}\\
U_{1,p} & 3(\xi^4+\oooo\xi^4+1)\cdot(1-\oooo\xi)\sqrt{-\xi}\\
E_{3,p} & (-3)(1+\xi)(1+\oooo\xi)(\xi+\oooo\xi-1)\cdot
(1-\xi)^{-1}\sqrt{-\xi} \\
Z_{1,p} & (-7)(\xi^2+\oooo\xi^2)\cdot(1-\oooo\xi)\sqrt{-\xi}\\
Q_{2,p} & (-6)(\xi+\oooo\xi+1)\cdot(1-\xi)^{-1}\sqrt{-\xi} \\
W_{1,p} & (-4)(\xi+\oooo\xi)\cdot(1-\oooo\xi)\sqrt{-\xi} \\
S_{1,p} & (-10)(\xi^2+\oooo\xi^2)\cdot(1-\oooo\xi) \sqrt{-\xi} 
\end{array}
\end{eqnarray}
Here observe that as in \eqref{5.57} and \eqref{5.61}
$(1-\oooo\xi)\sqrt{-\xi}>0$ and $(1-\xi)^{-1}\sqrt{-\xi}>0$.
In each of the eight cases we find
\begin{eqnarray}\label{5.63}
h_\xi(v_{1,\xi},v_{1,\xi})&>&0\qquad\textup{for }
\xi\not\in\{\zeta,\oooo\zeta\},\\
h_\xi(v_{1,\xi},v_{1,\xi})&<&0\qquad\textup{for }
\xi\in\{\zeta,\oooo\zeta\},\label{5.64}
\end{eqnarray}
and 
\begin{eqnarray}\label{5.65}
L(v_{1,\xi},\beta_1)=L(\frac{b_1}{t-\xi}(M_h)(\beta_1),\beta_1)\in\Z[\zeta]^*.
\end{eqnarray}

\medskip
{\bf Step 2:}
Define for each of the eight series 
\begin{eqnarray}\label{5.66}
b_5 :=\frac{b_1}{\Phi_m}\in\Z[t]\qquad\textup{unitary}.
\end{eqnarray}
Then
\begin{eqnarray*}
\begin{array}{c|c|c|c|c|c|c|c|c}
\textup{series}&W_{1,p}^\sharp & S_{1,p}^\sharp & U_{1,p} & E_{3,p} &
Z_{1,p} & Q_{2,p} & W_{1,p} & S_{1,p} \\ \hline
b_5 & 1 & \Phi_2 & 1 & \Phi_2 & \Phi_2 & \Phi_4\Phi_3 & 
\Phi_6\Phi_3\Phi_2 & \Phi_5\Phi_2 \end{array}
\end{eqnarray*}
and 
\begin{eqnarray}\label{5.67}
b_5(\xi)/b_5(\oooo\xi)\in\{\pm\xi^k\, |\, k\in\Z\}.
\end{eqnarray}

Define for each of the eight subseries with $m|p$ 
\begin{eqnarray}\label{5.68}
b_6 :=\frac{b_2}{\Phi_m}\in\Z[t]\qquad\textup{unitary}
\end{eqnarray}
and
\begin{eqnarray}
w(\xi)&:=& -\frac{h_\xi(v_{2,\xi},v_{2,\xi})}{h_\xi(v_{1,\xi},v_{1,\xi})}
=-\frac{\frac{b_2}{t-\oooo\xi}(\oooo\xi)\cdot L(v_{2,\xi},\beta_2)}
{\frac{b_1}{t-\oooo\xi}(\oooo\xi)\cdot L(v_{1,\xi},\beta_1)}\nonumber\\
&=&-\frac{b_6}{b_5}(\oooo\xi)\cdot 
\frac{L(v_{2,\xi},\beta_2)}{L(v_{1,\xi},\beta_1)}.\label{5.69}
\end{eqnarray}
Then 
\begin{eqnarray}\label{5.70}
b_5(\oooo\xi)w(\xi)= b_6(\oooo\xi)\cdot 
\frac{L(v_{2,\xi},\beta_2)}{L(v_{1,\xi},\beta_1)}\in\Z[\zeta].
\end{eqnarray}
It is in $\Z[\zeta]$ because of \eqref{5.65}.
The following table lists $w(\xi)$.
\begin{eqnarray}\label{5.71}
\begin{array}{c|l}
 & w(\xi) \\ \hline 
W_{1,p}^\sharp & (1+r)(+6)[(1-\xi)(1-\oooo\xi)(\xi+\oooo\xi)]^{-1} \\
S_{1,p}^\sharp & (1+r)(-2)[(\xi^2+\oooo\xi^2)(\xi^2+\oooo\xi^2-1)]^{-1}\\
U_{1,p} & (1+r)(-3)[(1-\xi)(1-\oooo\xi)(\xi^4+\oooo\xi^4+1)]^{-1}\\
E_{3,p} & (1+2r)(+3)(1-\xi)(1-\oooo\xi)
[(1+\xi)(1+\oooo\xi)(\xi+\oooo\xi-1)]^{-1} \\
Z_{1,p} & (1+2r)(+1)[\xi^2+\oooo\xi^2]^{-1}\\
Q_{2,p} & (1+2r)(+1)(1-\xi)(1-\oooo\xi)[\xi+\oooo\xi+1]^{-1}\\
W_{1,p} & (1+2r)(+\frac{3}{2})[\xi+\oooo\xi]^{-1} \\
S_{1,p} & (1+2r)(+\frac{1}{2})[\xi^2+\oooo\xi^2]^{-1}
\end{array}
\end{eqnarray}
The inequalities \eqref{5.57}\eqref{5.61}\eqref{5.63}\eqref{5.64} give
\begin{eqnarray}\label{5.72}
w(\xi)\left\{\begin{array}{ll}
<0& \textup{ for }\xi\not\in \{\zeta,\oooo\zeta\},\\
>0& \textup{ for }\xi\in\{\zeta,\oooo\zeta\}.\end{array} \right.
\end{eqnarray}

Using the $\Z[\zeta]$-basis $v_{1,\xi},v_{2,\xi}$ of 
$Ml(f)_{\xi,\Z[\zeta]}$,
the automorphism group $\Aut(Ml((f)_{\xi,\Z[\zeta]},h_\xi)$ can
be identified with the matrix group
\begin{eqnarray}
\{A(\xi)\in GL(2,\Z[\zeta])&|& \nonumber \\
\begin{pmatrix}-1&0\\0&w(\xi)\end{pmatrix}
&=& A(\xi)^t\cdot \begin{pmatrix}-1&0\\0&w(\xi)\end{pmatrix}\label{5.73}
\cdot \oooo{A(\xi)}\}.
\end{eqnarray}
The isomorphism is $A(\xi)\mapsto g$ with
\begin{eqnarray}\label{5.74}
g(v_{1,\xi},v_{2,\xi})&=&(v_{1,\xi},v_{2,\xi})\cdot A(\xi).
\end{eqnarray}
The inequalities \eqref{5.72} and 
theorem \ref{t3.2} tell that the matrix group 
in the case of $\xi=\zeta$ projects to an infinite Fuchsian group. 
Additionally, \ref{t3.2} tells  
that the elements of the matrix group for any $\xi$ 
can be represented by triples
$(a(\xi),c(\xi),\delta(\xi))\in\Z[\zeta]^2\times \{\pm\zeta^k\, |\, k\in\Z\}$
with 
\begin{eqnarray}\label{5.75}
a(\xi)a(\oooo{\xi})-1= w(\xi)\cdot c(\xi)c(\oooo{\xi}),
\end{eqnarray}
where
\begin{eqnarray}\label{5.76}
A(\xi) = \begin{pmatrix}
a(\xi)& w(\xi)\cdot c(\oooo{\xi})\cdot \delta(\xi)\\
c(\xi)&a(\oooo{\xi})\cdot\delta(\xi)\end{pmatrix} .
\end{eqnarray}
This gives a first approximation of $\Psi(G_\Z)$. It took into account
only the eigenspace $Ml(f)_{\xi,\Z[\zeta]}$ and the pairing $h_\xi$
which $L$ and complex conjugation induce on it. 

\medskip
{\bf Step 3:}  
Now \eqref{5.9} will be shown. 
We will use that the $B_j$ are Orlik blocks and
lemma \ref{t2.8} and \eqref{5.5} for $S_{1,10}$ and
\eqref{5.43} for all other singularities.

Let $g\in\ker\Psi\subset G_\Z$, i.e. 
$g|_{Ml(f)_\zeta}\in\C^*\cdot\id$.
Then $g|_{Ml(f)_\xi}\in\C^*\cdot\id$ for all $\xi$ with
$\Phi_m(\xi)=0$, and 
\begin{eqnarray}\label{5.78}
g((B_j)_{\Phi_m})=(B_j)_{\Phi_m}\quad\textup{ for }j=1,2.
\end{eqnarray}

Now $g(B_j)=B_j$ for $j=1,2$ follows in the case
$S_{1,10}$ from \eqref{5.5}. For all other singularities
$g(B_j)=B_j$ for $j=1,2$ follows with \eqref{5.43} 
(and \eqref{5.32} for $B_3$ in the case
$Z_{1,14r}$).

We want to apply lemma \ref{t2.8} to the Orlik blocks $B_1$
and $B_2$. One checks easily that all hypotheses are 
satisfied. In the case $Z_{1,14r}$ $B_3$ is glued to
$B_1$ by \eqref{5.32}. Therefore in all cases
\begin{eqnarray}\label{5.79}
g=(\varepsilon_1\cdot M_h^{k_1})|_{B_1}
\times (\varepsilon_2\cdot M_h^{k_2})|_{B_2}
\end{eqnarray}
for some $\varepsilon_1,\varepsilon_2\in\{\pm 1\}$
and $k_1,k_2\in\Z$. Now consider
\begin{eqnarray}\label{5.80}
\www g&:=& \varepsilon_2\cdot M_h^{-k_2}\circ g.
\end{eqnarray}
It satisfies 
\begin{eqnarray}
\www g|_{B_1}= \varepsilon_1\varepsilon_2\cdot
M_h^{k_1-k_2}|_{B_1},\quad \www g|_{B_2}=\id,\quad 
\www g|_{{Ml(f)}_\xi}\in\C^*\cdot\id,\nonumber \\
\textup{thus }\quad \www g|_{{Ml(f)_\xi}}=\id,\quad
\www g|_{{Ml(f)_{\Phi_m}}}=\id.\label{5.81}
\end{eqnarray}
Comparison with table \eqref{5.1} shows
\begin{eqnarray*}
\www g&=& \id \qquad\textup{for the first 5 series in }
\eqref{5.1},\\
\www g&=&\id \textup{ or }\www g=-M_h^{\frac{m}{2}(1+2r)}
\quad\textup{ for the last 3 series in }\eqref{5.1}.
\end{eqnarray*}
In any case, $\www g$ and $g$ are in 
$\{\pm M_h^k\, |\, k\in\Z\}$,
and thus $\ker\Psi=\{\pm M_h^k\, |\, k\in\Z\}$.

\medskip
{\bf Step 4:}
By step 2, $\Psi(G_\Z)$ is a subgroup of an infinite
Fuchsian group and therefore itself a Fuchsian group.
It rests to show that it is an infinite group.
By step 3, the kernel of $\Psi:G_\Z\to \Psi(G_\Z)$ 
is $\{\pm M_h^k\, |\, k\in\Z\}$, so it is finite.
Therefore it rests to show that $G_\Z$ is infinite.
We will see that the subgroup of elements $g\in G_\Z$ with
\begin{eqnarray}
g&=&\id \textup{ on any eigenspace }
Ml(f)_\lambda\textup{ with }\Phi_m(\lambda)\neq 0,\nonumber \\
\textup{i.e. }g&=&\id \textup{ on }(\www B_1)_{b_5}
\textup{ and on }(B_2)_{b_6}.\label{5.85}
\end{eqnarray}
is infinite. 

Consider an element $g\in G_\Z$ with \eqref{5.85}.
For all singularities except $S_{1,10}$ \eqref{5.4} holds.
For $S_{1,10}$ \eqref{5.85} implies $g(\gamma_4)=\pm\gamma_4$,
and then \eqref{5.36} gives $g\in\Aut(B_1\oplus B_2,L)$.
In the case of the series $Z_{1,14r}$, the element $g$ maps $B_1\oplus B_2$ 
to itself because $(B_1\oplus B_2)_\C$ contains 
$\ker \Phi_m(M_h)$. 
In any case, lemma \ref{t2.7} applies with 
$k=2, \Lambda^{(1)}=B_1,\Lambda^{(2)}=B_2,
e^{(1)}=\beta_1,e^{(2)}=\beta_2,
p_0=\Phi_m$. By \eqref{2.15} there are unique polynomials 
$p_{ij}\in \Z[t]_{< \deg b_i}$ for $i=1,2$ with
\begin{eqnarray}\label{5.86}
g(\beta_j)&=&p_{1j}(M_h)(\beta_1)+ p_{2j}(M_h)(\beta_2)
\end{eqnarray}
and 
\begin{eqnarray}\label{5.87}
\begin{array}{ll}
p_{11}=1+b_5\cdot q_{11},& p_{12}=b_5\cdot q_{12},\\
p_{21}=b_6\cdot q_{21},& p_{22}=1+b_6\cdot q_{22}
\end{array}
\end{eqnarray}
for suitable polynomials $q_{ij}\in\Z[t]_{<\varphi(m)}$.

$g$ restricts to an automorphism of the pair $(B_1\oplus B_2)_{\Phi_m},L)$.
By \eqref{2.16}, the matrix $A(\xi)$ from \eqref{5.74} in step 2
takes the form
\begin{eqnarray}\label{5.88}
A(\xi)=\begin{pmatrix} 
1+b_5(\xi)q_{11}(\xi) & b_6(\xi)q_{12}(\xi)\\
b_5(\xi)q_{21}(\xi) & 1+b_6(\xi)q_{22}(\xi)
\end{pmatrix} .
\end{eqnarray}
By step 2, this matrix $A(\xi)$ satisfies \eqref{5.75} and \eqref{5.76}.

Vice versa, any polynomials 
$q_{ij}\in \Z[t]_{<\varphi(m)}$ for $i=1,2$
such that the matrix in \eqref{5.88} satisfies \eqref{5.75} and 
\eqref{5.76}, give rise via \eqref{5.87} and \eqref{5.86}
to an element $g\in G_\Z$ with \eqref{5.85}.

We have to prove existence of infinitely many polynomials 
$q_{ij}\in \Z[t]_{<\varphi(m)}$
such that the matrix in \eqref{5.88} satisfies \eqref{5.75} and 
\eqref{5.76} and that $q_{12}(\xi)\neq 0$ and $q_{21}(\xi)\neq 0$.
We start by defining
\begin{eqnarray}\label{5.89}
w_0(\xi):= w(\xi)b_5(\xi)b_5(\oooo{\xi})\in\Z[\zeta]\cap\R
\end{eqnarray}
and asking for infinitely many 
solutions $a(\xi),f(\xi)\in\Z[\zeta]\cap \R$
of the Pell equation
\begin{eqnarray}\label{5.90}
a(\xi)^2-1=w_0(\xi)\cdot f(\xi)^2
\end{eqnarray}
with the additional condition
\begin{eqnarray}
w_0(\xi)&|& a(\xi)-1.\label{5.92}
\end{eqnarray}
Such solutions exist due to lemma \ref{t3.3}.
They give rise to the elements
\begin{eqnarray}\label{5.93}
q_{11}(\xi):= \frac{a(\xi)-1}{b_5(\xi)},&&
q_{12}(\xi):= f(\xi)\cdot\frac{w(\xi)b_5(\oooo\xi)}{b_6(\xi)},\\
q_{21}(\xi):= f(\xi), &&
q_{22}(\xi):= \frac{a(\xi)-1}{b_6(\xi)}.\label{5.94}
\end{eqnarray}
Here observe 
$$b_6(\xi)\, |\, w(\xi)b_5(\oooo\xi)\, |\, w_0(\xi)\, |\, 
a(\xi)-1,$$ 
see \eqref{5.70}, \eqref{5.67} and \eqref{5.65}.
These elements come from unique polynomials 
$q_{ij}\in\Z[t]_{<\varphi(m)}$.
These polynomials satisfy all desired properties.

\subsection{The series $W_{1,p}^\sharp$}\label{c5.1}

\begin{figure}[h]
\begin{tikzpicture}[line cap=round,line join=round,>=triangle 45,x=1.0cm,y=1.0cm]
\clip(-6.5,-1.6) rectangle (8,6);
\draw (-1,-1)-- (0,0);
\draw (0,3)-- (0,2);
\draw (0,1)-- (1,0);
\draw (0,1)-- (-1,-1);
\draw (-1,0)-- (0,1);
\draw (-1,0)-- (0,0);
\draw (0,0)-- (1,0);
\draw (1,0)-- (2,0);
\draw [dash pattern=on 4pt off 4pt] (0,3)-- (1,0);
\draw (0,1)-- (0,2);
\draw [dash pattern=on 4pt off 4pt] (-0.06,0)-- (0,1.01);
\draw [dash pattern=on 4pt off 4pt] (0.07,0)-- (0.07,1);
\draw (-1,0)-- (-2,0);
\draw (-2,0)-- (-3,0);
\draw (0,3)-- (2,0);
\draw (0,3)-- (-2,0);
\draw [dash pattern=on 4pt off 4pt] (-1,0)-- (0,3);
\draw [line width=1.2pt,dotted] (4,0)-- (4.5,0);
\draw [line width=1.2pt,dotted] (2,0)-- (2.95,0);
\draw [line width=1.2pt,dotted] (-5,0)-- (-4.08,0);
\draw [line width=1.2pt,dotted] (-3,0)-- (-3.66,0);
\begin{tiny}
\fill [color=black] (-1,-1) circle (2.0pt);
\draw[color=black] (-1.18,-0.77) node {3};
\fill [color=black] (-1,0) circle (2.0pt);
\draw[color=black] (-1.14,0.25) node {4};
\fill [color=black] (1,0) circle (2.0pt);
\draw[color=black] (1.14,0.27) node {5};
\fill [color=black] (-2,0) circle (2.0pt);
\draw[color=black] (-2.01,0.31) node {8};
\fill [color=black] (2,0) circle (2.0pt);
\draw[color=black] (2.35,-0.24) node {11+q resp. 12+q};
\fill [color=black] (0,1) circle (2.0pt);
\draw[color=black] (0.15,1.26) node {6};
\fill [color=black] (0,2) circle (2.0pt);
\draw[color=black] (0.15,2.27) node {1};
\fill [color=black] (0,3) circle (2.0pt);
\draw[color=black] (0.15,3.27) node {2};
\fill [color=black] (0,0) circle (2.0pt);
\draw[color=black] (0.11,-0.21) node {7};
\fill [color=black] (-3,0) circle (2.0pt);
\draw[color=black] (-2.86,0.27) node {9};
\fill [color=black] (-5,0) circle (2.0pt);
\draw[color=black] (-5.04,0.33) node {10+q resp. 11+q};
\fill [color=black] (4.5,0) circle (2.0pt);
\draw[color=black] (4.88,0.27) node {15+p};
\end{tiny}
\end{tikzpicture}
\caption[Figure 4.1]{The CDD of a
distinguished basis $e_1,\ldots,e_\mu$ for $W^\sharp_{1,2q-1}$ resp. $W^\sharp_{1,2q}$ from \cite[Tabelle 6 \& Abb. 16]{Eb81}}
\label{Fig:5.1}
\end{figure}
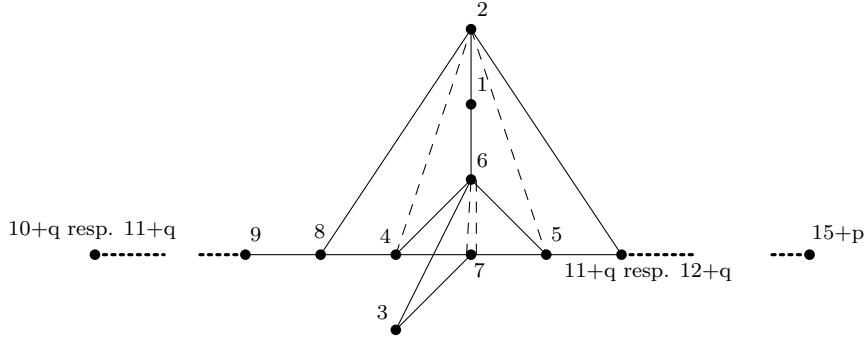

Here we only describe the case when $p=2q$ is even. 
But one can easily obtain the odd case $p=2q-1$ from that 
via replacing each $e_{\alpha+q}$ by $e_{\alpha-1+q}$ in the following lists. 
The monodromy acts on the distinguished basis 
$e_1,\ldots,e_\mu$ with the CDD in figure \ref{Fig:5.1}
as follows:
\begin{eqnarray*}
e_1&\mapsto& -e_1-e_2+e_3+e_4+e_5+e_6,\\
e_2&\mapsto& e_1+e_2+e_8+e_{12+q},\\
e_3&\mapsto& -e_1-e_3-e_6+e_7,\\
e_4&\mapsto& e_2-e_6+e_7+e_8,\\
e_5&\mapsto& e_2-e_6+e_7+e_{12+q},\\
e_6&\mapsto& e_1-2e_2+e_3+e_4+e_5+3e_6-2e_7,\\
e_7&\mapsto& -2e_2+e_3+e_4+e_5+2e_6-e_7,\\
e_{7+i}&\mapsto& e_{8+i}\qquad\textup{for } 1\leq i\leq 3+q,\\
e_{11+q}&\mapsto& -e_4-e_8-e_9-\ldots-e_{11+q},\\
e_{11+q+i}&\mapsto& e_{12+q+i}\qquad\textup{for } 1\leq i\leq 3+q,\\
e_{15+p}&\mapsto& -e_5-e_{12+q}-\ldots-e_{15+p}.
\end{eqnarray*}
By table \eqref{5.13} the generators of the Orlik blocks $B_1$ and $B_2$
are $\beta_1:=e_3$ and $\beta_2:=e_{8}$.
The monodromy acts on them as follows:
\begin{eqnarray}
e_3&\mapsto & -e_1-e_3-e_6+e_7\mapsto e_1+e_2-e_4-e_5-e_6 \nonumber \\ 
 &\mapsto& -e_1 \mapsto  e_{1}+e_{2}-e_{3}-e_{4}-e_{5}-e_{6},
\label{5.1.2}\\
e_{8}&\mapsto& e_{9}\mapsto \ldots \mapsto e_{11+q}
\mapsto -e_4-e_8-e_9- \ldots -e_{11+q}\nonumber\\
&\mapsto& -e_2+e_4+e_6-e_7 \mapsto -e_{12+q} \mapsto -e_{13+q} \mapsto \ldots \mapsto -e_{15+p} \nonumber\\
&\mapsto& e_{5}+e_{12+q}+\ldots+e_{15+p} \mapsto e_2-e_5-e_6+e_7 \mapsto e_{8}.\label{5.1.3}
\end{eqnarray}
Thus the characteristic polynomial 
of $M_h$ on $B_j$ is $b_j$, and the blocks are
\begin{eqnarray}\label{5.1.4}
B_1&=&\langle e_3,e_1,e_6-e_7,e_2-e_4-e_5-e_6\rangle,\\
B_2&=&\langle e_{8}, e_{9}, \ldots, e_{15+p}; e_{4}, e_5, -e_{2}+e_{6}-e_{7}\rangle.
\label{5.1.5}
\end{eqnarray}
This shows that $B_1$ and $B_2$ are primitive sublattices
with $B_1+B_2=B_1 \oplus B_2=Ml\left(f\right)$, i.e. $r_I=1$.

\subsection{The series $S_{1,p}^\sharp$}\label{c5.2}

\begin{figure}[h]
\begin{tikzpicture}[line cap=round,line join=round,>=triangle 45,x=1.0cm,y=1.0cm]
\clip(-5,-2.1) rectangle (9,5);
\draw (-1,-1)-- (0,0);
\draw (0,3)-- (0,2);
\draw (0,1)-- (1,0);
\draw (0,1)-- (-1,-1);
\draw (-1,0)-- (0,1);
\draw (-1,0)-- (0,0);
\draw (0,0)-- (1,0);
\draw (1,0)-- (2,0);
\draw [dash pattern=on 3pt off 3pt] (0,3)-- (1,0);
\draw (0,1)-- (0,2);
\draw [dash pattern=on 3pt off 3pt] (-0.06,0)-- (0,1.01);
\draw [dash pattern=on 3pt off 3pt] (0.07,0)-- (0.07,1);
\draw (-1,0)-- (-2,0);
\draw (-2,-2)-- (-1,-1);
\draw [line width=1.2pt,dotted] (-2,0)-- (-2.55,0);
\draw (0,3)-- (2,0);
\draw (0,3)-- (-2,0);
\draw [dash pattern=on 3pt off 3pt] (-1,0)-- (0,3);
\draw [line width=1.2pt,dotted] (4,0)-- (3.14,0);
\draw [line width=1.2pt,dotted] (2,0)-- (2.74,0);
\draw [line width=1.2pt,dotted] (-3.3,0)-- (-2.8,0);
\begin{tiny}
\fill [color=black] (-1,-1) circle (2.0pt);
\draw[color=black] (-1.13,-0.83) node {3};
\fill [color=black] (-1,0) circle (2.0pt);
\draw[color=black] (-1.11,0.19) node {4};
\fill [color=black] (1,0) circle (2.0pt);
\draw[color=black] (1.1,0.2) node {5};
\fill [color=black] (-2,0) circle (2.0pt);
\draw[color=black] (-2.01,0.23) node {9};
\fill [color=black] (2,0) circle (2.0pt);
\draw[color=black] (2.13,-0.18) node {11+q resp. 12+q};
\fill [color=black] (0,1) circle (2.0pt);
\draw[color=black] (0.11,1.19) node {6};
\fill [color=black] (0,2) circle (2.0pt);
\draw[color=black] (0.11,2.2) node {1};
\fill [color=black] (0,3) circle (2.0pt);
\draw[color=black] (0.11,3.2) node {2};
\fill [color=black] (4,0) circle (2.0pt);
\draw[color=black] (4.27,0.2) node {14+p};
\fill [color=black] (0,0) circle (2.0pt);
\draw[color=black] (0.1,-0.25) node {7};
\fill [color=black] (-2,-2) circle (2.0pt);
\draw[color=black] (-2.06,-1.77) node {8};
\fill [color=black] (-3.3,0) circle (2.0pt);
\draw[color=black] (-3.53,-0.16) node {10+q resp. 11+q};
\end{tiny}
\end{tikzpicture}
\caption[Figure 4.2]{The CDD of a
distinguished basis $e_1,\ldots,e_\mu$ for $S^\sharp_{1,2q-1}$ resp. $S^\sharp_{1,2q}$ from \cite[Tabelle 6 \& Abb. 16]{Eb81}}
\label{Fig:5.2}
\end{figure}
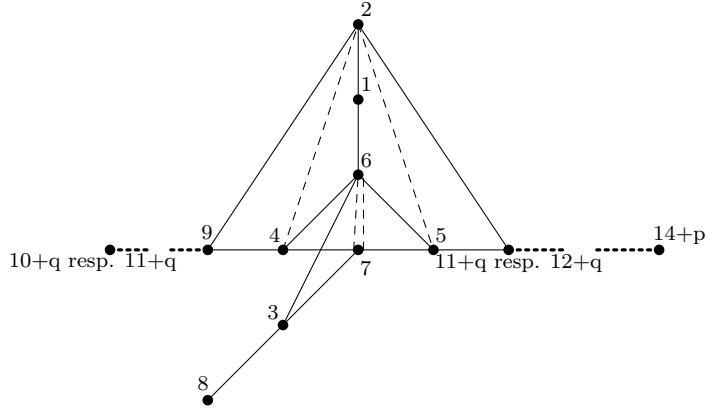

Again we only describe the case when $p=2q$ is even. 
But one can easily obtain the odd case
$p=2q-1$ from that via replacing each 
$e_{\alpha+q}$ by $e_{\alpha-1+q}$ in the following lists. 
The monodromy acts on the distinguished basis 
$e_1,\ldots,e_\mu$ with the CDD in figure \ref{Fig:5.2}
as follows:
\begin{eqnarray*}
e_1&\mapsto& -e_1-e_2+e_3+e_4+e_5+e_6,\\
e_2&\mapsto& e_1+e_2+e_9+e_{12+q},\\
e_3&\mapsto& -e_1-e_6+e_7+e_8,\\
e_4&\mapsto& e_2-e_6+e_7+e_9,\\
e_5&\mapsto& e_2-e_6+e_7+e_{12+q},\\
e_6&\mapsto& e_1-2e_2+e_3+e_4+e_5+3e_6-2e_7,\\
e_7&\mapsto& -2e_2+e_3+e_4+e_5+2e_6-e_7,\\
e_8&\mapsto& -e_3-e_8,\\
e_{8+i}&\mapsto& e_{9+i}\qquad\textup{for } 1\leq i\leq 2+q,\\
e_{11+q}&\mapsto& -e_4-e_9-e_{10}-\ldots-e_{11+q},\\
e_{11+q+i}&\mapsto& e_{12+q+i}\qquad\textup{for } 1\leq i\leq 2+q,\\
e_{14+p}&\mapsto& -e_5-e_{12+q}-e_{13+q}-\ldots-e_{14+p}.
\end{eqnarray*}
By table \eqref{5.13} the generators of the Orlik blocks $B_1$ and $B_2$
are $\beta_1:=e_8$ and $\beta_2:=e_{9}$.
The monodromy acts on them as follows:
\begin{eqnarray}
e_8&\mapsto & -e_3-e_8\mapsto e_1+e_3+e_6-e_7\nonumber \\
&\mapsto& -e_1-e_2+e_3+e_4+e_5+e_6+e_8 \nonumber \\
&\mapsto& -e_3-e_6+e_7 \mapsto  -e_8 ,\label{5.2.2}\\
e_9&\mapsto& e_{10}\mapsto \ldots \mapsto e_{11+q}
\mapsto -e_4-e_9-e_{10}- \ldots -e_{11+q}\nonumber\\
&\mapsto& -e_2+e_4+e_6-e_7 \mapsto -e_{12+q} \mapsto -e_{13+q} 
\mapsto \ldots \mapsto -e_{14+p} \nonumber\\
&\mapsto& e_{5}+e_{12+q}+\ldots+e_{14+p} \mapsto e_2-e_5-e_6+e_7 \mapsto e_{9}.
\label{5.2.3}
\end{eqnarray}
Thus the characteristic polynomial of $M_h$ on $B_j$ is $b_j$, 
and the blocks are
\begin{eqnarray}\label{5.2.4}
B_1&=&\langle e_8, e_3, e_6-e_7, e_1, -e_2+e_4+e_5+e_6 \rangle,\\
B_2&=&\langle e_{9}, e_{10}, \ldots, e_{14+p} ;e_{4},e_{5},-e_{2}+e_{6}-e_{7}\rangle.
\label{5.2.5}
\end{eqnarray}
This shows that $B_1$ and $B_2$ are primitive sublattices
with $B_1+B_2=B_1 \oplus B_2=Ml\left(f\right)$ and $r_I=1$.

\subsection{The series $U_{1,p}$}\label{c5.3}

\begin{figure}[h]
\begin{tikzpicture}[line cap=round,line join=round,>=triangle 45,x=1.0cm,y=1.0cm]
\clip(-6,-3.2) rectangle (11,4.5);
\draw (-1,-1)-- (0,0);
\draw (0,3)-- (0,2);
\draw (0,1)-- (1,0);
\draw (0,1)-- (-1,-1);
\draw (-1,0)-- (0,1);
\draw (-1,0)-- (0,0);
\draw (0,0)-- (1,0);
\draw (1,0)-- (2,0);
\draw [dash pattern=on 4pt off 4pt] (0,3)-- (1,0);
\draw (0,1)-- (0,2);
\draw [dash pattern=on 4pt off 4pt] (-0.06,0)-- (-0.06,1);
\draw [dash pattern=on 4pt off 4pt] (0.07,0)-- (0.07,1);
\draw (-1,0)-- (-2,0);
\draw (-2,-2)-- (-1,-1);
\draw (-2,0)-- (-3,0);
\draw (2,0)-- (3,0);
\draw (0,3)-- (2,0);
\draw (0,3)-- (-2,0);
\draw [dash pattern=on 4pt off 4pt] (-1,0)-- (0,3);
\draw (-3,-3)-- (-2,-2);
\draw [line width=1.2pt,dotted] (-3,0)-- (-3.72,0);
\draw [line width=1.2pt,dotted] (-4.53,0)-- (-4,0);
\draw [line width=1.2pt,dotted] (3,0)-- (3.66,0);
\draw [line width=1.2pt,dotted] (4.37,0)-- (4,0);
\begin{tiny}
\fill [color=black] (-1,-1) circle (2.0pt);
\draw[color=black] (-1.15,-0.81) node {3};
\fill [color=black] (-1,0) circle (2.0pt);
\draw[color=black] (-1.12,0.22) node {4};
\fill [color=black] (1,0) circle (2.0pt);
\draw[color=black] (1.12,0.23) node {5};
\fill [color=black] (-2,0) circle (2.0pt);
\draw[color=black] (-2.08,0.28) node {10};
\fill [color=black] (2,0) circle (2.0pt);
\draw[color=black] (2.31,0.23) node {12+q};
\fill [color=black] (0,1) circle (2.0pt);
\draw[color=black] (0.12,1.22) node {6};
\fill [color=black] (0,2) circle (2.0pt);
\draw[color=black] (0.12,2.23) node {1};
\fill [color=black] (0,3) circle (2.0pt);
\draw[color=black] (0.12,3.22) node {2};
\fill [color=black] (3,0) circle (2.0pt);
\draw[color=black] (3.32,0.23) node {13+q};
\fill [color=black] (4.37,0) circle (2.0pt);
\draw[color=black] (4.69,0.23) node {14+p};
\fill [color=black] (0,0) circle (2.0pt);
\draw[color=black] (0.1,-0.28) node {7};
\fill [color=black] (-2,-2) circle (2.0pt);
\draw[color=black] (-2.08,-1.72) node {8};
\fill [color=black] (-3,0) circle (2.0pt);
\draw[color=black] (-2.81,0.23) node {11};
\fill [color=black] (-4.53,0) circle (2.0pt);
\draw[color=black] (-4.42,0.26) node {11+q};
\fill [color=black] (-3,-3) circle (2.0pt);
\draw[color=black] (-3.04,-2.73) node {9};
\end{tiny}
\end{tikzpicture}
\caption[Figure 4.3]{The CDD of a
distinguished basis $e_1,\ldots,e_\mu$ for $U_{1,p}$
 from \cite[Tabelle 6 \& Abb. 16]{Eb81}}
\label{Fig:5.3}
\end{figure}
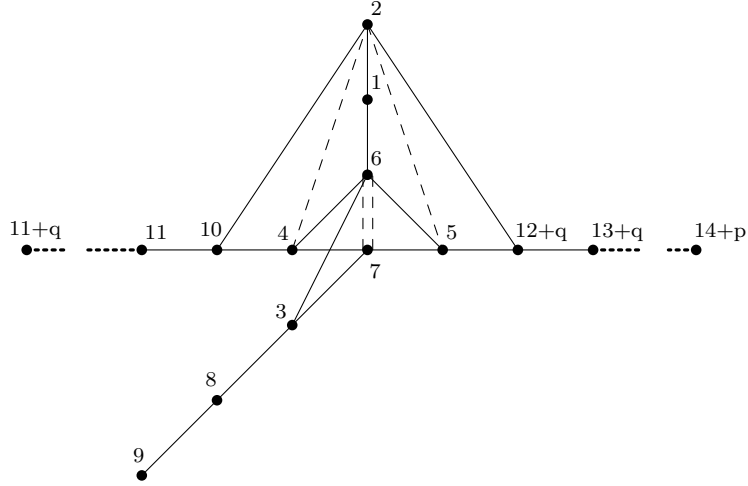

Here (and in all series except $W_{1,p}^\sharp$ and $S_{1,p}^\sharp$)
the list of the monodromy action 
on the distinguished basis 
$e_1,\ldots,e_\mu$ with the CDD in figure \ref{Fig:5.3}
includes both cases 
$p=2q$ and $p=2q-1$. It looks as follows: 
\begin{eqnarray*}
e_1&\mapsto& -e_1-e_2+e_3+e_4+e_5+e_6,\\
e_2&\mapsto& e_1+e_2+e_{10}+e_{12+q},\\
e_3&\mapsto& -e_1-e_6+e_7+e_8,\\
e_4&\mapsto& e_2-e_6+e_7+e_{10},\\
e_5&\mapsto& e_2-e_6+e_7+e_{12+q},\\
e_6&\mapsto& e_1-2e_2+e_3+e_4+e_5+3e_6-2e_7,\\
e_7&\mapsto& -2e_1+e_3+e_4+e_5+2e_6-e_7,\\
e_8&\mapsto& e_9,\\
e_9&\mapsto& -e_3-e_8-e_9,\\
e_{9+i}&\mapsto& e_{10+i}\qquad\textup{for } 1\leq i\leq 1+q,\\
e_{11+q}&\mapsto& -e_{4}-e_{10}-e_{11}-\ldots-e_{11+q},\\
e_{11+q+i}&\mapsto& e_{12+q+i} \qquad\textup{for } 1\leq i\leq 2+p-q,\\
e_{14+p}&\mapsto& -e_5-e_{12+q}-e_{13+q}-\ldots-e_{14+p}.
\end{eqnarray*}
By table \eqref{5.13} the generators of the Orlik blocks $B_1$ and $B_2$
are $\beta_1:=e_8$ and $\beta_2:=e_{10}$.
The monodromy acts on them as follows:
\begin{eqnarray}
e_8 &\mapsto & e_9\mapsto -e_3-e_8-e_9 \mapsto e_1+e_3+e_6-e_7 \nonumber \\
&\mapsto& -e_1-e_2+e_3+e_4+e_5+e_6+e_8\nonumber \\
&\mapsto& -e_6+e_7+e_8+e_9 \mapsto -e_1-e_3-e_6+e_7-e_8,\label{5.3.2}\\
e_{10}&\mapsto& e_{11}\mapsto \ldots \mapsto e_{11+q}
\mapsto -e_4-e_{10}-e_{11} - \ldots -e_{11+q}\nonumber\\
&\mapsto& -e_2+e_4+e_6-e_7 \mapsto -e_{12+q} \mapsto -e_{13+q} \mapsto \ldots \mapsto -e_{14+p} \nonumber\\
&\mapsto& e_{5}+e_{12+q}+\ldots+e_{14+p} \mapsto e_2-e_5-e_6+e_7 
\mapsto e_{10}.\label{5.3.3}
\end{eqnarray}
Thus the characteristic polynomial of $M_h$ on $B_j$ is $b_j$, 
and  the blocks are
\begin{eqnarray}\label{5.3.4}
B_1&=&\langle e_1, e_3, e_8, e_9, e_6-e_7, -e_2+e_4+e_5+e_6 \rangle,\\
B_2&=&\langle e_{10}, e_{11},\ldots,e_{14+p};e_{4}, e_{5}, -e_2+e_6-e_7  \rangle.
\label{5.3.5}
\end{eqnarray}

Again $B_1$ and $B_2$ are primitive sublattices
with $B_1+B_2=B_1 \oplus B_2=Ml\left(f\right)$ and $r_I=1$.

\subsection{The series $E_{3,p}$}\label{c5.4}

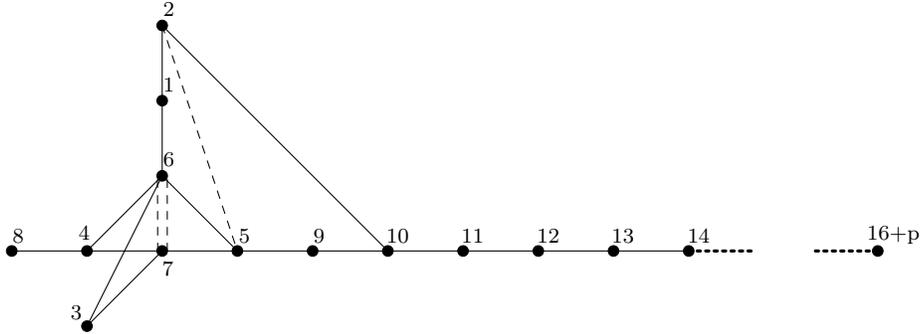
\begin{figure}[h]
\centering
\begin{tikzpicture}[line cap=round,line join=round,>=triangle 45,x=1.0cm,y=1.0cm]
\clip(-2.5,-1.4) rectangle (11,3.4););
\draw (-1.,-1.)-- (0.,0.);
\draw (0.,3.)-- (0.,2.);
\draw (0.,1.)-- (1.,0.);
\draw (0.,1.)-- (-1.,-1.);
\draw (-1.,0.)-- (0.,1.);
\draw (-1.,0.)-- (0.,0.);
\draw (0.,0.)-- (1.,0.);
\draw (1.,0.)-- (2.,0.);
\draw (-1.,0.)-- (-2.,0.);
\draw [dash pattern=on 3pt off 3pt] (0.,3.)-- (1.,0.);
\draw (0.,3.)-- (3.,0.);
\draw (2.,0.)-- (3.,0.);
\draw (0.,1.)-- (0.,2.);
\draw (3.,0.)-- (4.,0.);
\draw (4.,0.)-- (5.,0.);
\draw (5.,0.)-- (6.,0.);
\draw (6.,0.)-- (7.,0.);
\draw [dash pattern=on 3pt off 3pt] (-0.05921624320994279,0.)-- (-0.06285921894269199,1.0062273461178708);
\draw [dash pattern=on 3pt off 3pt] (0.06746025531853216,0.)-- (0.06732116058008084,1.0025079067029345);
\draw [line width=1.2pt,dotted] (7.,0.)-- (7.865703985214759,0.);
\draw [line width=1.2pt,dotted] (9.512977672125807,0.)-- (8.637046425911203,0.);
\begin{tiny}
\draw [fill=black] (-1.,-1.) circle (2.0pt);
\draw[color=black] (-1.142006890036921,-0.8124455505724937) node {3};
\draw [fill=black] (-1.,0.) circle (2.0pt);
\draw[color=black] (-1.0374180845187593,0.24651610529889312) node {4};
\draw [fill=black] (1.,0.) circle (2.0pt);
\draw[color=black] (1.0935788279137866,0.2072953032295825) node {5};
\draw [fill=black] (-2.,0.) circle (2.0pt);
\draw[color=black] (-1.9133493307333638,0.2072953032295825) node {8};
\draw [fill=black] (2.,0.) circle (2.0pt);
\draw[color=black] (2.0871724803363234,0.2072953032295825) node {9};
\draw [fill=black] (0.,1.) circle (2.0pt);
\draw[color=black] (0.08691157480147982,1.2139625563418885) node {6};
\draw [fill=black] (0.,2.) circle (2.0pt);
\draw[color=black] (0.08691157480147982,2.2075562087644243) node {1};
\draw [fill=black] (0.,3.) circle (2.0pt);
\draw[color=black] (0.08691157480147982,3.2142234618767302) node {2};
\draw [fill=black] (3.,0.) circle (2.0pt);
\draw[color=black] (3.1330605355179406,0.2072953032295825) node {10};
\draw [fill=black] (4.,0.) circle (2.0pt);
\draw[color=black] (4.126654187940478,0.2072953032295825) node {11};
\draw [fill=black] (5.,0.) circle (2.0pt);
\draw[color=black] (5.133321441052785,0.2072953032295825) node {12};
\draw [fill=black] (6.,0.) circle (2.0pt);
\draw[color=black] (6.126915093475321,0.2072953032295825) node {13};
\draw [fill=black] (7.,0.) circle (2.0pt);
\draw[color=black] (7.133582346587628,0.2072953032295825) node {14};
\draw [fill=black] (0.,0.) circle (2.0pt);
\draw[color=black] (0.0738379741117096,-0.23720712022260457) node {7};
\draw [fill=black] (9.512977672125807,0.) circle (2.0pt);
\draw[color=black] (9.722155283162131,0.2072953032295825) node {16+p};
\end{tiny}
\end{tikzpicture}
\caption[Figure 5.4]{The CDD of a
distinguished basis $e_1,\ldots,e_\mu$ for $E_{3,p}$
 from \cite[Tabelle 6 \& Abb. 16]{Eb81}}
\label{Fig:5.4}
	\end{figure}

Here the monodromy acts on the distinguished basis 
$e_1,\ldots,e_\mu$ with the CDD in figure \ref{Fig:5.4}
as follows:
\begin{eqnarray*}
e_1&\mapsto& e_3+e_4+e_5+e_6,\\
e_2&\mapsto& e_9+e_{10},\\
e_3&\mapsto& -e_1-e_3-e_6+e_7,\\
e_4&\mapsto& -e_1-e_6+e_7+e_8,\\
e_5&\mapsto& -e_1-e_6+e_7+e_9,\\
e_6&\mapsto& 2e_1-e_2+e_3+e_4+e_5+3e_6-2e_7,\\
e_7&\mapsto& e_1-e_2+e_3+e_4+e_5+2e_6-e_7,\\
e_8&\mapsto& -e_4-e_8,\\
e_9&\mapsto& e_1+e_2+e_{10},\\
e_{9+i}&\mapsto& e_{10+i}\qquad\textup{for } 1\leq i\leq 6+p,\\
e_{16+p}&\mapsto& -e_5-e_9-e_{10}-\ldots-e_{16+p}.
\end{eqnarray*}
By table \eqref{5.13} the generators of the Orlik blocks $B_1$ and $B_2$
are $\beta_1:=e_3$ and $\beta_2:=e_{10}$.
The monodromy acts on them as follows:
\begin{eqnarray}
e_3&\mapsto & -e_1-e_3-e_6+e_7\mapsto -e_4-e_5-e_6 \nonumber \\
&\mapsto& e_2-e_3-e_4-e_5-e_6-e_8-e_9\mapsto -e_5-e_7\nonumber \\
&\mapsto& e_2-e_3-e_4-e_5-e_6-e_9\mapsto -e_4-e_5-e_7-e_8\nonumber \\
&\mapsto& e_1+e_2-e_3-e_5-e_7-e_9\mapsto e_3+e_6-e_7\nonumber\\
&\mapsto& -e_3,
\label{5.4.2}\\
e_{10}&\mapsto& e_{11}\mapsto \ldots\mapsto e_{16+p}
\mapsto -e_5-\sum_{i=9}^{16+p}e_i \nonumber\\
&\mapsto& -e_2+e_5+e_6-e_7\mapsto -e_{10}.\label{5.4.3}
\end{eqnarray}
Thus the characteristic polynomial of $M_h$ on $B_j$ is $b_j$, 
and the blocks are
\begin{eqnarray}\label{5.4.4}
B_1&=&\langle e_1,e_3,e_4,e_8,e_6-e_7,e_5+e_6,e_2-e_9\rangle,\\
B_2&=&\langle e_{10},e_{11},\ldots,e_{16+p},e_5+e_9,e_2-e_6+e_7+e_9\rangle.
\label{5.4.5}
\end{eqnarray}
This shows that $B_1$ and $B_2$ are primitive sublattices
with $B_1+B_2=B_1\oplus B_2$. Furthermore 
$B_1\oplus B_2\supset\{2e_2\}$ and $B_1+B_2+\Z\cdot e_2=Ml(f)$.
This shows $[Ml(f):B_1\oplus B_2]=2=r_I$.

\subsection{The series $Z_{1,p}$}\label{c5.5}

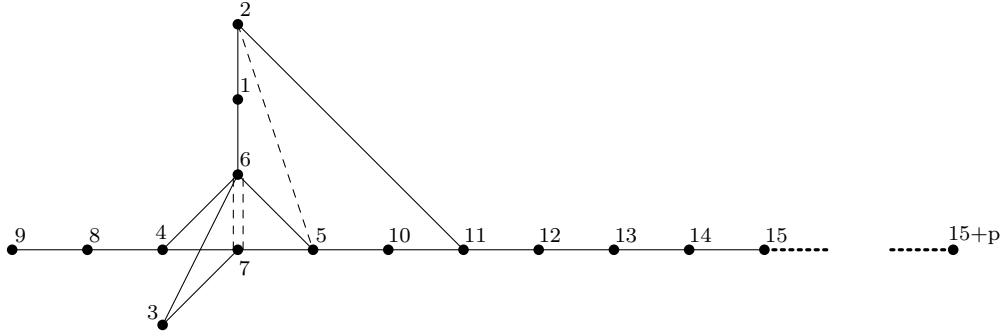
\begin{figure}[h]
\definecolor{ttqqqq}{rgb}{0,0,0}
\begin{tikzpicture}[line cap=round,line join=round,>=triangle 45,x=1.0cm,y=1.0cm]
\clip(-3.2,-1.3) rectangle (11,3.5);
\draw (-1,-1)-- (0,0);
\draw (0,3)-- (0,2);
\draw (0,1)-- (1,0);
\draw (0,1)-- (-1,-1);
\draw (-1,0)-- (0,1);
\draw (-1,0)-- (0,0);
\draw (0,0)-- (1,0);
\draw (1,0)-- (2,0);
\draw (-1,0)-- (-2,0);
\draw (-2,0)-- (-3,0);
\draw [dash pattern=on 3pt off 3pt] (0,3)-- (1,0);
\draw (0,3)-- (3,0);
\draw (2,0)-- (3,0);
\draw (0,1)-- (0,2);
\draw (3,0)-- (4,0);
\draw (4,0)-- (5,0);
\draw (5,0)-- (6,0);
\draw (6,0)-- (7,0);
\draw [dash pattern=on 3pt off 3pt] (-0.06,0)-- (-0.06,1.01);
\draw [dash pattern=on 3pt off 3pt] (0.07,0)-- (0.07,1);
\draw [line width=1.2pt,dotted] (7,0)-- (7.87,0);
\draw [line width=1.2pt,dotted] (9.51,0)-- (8.64,0);
\begin{tiny}
\fill [color=ttqqqq] (-1,-1) circle (2.0pt);
\draw[color=ttqqqq] (-1.13,-0.83) node {3};
\fill [color=ttqqqq] (-1,0) circle (2.0pt);
\draw[color=ttqqqq] (-1.02,0.23) node {4};
\fill [color=ttqqqq] (1,0) circle (2.0pt);
\draw[color=ttqqqq] (1.11,0.19) node {5};
\fill [color=ttqqqq] (-2,0) circle (2.0pt);
\draw[color=ttqqqq] (-1.9,0.19) node {8};
\fill [color=ttqqqq] (-3,0) circle (2.0pt);
\draw[color=ttqqqq] (-2.89,0.19) node {9};
\fill [color=ttqqqq] (2,0) circle (2.0pt);
\draw[color=ttqqqq] (2.15,0.19) node {10};
\fill [color=ttqqqq] (0,1) circle (2.0pt);
\draw[color=ttqqqq] (0.1,1.2) node {6};
\fill [color=ttqqqq] (0,2) circle (2.0pt);
\draw[color=ttqqqq] (0.1,2.19) node {1};
\fill [color=ttqqqq] (0,3) circle (2.0pt);
\draw[color=ttqqqq] (0.1,3.2) node {2};
\fill [color=ttqqqq] (3,0) circle (2.0pt);
\draw[color=ttqqqq] (3.16,0.19) node {11};
\fill [color=ttqqqq] (4,0) circle (2.0pt);
\draw[color=ttqqqq] (4.15,0.19) node {12};
\fill [color=ttqqqq] (5,0) circle (2.0pt);
\draw[color=ttqqqq] (5.16,0.19) node {13};
\fill [color=ttqqqq] (6,0) circle (2.0pt);
\draw[color=ttqqqq] (6.15,0.19) node {14};
\fill [color=ttqqqq] (7,0) circle (2.0pt);
\draw[color=ttqqqq] (7.16,0.19) node {15};
\fill [color=ttqqqq] (0,0) circle (2.0pt);
\draw[color=ttqqqq] (0.09,-0.25) node {7};
\fill [color=black] (9.51,0) circle (2.0pt);
\draw[color=black] (9.79,0.19) node {15+p};
\end{tiny}
\end{tikzpicture}
\caption[Figure 4.5]{The CDD of a
distinguished basis $e_1,\ldots,e_\mu$ for $Z_{1,p}$
 from \cite[Tabelle 6 \& Abb. 16]{Eb81}}
\label{Fig:5.5}
\end{figure}

Here the monodromy acts on the distinguished basis 
$e_1,\ldots,e_\mu$ with the CDD in figure \ref{Fig:5.5}
as follows:
\begin{eqnarray*}
e_1&\mapsto& e_3+e_4+e_5+e_6,\\
e_2&\mapsto& e_{10}+e_{11},\\
e_3&\mapsto& -e_1-e_3-e_6+e_7,\\
e_4&\mapsto& -e_1-e_6+e_7+e_8,\\
e_5&\mapsto& -e_1-e_6+e_7+e_{10},\\
e_6&\mapsto& 2e_1-e_2+e_3+e_4+e_5+3e_6-2e_7,\\
e_7&\mapsto& e_1-e_2+e_3+e_4+e_5+2e_6-e_7,\\
e_8&\mapsto& e_9,\\
e_9&\mapsto& -e_4-e_8-e_9,\\
e_{10}&\mapsto& e_1+e_2+e_{11},\\
e_{11+i}&\mapsto& e_{12+i}\qquad\textup{for } 1\leq i\leq 3+p,\\
e_{15+p}&\mapsto& -e_5-e_{10}-e_{11}-\ldots-e_{15+p}.
\end{eqnarray*}
Here there are three Orlik blocks $B_1,B_2$ and $B_3$. 
By table \eqref{5.13} their generators 
are $\beta_1:=e_8, \beta_2:=e_{11}$ and $\beta_3:=e_3+e_4-e_9$.
The monodromy acts on them as follows:
\begin{eqnarray}
e_8&\mapsto& e_{9} \mapsto -e_4-e_8-e_9 \mapsto e_1+e_4+e_6-e_7 \nonumber \\
&\mapsto&  e_3+e_4+e_5+e_6+e_8 \nonumber\\
&\mapsto& -e_1-e_2+e_4+e_5+e_7+e_8+e_9+e_{10} \nonumber \\
&\mapsto& -e_4-e_6+e_7\mapsto -e_8,
\label{5.5.2}\\
e_{11}&\mapsto& e_{12}\mapsto \ldots \mapsto e_{15+p}
\mapsto -e_5-\sum_{i=10}^{15+p}e_i \nonumber\\
&\mapsto& -e_2+e_5+e_6-e_7\mapsto -e_{11},\label{5.5.3}\\
e_{3}-e_4-e_9 &\mapsto& -e_3+e_4+e_9. \label{5.5.4}
\end{eqnarray}
Thus the characteristic polynomial of $M_h$ on $B_j$ is $b_j$, 
and the blocks are
\begin{eqnarray}
B_1&=&\langle e_8, e_9, e_4, e_1, e_6-e_7, e_3+e_5+e_6, \nonumber\\
&& -e_2+e_5+e_7+e_{10}\rangle,\label{5.5.5}\\
B_2&=&\langle e_{11}, e_{12}, \ldots, e_{15+p}; e_{5}+e_{10}, -e_{2}+e_{5}+e_{6}-e_{7}\rangle,\\
B_3&=&\langle e_{3}-e_{4}-e_{9}\rangle.
\label{5.5.6}
\end{eqnarray}
This shows that $B_1,B_2$ and $B_3$ are primitive sublattices
with $B_1+B_2+B_3=B_1\oplus B_2\oplus B_3$. Furthermore 
$B_1\oplus B_2\oplus B_3\supset\{2e_5\}$ and $B_1+B_2+B_3+\Z\cdot e_5=Ml(f)$.
This shows $[Ml(f):B_1\oplus B_2\oplus B_3]=2=r_I$.

\subsection{The series $Q_{2,p}$}\label{c5.6}

\begin{figure}[h]
\centering
\definecolor{ttqqqq}{rgb}{0,0,0}
\begin{tikzpicture}[line cap=round,line join=round,>=triangle 45,x=1.0cm,y=1.0cm]
\clip(-2.5,-2.2) rectangle (10,4);
\draw (-1,-1)-- (0,0);
\draw (0,3)-- (0,2);
\draw (0,1)-- (1,0);
\draw (0,1)-- (-1,-1);
\draw (-1,0)-- (0,1);
\draw (-1,0)-- (0,0);
\draw (0,0)-- (1,0);
\draw (1,0)-- (2,0);
\draw [dash pattern=on 4pt off 4pt] (0,3)-- (1,0);
\draw (0,3)-- (3,0);
\draw (2,0)-- (3,0);
\draw (0,1)-- (0,2);
\draw (3,0)-- (4,0);
\draw (4,0)-- (5,0);
\draw (5,0)-- (6,0);
\draw [dash pattern=on 4pt off 4pt] (-0.06,0)-- (-0.06,1.01);
\draw [dash pattern=on 4pt off 4pt] (0.07,0)-- (0.07,1);
\draw (-1,0)-- (-2,0);
\draw (-2,-2)-- (-1,-1);
\draw [line width=1.2pt,dotted] (6,0)-- (6.64,0);
\draw [line width=1.2pt,dotted] (7.96,0)-- (7.24,0);
\begin{tiny}
\fill [color=ttqqqq] (-1,-1) circle (2.0pt);
\draw[color=ttqqqq] (-1.15,-0.81) node {3};
\fill [color=ttqqqq] (-1,0) circle (2.0pt);
\draw[color=ttqqqq] (-1.03,0.28) node {4};
\fill [color=ttqqqq] (1,0) circle (2.0pt);
\draw[color=ttqqqq] (1.12,0.23) node {5};
\fill [color=ttqqqq] (-2,0) circle (2.0pt);
\draw[color=ttqqqq] (-1.88,0.23) node {9};
\fill [color=ttqqqq] (2,0) circle (2.0pt);
\draw[color=ttqqqq] (2.18,0.23) node {10};
\fill [color=ttqqqq] (0,1) circle (2.0pt);
\draw[color=ttqqqq] (0.12,1.22) node {6};
\fill [color=ttqqqq] (0,2) circle (2.0pt);
\draw[color=ttqqqq] (0.12,2.23) node {1};
\fill [color=ttqqqq] (0,3) circle (2.0pt);
\draw[color=ttqqqq] (0.12,3.22) node {2};
\fill [color=ttqqqq] (3,0) circle (2.0pt);
\draw[color=ttqqqq] (3.18,0.23) node {11};
\fill [color=ttqqqq] (4,0) circle (2.0pt);
\draw[color=ttqqqq] (4.18,0.23) node {12};
\fill [color=ttqqqq] (5,0) circle (2.0pt);
\draw[color=ttqqqq] (5.18,0.23) node {13};
\fill [color=ttqqqq] (6,0) circle (2.0pt);
\draw[color=ttqqqq] (6.17,0.23) node {14};
\fill [color=ttqqqq] (0,0) circle (2.0pt);
\draw[color=ttqqqq] (0.1,-0.28) node {7};
\fill [color=ttqqqq] (-2,-2) circle (2.0pt);
\draw[color=ttqqqq] (-2.08,-1.72) node {8};
\fill [color=ttqqqq] (7.96,0) circle (2.0pt);
\draw[color=ttqqqq] (8.01,0.28) node {14+p};
\end{tiny}
\end{tikzpicture}
\caption[Figure 5.6]{The CDD of a
distinguished basis $e_1,\ldots,e_\mu$ for $Q_{2,p}$
 from \cite[Tabelle 6 \& Abb. 16]{Eb81}}
\label{Fig:5.6}
\end{figure}
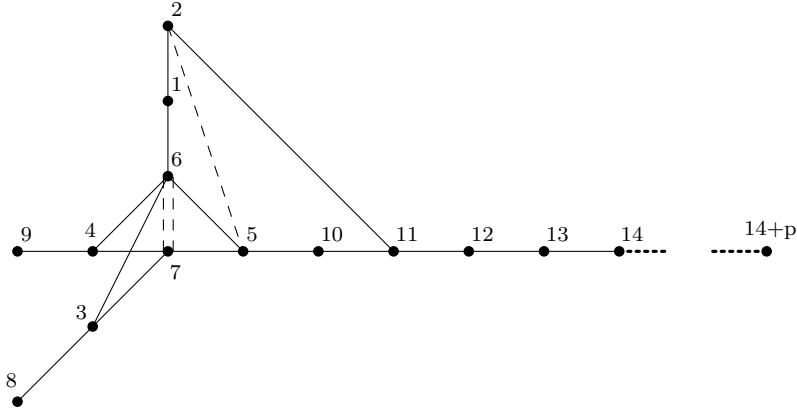

Here the monodromy acts on the distinguished basis 
$e_1,\ldots,e_\mu$ with the CDD in figure \ref{Fig:5.6}
as follows: 
\begin{eqnarray*}
e_1&\mapsto& e_3+e_4+e_5+e_6,\\
e_2&\mapsto& e_{10}+e_{11},\\
e_3&\mapsto& -e_1-e_6+e_7+e_8,\\
e_4&\mapsto& -e_1-e_6+e_7+e_9,\\
e_5&\mapsto& -e_1-e_6+e_7+e_{10},\\
e_6&\mapsto& 2e_1-e_2+e_3+e_4+e_5+3e_6-2e_7,\\
e_7&\mapsto& e_1-e_2+e_3+e_4+e_5+2e_6-e_7,\\
e_8&\mapsto& -e_3-e_8,\\
e_9&\mapsto& -e_4-e_9,\\
e_{10}&\mapsto& e_1+e_2+e_{11},\\
e_{10+i}&\mapsto& e_{11+i}\qquad\textup{for } 1\leq i\leq 3+p,\\
e_{14+p}&\mapsto& -e_5-e_{10}-e_{11}-\ldots-e_{14+p}.
\end{eqnarray*}
By table \eqref{5.13} the generators of the Orlik blocks $B_1$ and $B_2$
are $\beta_1:=e_8$ and $\beta_2:=e_{11}$.
The monodromy acts on them as follows:
\begin{eqnarray}
e_8 &\mapsto&  -e_3-e_8 \mapsto e_1+e_3+e_6-e_7 \mapsto e_3+e_4+e_5+e_6+e_8 \nonumber \\
&\mapsto& -e_1-e_2+e_4+e_5+e_7+e_9+e_{10} \mapsto-e_4-e_6+e_7 \nonumber \\
&\mapsto& -e_9 \mapsto e_4+e_9 \mapsto -e_1-e_4-e_6+e_7, 
\label{5.6.2}\\
e_{11}&\mapsto& e_{12}\mapsto \ldots \mapsto e_{14+p} \mapsto -e_5-\sum_{i=10}^{14+p}e_i \nonumber\\
&\mapsto& -e_2+e_5+e_6-e_7 \mapsto -e_{11}.\label{5.6.3}
\end{eqnarray}
Thus the characteristic polynomial of $M_h$ on $B_j$ is $b_j$, 
and the blocks are
\begin{eqnarray}
B_1&=&\langle e_8, e_3, e_9, e_4, e_1, e_6-e_7, \nonumber \\
&&e_5+e_6, -e_2+e_5+e_7+e_{10}\rangle, \label{5.6.4}\\
B_2&=&\langle e_{11}, e_{12}, \ldots, e_{14+p}; e_{5}+e_{10}, -e_{2}+e_{5}+e_{6}-e_{7}\rangle.
\label{5.6.5}
\end{eqnarray}
This shows that $B_1$ and $B_2$ are primitive sublattices
with $B_1+B_2=B_1\oplus B_2$. Furthermore 
$B_1\oplus B_2\supset\{2e_5\}$ and $B_1+B_2+\Z\cdot e_5=Ml(f)$.
This shows $[Ml(f):B_1\oplus B_2]=2=r_I$.

The proof of \eqref{5.4} for $Q_{2,p}$ was postponed to this subsection
and has to be given now. Recall the definition \eqref{5.37} of $b_4$
and recall $b_4=\Phi_4$ for $Q_{2,4s}$ and $b_4=1$ for the other $Q_{2,p}$.
The next aims are:
\begin{list}{}{}
\item[(i)] 
For $Q_{2,4s}$: To show for any 
$g\in G_\Z\cup\Aut(B_1\oplus B_2,L)$ 
\begin{eqnarray}\label{5.6.6} g:(B_1)_{b_4}\to (B_1)_{b_4}
\textup{ and }(B_2)_{b_4}\to (B_2)_{b_4}.
\end{eqnarray}
\item[(ii)]
For all $Q_{2,p}$: To find an element 
$\gamma_4\in (B_1)_{\Phi_4}$ with
\begin{eqnarray}\label{5.6.7}
B_1\oplus B_2 &=& \{a\in Ml(f)\, |\, L(a,\gamma_4)\equiv 0(2)\}\\
&=&\{a\in Ml(f)\, |\, L(a,M_h(\gamma_4))\equiv 0(2)\},\nonumber \\
g(\gamma_4)&\in&\{\pm\gamma_4,\pm M_h(\gamma_4)\}\qquad
\textup{for any }g\in G_\Z.\label{5.6.8}
\end{eqnarray}
\item[(iii)]
For all $Q_{2,p}$: To find an element $\gamma_5\in Ml(f)$ with
\begin{eqnarray}\label{5.6.9}
B_1+B_2+\Z\cdot \gamma_5=Ml(f)\\
\textup{and }\quad g(\gamma_5)\in Ml(f) \textup{ for any }
g\in \Aut(B_1\oplus B_2,L).\label{5.6.10}
\end{eqnarray}
\end{list}

For all $Q_{2,p}$ define
\begin{eqnarray}
\gamma_1&:=& \frac{b_1}{\Phi_4}(M_h)(\beta_1) = (\Phi_{12}\Phi_3)(M_h)(e_8)
\nonumber\\
&=& (t^6+t^5-t^3+t+1)(M_h)(e_8)\nonumber\\
&=& -2e_3-2e_4-e_5-2e_6+e_7-e_8-e_9.\label{5.6.11}
\end{eqnarray}
Obviously $M_h^2(\gamma_1)=-\gamma_1$.
By remark \ref{t2.6} (v), $(B_1)_{\Phi_4}$ is an Orlik block with
cyclic generator $\gamma_1$, so 
$(B_1)_{\Phi_4}=\Z\cdot \gamma_1\oplus \Z\cdot M_h(\gamma_1)$.
Calculate
\begin{eqnarray}\label{5.6.12}
M_h(\gamma_1)=2e_1+e_2-e_5+e_6-2e_7-e_8-e_9-e_{10}.
\end{eqnarray}

For $Q_{2,4s}$ define 
\begin{eqnarray}
\gamma_2&:=& \frac{b_2}{\Phi_4}(M_h)(\beta_2)
=\frac{t^{6+4s}+1}{t^2+1}(M_h)(e_{11}) \nonumber \\
&=& (t^{4+4s}-t^{2+4s}+t^{4s}-\ldots-t^2+1)(M_h)(e_{11})\nonumber\\
&=& -e_5-e_{10}+(-1)\sum_{j=1}^{2+2s}e_{10+2j}
+ (-2)\sum_{j=1}^{1+s} e_{9+4j}.\label{5.6.13}
\end{eqnarray}
Obviously $M_h^2(\gamma_2)=-\gamma_2$.
By remark \ref{t2.6} (v), $(B_2)_{\Phi_4}$ is an Orlik block with
cyclic generator $\gamma_2$, so 
$(B_2)_{\Phi_4}=\Z\cdot \gamma_2\oplus \Z\cdot M_h(\gamma_2)$.
Calculate
\begin{eqnarray}\label{5.6.14}
M_h(\gamma_2)=-e_2+e_5+e_6-e_7+
\sum_{j=1}^{2+2s}(-1)^{j+1}e_{10+2j}.
\end{eqnarray}

For $Q_{2,4s}$ define 
\begin{eqnarray}\label{5.6.15}
\gamma_3:= \frac{1}{2}(\gamma_1+M_h(\gamma_1)+\gamma_2+M_h(\gamma_2))
\end{eqnarray}
and observe
\begin{eqnarray}
\gamma_3 &=&e_1
-\sum_{j\in\{3,4,5,7,8,9,10\}}e_j-
\sum_{j=1}^{1+s}(e_{9+4j}+e_{10+4j})\nonumber \\
& \stackrel{!}{\in}& Ml(f).\label{5.6.16}
\end{eqnarray}
Together with $[Ml(f):B_1\oplus B_2]=2$ this shows \eqref{5.6.9}
and that $\gamma_1,M_h(\gamma_1),\gamma_3,M_h(\gamma_3)$
is a $\Z$-basis of $Ml(f)_{\Phi_4}$.
We want to calculate the matrices of $L$ with respect to the
basis $\gamma_1,M_h(\gamma_1),\gamma_2,M_h(\gamma_2)$ of 
$(B_1\oplus B_2)_{\Phi_4}$ and the basis 
$\gamma_1,M_h(\gamma_1),\gamma_3,M_h(\gamma_3)$
of $Ml(f)_{\Phi_4}$.
Essentially we need to calculate only the values 
$L(\gamma_1,\gamma_1)$ and $L(\gamma_2,\gamma_2)$, because of \eqref{5.3}
and because of the identities for any $a\in Ml(f)_{\Phi_4}$,
\begin{eqnarray}\label{5.6.17}
\begin{array}{lll}
L(a,M_h(a))&=&L(M_h(a),M_h^2(a))
=-L(M_h(a),a)\\
&=&L(a,a)=L(M_h(a),M_h(a)).\end{array}
\end{eqnarray}
Using $M_h^2(\gamma_j)=-\gamma_j$ and 
calculations similar to \eqref{2.12}, we find
\begin{eqnarray}\label{5.6.18}
L(\gamma_1,\gamma_1)
&=& L(\frac{b_1}{\Phi_4}(-M_h^{-1})(\gamma_1),e_8)
= 3\cdot L(\-M_h(\gamma_1),e_8)=3,\hspace*{1cm}\\
L(\gamma_2,\gamma_2)
&=& L(\frac{b_2}{\Phi_4}(M_h^{-1})(\gamma_2),e_{11})\nonumber\\
&=&(3+2s)\cdot L(\gamma_2,e_{11}) =3+2s,
\label{5.6.19}
\end{eqnarray}
thus 
\begin{eqnarray}\label{5.6.20}
L(\begin{pmatrix}\gamma_1\\M_h(\gamma_1)\\ \gamma_2\\M_h(\gamma_2)\end{pmatrix},
\begin{pmatrix}\gamma_1\\M_h(\gamma_1)\\ \gamma_2\\M_h(\gamma_2)\end{pmatrix}^t)
%(\gamma_1,M_h(\gamma_1),\gamma_2,M_h(\gamma_2))
=\begin{pmatrix}3 & 3 & 0 & 0 \\ -3 & 3 & 0 & 0 \\ 
0 & 0 & 3+2s & 3+2s \\ 0 & 0 & -(3+2s) & 3+2s \end{pmatrix}
\end{eqnarray}
and 
\begin{eqnarray}\label{5.6.21}
L(\begin{pmatrix}\gamma_1\\M_h(\gamma_1)\\ 
\gamma_3\\M_h(\gamma_3)\end{pmatrix},
\begin{pmatrix}\gamma_1\\M_h(\gamma_1)\\ 
\gamma_3\\M_h(\gamma_3)\end{pmatrix}^t)
%(\gamma_1,M_h(\gamma_1),\gamma_3,M_h(\gamma_3))
=\begin{pmatrix}3 & 3 & 3 & 0 \\ -3 & 3 & 0 & 3 \\ 
0 & 3 & 3+s & 3+s \\ -3 & 0 & -(3+s) & 3+s \end{pmatrix} .
\end{eqnarray} 
The quadratic form associated to the last matrix is 
\begin{eqnarray}\label{5.6.22}
\frac{3}{2}\cdot\left[ (x_1+x_3)^2+(x_1-x_4)^2+(x_2+x_3)^2
+(x_2+x_4)^2\right]\\
+s\cdot (x_3^2+x_4^2).\nonumber
\end{eqnarray}
This shows (first for $Q_{2,4s}$, but in fact for all $Q_{2,p}$)
\begin{eqnarray}\label{5.6.23}
\{a\in Ml(f)_{\Phi_4}\, |\, L(a,a)=3\}=\{\pm\gamma_1,\pm M_h(\gamma_1)\},
\end{eqnarray}
and because of $(B_1\oplus B_2)_{\Phi_4}\subset Ml(f)_{\Phi_4}$
\begin{eqnarray}\label{5.6.24}
\{a\in (B_1\oplus B_2)_{\Phi_4}\, |\, L(a,a)=3\}=\{\pm\gamma_1,\pm M_h(\gamma_1)\},
\end{eqnarray}
This implies that any $g\in G_\Z\cup\Aut(B_1\oplus B_2,L)$ maps
the set $\{\pm\gamma_1,\pm M_h(\gamma_1)\}$ to itself and thus
$(B_1)_{\Phi_4}$ to itself and thus the $L$-orthogonal sublattice
$(B_2)_{\Phi_4}$ to itself. This shows \eqref{5.6.6} and gives (i).

Define for all $Q_{2,p}$ 
\begin{eqnarray}\label{5.6.25}
\gamma_4&:=& \gamma_1+M_h(\gamma_1)\\
&=& 2e_1+e_2-2e_3-2e_4-2e_5-e_6-e_7-2e_8-2e_9-e_{10}.\nonumber
\end{eqnarray}
Observe 
\begin{eqnarray}\label{5.6.26}
M_h(\gamma_4)&=&-\gamma_1+M_h(\gamma_1)\\
&=&-2\gamma_1+\gamma_4.\label{5.6.26n}
\end{eqnarray}
\eqref{5.6.23} and \eqref{5.6.26} imply \eqref{5.6.8}.
\eqref{5.6.26n} implies the second equality in \eqref{5.6.7}.
One calculates
\begin{eqnarray}\label{5.6.27}
L(e_8,\gamma_4)=0.
\end{eqnarray}
This shows $L(e_8,M_h(\gamma_4))\equiv 0(2)$ 
(in fact, it is $=-2$).
The $M_h$-invariance of $L$ and the fact that $e_8$ is a cyclic generator
of the Orlik block $B_1$ give 
$B_1\subset\{a\in Ml(f)\, |\, L(a,\gamma_4)\equiv 0(2)\}$.
As \eqref{5.3} implies $L(B_2,\gamma_4)=0$, so
$B_1\oplus B_2\subset\{a\in Ml(f)\, |\, L(a,\gamma_4)\equiv 0(2)\}$.
Now $r_I=2$ and for example $L(e_2,\gamma_4)=-1\not\equiv 0(2)$
show \eqref{5.6.7} and (ii).
(ii) implies $G_\Z\subset \Aut(B_1\oplus B_2,L)$.

(iii) implies $\Aut(B_1\oplus B_2,L)\subset G_\Z$, but (iii)
has still to be proved.

We continue as in the final part of the proof of part (a) for the
other series. (i) holds. Lemma \ref{t2.8} can be applied.
Therefore \eqref{5.42} and \eqref{5.43} hold for $Q_{2,p}$.
The group $\Aut(B_1\oplus B_2,L)$ for $12\not| p $ is generated by
$M_h,-\id, M_h|_{B_1}\times \id|_{B_2}$ and $(-\id)|_{B_1}\times \id|_{B_2}$, and 
analogously for the group in \eqref{5.43} if $12|p$.

For $Q_{2,4s}$ we define $\gamma_5:=\gamma_3$. 
It satisfies \eqref{5.6.9}. If $12|4s$, it is in 
$(B_1)_{b_1/\Phi_m}+(B_3)_{b_2/\Phi_m}$, so we can work
with the group in \eqref{5.43}.
If $12\not| 4s$, we work with the group in \eqref{5.42}.
In both cases $\gamma_5$ satisfies \eqref{5.6.10}, because of
\begin{eqnarray}\label{5.6.28}
(M_h|_{B_1}\times \id|_{B_2})(\gamma_5)
&=&\gamma_5-M_h(\gamma_1)\in Ml(f),\\
((-\id)|_{B_1}\times \id|_{B_2})(\gamma_5)
&=& \gamma_5-(\gamma_1+M_h(\gamma_1))\in Ml(f).\label{5.6.29}
\end{eqnarray}
For other $Q_{2,p}$, we choose a different (rather simple)
$\gamma_5$,
\begin{eqnarray}\label{5.6.30}
\gamma_5&:=& e_{10} \\
&=& \frac{1}{2}(-e_2+e_6-e_7+e_{10})
-\frac{1}{2}(-e_2+e_6-e_7-e_{10}),\nonumber\\
\textup{with}&& -e_2+e_6-e_7+e_{10}\in B_1,\quad
-e_2+e_6-e_7-e_{10}\in B_2.\nonumber
\end{eqnarray}
Then \eqref{5.6.9} holds. And
\begin{eqnarray}\label{5.6.31}
(M_h|_{B_1}\times \id|_{B_2})(\gamma_5)
&=&e_1+e_2\in Ml(f),\\
((-\id)|_{B_1}\times \id|_{B_2})(\gamma_5)
&=& e_2-e_6+e_7\in Ml(f).\label{5.6.32}
\end{eqnarray}
In any case \eqref{5.6.9} and \eqref{5.6.10} and (iii) hold.
Thus $\Aut(B_1\oplus B_2,L)\subset G_\Z$, and \eqref{5.4}
is proved for $Q_{2,p}$.

\subsection{The series $W_{1,p}$}\label{c5.7}

\begin{figure}[h]
\begin{tikzpicture}[line cap=round,line join=round,>=triangle 45,x=1.0cm,y=1.0cm]
\clip(-6,-1.2) rectangle (11.5,6.4);
\draw (-1,-1)-- (0,0);
\draw (0,3)-- (0,2);
\draw (0,1)-- (1,0);
\draw (0,1)-- (-1,-1);
\draw (-1,0)-- (0,1);
\draw (-1,0)-- (0,0);
\draw (0,0)-- (1,0);
\draw (1,0)-- (2,0);
\draw [dash pattern=on 4pt off 4pt] (0,3)-- (1,0);
\draw (0,1)-- (0,2);
\draw (3,0)-- (4,0);
\draw [dash pattern=on 4pt off 4pt] (-0.06,0)-- (0,1.01);
\draw [dash pattern=on 4pt off 4pt] (0.07,0)-- (0.07,1);
\draw (-1,0)-- (-2,0);
\draw (-2,0)-- (-3,0);
\draw (-3,0)-- (-4,0);
\draw (2,0)-- (3,0);
\draw (0,3)-- (2,0);
\draw (0,3)-- (-2,0);
\draw [dash pattern=on 4pt off 4pt] (-1,0)-- (0,3);
\draw (0,3)-- (0,4);
\draw [line width=1.2pt,dotted] (0,4)-- (0,4.49);
\draw [line width=1.2pt,dotted] (0,5.25)-- (0,4.79);
\draw (-4,0)-- (-5,0);
\draw (4,0)-- (5,0);
\begin{tiny}
\fill [color=black] (-1,-1) circle (2.0pt);
\draw[color=black] (-1.18,-0.77) node {3};
\fill [color=black] (-1,0) circle (2.0pt);
\draw[color=black] (-1.14,0.25) node {4};
\fill [color=black] (1,0) circle (2.0pt);
\draw[color=black] (1.14,0.27) node {5};
\fill [color=black] (-2,0) circle (2.0pt);
\draw[color=black] (-2.01,0.31) node {8};
\fill [color=black] (2,0) circle (2.0pt);
\draw[color=black] (2.2,0.27) node {12};
\fill [color=black] (0,1) circle (2.0pt);
\draw[color=black] (0.15,1.26) node {6};
\fill [color=black] (0,2) circle (2.0pt);
\draw[color=black] (0.15,2.27) node {1};
\fill [color=black] (0,3) circle (2.0pt);
\draw[color=black] (0.15,3.27) node {2};
\fill [color=black] (3,0) circle (2.0pt);
\draw[color=black] (3.21,0.27) node {13};
\fill [color=black] (4,0) circle (2.0pt);
\draw[color=black] (4.21,0.27) node {14};
\fill [color=black] (0,0) circle (2.0pt);
\draw[color=black] (0.13,-0.33) node {7};
\fill [color=black] (-3,0) circle (2.0pt);
\draw[color=black] (-2.86,0.27) node {9};
\fill [color=black] (-4,0) circle (2.0pt);
\draw[color=black] (-3.78,0.27) node {10};
\fill [color=black] (0,4) circle (2.0pt);
\draw[color=black] (0.22,4.26) node {16};
\fill [color=black] (0,5.25) circle (2.0pt);
\draw[color=black] (0.38,5.51) node {15+p};
\fill [color=black] (-5,0) circle (2.0pt);
\draw[color=black] (-4.79,0.27) node {11};
\fill [color=black] (5,0) circle (2.0pt);
\draw[color=black] (5.22,0.27) node {15};
\end{tiny}
\end{tikzpicture}
\caption[Figure 4.7]{The CDD of a
distinguished basis $e_1,\ldots,e_\mu$ for $W_{1,p}$
 from \cite[Tabelle 6 \& Abb. 16]{Eb81}}
\label{Fig:5.7}
\end{figure}
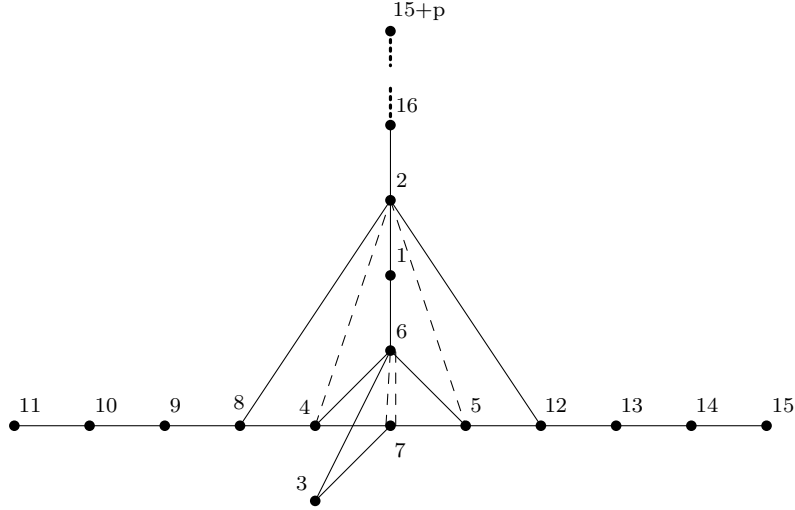

Here the monodromy acts on the distinguished basis 
$e_1,\ldots,e_\mu$ with the CDD in figure \ref{Fig:5.7}
as follows: 
\begin{eqnarray*}
e_1&\mapsto& -e_1-e_2+e_3+e_4+e_5+e_6,\\
e_2&\mapsto& 2e_1+2e_2+e_8+e_{12}+e_{16},\\
e_3&\mapsto& -e_1-e_3-e_6+e_7,\\
e_4&\mapsto& e_2-e_6+e_7+e_8,\\
e_5&\mapsto& e_2-e_6+e_7+e_{12},\\
e_6&\mapsto& e_1-2e_2+e_3+e_4+e_5+3e_6-2e_7,\\
e_7&\mapsto& -2e_2+e_3+e_4+e_5+2e_6-e_7,\\
e_8&\mapsto& e_9,\\
e_9&\mapsto& e_{10},\\
%\end{eqnarray*}
%\begin{eqnarray*}
e_{10}&\mapsto& e_{11},\\
e_{11}&\mapsto& -e_4-e_8-e_9-e_{10}-e_{11},\\
e_{12}&\mapsto& e_{13},\\
e_{13}&\mapsto& e_{14},\\
e_{14}&\mapsto& e_{15},\\
e_{15}&\mapsto& -e_5-e_{12}-e_{13}-e_{14}-e_{15},\\
e_{15+i}&\mapsto& e_{16+i}\qquad\textup{for } 1\leq i\leq p-1,\\
e_{15+p}&\mapsto& -e_1-e_2-e_{16}-e_{17}-\ldots-e_{15+p}.
\end{eqnarray*}
By table \eqref{5.13} the generators of the Orlik blocks $B_1$ and $B_2$
are $\beta_1:=e_3+e_9+e_{11}$ and $\beta_2:=e_{16}$.
The monodromy acts on them as follows:
\begin{eqnarray}
e_3+e_9+e_{11} &\mapsto&  -e_1-e_3-e_4-e_6+e_7-e_8-e_9-e_{11} \nonumber \\
&\mapsto&  e_1-e_5-e_7+e_{11}\nonumber \\
&\mapsto& -e_1-e_4-e_8-e_9-e_{10}-e_{11}-e_{12}\nonumber \\
&\mapsto& e_1-e_3-e_5-e_7-e_{13}\nonumber \\
&\mapsto& e_3+e_6-e_7-e_{12}-e_{14} \nonumber \\
&\mapsto& -e_3-e_{13}-e_{15}\nonumber \\
&\mapsto& e_1+e_3+e_5+e_6-e_7+e_{12}+e_{13}+e_{15}\nonumber \\
&\mapsto& -e_1+e_4+e_7 -e_{15}\nonumber \\
&\mapsto& e_1+e_5+e_8+e_{12}+e_{13}+e_{14}+e_{15} \label{5.7.2}\\
&\mapsto& -e_1+e_3+e_4+e_7+e_9 \nonumber\\
&\mapsto& -e_3-e_6+e_7+e_8+e_{10} \nonumber\\
&\mapsto& e_3+e_9+e_{11},\nonumber
\end{eqnarray}
\begin{eqnarray}
e_{16}&\mapsto& e_{17}\mapsto \ldots \mapsto e_{14+p} 
\mapsto e_{15+p} \nonumber \\
&\mapsto& -e_1-e_2-\sum_{i=16}^{15+p}e_i \nonumber\\
&\mapsto& -e_3-e_4-e_5-e_6-e_8-e_{12} \nonumber\\
&\mapsto& -e_4-e_5-e_7-e_8-e_9-e_{12}-e_{13}\nonumber\\
&\mapsto& -e_3-e_4-e_5-e_7-e_8-e_9-e_{10}-e_{12}-e_{13}-e_{14}\nonumber\\
&\mapsto& e_1-e_4-e_5+e_6-2e_7-\sum_{i=8}^{15}e_i \nonumber\\
&\mapsto& -e_2+e_4+e_5+2e_6-2e_7 \label{5.7.3}\\
&\mapsto& -e_{16}. \nonumber 
\end{eqnarray}
Thus the characteristic polynomial of $M_h$ on $B_j$ is $b_j$.  
Here the blocks $B_1$ and $B_2$ are
generated by the first $\deg b_1$ respectively $\deg b_2$
of the elements above. 
%and  the blocks are
%\begin{eqnarray}
%B_1&=&\langle e_3+e_9+e_{11}, -e_1-e_4-e_6+e_7-e_8,  e_1-e_5-e_7+e_{11},
%\nonumber\\
%&& -e_1-e_4-e_8-e_9-e_{10}-e_{11}-e_{12}, -e_3-e_{13}-e_{15},\nonumber\\
%&&  e_1-e_3-e_5-e_7-e_{13}, e_3+e_6-e_7-e_{12}-e_{14},\nonumber\\
%&& e_1+e_5+e_6-e_7+e_{12}, -e_1+e_4+e_7-e_{15} \rangle,\label{5.7.4}\\
%B_2&=&\langle e_{16}, \ldots, e_{15+p}; e_1+e_2, e_3+e_4+e_5+e_6+e_8+e_{12},
%\nonumber\\
%&& -e_3-e_6+e_7+e_9+e_{13}, e_3+e_4+e_{14},\nonumber\\ 
%&& -e_1-e_3-e_6+e_7+e_{11}+e_{15},\nonumber\\
%&& -e_2+e_4+e_5+2e_6-2e_7  \rangle.
%\label{5.7.5}
%\end{eqnarray}
Here $B_1+B_2=B_1\oplus B_2$ and $[Ml(f):B_1\oplus B_2]=2=r_I$
follow by the calculation of the determinant which expresses these
generators of $B_1$ and $B_2$ in the distinguished basis $e_1,\ldots,e_\mu$.
Then it also follows that $B_1$ and $B_2$ are primitive sublattices.

The proof of \eqref{5.4} for $W_{1,6s-3}$ was postponed to this subsection
and has to be given here. But the majority of the arguments was already
given in the proof of part (a). It rests to prove the following two points:

\begin{list}{}{}
\item[(i)]
\eqref{5.34} holds for $W_{1,3}$.
\item[(ii)]
In the case $W_{1,6s-3}$, any $g\in G_\Z\cup \Aut(B_1\oplus B_2,L)$
maps $(B_1)_{b_4}$ to itself and $(B_2)_{b_4}$ to itself.
Here $b_4=\Phi_6\Phi_2$.
\end{list}

For the rest of this subsection we restrict to $W_{6s-3}$.
Define for it
\begin{eqnarray}\label{5.7.4}
\delta_1 &:=& \frac{b_1}{\Phi_6\Phi_2}(M_h)(\beta_1)
=(\Phi_{12}\Phi_3)(M_h)(e_3+e_9+e_{11})\\
&=&\Phi_3(M_h)(e_9-e_{13})
= e_9+e_{10}+e_{11}-e_{13}-e_{14}-e_{15},\nonumber\\
\delta_2&:=& \frac{b_2}{\Phi_6\Phi_2}(M_h)(\beta_2)
=\frac{t^{6+p}+1}{t^3+1}(M_h)(e_{16})\label{5.7.5}\\
&=& (t^{3+p}-t^p+\ldots-t^3+1)(M_h)(e_{16})\nonumber\\
&=& e_1+e_2-\sum_{j\in\{3,4,5,7,8,9,10,12,13,14\}}e_j
+\sum_{j=1}^p e_{15+j}+\sum_{j=0}^{p/3-1}(-1)^j e_{16+3j}.\nonumber
\end{eqnarray}
$\delta_1$ and $\delta_2$ are cyclic generators of the Orlik blocks
$(B_1)_{\Phi_6\Phi_2}$ and $(B_2)_{\Phi_6\Phi_2}$, see remark \ref{t2.6} (v).
Thus $\delta_i,M_h(\delta_i)$ and $M_h^2(\delta_i)$ are a $\Z$-basis
of $(B_i)_{\Phi_6\Phi_2}$. One calculates

\begin{eqnarray}\label{5.7.6}
M_h(\delta_1)&=& -e_4-e_8-e_9+e_5+e_{12}+e_{13},\\
M_h^2(\delta_1)&=& -e_8-e_9-e_{10}+e_{12}+e_{13}+e_{14},\label{5.7.7}\\
M_h(\delta_2)&=& e_1+e_3+2e_6-2e_7%-e_9-e_{10}-e_{11}-e_{13}-e_{14}-e_{15}
-\sum_{j=9,10,11,13,14,15}e_j\nonumber\\
&+&\sum_{j=0}^{p/3-1}(-1)^j e_{17+3j},\label{5.7.8}\\
M_h^2(\delta_2)&=& -e_2+2e_4+2e_5+2e_6-e_7+e_8+e_9+e_{12}+e_{13}\nonumber\\
&+&\sum_{j=0}^{p/3-1}(-1)^je_{18+3j}.\label{5.7.9}
\end{eqnarray}

We need to calculate the $6\times 6$ matrix of values of $L$ for the
$\Z$-basis $\delta_1,M_h(\delta_1),M_h^2(\delta_1),\delta_2,
M_h(\delta_2),M_h^2(\delta_2)$ of $(B_1\oplus B_2)_{\Phi_6\Phi_2}$.
Because of \eqref{5.3}, it is block diagonal with two $3\times 3$ blocks.
Because $L$ is $M_h$-invariant and because of the identities
for any $a\in Ml(f)_{\Phi_6\Phi_2}$,
\begin{eqnarray*}
L(M_h(a),a)=-L(a,a),\quad L(M_h^2(a),a)=-L(a,M_h(a)),\\
L(a,M_h^2(a))=L(M_h(a),M_h^3(a))=-L(M_h(a),a)=L(a,a),
\end{eqnarray*}
each $3\times 3$ matrix is determined by two values. The matrices are
\begin{eqnarray}\label{5.7.10}
L(M_h^i(\delta_1),M_h^j(\delta_1))_{i,j=0,1,2} = 
\begin{pmatrix} 2 & 2 & 2 \\ -2 & 2 & 2 \\ -2 & -2 & 2 \end{pmatrix},\\
L(M_h^i(\delta_2),M_h^j(\delta_2))_{i,j=0,1,2} = 
\begin{pmatrix} 1+2s & 0 & 1+2s \\ -1-2s & 1+2s & 0 \\ 
0 & -1-2s & 1+2s \end{pmatrix}.\label{5.7.11}
\end{eqnarray}
Recall the definition $\www\gamma_2:=\frac{1}{2}(\gamma_1+\gamma_2)$
in \eqref{5.21}, and recall 
\begin{eqnarray}\label{5.7.12}
Ml(f)_{\Phi_2}=\Z\gamma_1\oplus\Z\www\gamma_2 
\stackrel{2:1}{\supset} \Z\gamma_1\oplus\Z\gamma_2=(B_1\oplus B_2)_{\Phi_2}.
\end{eqnarray}
Thus also
\begin{eqnarray}
Ml(f)_{\Phi_6\Phi_2}&=&
\langle \delta_1,M_h(\delta_1),M_h^2(\delta_1),\delta_2,M_h(\delta_2),
\www\gamma_2\rangle \nonumber\\
&\stackrel{2:1}{\supset}& (B_1\oplus B_2)_{\Phi_6\Phi_2},\label{5.7.13}
\end{eqnarray}
where
\begin{eqnarray*}
\www\gamma_2 =\frac{1}{2}(\gamma_1+\gamma_2)
= \frac{1}{2}(\delta_1-M_h(\delta_1)+M_h^2(\delta_1)
+\delta_2-M_h(\delta_2)+M_h^2(\delta_2)).
\end{eqnarray*}
The matrix of $L$ for the $\Z$-basis  
$\delta_1,M_h(\delta_1),M_h^2(\delta_1),\delta_2,M_h(\delta_2),\www\gamma_2$
of $Ml(f)_{\Phi_6\Phi_2}$ is
\begin{eqnarray}\label{5.7.14}
\begin{pmatrix} 2 & 2 & 2 & 0 & 0 & 1 \\
-2 & 2 & 2 & 0 & 0 & -1 \\
-2 & -2 & 2 & 0 & 0 & 1 \\ 
0 & 0 & 0 & 1+2s & 0 & 1+2s \\
0 & 0 & 0 & -1-2s & 1+2s & -1-2s \\ 
1 & -1 & 1 & 1+2s & -1-2s & 3+3s
\end{pmatrix}
\end{eqnarray}
The associated quadratic form 
$(x_1\ldots x_6)(\textup{matrix})\begin{pmatrix}x_1\\ \vdots\\ x_6\end{pmatrix}$
is 
\begin{eqnarray}\label{5.7.15}
&& \frac{1}{2}\left[ (2x_1+x_6)^2+(2x_2-x_6)^2+(2x_3+x_6)^2\right]\\
&+& \frac{1}{2}(1+2s)\left[(x_4-x_5+x_6)^2+(x_4+x_6)^2+(x_5-x_6)^2\right].
\nonumber
\end{eqnarray}
One finds
\begin{eqnarray}\label{5.7.16}
\{a\in Ml(f)_{\Phi_6\Phi_2}\, |\, L(a,a)=2\}
=\{\pm M_h^j(\delta_1)\, |\, j=0,1,2\},
\end{eqnarray}
and also 
\begin{eqnarray}\label{5.7.17}
\{a\in (B_1\oplus B_2)_{\Phi_6\Phi_2}\, |\, L(a,a)=2\}
=\{\pm M_h^j(\delta_1)\, |\, j=0,1,2\}.
\end{eqnarray}
Thus any $g\in G_\Z\cup \Aut(B_1\oplus B_2,L)$ 
maps $\delta_1$ to an element of
$\{\pm M_h^j(\delta_1)\, |\, j=0,1,2\}$. 
These are cyclic generators
of the Orlik block $(B_1)_{\Phi_6\Phi_2}$. 
Thus any $g\in G_\Z\cup \Aut(B_1\oplus B_2,L)$  maps $(B_1)_{\Phi_6\Phi_2}$
to itself. As $(B_2)_{\Phi_6\Phi_2}$ is the $L$-orthogonal sublattice
within $Ml(f)_{\Phi_6\Phi_2}$, such a $g$ maps also
$(B_2)_{\Phi_6\Phi_2}$ to itself.
This shows (ii) above.  Especially such a $g$ maps 
$(B_1)_{\Phi_2}$ to itself and its generator $\gamma_4=\gamma_1$ 
to $\pm\gamma_4$. This shows (i) above.

\subsection{The series $S_{1,p}$}\label{c5.8}

\begin{figure}[h]
\begin{tikzpicture}[line cap=round,line join=round,>=triangle 45,x=1.0cm,y=1.0cm]
\clip(-6.4,-2.3) rectangle (11,7);
\draw (-1,-1)-- (0,0);
\draw (0,3)-- (0,2);
\draw (0,1)-- (1,0);
\draw (0,1)-- (-1,-1);
\draw (-1,0)-- (0,1);
\draw (-1,0)-- (0,0);
\draw (0,0)-- (1,0);
\draw (1,0)-- (2,0);
\draw [dash pattern=on 4pt off 4pt] (0,3)-- (1,0);
\draw (0,1)-- (0,2);
\draw (3,0)-- (4,0);
\draw [dash pattern=on 4pt off 4pt] (-0.06,0)-- (0,1.01);
\draw [dash pattern=on 4pt off 4pt] (0.07,0)-- (0.07,1);
\draw (-1,0)-- (-2,0);
\draw (-2,-2)-- (-1,-1);
\draw (-2,0)-- (-3,0);
\draw (-3,0)-- (-4,0);
\draw (2,0)-- (3,0);
\draw (0,3)-- (2,0);
\draw (0,3)-- (-2,0);
\draw [dash pattern=on 4pt off 4pt] (-1,0)-- (0,3);
\draw (0,3)-- (0,4);
\draw [line width=1.2pt,dotted] (0,4)-- (0,4.49);
\draw [line width=1.2pt,dotted] (0,5.25)-- (0,4.79);
\begin{tiny}
\fill [color=black] (-1,-1) circle (2.0pt);
\draw[color=black] (-1.18,-0.77) node {3};
\fill [color=black] (-1,0) circle (2.0pt);
\draw[color=black] (-1.14,0.25) node {4};
\fill [color=black] (1,0) circle (2.0pt);
\draw[color=black] (1.14,0.27) node {5};
\fill [color=black] (-2,0) circle (2.0pt);
\draw[color=black] (-2.01,0.31) node {9};
\fill [color=black] (2,0) circle (2.0pt);
\draw[color=black] (2.2,0.27) node {12};
\fill [color=black] (0,1) circle (2.0pt);
\draw[color=black] (0.15,1.26) node {6};
\fill [color=black] (0,2) circle (2.0pt);
\draw[color=black] (0.15,2.27) node {1};
\fill [color=black] (0,3) circle (2.0pt);
\draw[color=black] (0.15,3.27) node {2};
\fill [color=black] (3,0) circle (2.0pt);
\draw[color=black] (3.21,0.27) node {13};
\fill [color=black] (4,0) circle (2.0pt);
\draw[color=black] (4.21,0.27) node {14};
\fill [color=black] (0,0) circle (2.0pt);
\draw[color=black] (0.13,-0.33) node {7};
\fill [color=black] (-2,-2) circle (2.0pt);
\draw[color=black] (-2.08,-1.68) node {8};
\fill [color=black] (-3,0) circle (2.0pt);
\draw[color=black] (-2.79,0.27) node {10};
\fill [color=black] (-4,0) circle (2.0pt);
\draw[color=black] (-3.78,0.27) node {11};
\fill [color=black] (0,4) circle (2.0pt);
\draw[color=black] (0.22,4.26) node {15};
\fill [color=black] (0,5.25) circle (2.0pt);
\draw[color=black] (0.38,5.51) node {14+p};
\end{tiny}
\end{tikzpicture}
\caption[Figure 4.8]{The CDD of a
distinguished basis $e_1,\ldots,e_\mu$ for $S_{1,p}$
 from \cite[Tabelle 6 \& Abb. 16]{Eb81}}
\label{Fig:5.8}
\end{figure}
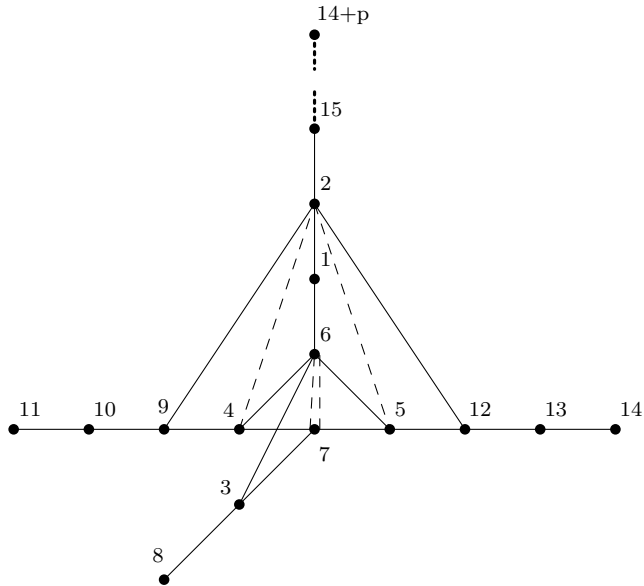

Here the monodromy acts on the distinguished basis 
$e_1,\ldots,e_\mu$ with the CDD in figure \ref{Fig:5.8}
as follows: 
\begin{eqnarray*}
e_1&\mapsto& -e_1-e_2+e_3+e_4+e_5+e_6,\\
e_2&\mapsto& 2e_{1}+2e_{2}+e_9+e_{12}+e_{15},\\
e_3&\mapsto& -e_{1}-e_6+e_7+e_8,\\
e_4&\mapsto& e_2-e_6+e_7+e_9,\\
e_5&\mapsto& e_2-e_6+e_7+e_{12},\\
e_6&\mapsto& e_1-2e_2+e_3+e_4+e_5+3e_6-2e_7,\\
e_7&\mapsto& -2e_2+e_3+e_4+e_5+2e_6-e_7,\\
e_8&\mapsto& -e_3-e_8,
\end{eqnarray*}
\begin{eqnarray*}
e_9&\mapsto& e_{10},\\
e_{10}&\mapsto& e_{11},\\
e_{11}&\mapsto& -e_4-e_9-e_{10}-e_{11},\\
e_{12}&\mapsto& e_{13},\\
e_{13}&\mapsto& e_{14},\\
e_{14}&\mapsto& -e_5-e_{12}-e_{13}-e_{14},\\
e_{14+i}&\mapsto& e_{15+i}\qquad\textup{for } 1\leq i\leq p-1,\\
e_{14+p}&\mapsto& -e_1-e_2-e_{15}-e_{16}-\ldots-e_{14+p}.
\end{eqnarray*}
By table \eqref{5.13} the generators of the Orlik blocks $B_1$ and $B_2$
are $\beta_1:=-e_8+e_{13}$ and $\beta_2:=e_{15}$.
The monodromy acts on them as follows:
\begin{eqnarray}
-e_8+e_{13} &\mapsto& e_3+e_8+e_{14} \nonumber \\
&\mapsto&  -e_1-e_3-e_5-e_6+e_7-e_{12}-e_{13}-e_{14} \nonumber \\
&\mapsto& e_{1}-e_{3}-e_{4}-e_{7}-e_{8} \nonumber \\
&\mapsto& e_{3}+e_{6}-e_{7}-e_{9}\nonumber \\
&\mapsto& e_{8}-e_{10} \nonumber \\
&\mapsto& -e_{3}-e_{8}-e_{11} \nonumber \\
&\mapsto& e_{1}+e_{3}+e_{4}+e_{6}-e_{7}+e_{9}+e_{10}+e_{11}\nonumber \\
&\mapsto& -e_{1}+e_{3}+e_{5}+e_{7}+e_{8}\nonumber \\
&\mapsto& -e_{3}-e_{6}+e_{7}+e_{12} \label{5.8.1}\\
&\mapsto& -e_8+e_{13},\nonumber
\end{eqnarray}
\begin{eqnarray}
e_{15}&\mapsto& e_{16}\mapsto \ldots \mapsto e_{14+p}
\mapsto -e_1-e_2-\sum_{i=15}^{14+p}e_i \nonumber\\
&\mapsto& -e_3-e_4-e_5-e_6-e_9-e_{12} \nonumber\\
&\mapsto& -e_3-e_4-e_5-e_7-e_8-e_9-e_{10}-e_{12}-e_{13}
\nonumber \\
&\mapsto& e_1-e_4-e_5+e_6-2e_7
-\sum_{j\in\{9,10,11,12,13,14\}}e_j\nonumber\\
&\mapsto& -e_2+e_4+e_5+2e_6-2e_7 \mapsto -e_{15}.\label{5.8.2}
\end{eqnarray}
Thus the characteristic polynomial of $M_h$ on $B_j$ is $b_j$. 
Here the blocks $B_1$ and $B_2$ are
generated by the first $\deg b_1$ respectively $\deg b_2$
of the elements above.
%\begin{eqnarray}\label{5.8.5}
%B_1&=&\langle -e_8+e_{13}, \rangle,\\
%B_2&=&\langle e_{15}, e_{16}, \ldots, e_{14+p}; e_1+e_2, %e_3+e_4+e_5+e_6+e_9+e_{12}, -e_6+e_7+e_8+e_{10}+e_{13}, e_1-e_3-e_6+e_7-e_8+e_{11}%+e_{14}, -e_2+e_4+e_5+2e_6-2e_7  \rangle.
%\label{5.5.6}
%\end{eqnarray}
Here $B_1+B_2=B_1\oplus B_2$ and $[Ml(f):B_1\oplus B_2]=2=r_I$
follow by the calculation of the determinant which expresses these
generators of $B_1$ and $B_2$ in the distinguished basis $e_1,\ldots,e_\mu$.
Then it also follows that $B_1$ and $B_2$ are primitive sublattices.

The proof of \eqref{5.5} for $S_{1,10}$ was postponed to this
section and has to be given here. From now on only $S_{1,10}$
is considered. \eqref{5.25} shows that 
$(Ml(f)_{\Phi_2},L)$ is an $A_2$-lattice with roots
$\{\pm \gamma_1,\pm\www\gamma_2,\pm(\www\gamma_2-\gamma_1)\}$.
Here $\gamma_1$ generates $(B_1)_{\Phi_2}$. 
We will show that $(B_1)_{\Phi_{10}}$ and $\pm\gamma_1$
satisfy the following special relationship:
\begin{eqnarray}
&&\left[\left((B_1)_{\Phi_{10}}+\Z\cdot a\right)_\Q\cap Ml(f) :
\left((B_1)_{\Phi_{10}}+\Z\cdot a\right)\right] \nonumber \\
&=& \left\{\begin{array}{ll}
5 & \textup{ if }a=\pm\gamma_1,\\
1 & \textup{ if }a\in\{\pm\www\gamma_2,
\pm(\www\gamma_2-\gamma_1)\}.
\end{array} \right.\label{5.8.3}
\end{eqnarray}
If $a=\pm\gamma_1$, then
\begin{eqnarray*}
\left((B_1)_{\Phi_{10}}+\Z\cdot a\right)_\Q\cap Ml(f)
&=& (B_1)_{\Phi_{10}\Phi_2} =\bigoplus_{j=0}^4
\Z\cdot (t^j\Phi_5)(M_h)(\beta_1),\\
(B_1)_{\Phi_{10}}+\Z\cdot a &=& 
(B_1)_{\Phi_{10}}+(B_1)_{\Phi_2} \\
= \bigoplus_{j=0}^3\Z\cdot (t^j\Phi_2\Phi_5)(M_h)(\beta_1)
&\oplus& \Z\cdot (\Phi_{10}\Phi_5)(M_h)(\beta_1),
\end{eqnarray*}
so the index is 
\begin{eqnarray*}
\left[\bigoplus_{j=0}^4\Z\cdot t^j \ :\ 
\bigoplus_{j=0}^3\Z\cdot t^j\Phi_2\oplus \Z\cdot \Phi_{10}
\right] = 5.
\end{eqnarray*}
Now recall that $(B_1)_{\Phi_{10}}$ is a primitive sublattice
of $Ml(f)$ and that
\begin{eqnarray*}
B_1\subset \bigoplus_{j=1}^{14}\Z\cdot e_j, \quad
\textup{ so } (B_1)_{\Phi_{10}}\subset \bigoplus_{j=1}^{14}
\Z\cdot e_j.
\end{eqnarray*}
Observe that 
\begin{eqnarray*}
\www\gamma_2 \equiv \www\gamma_2-\gamma_1
\equiv -\sum_{j=15}^{24}e_j\mod\sum_{j=1}^{14}\Z\cdot e_j.
\end{eqnarray*}
Because of the sum $-\sum_{j=15}^{24}e_j$
in $\www\gamma_2$ and in $\www\gamma_2-\gamma_1$, the 
sublattices $(B_1)_{\Phi_{10}}\oplus \Z\cdot\www\gamma_2$
and $(B_1)_{\Phi_{10}}\oplus\Z\cdot (\www\gamma_2-\gamma_1)$
are primitive in $Ml(f)$, so the index above is 1.
This shows \eqref{5.8.3}.

Now \eqref{5.5} is an easy consequence:
Consider an element $g\in G_\Z$ with $g((B_1)_{\Phi_{10}})
=(B_1)_{\Phi_{10}}$. It must map $\gamma_1$ to some root
of the $A_2$-lattice $(Ml(f)_{\Phi_{10}},L)$.
Because of \eqref{5.8.3}, the image must be $\pm\gamma_1$,
so $g((B_1)_{\Phi_2})=(B_1)_{\Phi_2}$.
Therefore $g((B_1)_{\Phi_{10}\Phi_2})=(B_1)_{\Phi_{10}\Phi_2}$
and by its $L$-orthogonality also
$g((B_2)_{\Phi_{10}\Phi_2})=(B_2)_{\Phi_{10}\Phi_2}$.

For $S_{1,10}$ $b_1=\Phi_{10}\Phi_5\Phi_2$ and
$b_2=\Phi_{30}\Phi_{10}\Phi_6\Phi_2$, so the 
eigenspaces with eigenvalues different from the roots of
$\Phi_{10}\Phi_2$ are one-dimensional and are either
in $(B_1)_\C$ or in $(B_2)_\C$. This implies \eqref{5.5}
for $S_{1,10}$.

This finishes the proof of theorem \ref{t5.1}.
\hfill $\Box$

\section{The group \texorpdfstring{$G_\Z$}{GZ} for the quadrangle
singularities}\label{c6}
\setcounter{equation}{0}

\noindent
The normal forms from \cite[\S 13]{AGV85} for the six 
families of quadrangle singularities 
will be listed below in section \ref{c10}. 
The quadrangle singularities can be seen as special 
0-th members of the eight bimodal series, with the two
series $W_{1,p}^\sharp$ and $W_{1,p}$ for $W_{1,0}$
and the two series $S_{1,p}^\sharp$ and $S_{1,p}$
for $S_{1,0}$.

The following table specializes the table \eqref{5.1}
to the case $p=0$. For $W_{1,0}$ and $S_{1,0}$, we have chosen
the specialization of the cases $W_{1,p}^\sharp$ and
$S_{1,p}^\sharp$, not $W_{1,p}$ and $S_{1,p}$.
The reason is that the Orlik blocks in theorem \ref{t5.1}
for $W_{1,p}^\sharp$ and $S_{1,p}^\sharp$ work also for
$W_{1,0}$ and $S_{1,0}$, but those for $W_{1,p}$ and $S_{1,p}$ 
work not for $W_{1,0}$ and $S_{1,0}$.
Again $b_1b_2$ respectively $b_1b_2b_3$ for $Z_{1,0}$
are the characteristic polynomials of the surface singularities.

\begin{eqnarray}\label{6.1}
\begin{array}{lllllll}
\textup{family} & \mu & b_1 & b_2 & b_3 & m & r_I\\   \hline 
W_{1,0} & 15 & \Phi_{12}& 
\Phi_{12}\Phi_6\Phi_4\Phi_3\Phi_2 & - & 12 & 1\\
S_{1,0} & 14 & \Phi_{10}\Phi_2 & 
\Phi_{10}\Phi_5\Phi_2 & - & 10 & 1\\
U_{1,0} & 14 & \Phi_9 & \Phi_9\Phi_3 & - & 9 & 1 \\
E_{3,0} & 16 & \Phi_{18}\Phi_2 & 
\Phi_{18}\Phi_6\Phi_2 & - & 18 & 2\\
Z_{1,0} & 15 & \Phi_{14}\Phi_2 & 
\Phi_{14}\Phi_2 & \Phi_2 & 14 & 2\\
Q_{2,0} & 14 & \Phi_{12}\Phi_4\Phi_3 & 
\Phi_{12}\Phi_4 & - & 12 & 2
\end{array}
\end{eqnarray}

The following theorem on the group $G_\Z$ 
has a strong similarity with the analogous theorem \ref{t5.1}
for the eight bimodal series.
And luckily, also large parts of the proof of theorem \ref{t5.1}
apply also to the case $p=0$.
We do not have \eqref{5.4} $G_\Z=\Aut(\bigoplus_{j\geq 1}B_j,L)$
for $E_{3,0},Z_{1,0},Q_{2,0}$. But we have an analogue
of the substitute \eqref{5.5} for $S_{1,10}$, the
formula \eqref{6.4}.
Contrary to theorem \ref{t5.1}, we need and give a precise
description of the induced Fuchsian group. The proof uses theorem
\ref{t3.6}. A part of the proof (a surjectivity) is 
postponed to section \ref{c10}.
For each family, denote $\zeta:=e^{2\pi i/m}\in S^1\subset\C$.

\begin{theorem}\label{t6.1}
For any surface singularity $f$ in any of the six families of
quadrangle singularities, the following holds.

(a) (See definition \ref{t2.3} for the notion {\rm Orlik block})
For all families except $Z_{1,0}$, there are Orlik blocks
$B_1,B_2\subset Ml(f)$, and for $Z_{1,0}$, there are
Orlik blocks $B_1,B_2,B_3\subset Ml(f)$ 
with the following properties.
The characteristic polynomial $p_{B_j}$ of the monodromy on $B_j$
is $b_j$. The sum $\sum_{j\geq 1}B_j$ 
is a direct sum $\bigoplus_{j\geq 1}B_j$, 
and it is a sublattice of $Ml(f)$ of full rank $\mu$ and of index $r_I$. 
Define
\begin{eqnarray}\label{6.2}
\www B_1:= \left\{\begin{array}{ll}
B_1 & \textup{for all series except }Z_{1,p},\\
B_1\oplus B_3 & \textup{for the series }Z_{1,p}.
\end{array}\right. 
\end{eqnarray}
Then 
\begin{eqnarray}\label{6.3}
L(\www B_1,B_2)&=&0=L(B_2,\www B_1),\\
g\in G_\Z\textup{ with }g((B_1)_{\Phi_{10}})=(B_1)_{\Phi_{10}}
&\Rightarrow& g(B_j)=B_j\textup{ for }j\geq 1.\label{6.4}
\end{eqnarray}

\medskip
(b) The eigenspace $Ml(f)_{\zeta}\subset Ml(f)_\C$ is 
2-dimensional.
The hermitian form $h_\zeta$ on it from lemma \ref{t2.2} (a) with 
$h_\zeta(a,b):=\sqrt{-\zeta}\cdot L(a,\oooo{b})$ for $a,b\in
Ml(f)_\zeta$ is nondegenerate and indefinite, so $\P (Ml(f)_{\zeta})\cong \P^1$ contains a half-plane 
\begin{eqnarray}\label{6.5}
\HH_\zeta:=\{\C\cdot a\, |\, a\in Ml(f)_{\zeta}\textup{ with }h_\zeta(a,a)<0\}
\subset\P(Ml(f)_\zeta).
\end{eqnarray}
Therefore the group $\Aut(Ml(f)_{\zeta},h_\zeta)/S^1\cdot\id$ is isomorphic
to $PSL(2,\R)$. The homomorphism
\begin{eqnarray}\label{6.6}
%G_{\Z,2}&:=& \{g\in G_\Z\, |\, g=\id \textup{ on }Ml(f)_{e(1/(9+p))}\},\\
\Psi:G_\Z\to \Aut(Ml(f)_{\zeta},h_\zeta)/S^1\cdot \id,
\quad g\mapsto g|_{Ml(f)_\zeta}\textup{mod}\, S^1\cdot \id,
\end{eqnarray}
is well-defined. $\Psi(G_\Z)$ 
is an infinite Fuchsian group acting on the half-plane 
$\HH_\zeta$. It is a triangle group of the same type as in
theorem \ref{t3.6}, so of the following type:
\begin{eqnarray}\label{6.7}
\begin{array}{l|l|l|l|l}
W_{1,0} & S_{1,0} & E_{3,0}\ \&\ U_{1,0} & Z_{1,0} & Q_{2,0} \\
(2,12,12) & (2,10,10) & (2,3,18) & (2,3,14) & (2,3,12) 
\end{array}
\end{eqnarray}
And 
\begin{eqnarray}\label{6.8}
\ker\Psi &=& \{\pm M_h^k\, |\, k\in\Z\} .
\end{eqnarray}
\end{theorem}

{\bf Proof:}
(a) We choose again (as in section \ref{c5})
for each of the six cases a distinguished basis with the
Coxeter-Dynkin diagram in \cite[Tabelle 6 and Abb. 16]{Eb81}.

The diagrams for $W_{1,p}^\sharp$ and $W_{1,p}$ specialize
both to the same diagram for $W_{1,0}$.
Though the description of the action of the monodromy on the
distinguished basis for $W_{1,p}^\sharp$ in \ref{c5.1} specializes
to $W_{1,0}$, but not the description for $W_{1,p}$ in \ref{c5.7}.
In the latter case $e_2\mapsto 2e_1+2e_2+e_8+e_{12}+e_{16}$,
but $e_{16}$ does not exist for $W_{1,0}$.
Therefore we work with the specialization to $p=0$ of the
formulas for $W_{1,p}^\sharp$ in subsection \ref{c5.1}.

The same applies to $S_{1,0}$. There we work with the
specialization to $p=0$ of the formulas for $S_{1,p}^\sharp$
in subsection \ref{c5.2}. 

The Orlik blocks $B_1$ and $B_2$ (and $B_3$ for $Z_{1,0}$)
are defined as in the proof of theorem \ref{t5.1},
there for $p>0$, now for $p=0$. By the same arguments,
the sum $\sum_{j\geq 1}B_j$ is a direct sum 
$\bigoplus_{j\geq 1}B_j$ and a sublattice of $Ml(f)$
of full rank $\mu$ and index $r_I$, and 
\eqref{6.3} holds. 

With  respect to part (a), it rests to show \eqref{6.4}.
In the cases $W_{1,0}$ and $U_{1,0}$, it is trivial
as $r_I=1$ and $b_1=\Phi_m$ and 
$B_1$ and $B_2$ are $L$-orthogonal.

In the cases $S_{1,0},E_{3,0},Z_{1,0}$ and $Q_{2,0}$,
the proof will be similar to the proof of \eqref{5.5}
for $S_{1,10}$ in subsection \ref{c5.8}. First we treat
$S_{1,0},E_{3,0}$ and $Z_{1,0}$ together, then we come
to $Q_{2,0}$. 

The following formulas in the proof of part (a) of theorem
\ref{t5.1} specialize to the cases $S_{1,0},E_{3,0}$
and $Z_{1,0}$: \eqref{5.10}--\eqref{5.26},
\eqref{5.28}, \eqref{5.33}, \eqref{5.35}.

The quadratic forms in \eqref{5.26} give now the following
variants of \eqref{5.27} and \eqref{5.29}:
\begin{eqnarray}\label{6.9}
\{a\in Ml(f)_{\Phi_2}\, |\, L(a,a)=5\} 
&=&\{\pm\gamma_1,\pm\gamma_2\}\quad\textup{for }S_{1,0},\\
\{a\in Ml(f)_{\Phi_2}\, |\, L(a,a)=6\} 
&=&\{\pm\gamma_1,\pm\www\gamma_2,\pm(\www\gamma_2-\gamma_1)\}
\quad\textup{for }E_{3,0},\nonumber\\ 
\{a\in Ml(f)_{\Phi_2}\, |\, L(a,a)=5\} 
&=&\{\pm(\gamma_1-3\gamma_2),\nonumber \\
&& \pm\www\gamma_2,
\pm(\www\gamma_2-\gamma_2)\}\quad\textup{for }Z_{1,0}.\nonumber
\end{eqnarray}
The first element (up to sign) of each of these three sets
generates in the corresponding case
$(B_1)_{\Phi_2}$. We claim that $(B_1)_{\Phi_m}$
and this first element satisfy the following special 
relationship. For $a$ in any of these three sets define
\begin{eqnarray}\label{6.10}
r(a):=[((B_1)_{\Phi_m}+\Z\cdot a)_\Q\cap Ml(f) :
((B_1)_{\Phi_m}+\Z\cdot a)]\in \Z_{\geq 1}.
\end{eqnarray}
Then we claim:
\begin{eqnarray}\label{6.11}
\begin{array}{ll|ll|ll|ll}
 & & S_{1,0} & & E_{3,0} & & Z_{1,0} \\ \hline 
a & r(a) & \pm\gamma_1 & 5 & \pm\gamma_1 & 3 &
\pm(\gamma_1-2\gamma_3) & 7 \\
a & r(a) & \pm\gamma_2 & 1 & 
\pm\www\gamma_2,\pm(\www\gamma_2-\gamma_1) & 1 &
\pm\www\gamma_2, \pm(\www\gamma_2-\gamma_2) & 1
\end{array}
\end{eqnarray}
The proof is the same as the proof of \eqref{5.8.3}
for $S_{1,10}$ in subsection \ref{c5.8}.
We use that for any unitary polynomial $p(t)\in\Z[t]$
\begin{eqnarray}\label{6.12}
\left[\bigoplus_{j=0}^{\deg p}\Z\cdot t^j :
\bigoplus_{j=0}^{\deg p-1}\Z\cdot t^j\Phi_2 \oplus
\Z\cdot p(t)\right] = |p(-1)|,
\end{eqnarray}
and
\begin{eqnarray}\label{6.13}
\Phi_{10}(-1)=5,\quad \Phi_{18}(-1)=3,\quad\Phi_{14}(-1)=7.
\end{eqnarray}
We also use
\begin{eqnarray}\label{6.14}
B_1\subset \sum_{j=1}^{m_1}\Z\cdot e_j\qquad
\textup{with }m_1:=8,9,10\textup{ \ for \ }S_{1,0}, E_{3,0},Z_{1,0}
\end{eqnarray}
and that the elements in the second line of \eqref{6.11} are
modulo $\sum_{j=1}^{m_1}\Z\cdot e_j$ 
\begin{eqnarray}\label{6.15}
S_{1,0}&:& \gamma_2\equiv e_9+e_{11}+e_{12}+e_{14},\\
E_{3,0}&:& \www\gamma_2\equiv e_{10}+e_{12}+e_{14}+e_{16},
\quad \www\gamma_2-\gamma_1\equiv \www\gamma_2, \nonumber\\
Z_{1,0}&:& \www\gamma_2\equiv e_{11}+e_{13}+e_{15},
\quad \www\gamma_2-\gamma_2\equiv -\www\gamma_2.\nonumber
%\begin{array}{l|l|l|}
%S_{1,0} & E_{3,0} & Z_{1,0}  \\
%\gamma_2\equiv \sum_{j\in\{9,11,12,14\}}e_j & 
%\www\gamma_2\equiv \sum_{j\in \{10,12,14,16\}}e_j &
%\www\gamma_2\equiv \sum_{j\in \{11,13,15\}}e_j \\
% & \www\gamma_2-\gamma_1\equiv \www\gamma_2 &
%\www\gamma_2-\gamma_2\equiv -\www\gamma_2
%\end{array}
\end{eqnarray}
Therefore $(B_1)_{\Phi_m}+\Z\cdot a$ for these elements $a$ 
is primitive in $Ml(f)$, and thus $r(a)=1$.

The derivation of \eqref{6.4} from \eqref{6.11} and 
\eqref{6.9} for $S_{1,0},E_{3,0}$ and $Z_{1,0}$
is almost the same as the derivation of \eqref{5.5} from 
\eqref{5.8.3} for $S_{1,10}$ in subsection \ref{c5.8}.

The only additional argument concerns $B_3=\Z\cdot \gamma_3$
in the case $Z_{1,0}$. 
Because of \eqref{5.28} any $g\in G_\Z$ maps $B_3$ to itself.
Because of $L(\gamma_1-2\gamma_3,\gamma_3)=1\neq 0$,
$B_3$ and $(B_1)_{\Phi_2}$ are glued together:
If $g=\varepsilon\cdot\id$ on $(B_1)_{\Phi_2}$ for some
$\varepsilon\in\{\pm 1\}$, then $g=\varepsilon\cdot\id$
on $B_3$.

Now we come to $Q_{2,0}$. The formulas
\eqref{5.6.2}--\eqref{5.6.5}, \eqref{5.6.7}--\eqref{5.6.8},
\eqref{5.6.11}--\eqref{5.6.22}, \eqref{5.6.25}--\eqref{5.6.27}
are also valid for $p=0$ respectively $s=0$.
The quadratic form in \eqref{5.6.22} now gives the following
variant of \eqref{5.6.23}:
\begin{eqnarray}
A := \{\gamma_1,\gamma_3,\gamma_1-\gamma_3+M_h(\gamma_3),
\gamma_1-M_h(\gamma_1)+M_h(\gamma_3)\},\label{6.16}\\
\{b\in Ml(f)_{\Phi_4}\, |\, L(b,b)=3\}
=\bigcup_{a\in A}\{\pm a,\pm M_h(a)\},\label{6.17}
\end{eqnarray}
so these are 16 elements which come in 4 sets of 4 elements
such that each set is $M_h$-invariant.
Recall that $M_h^2=-\id$ on $Ml(f)_{\Phi_4}$.
The set $\{\pm \gamma_1,\pm M_h(\gamma_1)\}$ generates
$(B_1)_{\Phi_4}$. 

We claim that $(B_1)_{\Phi_{12}}$ and this set satisfy the 
following special relationship. For $a\in A$ define
the index 
\begin{eqnarray}\label{6.18} 
r(a)&:=& \bigl[((B_1)_{\Phi_{12}}+\Z\cdot a+\Z\cdot M_h(a))_\Q
\cap Ml(f)   \\
&&   :\ ((B_1)_{\Phi_{12}}+\Z\cdot a+\Z\cdot M_h(a))\bigr]
\in\Z_{\geq 1}.\nonumber
\end{eqnarray}
Then we claim:
\begin{eqnarray}\label{6.19}
r(a) = \left\{\begin{array}{ll}
9 & \textup{ for }a=\gamma_1,\\
1 & \textup{ for }a\in\{\gamma_3,\gamma_1-M_h(\gamma_1)
+M_h(\gamma_3)\},\\
1\textup{ or }2 & \textup{ for }a=\gamma_1-\gamma_3
+M_h(\gamma_3).\end{array}\right. 
\end{eqnarray}
$r(\gamma_1)=9$ holds because of 
\begin{eqnarray}\label{6.20}
&&((B_1)_{\Phi_{12}}+\Z\cdot \gamma_1+
\Z\cdot M_h(\gamma_1))_\Q\cap Ml(f)\\
&=&(B_1)_{\Phi_{12}\Phi_4} 
= \bigoplus_{j=0}^5\Z\cdot (t^j\Phi_3)(M_h)(\beta_1),
\nonumber \\
&&(B_1)_{\Phi_{12}}+\Z\cdot \gamma_1+\Z\cdot M_h(\gamma_1)\\
&=&\bigoplus_{j=0}^3\Z\cdot (t^j\Phi_4\Phi_3)(M_h)(\beta_1)
\oplus\bigoplus_{j=0}^1\Z\cdot (t^j\Phi_{12}\Phi_3)(M_h)(\beta_1),
\nonumber
\end{eqnarray}
and  thus 
\begin{eqnarray}\label{6.22}
r(\gamma_1)=\left[\bigoplus_{j=0}^5 \Z\cdot t^j :
\bigoplus_{j=0}^3 \Z\cdot t^j\Phi_4\oplus 
\bigoplus_{j=0}^1\Z\cdot t^j\Phi_{12}
\right]=3\cdot 3.
\end{eqnarray}
For $a\in A-\{\gamma_1\}$, $r(a)\in\{1,2\}$ holds because of
\begin{eqnarray}\label{6.23}
B_1\subset \sum_{j=1}^{10}\Z\cdot e_j,
\end{eqnarray}
and because the elements $a$ and $M_h(a)$ for 
$a\in A-\{\gamma_1\}$ are modulo $\sum_{j=1}^{10}\Z\cdot e_j$
\begin{eqnarray}\label{6.24}
\gamma_1&\equiv& -e_{13}-e_{14}, \\ 
M_h(\gamma_1)&\equiv& e_{12}+e_{13}, \nonumber\\
\gamma_1-\gamma_3+M_h(\gamma_3)&\equiv& e_{12}+2e_{13}+e_{14}, 
\nonumber\\ 
M_h(\gamma_1-\gamma_3+M_h(\gamma_3))&\equiv& -e_{12}+e_{14},\nonumber\\
\gamma_1-M_h(\gamma_1)+M_h(\gamma_3) &\equiv& e_{12}+e_{13}, 
\nonumber\\ 
M_h(\gamma_1-M_h(\gamma_1)+M_h(\gamma_3))&\equiv& e_{13}+e_{14}.
\nonumber
\end{eqnarray}

The derivation of \eqref{6.4} for $Q_{2,0}$ from
\eqref{6.17} and \eqref{6.19} is a simple variant of the
derivation of \eqref{5.5} from \eqref{5.8.3} for
$S_{1,10}$ in subsection \ref{c5.8}:
Consider an element $g\in G_\Z$ with 
$g((B_1)_{\Phi_{12}}) = (B_1)_{\Phi_{12}}$. 
Because of \eqref{6.17}, it maps the set 
$\{\pm\gamma_1,\pm M_h(\gamma_1)\}$ to one of the four
sets on the right hand side of \eqref{6.17}.
Because of \eqref{6.19}, the image must be the set
$\{\pm\gamma_1,\pm M_h(\gamma_1)\}$ itself.
As this set generates $(B_1)_{\Phi_4}$, $g$ maps
$(B_1)_{\Phi_4}$ to itself.
Then $g$ maps the sets $(B_1)_{\Phi_{12}\Phi_4}$,
$B_1=(B_1)_{\Phi_{12}\Phi_4\Phi_3}$ and 
$B_2=(B_2)_{\Phi_{12}\Phi_4}$ to themselves.
This finishes the proof of part (a).

\medskip
(b) All the formulas and arguments in the proof of part (c)
of theorem \ref{t5.1} for the cases
$W_{1,12r}^\sharp, S_{1,10r}^\sharp, U_{1,9r}, E_{3,18r},
Z_{1,14r}$ and $Q_{2,12r}$ are also valid for $r=0$.

In step 3 now \eqref{6.4} is used instead of \eqref{5.4},
just as \eqref{5.5} for $S_{1,10}$. 
Therefore \eqref{6.7} holds and $\Psi(G_\Z)$ is an infinite
Fuchsian group.

By table \eqref{5.71}, the remarks \ref{t3.5} and theorem
\ref{t3.6}, $\Psi(G_\Z)$ is a subgroup of a triangle group
of the same type as in theorem \ref{t3.6}, for each case.
The proof of
theorem \ref{t10.1} will show that it is the full triangle group.
\hfill $\Box$

\section{Gauss-Manin connection and Brieskorn lattice}\label{c7}
\setcounter{equation}{0}

\noindent
The Gauss-Manin connection of isolated hypersurface 
singularities had been considered first by Brieskorn in 1970
\cite{Br70}.
Since then it had been described by many people in many
papers (K. Saito, Greuel, Pham, Varchenko, M. Saito,
Hertling, and others).
The following presentation will be short on the 
$\DD$-module foundations. But it will make the relations 
between the different pairings precise 
(more precise than anywhere in the literature).
And it will emphasize the computational aspects.
Other versions are in \cite{AGV88}, \cite{He93},
\cite{He95}, \cite{Ku98} and \cite{He02}.

Throughout most of this section, we consider 
a fixed isolated hypersurface singularity
$f:(\C^{n+1},0)\to(\C,0)$, its flat cohomology
bundle $\bigcup_{\tau\in\Delta^*}H^n(f^{-1}(\tau),\C)$,
and the space $H^\infty_\C$ of global flat multi-valued
sections (see section \ref{c4} for $H^\infty_\C$).
 
First we define the {\it elementary sections}
$es(A,\alpha)$, the spaces $C^{\alpha}$ which they generate,
and the $V$-filtration.

Any global flat multi-valued section $A\in H^\infty_\lambda$
and any choice of $\alpha\in\Q$ with $e^{-2\pi i \alpha}=\lambda$
leads to a holomorphic univalued section with specific 
growth condition at $0\in\Delta$, the 
{\it elementary section} $es(A,\alpha)$ with
\begin{eqnarray}\label{7.1}
es(A,\alpha)(\tau):=
e^{\log\tau(\alpha-\frac{N}{2\pi i})}\cdot A(\log \tau).
\end{eqnarray}
Recall that $N$ is the nilpotent part of the monodromy
$M_h$. 
Denote by $C^\alpha$ the $\C$-vector space of all
elementary sections with fixed $\alpha$ and $\lambda$.
The map 
\begin{eqnarray}\label{7.2}
\psi_\alpha:= es(.,\alpha):H^\infty_\lambda\to C^\alpha
\end{eqnarray}
is an isomorphism. The space 
$V^{mod}:=\bigoplus_{\alpha\in(-1,0]}\C\{\tau\}[\tau^{-1}]
\cdot C^\alpha$ is the space of all germs at 0 of the sheaf
of holomorphic sections on the flat cohomology bundle with 
moderate growth at 0. The {\it Kashiwara-Malgrange $V$-filtration}
is given by the subspaces
\begin{eqnarray}\label{7.3}
V^\alpha:=\bigoplus_{\beta\in[\alpha,\alpha+1)}\C\{\tau\}
\cdot C^\beta,\quad
V^{>\alpha}:=\bigoplus_{\beta\in(\alpha,\alpha+1]}\C\{\tau\}
\cdot C^\beta.
\end{eqnarray}
It is a decreasing filtration by free $\C\{\tau\}$-modules
of rank $\mu$ with 
$\Gr_V^\alpha=V^\alpha/V^{>\alpha}\cong C^\alpha$. And
\begin{eqnarray}
\tau:C^\alpha\to C^{\alpha+1}\textup{ bijective},&&
\tau\cdot es(A,\alpha)=es(A,\alpha+1),\nonumber\\
\ppp_{\tau}:C^\alpha\to C^{\alpha-1}\textup{ bijective}
&&\textup{if }\alpha\neq 0,\label{7.4}\\
\tau\ppp_\tau-\alpha:C^\alpha\to C^\alpha\textup{ nilpotent},&&
(\tau\ppp_\tau-\alpha)es(A,\alpha)=es(\frac{-N}{2\pi i}A,
\alpha).\nonumber
\end{eqnarray}
Therefore $\ppp_\tau^{-1}:V^{>-1}\to V^{>0}$ is an isomorphism,
and $V^{>-1}$ is a free $\C\{\{\ppp_\tau^{-1}\}\}$-module
of rank $\mu$. 

With the polarizing form $S$
(see \eqref{4.20}), we define a $\ppp_\tau^{-1}$-sesquilinear
pairing $K_f$ on $V^{>-1}$. Its restriction to
the Brieskorn lattice will be the restriction of
K. Saito's higher residue pairings to the Brieskorn lattice 
(which he defined on an extension of the Brieskorn lattice 
to a universal unfolding). 

\begin{lemma}\label{t7.1}
A unique pairing
\begin{eqnarray}\label{7.5}
K_f: V^{>-1}\times V^{>-1}\to \C\{\{\ppp_\tau^{-1}\}\}
\end{eqnarray}
with the properties in \eqref{7.6}--\eqref{7.9} exists.
In \eqref{7.6} and \eqref{7.7}
$A\in H^\infty_{e^{-2\pi i\alpha}},
B\in H^\infty_{e^{-2\pi i\beta}}$.
\begin{eqnarray}\label{7.6}
K_f(es(A,\alpha),es(B,\beta))&=&
\frac{1}{(2\pi i)^n}S(A,B)\cdot \ppp_\tau^{-1},\\
&&\textup{for }\alpha,\beta\in(-1,0),\alpha+\beta=-1,
\nonumber\\
K_f(es(A,\alpha),es(B,\beta))&=&
\frac{-1}{(2\pi i)^{n+1}}S(A,B)\cdot \ppp_\tau^{-2},
\label{7.7}\\
&&\textup{for }\alpha=\beta=0,
\nonumber\\
K_f:C^\alpha\times C^\beta\to 0 &&
\textup{for }\alpha,\beta\in\R_{>-1}, 
\alpha+\beta\notin\Z,\label{7.8}\\
\ppp_\tau^{-1}\cdot K_f(a,b)&=&
K_f(\ppp_\tau^{-1}a,b)=K_f(a,-\ppp_\tau^{-1}b)\label{7.9}\\
&&\textup{for }a,b\in V^{>-1}.\nonumber
\end{eqnarray}
It satisfies also (for $\alpha,\beta\in\R_{>-1}$)
\begin{eqnarray}
K_f:C^\alpha\times C^\beta\to 
\C\cdot \ppp_\tau^{-\alpha-\beta-2} && 
\textup{if }\alpha+\beta\in \Z,
\label{7.10}\\
K_f(\tau a,b)-K_f(a,\tau b)= [\tau,K_f(a,b)] &&
\textup{for } a,b\in V^{>-1},\label{7.11}
\end{eqnarray}
where $[\tau,\ppp_\tau^{-k}]=k\ppp_\tau^{-k-1}$.
If one writes $K_f(a,b)=\sum_{k\geq 1}K_f^{(-k)}(a,b)\cdot
\ppp_\tau^{-k}$ with $K_f^{(k)}(a,b)\in\C$, then
$K_f^{(-k)}$ is $(-1)^{k+n+1}$-symmetric.
\end{lemma}

{\bf Proof:}
It is clear that \eqref{7.6}--\eqref{7.9}
define a unique $\ppp_\tau^{-1}$-sesquilinear 
pairing on $V^{>-1}$. Its $\ppp_\tau^{-1}$-sesquilinearity
gives \eqref{7.10}.
One checks \eqref{7.11} with \eqref{7.4} and the infinitesimal
$N$-invariance of $S$.
The symmetry of the $K_f^{(k)}$ follows from the
symmetry of $S$ and the $\ppp_\tau^{-1}$-sesquilinearity
of $K_f$.
\hfill$\Box$ 

\begin{remark}\label{t7.2}
In the sections \ref{c9} and \ref{c10}, we will
prove the global Torelli conjecture for many families of 
marked bimodal surface singularities. We want to claim that it 
follows also for all suspensions of these families,
and also for the curve singularities, if the surface
singularities are themselves suspensions of curve 
singularities. 

The Milnor lattices of $f$ and $f+x_{n+1}^2$ are up to a sign
uniquely isomorphic. The normalized Seifert form $L^{hnor}$
and the group $G_\Z$ are the same, there is no problem.

But the Brieskorn lattices of $f$ and $f+x_{n+1}^2$
are not isomorphic.
In \cite{He93}, the second author had a lemma saying
that they are sufficiently similar and vary in the
same way in $\mu$-constant families.

Stronger and more elegant is the specialization
to $f+x_{n+1}^2$ of a Thom-Sebastiani formula.
But that requires to look at a Fourier-Laplace
transformation. In the present situation
of sections of moderate growth, this can be done in
a nice and explicit way. Lemma \ref{t7.3}, definition \ref{t7.4}
and theorem \ref{t7.5} do a good part of the work.
Theorem \ref{t7.9} gives a Thom-Sebastiani formula
for a Fourier-Laplace transform of the Brieskorn lattice. 
Theorem \ref{t7.7} states well known properties of the 
Brieskorn lattice.
\end{remark}

The pairing in lemma \ref{t7.3} had been considered
first by Pham \cite{Ph85}, see remark \ref{t7.6} (i).

\begin{lemma}\label{t7.3}
Let $\gamma_{-\pi}:H^n(f^{-1}(z),\C)\to H^n(f^{-1}(-z),\C)$
(respectively $\gamma_\pi$)
be the isomorphism by flat shift in mathematically negative
(respectively positive) direction. Define a pairing
\begin{eqnarray}\label{7.12}
P:H^n(f^{-1}(z),\C)\times H^n(f^{-1}(-z),\C)\to\C
\quad\textup{for }z\neq 0\\
\textup{by}\qquad 
P(a,b):= \frac{1}{(2\pi i)^{n+1}}\cdot L^{nor}(a,
\gamma_{-\pi}(b)).\nonumber
\end{eqnarray}
It is $(-1)^{n+1}$-symmetric and nondegenerate
and takes values in $(2\pi i)^{-(n+1)}\cdot\Z$ on
$H^n(f^{-1}(z),\Z)\times H^n(f^{-1}(-z),\Z)$. It is flat, 
i.e. it has constant values on pairs of flat sections
in the cohomology bundle.
\end{lemma}

{\bf Proof:}
The only property which might not be immediately
obvious, is the $(-1)^{n+1}$-symmetry. It follows
from the flatness, from $M_h\gamma_{-\pi}=\gamma_\pi$ and
\eqref{4.15}:
\begin{eqnarray}
(2\pi i)^{n+1}\cdot P(b,a)
&=&L^{nor}(b,\gamma_{-\pi}a)=(-1)^{n+1}L^{nor}(M_h\gamma_{-\pi}a,b) 
\nonumber\\
&=& (-1)^{n+1}L^{nor}(\gamma_\pi a,b)  
= (-1)^{n+1}L^{nor}(a,\gamma_{-\pi}b)\nonumber \\
&=&(2\pi i)^{n+1}\cdot (-1)^{n+1}\cdot P(a,b).\label{7.13}
\end{eqnarray}
\hfill$\Box$

\begin{definition}\label{t7.4} \cite[(7.47)]{He02}
For each $\alpha\in\R$ define the automorphism
\begin{eqnarray}
G^{(\alpha)}&:& H^\infty_{e^{-2\pi i\alpha}}
\to H^\infty_{e^{-2\pi i\alpha}},\nonumber\\
G^{(\alpha)}&:=& \sum_{k\geq 0}\frac{1}{k!}\Gamma^{(k)}(\alpha)
\cdot \left(\frac{-N}{2\pi i}\right)^k 
={}''\Gamma\left(\alpha\cdot\id+
\frac{-N}{2\pi i}\right){}'' .\label{7.14}
\end{eqnarray}
Define the automorphism
\begin{eqnarray}\label{7.15}
G:=\sum_{\alpha\in (-1,0]}G^{(\alpha)}:H^\infty_\C\to
H^\infty_\C.
\end{eqnarray}
\end{definition}

The following theorem was first formulated in
\cite[Proposition 7.7]{He03}.
A detailed proof is in \cite[Theorem 5.2]{BH17}.
The most difficult part is the proof of \eqref{7.21}.

\begin{theorem}\label{t7.5}
(a) Let $\tau$ and $z$ both be coordinates on $\C$. For 
$\alpha>0$ and $A\in H^\infty_{e^{-2\pi i\alpha}}$,
the Fourier-Laplace transformation $FL$ with
\begin{eqnarray}\label{7.16}
FL(es(A,\alpha-1)(\tau))(z):=
\int_0^{\infty\cdot z}e^{-\tau/z}\cdot es(A,\alpha-1)(\tau)d\tau
\end{eqnarray}
is well defined and maps the elementary section 
$es(A,\alpha-1)(\tau)$ in $\tau$ to the elementary section
\begin{eqnarray}\label{7.17}
FL(es(A,\alpha-1)(\tau))(z)=es(G^{(\alpha)}A,\alpha)(z)
\end{eqnarray}
in $z$.

\medskip
(b) It extends to a well defined isomorphism
\begin{eqnarray}\label{7.18}
FL:\sum_{\alpha\in(-1,0]}\C\{\ppp_\tau^{-1}\}\cdot C^\alpha_\tau 
\to V^{>0}_z .
\end{eqnarray}
Here the indices $\tau$ at $C^\alpha$ and $z$ at 
$V^{>0}$ indicate that the coordinate $\tau$ respectively
$z$ has to be used.
It satisfies for $a,b\in \sum_{\alpha\in(-1,0]}
\C\{\ppp_\tau^{-1}\}\cdot C^\alpha_\tau$
\begin{eqnarray}\label{7.19}
FL(\ppp_\tau^{-1}a)&=&z\cdot FL(a),\\
FL(\tau\cdot a)&=&z^2\ppp_z FL(a),\label{7.20}\\
P(FL(a),FL(b))&=&\sum_{k\geq 1}c_kz^l\quad\textup{if }
K_f(a,b)=\sum_{k\geq 1}c_k\ppp_\tau^{-k}.\label{7.21}
\end{eqnarray}
\end{theorem}

\begin{remarks}\label{t7.6}
(i) Pham \cite{Ph85} defined the pairing $P$ in lemma \ref{t7.3}
starting with an intersection form for Lefschetz thimbles.
In our situation, $H_n(f^{-1}(z),\Z)$ for $z\in\Delta^*$
is canonically isomorphic to the $\Z$-module generated 
by Lefschetz thimbles above the straight path from $0$ to $z$. 
And it is easy to see that the pairing
\begin{eqnarray}\label{7.22}
(-1)^{n(n+1)/2}\cdot L^{hnor}(.,\gamma_{-\pi}):\\
H_n(f^{-1}(z),\Z)\times H_n(f^{-1}(-z),\Z)\to\Z \nonumber
\end{eqnarray}
for $z\in\Delta^*$ 
is the intersection form for Lefschetz thimbles \cite{He05}.
This formula connects lemma \ref{t7.3} with
Pham's definition.

\medskip
(ii) Neither Pham nor K. Saito knew the definition
\ref{t7.2} of $K_f$ with the polarizing form $S$.
Pham had the version of \eqref{7.21} with 
K. Saito's higher residue pairings \cite{SaK83} instead of $K_f$.
He did not consider explicitly the automorphisms
$G^{(\alpha)}$ and \eqref{7.17}.

\medskip
(iii) Because of \eqref{7.19},
we have to consider on the left hand side of \eqref{7.18}
and in \eqref{7.19}--\eqref{7.21}
the subspace $\sum_{\alpha\in(-1,0]}
\C\{\ppp_\tau^{-1}\}\cdot C^\alpha_\tau$ of $V^{>-1}_\tau$.
The convergence condition is stronger.
\end{remarks}

Now we come to the Brieskorn lattice.
It is a free $\C\{\tau\}$-module $H_0''(f)\subset V^{>-1}$ 
of rank $\mu$ which had first been studied by Brieskorn
\cite{Br70}. The name {\it Brieskorn lattice}
is due to \cite{SaM89}, the notation $H_0''(f)$ is from
\cite{Br70}. 
The Brieskorn lattice is generated by germs of sections
$s[\omega]$ from holomorphic $(n+1)$-forms 
$\omega\in\Omega^{n+1}_X$: Integrating the 
Gelfand-Leray form $\frac{\omega}{df}|_{f^{-1}(\tau)}$
over cycles in $H_n(f^{-1}(\tau),\C)$ gives a holomorphic
section $s[\omega]$ in the cohomology bundle, whose germ
$s[\omega]_0$ at 0 is in fact in $V^{>-1}$
(this was proved first by Malgrange).
The following theorem collects well known properties of the
Brieskorn lattice. Afterwards we make comments on their proofs.
See also \cite{He02}.

\begin{theorem}\label{t7.7}
Algebraic properties:
\begin{eqnarray}\label{7.23}
H_0''(f)&\cong&\Omega^{n+1}_{\C^{n+1},0}/df\land 
d\Omega^{n-1}_{\C^{n+1},0},\\
\ppp_\tau^{-1}:H_0''(f)&\stackrel{\cong}{\longrightarrow} & 
H_0'(f)\subset H_0''(f)\nonumber\\
\textup{with }H_0'(f) &\cong & 
df\land \Omega^{n}_{\C^{n+1},0}/df\land 
d\Omega^{n-1}_{\C^{n+1},0},\nonumber\\
\textup{and }\ppp_\tau:s[df\land\eta]_0&\mapsto& 
s[d\eta]_0\label{7.24}.
\end{eqnarray}
Compatibility with $K_f$: $K_f$  
is the restriction to $H_0''(f)$ of 
K. Saito's higher residue pairings.
It satisfies
\begin{eqnarray}\label{7.25}
K_f: H_0''(f)\times H_0''(f)\to \ppp_\tau^{-n-1}\cdot
\C\{\{\ppp_\tau^{-1}\}\}.
\end{eqnarray}
The leading part 
\begin{eqnarray}\label{7.26}
K_f^{(-n-1)}:H_0''(f)/H_0'(f)\times
H_0''(f)/H_0'(f)\to\C
\end{eqnarray}
is symmetric (lemma \ref{t7.2}) and nondegenerate.
It is Grothendieck's residue pairing on
$\Omega^{n+1}_{\C^{n+1},0}/df\land \Omega^n_{\C^{n+1},0}$.

Relation to Steenbrink's Hodge filtration
$F^\bullet H^\infty_\C$: For $\lambda=e^{-2\pi i\alpha}$
with $\alpha\in (-1,0]$,
\begin{eqnarray}\label{7.27}
F^p_{St}H^\infty_\lambda&=&
\psi_\alpha^{-1}
\Bigl(\ppp_\tau^{n-p}\Gr_V^{n-p+\alpha}H_0''(f)\Bigr).
\end{eqnarray}
Define the unordered tuple 
$\Sp(f)=\sum_{i=1}^\mu(\alpha_i)=\sum_{\alpha\in\Q}
d(\alpha)\cdot (\alpha)\in \Z_{\geq 0}[\Q]$ 
of spectral numbers $\alpha_1,\ldots,\alpha_\mu\in\Q$ by  
\begin{eqnarray}\label{7.28}
d(\alpha)&:=&\dim \Gr_V^\alpha H_0''-\dim\Gr_V^\alpha H_0'.
\end{eqnarray}
Number them such that $\alpha_1\leq \ldots \leq \alpha_\mu$.
Then they satisfy the symmetry
\begin{eqnarray}\label{7.29}
\alpha_{i}+\alpha_{\mu+1-i}=n-1
\end{eqnarray}
and
\begin{eqnarray}\label{7.30}
-1<\alpha_1\leq \ldots\leq \alpha_\mu<n,\\
V^{>-1}\supset H_0''\supset V^{n-1},\nonumber\\
0=F^{n+1}H^\infty,\ F^0H^\infty_{\neq 1}=H^\infty_{\neq 1},
\ F^1H^\infty_1=H^\infty_1.\nonumber
\end{eqnarray}
\end{theorem}

The algebraic properties had been proved by Brieskorn
\cite{Br70} with some help by Sebastiani.
That $K_f$ is the restriction to $H_0''(f)$ of K. Saito's
higher residue pairings \cite{SaK83} follows from 
\eqref{7.21} and Pham's identification 
of $P$ with the Fourier-Laplace transform of 
K. Saito's higher residue pairings \cite{Ph85}.
See \cite{He02} for an alternative reasoning.
Then \eqref{7.25} and the properties of \eqref{7.26}
follow from K. Saito's work.

Steenbrink defined the Hodge filtration $F^\bullet_{St}$ 
first using resolution of singularities \cite{St77}.
Then Varchenko \cite{Va80-1} constructed a closely related
Hodge filtration $F^\bullet_{Va}$ from the Brieskorn lattice
$H_0''(f)$. Scherk and Steenbrink
\cite{SS85} (and also M. Saito) modified this 
construction to recover $F^\bullet_{St}$. This is \eqref{7.27}.
Then \eqref{7.29} and \eqref{7.30} follow from properties
of the Hodge filtration. Though $V^{>-1}\supset H_0''$
was proved before by Malgrange.

\begin{remark}\label{t7.8}
The Fourier-Laplace transformation $FL$ is defined
on any sum of elementary sections with the stronger
convergence condition in \eqref{7.17}. Therefore it is not
defined on arbitrary elements of $H_0''$. But 
because of \eqref{7.30},
\begin{eqnarray}\label{7.31}
H_0''=(H_0''\cap \bigoplus_{-1<\alpha<n-1}C^\alpha_\tau)
\oplus V^{n-1}_\tau,
\end{eqnarray}
and the elements of the first summand are finite sums of
elementary sections. Therefore the space
\begin{eqnarray}\label{7.32}
FL(H_0''\cap \bigoplus_{-1<\alpha<n-1}C^\alpha_\tau)\oplus
V^n_z
\end{eqnarray}
is a well-defined free $\C\{z\}$-module of rank $\mu$.
For simplicity we call it $FL(H_0'')$, although that
is not completely correct. It satisfies 
\begin{eqnarray}\label{7.33}
z^2\ppp_z: FL(H_0'')&\to& FL(H_0''),\\
\textup{and }
P:FL(H_0'')\times FL(H_0'')&\to& z^{n+1}\cdot\C\{z\},\label{7.34}
\end{eqnarray}
and the leading part of $P$ is a symmetric and nondegenerate
pairing on $FL(H_0'')/z\cdot FL(H_0'')$, all of this
because of \eqref{7.19}--\eqref{7.21}, 
\eqref{7.24}--\eqref{7.26}.
It thus satisfies all properties of a TERP-structure
\cite[definition 2.12]{He02}.
Because of the $\Z$-lattice $H^\infty_\Z$
and the $\Z$-lattice bundle in the cohomology,
we can even call it a TEZP-structure.
More precisely, we denote as TEZP structure the following
tuple.
\begin{eqnarray}\label{7.35}
TEZP(f):=(H^\infty_\Z,L^{nor},V^{mod}_z,P,FL(H_0''))(f).
\end{eqnarray}
Here $V^{mod}_z$ comes equipped with the actions of 
$z,\ppp_z^{-1}$ and $z\ppp_z$.
We formulated theorem \ref{t7.5} and introduced
$FL(H_0'')$ because of the following
Thom-Sebastiani result.
\end{remark}

\begin{theorem}\cite{SS85}\cite[Theorem 6.4]{BH17}\label{t7.9}
Consider besides $f(x_0,\ldots,x_n)$ a second
singularity $g(x_{n+1},\ldots,x_{n+m+1})$. Then
\begin{eqnarray}\label{7.36}
TEZP(f+g)\cong TEZP(f)\otimes TEZP(g).
\end{eqnarray}
\end{theorem}

\begin{remarks}\label{t7.10}
(i) The isomorphism for the data
$(H^\infty_\Z,L^{nor})$ is the classical Thom-Sebastiani result
in \eqref{4.8} and \eqref{4.10}.
The isomorphism for $P$ follows from its definition
with $L^{nor}$. The isomorphism for $V^{mod}_z$ is trivial.
The isomorphism for $H_0''$ was essentially proved in 
\cite[(8.7) Lemma]{SS85}.
Though Scherk and Steenbrink did not make the compatibility with
the topological Thom-Sebastiani isomorphism
between the cohomology bundles precise,
and they avoided the use of the Fourier-Laplace transformation.
They obtained a $\ppp_\tau^{-1}$-linear isomorphism
$H_0''(f+g)\cong H_0''(f)\otimes H_0''(g)$.

\medskip
(ii) This is still fine. But then they mixed
$\ppp_\tau^{-1}$-linearity and $\tau$-linearity and
went with this isomorphism directly into the defining
formula \eqref{7.27} of $F^{\bullet}_{St}$.
This lead them to a wrong Thom-Sebastiani formula
for $F^\bullet_{St}$ in \cite[Theorems (8.2) and (8.11)]{SS85}.
But the true Thom-Sebastiani formula is quite close
\cite[Corollary 6.5]{BH17}.
One has to replace in \cite[Theorems (8.2) and (8.11)]{SS85}
$F^\bullet_{St}$ by $G(F^{\bullet}_{St})$.
This follows immediately from \eqref{7.27}
and \eqref{7.36}. 
Of course, in the case $N=0$, the isomorphism
$G$ in definition \ref{t7.4} is just a rescaling,
and then $G(F^\bullet_{St})=F^\bullet_{St}$,
so then their Thom-Sebastiani formula is correct.

\medskip
(iii) As a corollary of theorem \ref{t7.9}, we obtain
for a suspension of $f$
\begin{eqnarray}\label{7.40}
TEZP(f+x_{n+1}^2)\cong TEZP(f)\otimes TEZP(x_{n+1}^2).
\end{eqnarray}
This allows us to consider in the sections \ref{c9}
and \ref{c10} only the surface singularities.
More generally, it implies the corollary \ref{t8.14}.
This corollary is the reason why we introduced
$FL(H_0''(f))$. Formula \eqref{7.36} and this
corollary are more elegant and general than the arguments
with which suspensions were treated in \cite{He93},
\cite{He95}, \cite{He11} and \cite{GH17}.

\medskip
(iv) The Thom-Sebastiani formula for $F^\bullet_{St}$
expresses in the case of a suspension 
$F^\bullet_{St}(f+x_{n+1}^2)$ in terms of $F^\bullet_{St}(f)$.
It is made explicit in \cite[Theorem 4.6]{BH17}.
It can be seen as a square root of a Tate twist, because
$F^\bullet_{St}(f)$ and 
$F^\bullet_{St}(f+x_{n+1}^2+x_{n+2}^2)$ are simply 
related by a Tate twist.
$f$ and $f+x_{n+1}^2+x_{n+2}^2$ have the same polarizing form
$S$ by \eqref{4.20} and \eqref{4.21}, because
$M_h(f)=M_h(f+x_{n+1}^2+x_{n+2}^2)$. But the polarizing form
of $f+x_{n+1}^2$ is quite different, because of 
$M_h(f+x_{n+1}^2)=-M_h(f)$ and \eqref{4.20} and \eqref{4.21}.
The formula in \cite[Theorem 4.6]{BH17} which expresses
$F^\bullet_{St}(f+x_{n+1}^2)$ in terms of $F^\bullet_{St}(f)$
involves the $G^{(\alpha)}$ from definition \ref{t7.4}
and is compatible with the isotropy condition \eqref{4.18}
and (the generalization in the case $N\neq 0$ of)
the positivity condition \eqref{4.19}.
\end{remarks}

Fix for a moment a reference singularity $f_0$. 
In \cite{He99} a classifying space $D_{PMHS}(f_0)$
and a classifying space $D_{BL}(f_0)$ are constructed.
$D_{PMHS}$ is a classifying space for $M_s$-invariant
Hodge filtrations $F^\bullet$ on $H^\infty_\C(f_0)$
such that $(H^\infty_{\neq 1},H^\infty_{\neq 1,\Z},F^\bullet,
W,-N,S)$ and $(H^\infty_1,H^\infty_{1,\Z},F^\bullet,W,-N,S)$
are polarized mixed Hodge structures of weight $n$ 
respectively $n+1$ with the same Hodge numbers as
$F^\bullet_{St}(f_0)$. 

And $D_{BL}$ is a classifying space for subspaces
$\LL_0\subset V^{>-1}_\tau$ with the following properties:

\begin{list}{}{}
\item[($\alpha$)] 
$\LL_0$ is a free $\C\{\tau\}$-module of rank $\mu$.
\item[($\beta$)]
$\LL_0$ is a free $\C\{\{\ppp_\tau^{-1}\}\}$-module
of rank $\mu$. 
\item[($\gamma$)]
The filtration $F^\bullet$ in $H^\infty_\C(f_0)$
which is constructed by formula \eqref{7.27} with
$\LL_0$ instead of $H_0''(f_0)$ is in $D_{PMHS}$.
\item[($\delta$)]
It satisfies $K_f(\LL_0,\LL_0)\subset \ppp_\tau^{-n-1}\cdot
\C\{\{\ppp_\tau^{-1}\}\}$. 
\end{list}

\begin{theorem}\label{t7.11}
Fix a reference singularity $f(x_0,\ldots,x_n)$.

(a) \cite[ch. 2]{He99} $D_{PMHS}(f_0)$ is a real homogeneous
space and a complex manifold.
It is a locally trivial bundle over a product $D_{PHS}$
of classifying spaces for pure polarized Hodge structures.
The fibers carry an affine algebraic structure and are isomorphic
to $C^{N_{PMHS}}$ for some $N_{PMHS}\in\Z_{\geq 0}$.
The group $G_\Z(f_0)$ acts properly discontinuously 
on $D_{PMHS}$.

\medskip
(b) \cite[ch. 5]{He99} $D_{BL}(f_0)$ is a complex
manifold and a locally trivial bundle over $D_{PMHS}$.
The fibers have a natural $\C^*$-action with negative weights
and are affine algebraic manifolds and are isomorphic
to $\C^{N_{BL}}$ for some $N_{BL}\in \Z_{\geq 0}$.
The group $G_\Z(f_0)$ acts properly discontinuously on $D_{BL}$.

\medskip
(c) $D_{PMHS}(f_0)$ and $D_{PMHS}(f_0+x_{n+1}^2)$ are 
canonically isomorphic.
$D_{BL}(f_0)$ and $D_{BL}(f_0+x_{n+1}^2)$ are 
canonically isomorphic.
\end{theorem}

Part (c) is not formulated in \cite{He99}.
The isomorphism $D_{BL}(f_0)\to D_{BL}(f_0+x_{n+1}^2)$
is given by the generalization of \eqref{7.36}, namely the map
\begin{eqnarray}\label{7.41}
\LL_0\mapsto FL^{-1}\Bigl(FL(\LL_0)\otimes 
FL(H_0''(x_{n+1}^2))\Bigr).
\end{eqnarray}
The isomorphism $D_{PMHS}(f_0)\to D_{PMHS}(f_0+x_{n+1}^2)$
is obtained by applying $\Gr_V^\bullet$.
It follows also from \cite[Theorem 4.6]{BH17}.

\bigskip
In the sections \ref{c9} and \ref{c10},
$\mu$-constant families of singularities in two parameters
will be studied. The following definition and theorem
treat a more general situation. It had been considered
especially in \cite{Va80-2} \cite{AGV88} \cite{SaM91} 
\cite{He93} \cite{Ku98}.

\begin{definition}\label{t7.12}
A holomorphic $\mu$-constant family of singularities
consists of a number $\mu\in\Z_{\geq 1}$, 
a complex manifold $T$, an open neighborhood $X\subset
\C^{n+1}\times T$ of $\{0\}\times T$ and a holomorphic
function $F:X\to \C$ such that 
$F_t:=F|_{X_t}$ with $X_t:=X\cap\C^{n+1}\times\{t\}$
for any $t\in T$ has an isolated singularity at 0 with 
Milnor number $\mu$.
\end{definition}

\begin{theorem}\label{t7.13}
Consider a holomorphic $\mu$-constant family as in 
definition \ref{t7.12}.

\medskip
(a) The Milnor lattices $(Ml(F_t),L)$ with Seifert forms
for $t\in T$ are locally canonically isomorphic.
They glue to a local system $\bigcup_{t\in T}Ml(F_t)$
of free $\Z$-modules of rank $\mu$.

\medskip
(b) Therefore also the spaces $C^\alpha(F_t)$,
$V^{mod}_\tau(F_t)$, $V^\alpha_\tau(F_t)$ are locally
canonically isomorphic and glue to local systems.

\medskip
(c) But the Brieskorn lattices 
$H_0''(F_t)\subset V^{>-1}_\tau(F_t)$
vary holomorphically. For $\omega\in \Omega^{n+1}_{X/T}$,
$s[\omega]_0(t):=s[\omega|_{X_t}]_0\in H_0''(F_t)$.
Let $\xi$ be a holomorphic vector field on $T$.
Its canonical lifts to $\C\times T$ (with coordinate $\tau$
on $\C$) and $X$ are also denoted $\xi$.
The covariant derivative of $s[\omega]_0(t)$ by $\xi $ is
\begin{eqnarray}\label{7.42}
\xi\ s[\omega]_0(t)=
s[\Lie_\xi\omega]_0(t) +(-\ppp_\tau)s[\xi(F)\cdot\omega]_0(t).
\end{eqnarray}

(d) All germs $F_t$ have the same spectrum.
\end{theorem}

\begin{remarks}\label{t7.14}
(i) Part (a) is less trivial than one might expect,
as it is not clear whether $\varepsilon(t)$ and $\delta(t)$
in the definition of a Milnor fibration $F_t:X(\varepsilon(t),
\delta(t))\to \Delta_{\delta(t)}$ can be chosen as 
continuous functions in $t$. But lemma \ref{t2.2}
in \cite{LR73} saves the situation.
See \cite{Va80-2} \cite{He93} \cite{Ku98} \cite{He11}
for details.

\medskip 
(ii) Part (b) follows from part (a).
Formula \eqref{7.42} is well known, see e.g.
\cite{Va80-2} \cite{AGV88} \cite{He93} \cite{Ku98}.
Part (d) is proved in \cite{Va82}.

\medskip 
(iii) The bundle $\bigcup_{t\in T}H_0''(F_t)\subset
\bigcup_{t\in T}V^{>-1}_\tau(F_t)$ can be seen as a 
germ along $\{0\}\times T$ on $(\C,0)\times T$
of a holomorphic rank $\mu$ bundle.

$s[\omega]_0$ for $\omega\in \Omega^{n+1}_{X/T}$
is a holomorphic section in this bundle.

But in theorem \ref{t9.6} and theorem \ref{10.6} we will be
imprecise and consider $s[\omega]_0$ as a possibly
multi-valued holomorphic map 
$s[\omega]_0:T\to V^{>-1}_\tau(F_{t^0})$
for a reference singularity $F_{t^0}$.

\medskip 
(iv) $s[\omega]_0$ is a sum 
$s[\omega]_0=\sum_{\alpha>-1} s(\omega,\alpha)$
of holomorphic families $s(\omega,\alpha)(t)\in C^{\alpha}(F_t)$,
$t\in T$, of elementary sections. For each $t\in T$,
\begin{eqnarray}\label{7.43}
\alpha(s[\omega]_0(t)):=
\alpha(\omega|_{X_t}):=\min(\alpha\, |\, s(\omega,\alpha)(t)
\neq 0)
\end{eqnarray}
is the {\it order} of $s[\omega]_0(t)$,
and $s(\omega,\alpha(\omega|_{X_t}))(t)$ is its 
{\it principal part}.
The order is upper semicontinuous in $t$.

\medskip
(v) A notation: $\omega_0:=dx_0\ldots dx_n$.
\end{remarks}

All bimodal series singularities in table \eqref{9.1}
except $W_{1,p}^\sharp$ (see remark \ref{t9.5} 
for $W_{1,p}^\sharp$) are Newton nondegenerate.
All quadrangle singularities in table \eqref{10.1}
are semiquasihomogeneous.
For such singularities there are useful results for the 
computation of the order $\alpha(\omega|_{X_t})$,
which we describe in the following.
We start with a definition of Kouchnirenko.

\begin{definition}\label{t7.15}
Let $f:(\C^{n+1},0)\to(\C,0)$ be a singularity.

(a) \cite{Ko76} Write $f=\sum_{i\in \Z_{\geq 0}^{n+1}}a_ix^i$
and define 
\begin{eqnarray}\label{7.44}
\supp(f)&:=&\{i\in\Z_{\geq 0}^{n+1}\, |\, a_i\neq 0\},\\
\Gamma_+(f)&:=& \Bigl(\textup{convex hull of }
\bigcup_{i\in\supp(f)}(i+\R_{\geq 0}^{n+1})\Bigr)\subset\R^{n+1},
\nonumber\\
\Gamma_{com}(f)&:=& \{\sigma\, |\, \sigma\textup{ is a compact
face of }\Gamma_+(f)\},\nonumber\\
\Gamma_{com,n}(f)&:=& \{\sigma\in \Gamma_{com}(f)\,  |\, 
\dim\sigma=n\},\nonumber\\
l_\sigma&:&\R^{n+1}\to\R\quad 
\textup{ for }\sigma\in\Gamma_{com,n}(f)
\nonumber
\end{eqnarray}
as the linear function with $\sigma\subset l_\sigma^{-1}(1)$.

\medskip
(b) \cite{SaM88}\cite{KV85}
The Newton order $\nu:\C\{x_0,\ldots,x_n\}\to \Q_{\geq 0}\cup
\{\infty\}$ is
\begin{eqnarray}\label{7.45}
\nu(\sum_ib_ix^i):=\min(l_\sigma(i)\, |\, \textup{all }
i\textup{ with }b_i\neq 0,
\textup{ all }\sigma\in \Gamma_{com,n}(f)\}.
\end{eqnarray}
The Newton order $\nu:\Omega^{n+1}_{\C^{n+1},0}\to
\Q_{>0}\cup\{\infty\}$ is 
\begin{eqnarray}\label{7.46}
\nu((\sum_ib_ix^i)\cdot \omega_0):= \nu((\sum_ib_ix^i)x_0\ldots x_n).
\end{eqnarray}
The Newton order $\oooo\nu:H_0''(f)\to \Q_{>0}\cup\{\infty\}$ is
\begin{eqnarray}\label{7.47}
\oooo{\nu}:=\max(\nu(\eta)\, |\, \eta\equiv\omega\mod
df\land d\Omega^{n-1}_{\C^{n+1},0}).
\end{eqnarray}

(c) \cite{Ko76} For $\sigma\in\Gamma_{com}(f)$ define
$f_\sigma:=\sum_{i\in\sigma}a_ix^i$. 
The singularity $f$ is {\it Newton nondegenerate} if 
for each $\sigma\in \Gamma_{com}(f)$ 
the Jacobi ideal $J(f_\sigma)$ of $f_\sigma$ has no zero
in $(\C^*)^{n+1}$. It is {\it convenient} if $f$ contains
for each index $j\in\{0,\ldots,n\}$ a monomial
$x_j^{m_j}$ for some $m_j\geq 2$.
\end{definition}

The following theorem was proved in 1983 by M. Saito \cite{SaM88}.
The proof shortly afterwards by Khovanskii and Varchenko
\cite{KV85} is completely different.

\begin{theorem}\label{t7.16}
Let $f$ be a Newton nondegenerate and convenient singularity.
For any $\omega\in \Omega^{n+1}_{\C^{n+1},0}$,
its order $\alpha(\omega)$ (defined in remark \ref{t7.14} (iv))
is $\alpha(\omega)=\oooo{\nu}(\omega)-1$.
\end{theorem}

The following corollary is an easy consequence.
It is proved in \cite[Satz 1.10]{He93}.

\begin{corollary}\label{t7.17}
Let $f$ be a Newton nondegenerate and convenient singularity.
Define
\begin{eqnarray}\label{7.48}
s(f)&:=& \min\left(\nu(\frac{\ppp f}{\ppp x_j}\cdot 
\omega_0)-1\, |\, j\in\{0,\ldots,n\}\right)>0,\\
I(f)&:=& \{i\in \Z_{\geq 0}^{n+1}\, |\, 
\nu(x^i\omega_0)-1<s(f)\}.\label{7.49}
\end{eqnarray}
Then for $i\in I(f)$
\begin{eqnarray}\label{7.50}
\alpha(x^i\omega_0)=\nu(x^i\omega_0)-1,
\end{eqnarray}
the numbers $\alpha(x^i\omega_0)$, $i\in I(f)$, are the
spectral numbers in the interval $(-1,s(f))$,
and 
\begin{eqnarray}\label{7.51}
\alpha((\sum_ib_ix^i)\cdot\omega_0) =\left\{
\begin{array}{l}
\min(\alpha(x^i\omega_0)\, |\, i\in I(f),b_i\neq 0)\\
\hspace*{1cm} \textup{if an }i\in I(f)\textup{ with }
b_i\neq 0\textup{ exists},\\
\geq s(f)\qquad \textup{else}.\end{array}\right.
\end{eqnarray}
\end{corollary}

\begin{remarks}\label{t7.18}
(i) We expect that theorem \ref{t7.16} holds also without
the condition that $f$ is convenient.
This would be desirable as many normal forms of singularities
are Newton nondegenerate, but not convenient.

\medskip
(ii) A singularity is $(\mu+1)$-determined, i.e.
$f+g\sim_\RR f$ for any $g\in {\bf m}^{\mu+1}$,
where ${\bf m}$ is the maximal ideal in $\C\{x\}$ \cite{Ma68}.
If $f$ is Newton nondegenerate, then $f+\sum_{j=0}^nc_jx_j^{m_j}$
for arbitrary $m_j\geq \mu+1$ and sufficiently generic
$c_j\in \C^*$ is Newton nondegenerate and convenient
and right equivalent to $f$. 

Furthermore, because of ${\bf m}^\mu\subset J(f)$ and the
Artin approximation theorem, one can choose a coordinate
change $\varphi$ with $f+\sum_{j=0}^nc_jx_j^{m_j}
=f\circ\varphi$ such that all 
$\varphi_j-x_j\in {\bf m}^{\min(m_k)-\mu}$. 
Unfortunately, this is not sufficient for a generalization 
of theorem \ref{t7.16} to the case where $f$ is not 
convenient.

\medskip
(iii) We claim that the calculations in the proof of theorem
\ref{t9.6} can be carried out with almost no change
(but with additional terms) for $f+\sum_{j=0}^nc_jx_j^{m_j}$
with large $m_j$ and that they give essentially the
same results. With this claim, we justify that we calculate
in the proof of theorem \ref{t9.6} with the normal forms $f$ 
in table \eqref{9.1} which are almost all not convenient,
but that we apply theorem \ref{t7.16} and corollary \ref{t7.17}.

\medskip
(iv) Theorem \ref{t7.16} holds without the condition that
$f$ is convenient if $f$ is semiquasihomogeneous.
That is the case when there is only one compact face 
of dimension $n$.
\end{remarks}

\begin{definition}\label{t7.19}
(a) A singularity $f$ is {\it semiquasihomogeneous} with
weights $w_0,\ldots,w_n\in\Q_{>0}$ if 
\begin{eqnarray}\label{7.52}
f=\sum_{i\in \Z_{\geq 0}^{n+1}}a_ix^i\textup{ with }
\deg_wx^i\geq 1\textup{ for all }i\textup{ with }a_i\neq 0,
\end{eqnarray}
and the quasihomogeneous polynomial
\begin{eqnarray}\label{7.53}
f_{qh}:= \sum_{i:\, \deg_wx^i=1}a_ix^i
\end{eqnarray}
has an isolated singularity at 0.

\medskip
(b) A singularity $f$ is {\it quasihomogeneous} if it
is semiquasihomogeneous with $f=f_{qh}$.
\end{definition}

A quasihomogeneous singularity $f$ satisfies the Euler equation
\begin{eqnarray}\label{7.54}
f=\sum_{j=0}^n w_jx_j\frac{\ppp f}{\ppp x_j}.
\end{eqnarray}
This equation and \eqref{7.24} and elementary calculations
in \cite{Br70} imply part (a) of the following lemma.

\begin{lemma}\label{t7.20}
(a) Let $f$ be a quasihomogeneous singularity with
weights $(w_0,\ldots,w_n)$. If $\omega=x^i\omega_0$ is a
monomial differential form then
\begin{eqnarray}
\textup{either}&&s[\omega]_0=0\nonumber\\
\textup{or}&& \alpha(\omega)=\deg_w(x^ix_0\ldots x_n)-1
\textup{ and }s[\omega]_0=s(\omega,\alpha(\omega)).\hspace*{1cm}
\label{7.55}
\end{eqnarray}

(b) Let $f$ be a semiquasihomogeneous singularity
with weights $(w_0,\ldots,w_n)$ and $f\neq f_{qh}$.
The 1-parameter family $f_{qh}+t\cdot (f-f_{qh})$
is a $\mu$-constant family. 
If $\omega=x^i\omega_0$ is a monomial differential form then
\begin{eqnarray}\label{7.56}
\alpha(\omega)&\geq & \deg_w(x^ix_0\ldots x_n)-1,\\
s(\omega&,&\deg_w(x^ix_0\ldots x_n)-1)(t)= 
s[\omega]_0(0),\nonumber\\
s(\omega,\alpha)(t)&=&
\sum_{k\geq 0}\frac{1}{k!}\cdot t^k\cdot (-\ppp_\tau)^k
s((f-f_{qh})^k\cdot \omega,\alpha+k)(0).\nonumber
\end{eqnarray}
The last expression is polynomial in $t$
because $\alpha((f-f_{qh})^k\omega)>\alpha+k$
for large $k$.
\end{lemma}

{\bf Proof of part (b):} 
In \cite[ch. 12]{AGV85} it is shown that $f_{qh}+t(f-f_{qh})$
is a $\mu$-constant family. The other assertions follow
with theorem \ref{t7.13} (c) and part (a) of lemma \ref{t7.20}.
\hfill $\Box$

\section{Review on marked singularities, their moduli spaces, \texorpdfstring{$\mu$}{mu}-constant monodromy groups and Torelli conjectures}
\label{c8}
\setcounter{equation}{0}

\noindent
This paper and the paper \cite{GH17} 
complete the study of the data in the title of this section
for the singularities of modality $\leq 2$.
These data were introduced in \cite{He11}. 
Here we review them. We start with the notions 
{\it marked singularity}
and  {\it strongly marked singularity}. 

\begin{definition}\label{t8.1}
Fix one reference singularity $f_0$.

(a)
Then a strong marking for any singularity $f$
in the $\mu$-homotopy class of $f_0$ 
(i.e. there is a 1-parameter family of singularities with
constant Milnor number connecting $f$ and $f_0$)
is an isomorphism $\rho:(Ml(f),L)\to (Ml(f_0),L)$.

(b)
The pair $(f,\rho)$ is a {\it strongly marked singularity}.
Two strongly marked singularities $(f_1,\rho_1)$ and
$(f_2,\rho_2)$ are right equivalent (notation: $\sim_\RR$) 
if a coordinate change
$\varphi:(\C^{n+1},0)\to(\C^{n+1},0)$ with
$$f_1=f_2\circ \varphi\quad\textup{and}\quad
\rho_1=\rho_2\circ\varphi_{hom}$$
exists, where $\varphi_{hom}:(Ml(f_1),L)\to (Ml(f_2),L)$ 
is the induced isomorphism.

(c)
The notion of a marked singularity is slightly weaker.
If $f$ and $\rho$ are as above, then the pair ($f,\pm \rho)$
is a {\it marked singularity} (writing $\pm\rho$, the set
$\{\rho,-\rho\}$ is meant, neither $\rho$ nor 
$-\rho$ is preferred). 

(d)
Two marked singularities $(f_1,\pm\rho_1)$ and $(f_2,\pm\rho_2)$ 
are right equivalent (notation: $\sim_\RR$) 
if a coordinate change $\varphi$ with 
$$f_1=f_2\circ \varphi\quad
\textup{and}\quad\rho_1=\varepsilon\rho_2\circ\varphi_{hom}
\qquad\textup{for some }\varepsilon\in\{\pm 1\}$$
exists.
\end{definition}

\begin{remarks}\label{t8.2}
(i) The notion of a marked singularity behaves better than the
notion of a strongly marked singularity, because it is not
known whether all $\mu$-homotopy families of singularities
satisfy one of the following two properties:
\begin{eqnarray}\label{8.1}
\textup{Assumption (8.1):}&&\textup{Any singularity in the }
\mu\textup{-homotopy}\\
&&\textup{class of }f_0\textup{ has multiplicity }\geq 3.\nonumber\\
\textup{Assumption (8.2):}&&\textup{Any singularity in the }
\mu\textup{-homotopy} \label{8.2}\\
&&\textup{class of }f_0\textup{ has multiplicity }2.\nonumber
\end{eqnarray}
We expect that always one of two assumptions holds.
For curve singularities and singularities right equivalent
to semiquasihomogeneous singularities and all singularities
with modality $\leq 2$ this is true, but
in general it is not known. In a $\mu$-homotopy family
where neither of the two assumptions holds, strong marking
behaves badly, see (ii).

\medskip
(ii) If $\textup{mult}(f)=2$ then $(f,\rho)\sim_\RR (f,-\rho)$,
which is easy to see. If $\textup{mult}(f)\geq 3$, then
$(f,\rho)\not\sim_\RR(f,-\rho)$, whose proof in \cite{He11}
is quite intricate. These properties imply that the
moduli space for strongly marked singularities discussed below
is not Hausdorff in the case of a $\mu$-homotopy class 
which satisfies neither one of the assumptions \eqref{8.1} 
or \eqref{8.2}.
\end{remarks}

In \cite{He02} a moduli space $M_\mu(f_0)$ was constructed
for the $\mu$-homotopy class of any singularity $f_0$. 
As a set it is simply
the set of right equivalence classes of singularities in the
$\mu$-homotopy class of $f_0$. But in \cite{He02}
it is constructed as an analytic geometric quotient,
and it is shown that it is locally isomorphic to the 
$\mu$-constant stratum of a singularity modulo the action of
a finite group. The $\mu$-constant stratum of a singularity
is the germ $(S_\mu,0)\subset (M,0)$ within the germ of the
base space of a universal unfolding $F$ of $f$, such that for
a suitable representative 
\begin{eqnarray}\label{8.3}
S_\mu=\{t\in M\, |\, F_t\textup{ has only one singularity }x_0
\textup{ and }F_t(x_0)=0\}.
\end{eqnarray}
It comes equipped with a canonical complex structure,
and $M_\mu$ inherits a canonical complex structure, 
see the chapters 12 and 13 in \cite{He02}.

In \cite{He11} analogous results for marked singularities 
were proved. A better property is that $M^{mar}_\mu$ is
locally isomorphic to a $\mu$-constant stratum without
dividing out a finite group action. 
Therefore one can consider it as a {\it global $\mu$-constant
stratum} or as a {\it Teichm\"uller space for singularities}.
The following theorem
collects results from \cite[theorem 4.3]{He11}.

\begin{theorem}\label{t8.3}
Fix one reference singularity $f_0$. Define the sets
\begin{eqnarray}\label{8.4}
M^{smar}_\mu(f_0) &:=& 
\{\textup{strongly marked }(f,\rho)\, |\, \\
&&f\textup{ in the }\mu\textup{-homotopy class of }f_0\}/
\sim_\RR,  \nonumber \\
M^{mar}_\mu(f_0) &:=& 
\{\textup{marked }(f,\pm\rho)\, |\,  
\label{8.5} \\
&&f\textup{ in the }\mu\textup{-homotopy class of }f_0\}/
\sim_\RR.  \nonumber
\end{eqnarray}

(a) $M^{mar}_\mu(f_0)$ carries a natural canonical complex structure.
It can be constructed with the underlying reduced complex
structure as an analytic geometric quotient 
(see \cite[theorem 4.3]{He11} for details).

(b) The germ $(M^{mar}_\mu(f_0),[(f,\pm\rho)])$ with its canonical
complex structure is isomorphic to the $\mu$-constant stratum 
of $f$ with its canonical complex structure 
(see \cite[chapter 12]{He02} for the definition of that).

(c) For any $\psi\in G_\Z(f_0)=:G_\Z$, the map
$$\psi_{mar}:M_\mu^{mar}\to M_\mu^{mar},\quad
[(f,\pm\rho)]\to [(f,\pm\psi\circ\rho)]$$
is an automorphism of $M_\mu^{mar}$.  The action
$$G_\Z\times M_\mu^{mar}\to M_\mu^{mar},\quad 
(\psi,[(f,\pm\rho)]\mapsto \psi_{mar}([(f,\pm\rho)])$$
is a group action from the left.

(d) The action of $G_\Z$ on $M_\mu^{mar}$ is properly discontinuous.
The quotient $M_\mu^{mar}/G_\Z$ is the moduli space $M_\mu$ 
for right equivalence classes in the 
$\mu$-homotopy class of $f_0$, with its canonical complex structure.
Especially, $[(f_1,\pm\rho_1)]$ and $[(f_2,\pm\rho_2)]$ are in one
$G_\Z$-orbit if and only if $f_1$ and $f_2$ are right equivalent.

(e) If assumption \eqref{8.1} or \eqref{8.2} holds then 
(a) to (d) 
are also true for $M_\mu^{smar}$ and $\psi_{smar}$
with $\psi_{smar}([(f,\rho)]):=[(f,\psi\circ\rho)]$.
If neither \eqref{8.1} nor \eqref{8.2} holds then the natural
topology on $M^{smar}_\mu$ is not Hausdorff.
\end{theorem}

We stick to the situation in theorem \ref{t8.3} and
define two subgroups of $G_\Z(f_0)$.
The definitions in 
\cite[definition 3.1]{He11} are different, they use 
$\mu$-constant families. The following definitions are
a part of theorem 4.4 in \cite{He11}.

\begin{definition}\label{t8.4}
Let $(M^{mar}_\mu)^0$ be the topological component
of $M^{mar}_\mu$ (with its reduced complex structure) which
contains $[(f_0,\pm\id)]$. Then
\begin{eqnarray}\label{8.6}
G^{mar}(f_0)&:=& \{\psi\in G_\Z\, |\, \psi\textup{ maps }
(M_\mu^{mar})^0 \textup{ to itself}\}\subset G_\Z(f_0).
\end{eqnarray}
If assumption \eqref{8.1} or \eqref{8.2} holds, 
$(M^{smar}_\mu)^0$ and 
$G^{smar}(f_0)\subset G_\Z(f_0)$ are defined analogously.
\end{definition}

The following theorem is also proved in \cite{He11}.

\begin{theorem}\label{t8.5}
(a) In the situation above, the map
\begin{eqnarray*}
G_\Z/G^{mar}(f_0)&\to& \{\textup{topological components of }
M_\mu^{mar}\}\\
\psi\cdot G^{mar}(f_0)&\mapsto& \textup{the component }
\psi_{mar}((M_\mu^{mar})^0)
\end{eqnarray*}
is a bijection.

(b) If assumption \eqref{8.1} or \eqref{8.2} holds then (a) 
is also true for $M_\mu^{smar}$ and $G^{smar}(f_0)$.

(c) $-\id\in G_\Z$ acts trivially on $M_\mu^{mar}(f_0)$. 
Suppose that assumption \eqref{8.2} holds and that
$f_0=g_0(x_0,\ldots,x_{n-1})+x_n^2$. Then $-\id$
acts trivially on $M^{smar}_\mu(f_0)$ and 
\begin{eqnarray}\label{8.7}
\begin{array}{ccccc}
M_\mu^{smar}(f_0)&=& M_\mu^{mar}(f_0) &=& M_\mu^{mar}(g_0),\\
G^{smar}(f_0)&=& G^{mar}(f_0) &=& G^{mar}(g_0).
\end{array}
\end{eqnarray}
Suppose additionally that assumption \eqref{8.1} holds for $g_0$
(instead of $f_0$). 
Then $\{\pm\id\}$ acts freely on $M_\mu^{smar}(g_0)$,
and the quotient map
$$M_\mu^{smar}(g_0) \stackrel{/\{\pm\id\}}{\longrightarrow}
M_\mu^{mar}(g_0),\quad [(f,\rho)]\mapsto [(f,\pm\rho)]$$
is a double covering.
\end{theorem}

The following conjecture was formulated as
conjecture 3.2 in \cite{He11}. 

\begin{conjecture}\label{t8.6}\cite[Conjecture 3.2]{He11}
(a) Fix a singularity $f_0$. Then $M^{mar}_\mu$ is
connected. Equivalently (in view of theorem \ref{t8.5} (a)):
$G^{mar}(f_0)=G_\Z.$

(b) If the $\mu$-homotopy class of $f_0$ satisfies 
assumption \eqref{8.1}, then $-\id\notin G^{smar}(f_0)$.
\end{conjecture}

The study of the singularities with modality $\leq 2$
in \cite{He11}\cite{GH17} and this paper gives:
Part (b) is true for all singularities with modality
$\leq 2$. Part (a) is true for almost all singularities
with modality $\leq 2$, but not for all.
The exceptions are the subseries for $p=m\cdot r$
of the eight bimodal series. 
This is a part of theorem \ref{t9.1}.
Now we expect that part (a) will be wrong for
many singularities. 

Using the other definition of $G^{mar}$ in \cite{He11},
part (a) says that up to $\pm\id$, any element of $G_\Z$
can be realized as  transversal monodromy of a
$\mu$-constant family with parameter space $S^1$.
As it is wrong for some singularities and probably 
for many more, part (a) of conjecture \ref{t8.6} 
has to be replaced now by the
question whether the subgroup $G^{mar}$ of $G_\Z$ 
can be described in a nice conceptual way.

In order to understand the stabilizers
$\textup{Stab}_{G_\Z}([(f,\rho)])$ and 
$\textup{Stab}_{G_\Z}([(f,\pm\rho)])$
of points $[(f,\rho)]\in M^{smar}_\mu(f_0)$ and 
$[(f,\pm\rho)]\in M^{mar}_\mu(f_0)$,
we have to look at
the {\it symmetries} of a single singularity.
These had been discussed in \cite[chapter 13.2]{He02}.
The discussion had been taken up again in \cite{He11}.

\begin{definition}\label{t8.7}
Let $f_0=f_0(x_0,\ldots,x_n)$ be a reference singularity 
and let $f$ be any singularity in the $\mu$-homotopy class of $f_0$.
If $\rho$ is a marking, then $G_\Z(f)=\rho^{-1}\circ G_\Z\circ\rho$.

We define
\begin{eqnarray}\label{8.8}
\RR &:=&\{\varphi:(\C^{n+1},0)\to (\C^{n+1},0)\ 
\textup{ biholomorphic}\},\\
\RR^f &:=& \{\varphi\in\RR\, |\, f\circ \varphi=f\},\label{8.9}\\
R_f&:=& j_1\RR^f/(j_1\RR^f)^0,\label{8.10}\\
G^{smar}_\RR(f)&:=&\{\varphi_{hom}\, |\, \varphi\in\RR^f\}
\subset G_\Z(f),\label{8.11}\\
G^{mar}_\RR(f)&:=& \{\pm\psi\, |\, \psi\in G^{smar}_\RR(f)\},
\label{8.12}\\
G^{smar,gen}_\RR(f_0)&:=& \bigcap_{[(f,\rho)]\in M_\mu^{smar}}
\rho^{-1}\circ G^{smar}_\RR(f)\circ\rho\subset G_\Z.\label{8.13}
\end{eqnarray}
\end{definition}

Again, the definition of $G^{smar}_\RR$ is different from
the definition in \cite[definition 3.1]{He11}.
The characterization in \eqref{8.11} is 
\cite[theorem 3.3. (e)]{He11}.
$R_f$ is the finite group of components of the group
$j_1\RR^f$ of 1-jets of coordinate changes which leave $f$
invariant.
The following theorem collects results from several theorems
in \cite{He11}.

\begin{theorem}\label{t8.8}
Consider the data in definition \ref{t8.7}.

(a) If $\textup{mult}(f)\geq 3$ then $j_1\RR^f= R_f$.

(b) The homomorphism $()_{hom}:\RR^f\to G_\Z(f)$
factors through $R_f$. Its image is $(R_f)_{hom}=G^{smar}_\RR(f)
\subset G_\Z(f)$.

(c) The homomorphism $()_{hom}:R_f\to G^{smar}_\RR(f)$ 
is an isomorphism.

(d) 
\begin{eqnarray}\label{8.14}
-\id\notin G^{smar}_\RR(f)\iff \mult f\geq 3.
\end{eqnarray}
Equivalently: $G^{mar}_\RR(f)=G^{smar}_\RR(f)$ if $\mult f=2$, and
$G^{mar}_\RR(f)=G^{smar}_\RR(f)\times\{\pm\id\}$ if $\mult f\geq 3$.

(e) $G^{mar}_\RR(f)=G^{mar}_\RR(f+x_{n+1}^2)$.

(f) $M_h\in G^{smar}(f)$. If $f$ is quasihomogeneous then 
$M_h\in G^{smar}_\RR(f)$.

(g) For any $[(f,\rho)]\in M_\mu^{smar}$
\begin{eqnarray}\label{8.15}
\Stab_{G_\Z}([(f,\rho)]) &=& \rho\circ G^{smar}_\RR(f)\circ 
\rho^{-1},\\
\Stab_{G_\Z}([(f,\pm\rho)]) 
&=& \rho\circ G^{mar}_\RR(f)\circ \rho^{-1}.\label{8.16}
\end{eqnarray}
(\eqref{8.15} does not require assumption \eqref{8.1} or \eqref{8.2}).
As $G_\Z$ acts properly discontinuously on $M^{mar}_\mu(f_0)$,
$G^{smar}_\RR(f)$ and $G^{mar}_\RR(f)$ are finite.
(But this follows already from the finiteness of $R_f$ and (b).)
\end{theorem}

The group $G^{smar,gen}_\RR(f_0)$ in \eqref{8.13}
had not been considered in \cite{He11}.
Usually it is very small. It is useful because of the 
following elementary fact.

\begin{lemma}\label{t8.9}
Let $T$ be the parameter space of a $\mu$-constant
family as in definition \ref{t7.12}.
The transversal monodromy of it is the
representation $\pi_1(T,t^0)\to G_\Z(F_{t^0})$
which comes from the local system
$\bigcup_{t\in T}Ml(F_t)$. 

If its image is in $G^{smar,gen}_\RR(F_{t^0})$, then
there is a natural map $T\to M_\mu^{smar}(F_{t^0})$.
\end{lemma}

{\bf Proof:}
The trivial strong marking $+\id$ for $F_{t^0}$
induces along any path strong markings of other singularities
$F_t$. Two paths which meet at a point $t$, might not 
induce the same strong marking of $F_t$, but the two markings
differ only by an element of $G^{smar}_\RR(F_t)$.
Therefore they induce the same right equivalence class
of a marked singularity. \hfill$\Box$ 

\bigskip
Finally, we come to the Brieskorn lattices of 
marked singularities and Torelli problems.
After fixing a reference singularity $f_0$, 
a marked singularity $(f,\pm\rho)$  comes equipped with a 
{\it marked Brieskorn lattice} $BL(f,\pm\rho)$.
The classifying space $D_{BL}(f_0)$ in theorem \ref{t7.11}
is a classifying space for marked Brieskorn lattices.
Theorem \ref{t7.13} implies part (a) of the following theorem.

\begin{theorem}\label{t8.10}
Fix one reference singularity $f_0$.

(a) There is a natural holomorphic period map
\begin{eqnarray}\label{8.17}
BL:M^{mar}_\mu(f_0)\to D_{BL}(f_0).
\end{eqnarray}
It is $G_\Z$-equivariant.

(b) \cite[theorem 12.8]{He02}
It is an immersion, here the reduced complex structure on 
$M^{mar}_\mu(f_0)$ is considered.
\end{theorem}

The second author conjectured part (b) of the following
global Torelli conjecture in \cite{He93}, part (c) in \cite{He02} 
and part (a) in \cite{He11}.

\begin{conjecture}\label{t8.11}
Fix one reference singularity $f_0$. 

(a) The period
map $BL:M^{mar}_\mu\to D_{BL}$ is injective.

(b) The period map
$LBL:M_\mu=M^{mar}_\mu/G_\Z\to D_{BL}/G_\Z$ is injective.

(c) For any singularity $f$ in the $\mu$-homotopy class
of $f_0$ and any marking $\rho$,
\begin{eqnarray}\label{8.19}
\textup{Stab}_{G_\Z}([(f,\pm\rho)])
=\textup{Stab}_{G_\Z}(BL([(f,\pm\rho)]))
\end{eqnarray}
(only $\subset$ and the finiteness of both groups are clear).
\end{conjecture}

The second author has a long-going project on Torelli type
conjectures. Already in \cite{He93}, part (b) was 
proved for all simple and unimodal singularities and 
almost all bimodal singularities 
(all except 3 subseries of the 8 bimodal series).
This was possible without the general construction of
$M_\mu$ and $D_{BL}$, which came later in \cite{He02}
and \cite{He99}. In the concrete cases considered in \cite{He93},
it is easy to identify a posteriori the spaces $M_\mu$
and $D_{BL}$. 

The following lemma from \cite{He11} clarifies the logic
between the parts (a), (b) and (c) of conjecture \ref{t8.11}.

\begin{lemma}\label{t8.12}
In conjecture \ref{t8.11}, (a) $\iff$ (b) and (c).
\end{lemma}

Part (a) of conjecture \ref{t8.11} was proved in 
\cite{He11} for the simple and those 22 of the 28 exceptional
unimodal and bimodal singularities, where all eigenvalues
of the monodromy have multiplicity one.
In \cite{GH17} part (a) was proved for the remaining unimodal 
and the remaining exceptional bimodal singularities.
In the sections \ref{c9} and \ref{c10}, part (a)
will be proved for the remaining bimodal singularities,
namely the bimodal series singularities and the 
quadrangle singularities. 

As part (b) had been proved for almost all singularities
with modality $\leq 2$, the main work in \cite{GH17}
and here is the good control of the group $G_\Z$.
But that is surprisingly difficult. 
In the case of the bimodal singularities in this paper,
also the control of the Gauss-Manin connection side
had to be improved: We provide better information on
the transversal monodromy of the studied families
than in \cite{He93}. Due to this improvement, also
the annoying gap of 3 subseries of the 8 bimodal series,
where part (b) was not proved in \cite{He93},
could be closed here.

\begin{remark}\label{t8.13}
In the sections \ref{c9} and \ref{c10}, we will
restrict to consider surface singularities, i.e. 
singularities in $3$ variables. 
This is justified by the following corollary.
It is an application for suspensions of the
Thom-Sebastiani formula for the Fourier-Laplace transforms
of Brieskorn lattices in theorem \ref{t7.9}.
This is elegant, but the preparations in section 
\ref{c7} were heavy. In the earlier papers
\cite{He93}\cite{He11}\cite{GH17}, we had dealt
with this problem in a less conceptual, but leaner way,
sometimes with extra calculations for curve singularities.
\end{remark}

\begin{corollary}\label{t8.14}
Consider the $\mu$-homotopy class of a reference singularity
$f_0(x_0,\ldots,x_n)$ which satisfies assumption \eqref{8.1}
and such that for any $m\geq 1$ the $\mu$-homotopy class
of $f_0+\sum_{j=n+1}^{n+m}x_j^2$ satisfies 
assumption \eqref{8.2}.

Fix a number $m\geq 1$. 
The global Torelli conjecture \ref{t8.11} (a)
holds for $f_0$ if any only if it holds 
for the reference singularity $f_0+\sum_{j=n+1}^{n+m}x_j^2$
\end{corollary}

{\bf Proof:}
By \eqref{8.7}, $M_\mu^{mar}(f_0)$ and 
$M_\mu^{mar}(f_0+\sum_{j=n+1}^{n+m}x_j^2)$ are canonically
isomorphic. By theorem \ref{t7.11} (c),
the classifying spaces $D_{BL}(f_0)$ and
$D_{BL}(f_0+\sum_{j=n+1}^{n+m}x_j^2)$ are canonically isomorphic.
It rests to see that these isomorphisms are compatible
with the period maps $BL$ for $f_0$ and for 
$f_0+\sum_{j=n+1}^{n+m}x_j^2$. This is also rather 
clear from the formula \eqref{7.40} for the
TEZP-structure of a suspension.
\hfill$\Box$

\section{Period maps and Torelli results for the
bimodal series and 
\texorpdfstring{$G_\Z\supsetneqq G^{mar}$}{GZ/Gmar} 
for the subseries}\label{c9}
\setcounter{equation}{0}

\noindent
In this section we will prove for the bimodal series the strong
global Torelli conjecture \ref{t8.11} (a), 
the conjecture \ref{t8.6} (b)
$-\id\notin G^{smar}$ and for the singularities 
with $m\not|p$ the conjecture \ref{t8.6} (a)
$G_\Z=G^{mar}$. But for the singularities in the subseries
with $m|p$, we will see $G_\Z\supsetneqq G^{mar}$,
$|G_\Z|=\infty,|G^{mar}|<\infty$. 
Theorem \ref{t9.1} states these results in more detail.

The singularities in the eight bimodal series 
$W_{1,p}^\sharp$, $S_{1,p}^\sharp$, $U_{1,p}$,
$E_{3,p}$, $Z_{1,p}$, $Q_{2,p}$, $W_{1,p}$ and $S_{1,p}$ have
as surface singularities the normal forms in 
table \eqref{9.1} \cite[15.1]{AGV85}.
Here $p\geq 1$ and $q\geq 1$, 
and the parameters $(t_1,t_2)$ are in $T:=(\C-\{0\})\times \C$.

\begin{eqnarray}\label{9.1}
\begin{array}{ll}
W_{1,2q-1}^\sharp
& (x^2+y^3)^2+(t_1+t_2y)xy^{4+q}+z^2\\
W_{1,2q}^\sharp
& (x^2+y^3)^2+(t_1+t_2y)x^2y^{3+q}+z^2\\
S_{1,2q-1}^\sharp
& x^2z+y^3z+yz^2+(t_1+t_2y)xy^{3+q}\\
S_{1,2q}^\sharp
& x^2z+y^3z+yz^2+(t_1+t_2y)x^2y^{2+q}\\
U_{1,2q-1}
& x^3+xz^2+xy^3+(t_1+t_2y)y^{1+q}z^2\\
U_{1,2q}
& x^3+xz^2+xy^3+(t_1+t_2y)y^{3+q}z\\
E_{3,p}
& x^3+x^2y^3+(t_1+t_2y)y^{9+p}+z^2\\
Z_{1,p}
& x^3y+x^2y^3+(t_1+t_2y)y^{7+p}+z^2\\
Q_{2,p}
& x^3+yz^2+x^2y^2+(t_1+t_2y)y^{6+p}\\
W_{1,p}
& x^4+x^2y^3+(t_1+t_2y)y^{6+p}+z^2\\
S_{1,p}
& x^2z+yz^2+x^2y^2+(t_1+t_2y)y^{5+p}
\end{array}
\end{eqnarray}

Recall that table \eqref{5.1} lists for these singularities
the Milnor number $\mu$, the characteristic polynomials
$b_j$, $j\geq 1$, of the monodromy on the Orlik blocks
$B_j$ in theorem \ref{t5.1}, the order $m$ of the monodromy
on $B_1$ and the index $r_I=[Ml(f):\bigoplus_{j\geq 1}B_j]$.
The order of the monodromy on $B_2$ is
\begin{eqnarray}\label{9.2}
m+r_I\cdot p=:m_2.
\end{eqnarray}
We will need the space $T^{cov}:=(\C-\{0\})\times \C$
and the $m_2$-fold covering
\begin{eqnarray}\label{9.3}
c_T:T^{cov}\to T,\qquad (\tau_1,t_2)\mapsto (\tau_1^{m_2},t_2).
\end{eqnarray}
For each 2-parameter family of singularities in table \eqref{9.1},
we choose $f_0:= f_{(1,0)}$ as reference singularity.
In the following, we will write $M_\mu^{mar}$,
$(M_\mu^{mar})^0$, $G_\Z$, $G^{mar}$, $Ml$, $H^\infty$
and $C^\alpha$ for 
$M_\mu^{mar}(f_0)$,$(M_\mu^{mar}(f_0))^0$, $G_\Z(f_0)$, 
$G^{mar}(f_0)$, $Ml(f_0)$, $H^\infty(f_0)$ and $C^\alpha(f_0)$.

We denote by $M_T\in G_\Z$ the monodromy of the homology bundle
$\bigcup_{(t_1,t_2)\in T}Ml(f_{(t_1,t_2)})\to T$
along the cycle $\{(e^{2\pi i s},0)\, |\, s\in [0,1]\}$.
We call $M_T$ the {\it transversal monodromy}.
By the other definition of $G^{mar}$ in \cite{He11}, $M_T\in G^{mar}$.
As always, $\zeta:=e^{2\pi i /m}$.

\begin{theorem}\label{t9.1}
Consider a family of bimodal series singularities in table 
\eqref{9.1}

(a) $M_T^{m_2}=\id$. Therefore the pull back to $T^{cov}$
with $c_T$ of the family of singularities over $T$ 
has trivial transversal monodromy.
Thus the strong marking $+\id$ for $f_{(1,0)}$ induces
a well defined strong marking for each singularity of
this family over $T^{cov}$. This gives a map
$T^{cov}\to (M_\mu^{smar})^0$ and a map 
$T^{cov}\to (M_\mu^{mar})^0$.

\medskip
(b) Both maps are isomorphisms.
And $-\id\notin G^{smar}$, where $G^{smar}$ is the
group for the singularities of multiplicity $\geq 3$,
namely the curve singularities
$W_{1,p}^\sharp, E_{3,p}, Z_{1,p}, W_{1,p}$
and the surface singularities $S_{1,p}^\sharp, U_{1,p},
Q_{2,p}, S_{1,p}$. So, conjecture \ref{t8.6} (b) is true.

\medskip
(c) The period map $BL:M_\mu^{mar}\to D_{BL}$ is an embedding.
So, the strong global Torelli conjecture \ref{t8.11} (a) is true.

\medskip
(d) If $m\not|p$ then $G_\Z=G^{mar}$.
So, here conjecture \ref{t8.6} (a) is true.

\medskip
(e) In the case of the subseries with $m|p$, 
$G_\Z\supsetneqq G^{mar}$. So, here conjecture \ref{t8.6} (a)
is wrong. More precisely, $G^{mar}$ and $G_\Z$ are as follows.
$M_T$ has on the 2-dimensional $\C$-vector space
$Ml_\zeta$ the eigenvalues $1$ and $\oooo\zeta$.
Let $Ml_{\zeta,1}$ be the 1-dimensional 
eigenspace of $M_T$ on $Ml_\zeta$
with eigenvalue $1$. Then $|G_\Z|=\infty$ and
$|G^{mar}|<\infty$ and 
\begin{eqnarray}\label{9.4}
G^{mar}&=& \{g\in G_\Z\, |\, g(Ml_{\zeta,1})=Ml_{\zeta,1}\}.
\end{eqnarray}
$\Psi(G_\Z)$ is an infinite Fuchsian group by theorem \ref{t5.1}
(c). $\Psi(G^{mar})$ is the finite subgroup of elliptic 
elements which fix the point
$[Ml_{\zeta,1}]\in \HH_\zeta$
($\HH_\zeta$ was defined in \eqref{5.7}).
And $M_\mu^{mar}$ consists of infinitely many copies of $T^{cov}$.
\end{theorem}

Theorem \ref{t9.1} will be proved in this section in several 
steps. It builds on two hard results.
The first and more difficult one is theorem \ref{t5.1} 
on $G_\Z$. The second one is easier, but still rather
technical. It is the calculation of the multi-valued
period map  $T\to D_{BL}$.
The results are fixed in theorem \ref{t9.6}.

But we prefer to present the nice geometry before the 
technical details. Therefore we will now explain everything
what can be understood without going into the details
of the Gauss-Manin connection and theorem \ref{t9.6}.
Afterwards we will come to the Gauss-Manin connection
and theorem \ref{t9.6}.

Define
\begin{eqnarray}\label{9.5}
\alpha_1:=\frac{-1}{m} < \beta_1:=\frac{-1}{m_2} < 0 
<\alpha_2:=\frac{1}{m_2} < \beta_2:=\frac{1}{m}
\end{eqnarray}
and recall that $\psi_\alpha:H^\infty\to C^\alpha$,
$A\mapsto es(A,\alpha)$, is an isomorphism.
Therefore and because of table \eqref{5.1} 
\begin{eqnarray}\label{9.6}
\dim C^{\beta_1}=\dim C^{\alpha_2}=1,\\
\dim C^{\alpha_1}=\dim C^{\beta_2}=\left\{\begin{array}{ll}
1&\textup{ if }m\not|p,\\
2&\textup{ if }m|p.\end{array}\right.  \nonumber
\end{eqnarray}

For the cases with $m\not|p$, define the 2-dimensional space
\begin{eqnarray}
D_{BL}^{sub}&:=&
\{\C\cdot (v_1+v_2+v_4)\, |\, v_1\in C^{\alpha_1}-\{0\},
v_2\in C^{\beta_1}-\{0\},v_4\in C^{\beta_2}\}\nonumber \\
&=& \{\C\cdot (v_1^0+\rho_1 v_2^0+\rho_2 v_4^0)\, |\, 
(\rho_1,\rho_2)\in(\C-\{0\})\times\C\}  \label{9.7}\\
&&\textup{for some generators }v_1^0,v_2^0,v_4^0\textup{ of }
C^{\alpha_1},C^{\beta_1},C^{\beta_2} \nonumber \\
&\cong& (\C-\{0\})\times \C.\nonumber
\end{eqnarray}

For the cases with $m|p$, the polarizing form $S$
defines an indefinite hermitian form 
$((a,b)\mapsto S(a,\oooo b))$ on $H^\infty_\zeta$.
This follows from the corresponding statement for 
$h_\zeta$ on $Ml_\zeta$ in theorem \ref{t5.1},
from lemma \ref{t2.2} (b) and from the relation between
Seifert form $L$ and polarizing form $S$,
see \eqref{4.20}.
Thus we get a half-plane
\begin{eqnarray}
\HH(C^{\alpha_1})&:=& \{\C\cdot v\, |\, 
v\in C^{\alpha_1}\textup{ with }S(\psi^{-1}_{\alpha_1}(v),
\oooo{\psi_{\alpha_1}^{-1}(v)})<0\}\nonumber\\
&\subset& \P(C^\alpha).\label{9.8}
\end{eqnarray}
Now define for the cases with $m|p$ the 3-dimensional space
\begin{eqnarray}
D_{BL}^{sub}&:=&
\{\C\cdot (v_1+v_2+v_4)\, |\, v_1\in C^{\alpha_1}-\{0\}
\textup{ with }[\C\cdot v_1]\in \HH(C^{\alpha_1}),\nonumber\\
&& \hspace*{1cm} 
v_2\in C^{\beta_1}-\{0\},v_4\in \C\cdot\psi_{\beta_2}
(\oooo{\psi_{\alpha_1}^{-1}(v_1)})\subset C^{\beta_2}\}
\label{9.9} \\
&\cong & \HH(C^{\alpha_1})\times (\C-\{0\})\times \C.
\nonumber
\end{eqnarray}

\begin{theorem}\label{t9.2}
(a) $D_{BL}^{sub}$ embeds canonically into $D_{BL}$.

\medskip
(b) For suitable $v_1^0\in C^{\alpha_1}-\{0\},
v_2^0\in C^{\beta_1}-\{0\}$ and for
$v_4^0:=\psi_{\beta_2}(\oooo{\psi_{\alpha_1}^{-1}(v_1^0)})
\in C^{\beta_2}-\{0\}$, the multi-valued period map
$BL_T:T\to D_{BL}$ has its image in $D_{BL}^{sub}$ and
takes the form
\begin{eqnarray}\label{9.10}
(t_1,t_2)\mapsto \C\cdot \left( v_1^0+ t_1^{1/m_2}\cdot v_2^0
+\left(\frac{t_2}{t_1}+r(t_1)\right)v_4^0ß\right)
\end{eqnarray}
with 
\begin{eqnarray}\label{9.11}
r(t_1)=\left\{
\begin{array}{ll}
0 & \textup{in the cases }(r_I=1\ \&\ p\geq 3),\\
& \textup{the cases }(r_I=2\ \&\ p\geq 2)\\
& \textup{and the case }U_{1,2},\\
c_T\cdot t_1 & \textup{in the cases }(r_I=2\ \&\ p=1)\\
& \textup{and the cases }W_{1,2}^\sharp\textup{ and }
S_{1,2}^\sharp,\\
c_T\cdot t_1^2 & \textup{in the cases }(r_I=1\ \&\ p=1),
\end{array} \right. 
\end{eqnarray}
for a suitable constant $c_T\in\C$. In the cases with $m|p$,
the transversal monodromy $M_T$ has on $C^{\alpha_1}$
the eigenvalues $1$ and $\oooo\zeta$,
and $\C\cdot v_1^0$ is the eigenspace with eigenvalue $1$.
The class $[\C\cdot v_1^0]$ is in $\HH(C^{\alpha_1})$.

\medskip
(c) The induced period map $BL_{T^{cov}}:T^{cov}\to D_{BL}^{sub}$
is an isomorphism if $m\not|p$ and an isomorphism
to the fiber above $[\C\cdot v_1^0]\in\HH(C^{\alpha_1})$
of the projection $D_{BL}^{sub}\to \HH(C^{\alpha_1})$
if $m|p$.

\medskip
(d) In the case of the subseries $U_{1,9r}$, $G^{mar}$
contains an element $g_3$ such that $\Psi(g_3)$ is elliptic
of order $18$ (for all subseries with $p=m\cdot r$, 
$\Psi(M_T)$ is elliptic of order $m$, for $U_{1,9r}$ $m=9$).

\medskip
(e) $f_{(t_1,t_2)}$ and $f_{(\www t_1,\www t_2)}$ are right
equivalent
\begin{eqnarray}\label{9.12}
\iff\left\{
\begin{array}{l}
\exists\ k\in\Z\textup{ with }
(\www t_1,\www t_2)=(\zeta^{r_Ipk}\cdot t_1,
\zeta^{(r_Ip+2)k}\cdot t_2)\\
\hspace*{1cm}\textup{for all 8 series except }U_{1,2q},\\
\exists\ k\in\Z\textup{ and }\varepsilon\in\{\pm 1\}
\textup{ with }\\
\hspace*{1cm}(\www t_1,\www t_2)
=(\varepsilon\zeta^{r_Ipk}\cdot t_1,
\varepsilon\zeta^{(r_Ip+2)k}\cdot t_2)\textup{ for }U_{1,2q}.
\end{array}  \right. 
\end{eqnarray}
\end{theorem}

The parts (a), (b) and (d) of theorem \ref{t9.2} will be 
proved after theorem \ref{t9.6}.

\bigskip
{\bf Proof of theorem \ref{t9.2} (c) and (e):}

(c) This follows immediately from \eqref{9.10}.

(e) First we prove $\Leftarrow$.
We give explicit coordinate changes. A case by case
comparison with the normal forms in table \eqref{9.1}
shows that the following equality \eqref{9.13} holds.
Here $(\delta_1,\delta_2,\delta_3)$ are as in table \eqref{9.14},
and $k\in\Z$.
\begin{eqnarray}\label{9.13}
f_{(t_1,t_2)}(x\cdot \zeta^{\delta_1\cdot k},
y\cdot \zeta^{\delta_2\cdot k},z\cdot \zeta^{\delta_3\cdot k})
=f_{(t_1\cdot\zeta^{r_Ipk},t_2\cdot \zeta^{(r_Ip+2)k})}
(x,y,z).
\end{eqnarray}
\begin{eqnarray}\label{9.14}
\begin{array}{llll}
 & \delta_1 & \delta_2 & \delta_3 \\
W_{1,p}^\sharp\textup{ and }W_{1,p} & 3 & 2 & 0 \\
S_{1,p}^\sharp\textup{ and }S_{1,p} & 3 & 2 & 4 \\
U_{1,p} & 3 & 2 & 3 \\
E_{3,p} & 6 & 2 & 0 \\
Z_{1,0} & 4 & 2 & 0 \\
Q_{2,p} & 4 & 2 & 5 
\end{array}
\end{eqnarray}
In the case $U_{1,2q}$ we have additionally
\begin{eqnarray}\label{9.15}
f_{(t_1,t_2)}(x,y,-z) & =& f_{(-t_1,-t_2)}(x,y,z).
\end{eqnarray}
This shows $\Leftarrow$.

Now we prove $\Rightarrow$. 
Let $f_{(t_1,t_2)}$ and $f_{(\www t_1,\www t_2)}$ be 
right equivalent. Then $BL_T(t_1,t_2)$ and 
$BL_T(\www t_1,\www t_2)$ are isomorphic, so a $g\in G_\Z$
with $g(BL_T(t_1,t_2))=BL_T(\www t_1,\www t_2)$ exists.
We claim that $v_1^0,v_2^0$ and $v_4^0$ are eigenvectors
of $g$ with some eigenvalues $\lambda_1,\lambda_2$
and $\oooo{\lambda_1}$.
For $v_2^0$ this is trivial as $\dim C^{\beta_1}=1$,
for $v_1^0$ in the case $m\not|p$ also.
In the case $m|p$, it follows for $v_1^0$ from \eqref{9.10}.
For $v_4^0$ use 
$v_4^0=\psi_{\beta_2}(\oooo{\psi_{\alpha_1}^{-1}(v_1^0)})$.
We claim also
\begin{eqnarray}\label{9.16}
\lambda_1\in\Eiw(\zeta),\quad \lambda_2\in \Eiw(e^{2\pi i/m_2}).
\end{eqnarray}
For $\lambda_2$ this is a consequence of the following
three facts and of theorem \ref{t2.9} (a)\&(b).

\begin{list}{}{}
\item[(i)]
The 1-dimensional eigenspace $Ml_{e^{2\pi i/m_2}}$ is already
defined over $\Q(e^{2\pi i/m_2})$. Therefore
$\lambda_2\in \Q(e^{2\pi i/m_2})$.
\item[(ii)]
$|\lambda_2|=1$ because $L$ pairs
$Ml_{e^{2\pi i/m_2}}$ and $Ml_{e^{-2\pi i/m_2}}$.
\item[(iii)]
$\lambda_2$ is an algebraic integer because $g\in G_\Z$.
\end{list}

If $m\not|p$, the same reasoning applies also to $\lambda_1$.
Suppose for a moment $m|p$. 

By part (b), the transversal monodromy $M_T$ acts on 
$C^{\alpha_1}$ and on $H^\infty_\zeta$ with eigenvalues 
$1$ and $\zeta$, and the 1-dimensional eigenspaces
with eigenvalue $1$ are 
$\C\cdot v_1^0$ and $\C\cdot\psi_{\alpha_1}^{-1}(v_1^0)$.
Therefore $\C\cdot\psi_{\alpha_1}^{-1}(v_1^0)$ 
is already defined over $\Q(\zeta)$, i.e. 
$\C\cdot\psi_{\alpha_1}^{-1}(v_1^0)\cap H^\infty_{\Q(\zeta)}$
is a 1-dimensional $\Q(\zeta)$-vector space.
This implies (i) $\lambda_1\in\Q(\zeta)$.
(ii) $|\lambda_1|=1$ holds because $v_1^0\in\HH(C^{\alpha_1})$.
And (iii) ($\lambda_1$ is an algebraic integer) holds anyway.
Again with theorem \ref{t2.9} (a)\&(b) we conclude
$\lambda_1\in \Eiw(\zeta)$. Now \eqref{9.16}
is proved in all cases.

The equality $g(BL_T(t_1,t_2))=BL_T(\www t_1,\www t_2)$
becomes
\begin{eqnarray}
&&\C\cdot
\left(\lambda_1\cdot v_1^0+\lambda_2\cdot t_1^{1/m_2}\cdot 
v_2^0 +\oooo{\lambda_1}\left(\frac{t_2}{t_1}+r(t_1)\right)
\cdot v_4^0\right) \nonumber \\
&=& \C\cdot\left( v_1^0 + \www t_1^{1/m_2}\cdot v_2^0
+\left(\frac{\www t_2}{\www t_1}+r(\www t_1)\right)
\cdot v_4^0\right), \nonumber \\
\textup{so}&& \www t_1^{1/m_2} = \lambda_2\oooo{\lambda_1}
\cdot t_1^{1/m_2},\quad 
\frac{\www t_2}{\www t_1}+r(\www t_1) =
\oooo{\lambda_1}^2\left(\frac{t_2}{t_1}+r(t_1)\right),
\nonumber \\
\textup{so}&& \www t_1=\lambda_2^{m_2}\oooo{\lambda_1}^{m_2}
\cdot t_1,\nonumber\\
\textup{and}&& 
\www t_2 = \oooo{\lambda_1}^2\cdot\frac{\www t_1}{t_1}\cdot t_2
+ \www t_1\cdot (\oooo{\lambda_1}^2\cdot r(t_1)-r(\www t_1)).
\label{9.17}
\end{eqnarray}
Because of \eqref{9.16}, we can write $\lambda_1$
and $\lambda_2$ as follows, here $k,l\in\Z$ and 
$\varepsilon_1,\varepsilon_2\in\{\pm 1\}$.
\begin{eqnarray}\label{9.18}
\begin{array}{lll}
 & \lambda_1 & \lambda_2 \\
\textup{All cases with }m\equiv 0(2),\ m_2\equiv 0(2) & 
\oooo{\zeta}^k & e^{2\pi il/m_2} \\
\textup{The cases }W_{1,2q-1}^\sharp\textup{ and }
S_{1,2q-1}^\sharp & 
\varepsilon_2\cdot \oooo{\zeta}^k & 
\varepsilon_2\cdot e^{2\pi il/m_2} \\
\textup{The cases }U_{1,2q-1} &
\varepsilon_1\cdot \oooo{\zeta}^k &
e^{2\pi il/m_2}\\
\textup{The cases }U_{1,2q} &
\varepsilon_1\cdot\oooo{\zeta}^k &
\varepsilon_2\cdot e^{2\pi i l/m_2}
\end{array}
\end{eqnarray}
One checks that \eqref{9.17} boils down to
\begin{eqnarray}\label{9.19}
\www t_1 = \zeta^{r_Ipk}\cdot t_1,\quad 
\www t_2 = \zeta^{(r_Ip+2)k}\cdot t_2,
\end{eqnarray}
in all cases except $U_{1,2q}$. In the cases $U_{1,2q}$,
it boils down to
\begin{eqnarray}\label{9.20}
\www t_1 = \varepsilon_1\varepsilon_2\cdot \zeta^{pk}\cdot t_1,
\quad \www t_2 = \varepsilon_1\varepsilon_2\cdot 
\zeta^{(p+2)k}.
\end{eqnarray}
This finishes the proof of $\Rightarrow$ and
the proof of theorem \ref{t9.2} (e).
\hfill $\Box$

\bigskip
The statements in theorem \ref{t9.1} on the transversal
monodromy ($M_T^{m_2}=\id$, $M_T$ has the eigenvalues
$1$ and $\zeta$ on $Ml_\zeta$) will be proved after
theorem \ref{t9.6}. The rest of theorem \ref{t9.1}
will be proved now.

\bigskip
{\bf Proof of theorem \ref{t9.1}}
(without the statements on $M_T$):

(a) This is clear.

(b) Consider the maps

\begin{eqnarray}\label{9.21}
\begin{xy}
  \xymatrix{
      T^{cov}  \ar[r] \ar[d]_{\cong}    &  (M_\mu^{smar})^0 \ar[d]^{BL}  \\
     D_{BL}^{sub}   \ar@{^{(}->}[r]            &   D_{BL}    
  }
\end{xy}
\end{eqnarray}

As $T^{cov}\hookrightarrow D_{BL}^{sub}\hookrightarrow D_{BL}$
is an embedding, $T^{cov}\to (M_\mu^{smar})^0$ is an 
embedding.

Both spaces $T^{cov}$ and $(M_\mu^{smar})^0$ are locally
$\mu$-constant strata of universal unfoldings
and are therefore smooth of dimension 2.
$D_{BL}^{sub}$ is almost closed in $D_{BL}$.
Its closure consists of itself and the
space $\{\C\cdot(v_1+v_4)\, |\, v_1\textup{ and }v_4
\textup{ as in \eqref{9.7} or \eqref{9.9}}\}$
(so $v_2=0$). No $g\in G_\Z$ maps a point of this space
to a point of $D_{BL}^{sub}$. 
And $T^{cov}$ contains representatives of any right
equivalence class in the $\mu$-homotopy family.
Therefore the image of $(M_\mu^{smar})^0$ in $D_{BL}$
cannot be bigger than $D_{BL}^{sub}$.
Thus $T^{cov}\cong (M_\mu^{smar})^0$.

In the case of singularities of multiplicity 2,
$M_\mu^{smar}\cong M_\mu^{mar}$ holds anyway
by theorem \ref{t8.5} (c), and then also 
$(M_\mu^{smar})^0\cong (M_\mu^{mar})^0$ holds.

Consider the case of singularities of multiplicity $\geq 3$.
Then $-\id\in G_\Z$ acts nontrivially on $M_\mu^{smar}$ by 
theorem \ref{t8.5} (c). It acts trivially on $D_{BL}$.
The map $(M_\mu^{smar})^0\to D_{BL}$ is an embedding.
Therefore $-\id\in G_\Z$ does not act on 
$(M_\mu^{smar})^0$, therefore $-\id\notin G^{smar}$.
Then $(M_\mu^{smar})^0\to (M_\mu^{mar})^0$ is an
isomorphism by theorem \ref{t8.5} (c).

\medskip
(c) for $m\not|p$ and (d):
$(M_\mu^{mar})^0\stackrel{\cong}{\longrightarrow}
T^{cov}\stackrel{\cong}{\longrightarrow}
D_{BL}^{sub}\hookrightarrow D_{BL}$
is an embedding. $G_\Z=G^{mar}$ would imply
$M_\mu^{mar}=(M_\mu^{mar})^0$.
Therefore it is sufficient to prove $G_\Z=G^{mar}$.

Let $g_1\in G_\Z$. It acts on $D_{BL}^{sub}$. 
By the proof of theorem \ref{t9.2} (e), the map
\begin{eqnarray}\label{9.22}
(M_\mu^{mar})^0/G^{mar} \to D_{BL}^{sub}/G_\Z
\end{eqnarray}
is an isomorphism.
Therefore an element $g_2\in G^{mar}$ exists which acts 
in the same way on $D_{BL}^{sub}$ as $g_1$.
Consider $g_3:=g_1\circ g_2^{-1}$. It acts trivially
on $D_{BL}^{sub}$. It has eigenvalues $\lambda_1$,
$\lambda_2$ and $\oooo{\lambda_1}$ on
$C^{\alpha_1},C^{\beta_1}$ and $C^{\beta_2}$.
Therefore
\begin{eqnarray}
&&\C(v_1+v_2+v_4)=\C(\lambda_1\cdot v_1+\lambda_2\cdot v_2+
\oooo{\lambda_1}\cdot v_4)\nonumber\\
&&\hspace*{1cm}\textup{ for any }\C(v_1+v_2+v_4)\in D_{BL}^{sub},
\nonumber\\
&&\textup{thus }
\lambda_2\oooo{\lambda_1}=1,\quad 
\oooo{\lambda_1}^2=\id,\quad\textup{so }
\lambda_1=\lambda_2\in\{\pm 1\},\nonumber \\
&&\textup{and }g_3=\lambda_1\cdot\id \textup{ on }
Ml_\zeta\oplus Ml_{e^{2\pi i/m_2}}.\label{9.23}
\end{eqnarray}
$G_\Z$ was determined in theorem \ref{t5.1} (b).
It contains very few automorphisms $g_3$ with \eqref{9.23}.
Formula \eqref{5.6} and table \eqref{5.1} show that the
group $\{g\in G_\Z\, |\, g=\pm\id\textup{ on }
Ml_\zeta\oplus Ml_{e^{2\pi i/m_2}}\}$ is as follows:
\begin{eqnarray}\label{9.24}
&&\{\pm\id\} \qquad \textup{in the cases }
W_{1,2q-1}^\sharp, S_{1,2q-1}^\sharp,U_{1,2q},E_{3,p},Z_{1,p},\\
&& \{\pm\id,\pm (\id|_{B_1}\times (-M_h^{m_2/2})|_{B_2})\}
\quad \textup{in the cases }
W_{1,2q}^\sharp, S_{1,2q}^\sharp, U_{1,2q-1},\nonumber\\
&& \{\pm\id, \pm((-M_h^{m/2})|_{B_1}\times \id|_{B_2})\}
\qquad \textup{in the cases }
Q_{2,p}, W_{1,p}, S_{1,p}.\nonumber
\end{eqnarray}

{\bf Claim:}
\begin{eqnarray}\label{9.25}
\{g\in G_\Z\, |\, g=\pm\id\textup{ on }
Ml_\zeta\oplus Ml_{e^{2\pi i/m_2}}\}
=G_\RR^{mar}.
\end{eqnarray}
This claim shows $g_3\in G_\RR^{mar}$ and 
$g_1\in G^{mar}$, so that $G_\Z=G^{mar}$.

The inclusion $\supset$ in \eqref{9.25} holds because of the following: 
Any element  of $G^{mar}_\RR=G^{mar}_\RR(f_{(1,0)})$ acts on 
$D_{BL}^{sub}$ with $BL_T(1,0)$ as fixed point. 
The proof of theorem \ref{t9.2} (e) shows that it acts
then trivially on $D_{BL}^{sub}$.

The group $G^{mar}_\RR$ contains $\pm\id$.
In order to prove equality in \eqref{9.25} 
for the cases in the second
and third line of \eqref{9.24},
it is sufficient to show that $G^{mar}_\RR$ contains
more elements than $\pm\id$.
Equivalent is that $G^{smar}_\RR(f) $ for a generic
singularity $f$ with multiplicity $\geq 3$ contains
one other element than $+\id$. 
The following table lists coordinate changes which 
give such an element.
\begin{eqnarray}\label{9.26}
\begin{array}{ll}
W_{1,2q}^\sharp &
(x,y)\mapsto (-x,y) \\
S_{1,2q}^\sharp &
(x,y,z)\mapsto (-x,y,z) \\
U_{1,2q-1} &
(x,y,z)\mapsto (x,y,-z) \\
Q_{2,p} &
(x,y,z)\mapsto (x,y,-z) \\
W_{1,p} & 
(x,y)\mapsto (-x,y) \\
S_{1,p} &
(x,y,z)\mapsto (-x,y,z)
\end{array}
\end{eqnarray}
This proves the claim and finishes the proof
of (c) for $m\not|p$ and (d).

\medskip
(c) for $m|p$ and (e): 
First we prove \eqref{9.4}.

$\Psi(M_T)$ is an elliptic element with fixed point
$[Ml_{\zeta,1}]\in \HH_\zeta$ and angle
$\frac{2\pi }{m}=\arg(\frac{\zeta}{1})$.
All elements of $G^{mar}$, including $M_T$, act on
$\HH(C^{\alpha_1})$ as elliptic elements with fixed
point $[\C\cdot v_1^0]$, because all elements
in $G^{mar}$ act on $(M_\mu^{mar})^0$ and on its image
$BL_{T^{cov}}((M_\mu^{mar})^0)\subset D_{BL}^{sub}$.
Therefore all elements of $G^{mar}$ act on $\HH_\zeta$
as elliptic elements with fixed point
$[Ml_{\zeta,1}]$. This shows $\subset$ in \eqref{9.4}.

Now let $g_1\in \{g\in G_\Z\, |\, g(Ml_{\zeta,1})=Ml_{\zeta,1}\}$.
It has an eigenvalue $\lambda_1$ on $Ml_{\zeta,1}$
and an eigenvalue $\lambda_2$ on the other eigenspace
within $Ml_{\zeta}$
(which is the $h_\zeta$-orthogonal subspace of $Ml_\zeta$).
By \eqref{9.16} $\lambda_1$ and $\lambda_2\in\Eiw(\zeta)$.
Therefore $\Psi(g_1)$ is an elliptic element with 
fixed point $[Ml_{\zeta,1}]\in\HH_\zeta$ and
angle $\arg\frac{\lambda_2}{\lambda_1}$.

In all cases except possibly $U_{1,9r}$, the product
$g_2=g_1\circ M_T^{k}$ for a suitable $k\in\Z$
acts trivially on $\HH_\zeta$.
In the cases $U_{1,9r}$, the product $g_2=g_1\circ g_3^k$
for $g_3\in G^{mar}$ as in theorem \ref{t9.2} (d)
does the same.

Formula \eqref{5.9} in theorem \ref{t5.1} (c) applies to
$g_2$ and shows $g_2\in\{\pm M_h^k\, |\, k\in\Z\}$.
Therefore $g_2\in G^{mar}$ and $g_1\in G^{mar}$.
This shows $\supset$ in \eqref{9.4},
so \eqref{9.4} is now proved.

Especially, $\Psi(G^{mar})$ and $G^{mar}$ are finite.
By theorem \ref{t5.1} (c), $\Psi(G_\Z)$ and $G_\Z$
are infinite. Therefore $G_\Z\supsetneqq G^{mar}$.

By theorem \ref{t8.5} (a), $M_\mu^{mar}$ consists of infinitely
many copies of $(M_\mu^{mar})^0$. 

If two different copies would have intersecting images in
$D_{BL}$ under the period map $BL$, 
the images would coincide, and there would
be a copy different from $(M_\mu^{mar})^0$ with the same
image in $D_{BL}$ as $(M_\mu^{mar})^0$. 
An element $g_3\in G_\Z$ which maps $(M_\mu^{mar})^0$ to
this copy would be in
$\{g\in G_\Z\, |\, g(Ml_{\zeta,1})=Ml_{\zeta,1}\}-G^{mar}
=\emptyset$, a contradiction.
Therefore $BL:M_\mu^{mar}\to D_{BL}^{sub}$
is an embedding. \hfill $\Box$

\begin{remarks}\label{t9.3}
(i) The arithmetic triangle group of type $(2,3,14)$ 
for $Z_{1,0}$ in 
theorem \ref{t3.6} contains elliptic elements of order $3$
although $\arg\zeta=\frac{2\pi }{14}$ 
and the matrices defining these
elliptic elements are in $GL(2,\Z[\zeta])$.
The eigenspaces in $M(2\times 1,\C)$
of these matrices 
are not defined over $\Q(\zeta)$, but only over
$\Q(e^{2\pi i/3},\zeta)$. 
This example shows that \eqref{9.16} in the case $m|p$
and the arguments proving it are nontrivial.

\medskip
(ii) In 1993, the second author worked on the Torelli conjecture
for the unmarked bimodal series singularities.
He missed to consider $M_T$ carefully and thus was not sure
which elliptic elements fix $[\C\cdot v_1^0]\in\HH(C^{\alpha_1})$.
Therefore he could not prove the Torelli conjecture for
the unmarked singularities in the subseries
$S_{1,10r}^\sharp, S_{1,10r}$ and $Z_{1,14r}$.
Now theorem \ref{t9.1} gives the marked and unmarked
Torelli theorem for all bimodal series singularities.
\end{remarks}

Now we come to the spectral numbers and the classifying space $D_{BL}$.

\begin{lemma}\label{t9.4}
Consider a family of bimodal series singulariteis in table \eqref{9.1}.

(a) The spectral numbers $\alpha_1,...,\alpha_\mu$ with
$\alpha_1\leq ...\leq \alpha_\mu$ satisfy
\begin{eqnarray}\label{9.27}
\alpha_1=\frac{-1}{m}<\alpha_2=\frac{1}{m_2} <\alpha_3\leq ...
\leq \alpha_{\mu-2}\\
<\alpha_{\mu-1} =1-\frac{1}{m_2}<\alpha_\mu=1+\frac{1}{m}
\nonumber
\end{eqnarray}
and are uniquely determined by this and the characteristic 
polynomial $\prod_{j\geq 1}b_j$ of the monodromy with $b_j$ 
as in table \eqref{5.1}.

\medskip
(b) Recall from \eqref{9.5} $\beta_1=\frac{-1}{m_2}=-\alpha_2$ 
and $\beta_2=\frac{1}{m}=-\alpha_1$. Then
\begin{eqnarray}\label{9.28}
\dim C^{\alpha_1}=\left\{\begin{array}{ll}
1 & \textup{ if }m\not| p,\\
2 & \textup{ if }m|p,\end{array}\right. 
\end{eqnarray}
\begin{eqnarray}\label{9.29}
\dim C^{\beta}=\left\{\begin{array}{ll}
1 & \textup{for }\beta\in(\alpha_1,\beta_2)\cap\frac{1}{m_2}(\Z-\{0\})
\textup{ if }r_I=1,\\
 & \textup{and for }\beta\in(\alpha_1,\beta_2)\cap
 (\frac{1}{m_2}+\frac{2}{m_2}\Z)\textup{ if }r_I=2,\\
0 & \textup{for other }\beta\in(\alpha_1,\beta_2).\end{array}\right. 
\end{eqnarray}
The following two pictures illustrate this for $2m<p<3m$,  
the first for $r_I=1$, the second for $r_I=2$.

%\bigskip
%\includegraphics[width=1.0\textwidth]{pic-kap9-1.jpg}

\begin{scriptsize}
\begin{tikzpicture}[line cap=round,line join=round,>=triangle 45,x=1.0cm,y=1.0cm]
\clip(-7.3,-0.6) rectangle (4.3,1.2);
\fill[line width=0.8pt,color=ttqqqq,fill=ttqqqq,pattern=north east lines,pattern color=ttqqqq] (-6.2,0.) -- (-6.2,1.) -- (-6.,1.) -- (-6.,0.) -- cycle;
\fill[line width=0.8pt,color=ttqqqq, fill opacity=0] (-6.6,0.) -- (-6.6,1.) -- (-6.8,1.) -- (-6.8,0.) -- cycle;
\fill[line width=0.8pt,color=ttqqqq, fill opacity=0] (-5.2,0.) -- (-5.,0.) -- (-5.,1.) -- (-5.2,1.) -- cycle;
\fill[line width=0.8pt,color=ttqqqq, fill opacity=0] (-4.,1.) -- (-3.8,1.) -- (-3.8,0.) -- (-4.,0.) -- cycle;
\fill[line width=0.8pt,color=ttqqqq, fill opacity=0] (-2.8,1.) -- (-2.6,1.) -- (-2.6,0.) -- (-2.8,0.) -- cycle;
\fill[line width=0.8pt,color=ttqqqq,fill=ttqqqq,pattern=north east lines,pattern color=ttqqqq] (-0.6,1.) -- (-0.4,1.) -- (-0.4,0.) -- (-0.6,0.) -- cycle;
\fill[line width=0.8pt,color=ttqqqq,fill=ttqqqq,pattern=north east lines,pattern color=ttqqqq] (0.6,1.) -- (0.8,1.) -- (0.8,0.) -- (0.6,0.) -- cycle;
\fill[line width=0.8pt,color=ttqqqq,fill=ttqqqq,pattern=north east lines,pattern color=ttqqqq] (1.8,0.) -- (2.,0.) -- (2.,1.) -- (1.8,1.) -- cycle;
\fill[line width=0.8pt,color=ttqqqq, fill opacity=0] (2.8,0.) -- (2.8,1.) -- (3.,1.) -- (3.,0.) -- cycle;
\fill[line width=0.8pt,color=ttqqqq,fill=ttqqqq,pattern=north east lines,pattern color=ttqqqq] (3.4,0.) -- (3.4,1.) -- (3.6,1.) -- (3.6,0.) -- cycle;
\draw [line width=1.2pt] (-8.,0.)-- (8.,0.);
\draw [line width=0.8pt,color=ttqqqq] (-6.2,0.)-- (-6.2,1.);
\draw [line width=0.8pt,color=ttqqqq] (-6.2,1.)-- (-6.,1.);
\draw [line width=0.8pt,color=ttqqqq] (-6.,1.)-- (-6.,0.);
\draw [line width=0.8pt,color=ttqqqq] (-6.,0.)-- (-6.2,0.);
\draw [line width=0.8pt,color=ttqqqq] (-6.6,0.)-- (-6.6,1.);
\draw [line width=0.8pt,color=ttqqqq] (-6.6,1.)-- (-6.8,1.);
\draw [line width=0.8pt,color=ttqqqq] (-6.8,1.)-- (-6.8,0.);
\draw [line width=0.8pt,color=ttqqqq] (-6.8,0.)-- (-6.6,0.);
\draw [line width=0.8pt,color=ttqqqq] (-5.2,0.)-- (-5.,0.);
\draw [line width=0.8pt,color=ttqqqq] (-5.,0.)-- (-5.,1.);
\draw [line width=0.8pt,color=ttqqqq] (-5.,1.)-- (-5.2,1.);
\draw [line width=0.8pt,color=ttqqqq] (-5.2,1.)-- (-5.2,0.);
\draw [line width=0.8pt,color=ttqqqq] (-4.,1.)-- (-3.8,1.);
\draw [line width=0.8pt,color=ttqqqq] (-3.8,1.)-- (-3.8,0.);
\draw [line width=0.8pt,color=ttqqqq] (-3.8,0.)-- (-4.,0.);
\draw [line width=0.8pt,color=ttqqqq] (-4.,0.)-- (-4.,1.);
\draw (-6.6,0.012) node[anchor=north west] {$\alpha_{1}=\frac{-1}{m}$};
\draw (-3.4,0.012) node[anchor=north west] {$\beta_{1}=\frac{-1}{m_2}$};
\draw (1.65,0.012) node[anchor=north west] {$\alpha_{4}$};
\draw (0.45,0.012) node[anchor=north west] {$\alpha_{3}$};
\draw [line width=0.8pt,color=ttqqqq] (-2.8,1.)-- (-2.6,1.);
\draw [line width=0.8pt,color=ttqqqq] (-2.6,1.)-- (-2.6,0.);
\draw [line width=0.8pt,color=ttqqqq] (-2.6,0.)-- (-2.8,0.);
\draw [line width=0.8pt,color=ttqqqq] (-2.8,0.)-- (-2.8,1.);
\draw [line width=0.8pt,color=ttqqqq] (-0.6,1.)-- (-0.4,1.);
\draw [line width=0.8pt,color=ttqqqq] (-0.4,1.)-- (-0.4,0.);
\draw [line width=0.8pt,color=ttqqqq] (-0.4,0.)-- (-0.6,0.);
\draw [line width=0.8pt,color=ttqqqq] (-0.6,0.)-- (-0.6,1.);
\draw [line width=0.8pt,color=ttqqqq] (0.6,1.)-- (0.8,1.);
\draw [line width=0.8pt,color=ttqqqq] (0.8,1.)-- (0.8,0.);
\draw [line width=0.8pt,color=ttqqqq] (0.8,0.)-- (0.6,0.);
\draw [line width=0.8pt,color=ttqqqq] (0.6,0.)-- (0.6,1.);
\draw [line width=0.8pt,color=ttqqqq] (1.8,0.)-- (2.,0.);
\draw [line width=0.8pt,color=ttqqqq] (2.,0.)-- (2.,1.);
\draw [line width=0.8pt,color=ttqqqq] (2.,1.)-- (1.8,1.);
\draw [line width=0.8pt,color=ttqqqq] (1.8,1.)-- (1.8,0.);
\draw [line width=0.8pt,color=ttqqqq] (2.8,0.)-- (2.8,1.);
\draw [line width=0.8pt,color=ttqqqq] (2.8,1.)-- (3.,1.);
\draw [line width=0.8pt,color=ttqqqq] (3.,1.)-- (3.,0.);
\draw [line width=0.8pt,color=ttqqqq] (3.,0.)-- (2.8,0.);
\draw [line width=0.8pt,color=ttqqqq] (3.4,0.)-- (3.4,1.);
\draw [line width=0.8pt,color=ttqqqq] (3.4,1.)-- (3.6,1.);
\draw [line width=0.8pt,color=ttqqqq] (3.6,1.)-- (3.6,0.);
\draw [line width=0.8pt,color=ttqqqq] (3.6,0.)-- (3.4,0.);
\draw (2.45,0.012) node[anchor=north west] {$\beta_{2}=\frac{1}{m}$};
\draw (-1.1,0.012) node[anchor=north west] {$\alpha_{2}=\frac{1}{m_2}$};
\begin{scriptsize}
\draw [color=ttqqqq] (-1.6,0.)-- ++(-2.5pt,0 pt) -- ++(5.0pt,0 pt) ++(-2.5pt,-2.5pt) -- ++(0 pt,5.0pt);
\draw[color=ttqqqq] (-1.604163351327485,0.2418611996401875) node {0};
\end{scriptsize}
\end{tikzpicture}
\end{scriptsize}

%\bigskip
%\includegraphics[width=1.0\textwidth]{pic-kap9-2.jpg}

\begin{scriptsize}
\begin{tikzpicture}[line cap=round,line join=round,>=triangle 45,x=1.2cm,y=1.0cm]
\clip(-6.5,-0.6) rectangle (3.5,1.2);
\fill[line width=0.8pt,color=ttqqqq,fill=ttqqqq,pattern=north east lines,pattern color=ttqqqq] (-5.6,0.) -- (-5.6,1.) -- (-5.4,1.) -- (-5.4,0.) -- cycle;
\fill[line width=0.8pt,color=ttqqqq,fill opacity=0] (-6.,0.) -- (-6.,1.) -- (-6.2,1.) -- (-6.2,0.) -- cycle;
\fill[line width=0.8pt,color=ttqqqq,fill opacity=0] (-4.6,0.) -- (-4.4,0.) -- (-4.4,1.) -- (-4.6,1.) -- cycle;
\fill[line width=0.8pt,color=ttqqqq,fill opacity=0] (-3.4,1.) -- (-3.2,1.) -- (-3.2,0.) -- (-3.4,0.) -- cycle;
\fill[line width=0.8pt,color=ttqqqq,fill opacity=0] (-2.2,1.) -- (-2.,1.) -- (-2.,0.) -- (-2.2,0.) -- cycle;
\fill[line width=0.8pt,color=ttqqqq,fill=ttqqqq,pattern=north east lines,pattern color=ttqqqq] (-1.2,1.) -- (-1.,1.) -- (-1.,0.) -- (-1.2,0.) -- cycle;
\fill[line width=0.8pt,color=ttqqqq,fill=ttqqqq,pattern=north east lines,pattern color=ttqqqq] (-0.004842576466940085,0.995967724566367) -- (0.19515742353305898,0.995967724566367) -- (0.19515742353305898,-0.004032275433632243) -- (-0.004842576466940085,-0.004032275433632243) -- cycle;
\fill[line width=0.8pt,color=ttqqqq,fill=ttqqqq,pattern=north east lines,pattern color=ttqqqq] (1.2,0.) -- (1.4,0.) -- (1.4,1.) -- (1.2,1.) -- cycle;
\fill[line width=0.8pt,color=ttqqqq,fill opacity=0] (2.2,0.) -- (2.2,1.) -- (2.4,1.) -- (2.4,0.) -- cycle;
\fill[line width=0.8pt,color=ttqqqq,fill=ttqqqq,pattern=north east lines,pattern color=ttqqqq] (2.8,0.) -- (2.8,1.) -- (3.,1.) -- (3.,0.) -- cycle;
\draw [line width=1.2pt] (-8.,0.)-- (8.,0.);
\draw [line width=0.8pt,color=ttqqqq] (-5.6,0.)-- (-5.6,1.);
\draw [line width=0.8pt,color=ttqqqq] (-5.6,1.)-- (-5.4,1.);
\draw [line width=0.8pt,color=ttqqqq] (-5.4,1.)-- (-5.4,0.);
\draw [line width=0.8pt,color=ttqqqq] (-5.4,0.)-- (-5.6,0.);
\draw [line width=0.8pt,color=ttqqqq] (-6.,0.)-- (-6.,1.);
\draw [line width=0.8pt,color=ttqqqq] (-6.,1.)-- (-6.2,1.);
\draw [line width=0.8pt,color=ttqqqq] (-6.2,1.)-- (-6.2,0.);
\draw [line width=0.8pt,color=ttqqqq] (-6.2,0.)-- (-6.,0.);
\draw [line width=0.8pt,color=ttqqqq] (-4.6,0.)-- (-4.4,0.);
\draw [line width=0.8pt,color=ttqqqq] (-4.4,0.)-- (-4.4,1.);
\draw [line width=0.8pt,color=ttqqqq] (-4.4,1.)-- (-4.6,1.);
\draw [line width=0.8pt,color=ttqqqq] (-4.6,1.)-- (-4.6,0.);
\draw [line width=0.8pt,color=ttqqqq] (-3.4,1.)-- (-3.2,1.);
\draw [line width=0.8pt,color=ttqqqq] (-3.2,1.)-- (-3.2,0.);
\draw [line width=0.8pt,color=ttqqqq] (-3.2,0.)-- (-3.4,0.);
\draw [line width=0.8pt,color=ttqqqq] (-3.4,0.)-- (-3.4,1.);
\draw (-5.95,0.0) node[anchor=north west] {$\alpha_{1}=\frac{-1}{m}$};
\draw (-2.75,0.0) node[anchor=north west] {$\beta_{1}=\frac{-1}{m_2}$};
\draw (1.13,0.0) node[anchor=north west] {$\alpha_{4}$};
\draw (-0.08,0.0) node[anchor=north west] {$\alpha_{3}$};
\draw [line width=0.8pt,color=ttqqqq] (-2.2,1.)-- (-2.,1.);
\draw [line width=0.8pt,color=ttqqqq] (-2.,1.)-- (-2.,0.);
\draw [line width=0.8pt,color=ttqqqq] (-2.,0.)-- (-2.2,0.);
\draw [line width=0.8pt,color=ttqqqq] (-2.2,0.)-- (-2.2,1.);
\draw [line width=0.8pt,color=ttqqqq] (-1.2,1.)-- (-1.,1.);
\draw [line width=0.8pt,color=ttqqqq] (-1.,1.)-- (-1.,0.);
\draw [line width=0.8pt,color=ttqqqq] (-1.,0.)-- (-1.2,0.);
\draw [line width=0.8pt,color=ttqqqq] (-1.2,0.)-- (-1.2,1.);
\draw [line width=0.8pt,color=ttqqqq] (-0.00484257,0.995)-- (0.195157,0.9959677);
\draw [line width=0.8pt,color=ttqqqq] (0.195157,0.995)-- (0.195157,-0.004);
\draw [line width=0.8pt,color=ttqqqq] (0.195157,-0.004)-- (-0.00484,-0.004);
\draw [line width=0.8pt,color=ttqqqq] (-0.00484,-0.004)-- (-0.00484,0.9959677);
\draw [line width=0.8pt,color=ttqqqq] (1.2,0.)-- (1.4,0.);
\draw [line width=0.8pt,color=ttqqqq] (1.4,0.)-- (1.4,1.);
\draw [line width=0.8pt,color=ttqqqq] (1.4,1.)-- (1.2,1.);
\draw [line width=0.8pt,color=ttqqqq] (1.2,1.)-- (1.2,0.);
\draw [line width=0.8pt,color=ttqqqq] (2.2,0.)-- (2.2,1.);
\draw [line width=0.8pt,color=ttqqqq] (2.2,1.)-- (2.4,1.);
\draw [line width=0.8pt,color=ttqqqq] (2.4,1.)-- (2.4,0.);
\draw [line width=0.8pt,color=ttqqqq] (2.4,0.)-- (2.2,0.);
\draw [line width=0.8pt,color=ttqqqq] (2.8,0.)-- (2.8,1.);
\draw [line width=0.8pt,color=ttqqqq] (2.8,1.)-- (3.,1.);
\draw [line width=0.8pt,color=ttqqqq] (3.,1.)-- (3.,0.);
\draw [line width=0.8pt,color=ttqqqq] (3.,0.)-- (2.8,0.);
\draw (1.8,0.0) node[anchor=north west] {$\beta_{2}=\frac{1}{m}$};
\draw (-1.5176964072373837,0.0) node[anchor=north west] {$\alpha_{2}=\frac{1}{m_2}$};
\begin{scriptsize}
\draw [color=ttqqqq] (-1.6,0.)-- ++(-2.5pt,0 pt) -- ++(5.0pt,0 pt) ++(-2.5pt,-2.5pt) -- ++(0 pt,5.0pt);
\draw[color=ttqqqq] (-1.597,0.297) node {0};
\end{scriptsize}
\end{tikzpicture}
\end{scriptsize}

\medskip
(c) Denote by $(*)$ the condition
\begin{eqnarray*}
(*):\hspace*{1cm} \beta\in(\alpha_1,0)\textup{ with }C^\beta\neq\{0\}
\quad  (\textup{then }\dim C^\beta=1).
\end{eqnarray*}
If $m\not|p$ the classfying space $D_{BL}$ in \cite{He99} is 
\begin{eqnarray}\label{9.30}
D_{BL}&=& \{\C\cdot(v_1+\sum_{\beta:(*)}v_{(\beta)}+v_2)\, |\, \\
&& v_1\in C^{\alpha_1}-\{0\},v_{(\beta)}\in C^\beta,v_2\in C^{\beta_2}\}
\nonumber \\
&\cong & \C^{N_{BL}}\qquad \textup{with }N_{BL}:=|\{\beta:(*)\}|+1.
\nonumber
\end{eqnarray}
In \eqref{9.8} $\HH(C^{\alpha_1})$ was defined for $m|p$.
If $m|p$ then $D_{BL}$ is 
\begin{eqnarray}\label{9.31}
D_{BL}&=& \{\C\cdot(v_1+\sum_{\beta:(*)}v_{(\beta)}+v_2)\, |\,  \\
&& v_1\in C^{\alpha_1}-\{0\}\textup{ with }[\C\cdot v_1]\in \HH(C^{\alpha_1}),
\nonumber \\
&& v_{(\beta)}\in C^\beta,v_2\in\C\cdot\psi_{\beta_2}
(\oooo{\psi_{\alpha_1}^{-1}(v_1)})\subset C^{\beta_2}\}\nonumber\\
&\cong& \HH(C^{\alpha_1})\times \C^{N_{BL}}
\qquad \textup{with }N_{BL}:=|\{\beta:(*)\}|+1.
\nonumber
\end{eqnarray}
\end{lemma}

{\bf Proof:}
(a) The spectral numbers are well known \cite[13.3.4, p. 389]{AGV88}.
They also follow from corollary \ref{t7.17} and the proof
of theorem \ref{t9.6}.

(b) \eqref{9.28} follows from $\dim C^{\alpha_1}=\dim Ml_\zeta$
and $\Phi_m\not| b_2 \iff m\not|p$.
\eqref{9.29} follows from the values of $b_j$ 
in table \eqref{5.1}.

(c) The spectral numbers and the numbers $\beta$ with 
$C^\beta\neq\{0\}$ give for each $\LL_0\in D_{BL}$
\begin{eqnarray}\label{9.32}
\LL_0 = \C\cdot \sigma_1 \oplus \LL_0\cap 
\bigoplus_{\beta: \,\alpha_2\leq \beta\leq\beta_2}C^\beta
\oplus V^{>\beta_2}
\end{eqnarray}
where 
\begin{eqnarray}\label{9.33}
\alpha(\sigma_1)=\alpha_1,\quad
\sigma_1\in C^{\alpha_1}\oplus\bigoplus_{\beta:(*)}C^\beta 
\oplus C^{\beta_2}.
\end{eqnarray}
Here observe that for $\beta$ with $\alpha_2\leq\beta<\beta_2$
and $C^\beta\neq\{0\}$, the space $C^\beta$ is one-dimensional
and is generated by the principal part of a section in $\LL_0$.

If $m\not|p$ then $\dim C^{\beta_2}=1$ and $C^{\beta_2}$ 
is not generated by the principal part of a section in $\LL_0$.
If $m|p$ then $\dim C^{\beta_2}=2$ and the one-dimensional
subspace $\{v\in C^{\beta_2}\, |\, 
K_f^{(-2)}(v,s(\sigma_1,\alpha_1))=0\}\subset C^{\beta_2}$
is in $\LL_0$, because then $\beta_2$ is a spectral number
with multiplicity 1. And then the principal part 
$s(\sigma_1,\alpha_1)$ must be compatible with a polarized
Hodge structure of weight $2$ on $H^\infty_\zeta
\oplus H^\infty_{\oooo\zeta}$. 
This amounts to 
$[\C\cdot s(\sigma_1,\alpha_1)]\in\HH(C^{\alpha_1})$.
Especially then
\begin{eqnarray}\label{9.34}
C^{\beta_2}= \C\cdot \psi_{\beta_2}(\oooo{\psi_{\alpha_1}^{-1}
s(\sigma_1,\alpha_1)})\oplus 
\{v\in C^{\beta_2}\, |\, K_f^{(-2)}(v,s(\sigma_1,\alpha_1))=0\},
\end{eqnarray}
and $\sigma_1$ can be chosen with 
\begin{eqnarray}\label{9.35}
\alpha(\sigma_1)=\alpha_1,\quad 
\sigma_1\in C^{\alpha_1}\oplus\bigoplus_{\beta:(*)}C^\beta 
\oplus \C\cdot \psi_{\beta_2}(\oooo{\psi_{\alpha_1}^{-1}
s(\sigma_1,\alpha_1)}).
\end{eqnarray}
$\sigma_1$ is (up to rescaling) 
uniquely determined by \eqref{9.33} if $m\not|p$
and by \eqref{9.35} if $m|p$. And it can be chosen freely
with \eqref{9.33} respectively with \eqref{9.35} and 
$[\C\cdot s(\sigma_1,\alpha_1)]\in\HH(C^{\alpha_1})$.
The condition $(\delta)$  $K_f^{(-2)}(\LL_0,\LL_0)=0$ on $D_{BL}$
directly before theorem \ref{t7.11} implies that
$\LL_0\cap \bigoplus_{\alpha_2\leq \beta\leq \beta_2}C^\beta$
is uniquely determined by $\sigma_1$. 
Therefore $\LL_0$ is uniquely determined by $\sigma_1$. 
Therefore $D_{BL}$ is as stated in \eqref{9.30} and \eqref{9.31}.
\hfill$\Box$

\begin{remarks}\label{t9.5}
(i) All the normal forms in table \eqref{9.1} except $W_{1,p}^\sharp$
are Newton nondegenerate. But also the normal form
$f_p(x,y,\www z)$ for $W_{1,p}^\sharp$ in table \eqref{9.1}
can be made easily Newton nondegenerate with the coordinate change
$\www z = z + i(x^2+y^3)$. Then
\begin{eqnarray}\label{9.36}
f_p(x,y,z+i(x^2+y^3))&=& z^2+2ix^2z+2iy^3z \\
&+&\left\{\begin{array}{ll}
(t_1+t_2y)xy^{4+q} & \textup{ if }p=2q-1,\\
(t_1+t_2y)x^2y^{3+q} & \textup{ if }p=2q.
\end{array}\right.  \nonumber
\end{eqnarray}

(ii) The Newton boundaries of the normal forms in table \eqref{9.1}
except for $W_{1,p}^\sharp$ and of the normal form in \eqref{9.36}
for $W_{1,p}^\sharp$ have each two compact $n$-dimensional faces
$\sigma_1$ and $\sigma_2$. The following table lists
the corresponding linear forms $l_{\sigma_j}$ 
and the value $s(f)$ from corollary \ref{t7.17}.
A linear form is encoded by the values
$(l_{\sigma_j}(x),l_{\sigma_j}(y),l_{\sigma_j}(z))$.

\begin{eqnarray}\label{9.37}
\begin{array}{llll}
W_{1,p}^\sharp & \sigma_1: \frac{1}{12}(3,2,6)  &
\sigma_2: \frac{1}{12+p}(3,2,6+p) &
\frac{5}{12+p} \\
S_{1,p}^\sharp & \sigma_1: \frac{1}{10}(3,2,4)  &
\sigma_2: \frac{1}{10+p}(3,2,4+p) &
\frac{5}{10+p} \\
U_{1,p} & \sigma_1: \frac{1}{9}(3,2,3)  &
\sigma_2: \frac{1}{9+p}(3+p,2,3) &
\frac{5}{9+p} \\
E_{3,p} & \sigma_1: \frac{1}{18}(6,2,9)  &
\sigma_2: \frac{1}{2(9+p)}(6+p,2,9+p) &
\frac{4}{9} \\
Z_{1,p} & \sigma_1: \frac{1}{14}(4,2,7)  &
\sigma_2: \frac{1}{2(7+p)}(4+p,2,7+p) &
\frac{3}{7} \\
Q_{2,p} & \sigma_1: \frac{1}{12}(4,2,5)  &
\sigma_2: \frac{1}{2(6+p)}(4+p,2,5+p) &
\frac{1}{2} \\
W_{1,p} & \sigma_1: \frac{1}{12}(3,2,6)  &
\sigma_2: \frac{1}{2(6+p)}(3+p,2,6+p) &
\frac{5}{12} \\
S_{1,p} & \sigma_1: \frac{1}{10}(3,2,4)  &
\sigma_2: \frac{1}{2(5+p)}(3+p,2,4+p) &
\frac{1}{2} 
\end{array}
\end{eqnarray}
\end{remarks}

\begin{theorem}\label{t9.6}
Consider the normal form in \eqref{9.36} for $W_{1,p}^\sharp$
and the normal forms in table \eqref{9.1} for the other 
seven series. Recall the notation $\omega_0:=dxdydz$
from remark \ref{t7.14} (v). Define
\begin{eqnarray*}
b_1 &:=& s(\omega_0,\alpha_1)(1,0) \in C^{\alpha_1},\\
b_2 &:=& s(\omega_0,\beta_1)(1,0) \in C^{\beta_1},\\
b_3 &:=& s(y\omega_0,\alpha_2)(1,0) \in C^{\alpha_2},\\
b_4 &:=& s(y\omega_0,\beta_2)(1,0) \in C^{\beta_2}.
\end{eqnarray*}
If $m|p$, choose  $b_5\in C^{\beta_2}$ with
$\C\cdot b_5 =\{v\in C^{\beta_2}\, |\, 
K_f^{(-2)}(b_1,v)=0\}$.

\medskip
(a) All $b_j\neq 0$. And $K_f^{(-2)}(b_1+b_2,b_3+b_4)=0$.
If $m|p$ then $C^{\beta_2}=\C\cdot b_4\oplus \C\cdot b_5$.

\medskip
(b) We write $t=(t_1,t_2)$. 
Recall the notation $\alpha(s[\omega]_0(t))=\min(\alpha\, |\, 
s(\omega,\alpha)(t)\neq 0)$ from remark \ref{7.14} (iv).
\begin{eqnarray}\label{9.38}
\alpha(s[\omega_0]_0(t))
&=&\alpha_1, \\
s(\omega_0,\alpha_1)(t)
&=& b_1, \label{9.39} \\
s(\omega_0,\beta)(t)
&=& 0\quad\textup{for }\alpha_1<\beta<\beta_1,\label{9.40} \\
s(\omega_0,\beta_1)(t)
&=& t_1^{1/m_2}\cdot b_2,\label{9.41} \\
s(\omega_0,\alpha_2)(t)
&=& \frac{t_2}{t_1}\cdot\frac{-1}{m_2}\cdot 
t_1^{-1/m_2}\cdot b_3 + s(\omega,\alpha_2)(t_1,0),\label{9.42}\\
s(\omega_0,\beta_2)(t)
&&\left\{ \begin{array}{ll} =s(\omega_0,\beta_2)(t_1,0)
&\textup{ if }m\not|p,\\
\in s(\omega_0,\beta_2)(t_1,0)+\C\cdot b_5
&\textup{ if }m|p,\end{array}\right. \label{9.43}
\end{eqnarray}
\noindent
with
\begin{eqnarray}\label{9.44}
\begin{array}{l|l|l}
 & s(\omega_0,\alpha_2)(t_1,0) & s(\omega_0,\beta_2)(t_1,0) \\
 \hline 
(r_I=2\, \&\,  p\geq 2)\textup{ or } & & \\
(r_I=1\, \&\,  p\geq 3) \textup{ or }U_{1,2} & 0 & 0 \\ \hline 
W_{1,1}^\sharp, S_{1,1}^\sharp,U_{1,1} & 
c_1\cdot t_1^{2-1/m_2}\cdot b_3 &
c_2\cdot t_1^2\cdot b_4 \\ \hline 
W_{1,2}^\sharp, S_{1,2}^\sharp, E_{3,1} & & \\
Z_{1,1},Q_{2,1},W_{1,1},S_{1,1} & 
c_1\cdot t_1^{1-1/m_2}\cdot b_3 &
c_2\cdot t_1\cdot b_4
\end{array}
\end{eqnarray}
for some values $c_1,c_2\in\C$.

\begin{eqnarray}\label{9.45}
\alpha(s[y\omega_0]_0(t)
&=& \alpha_2, \\
s(y\omega_0,\alpha_2)(t) 
&=& t_1^{-1/m_2}\cdot b_3,\label{9.46} \\
s(y\omega_0,\beta_2)(t)
&&\left\{ \begin{array}{ll}
=b_4 & \textup{ if }m\not|p\textup{ or }t_2=0,\\
\in b_4+\C\cdot b_5 & \textup{ if }m|p,
\end{array}\right . \label{9.47}
\end{eqnarray}
\begin{eqnarray}\label{9.48}
s(\sigma,\beta_2)(t)  \left\{
\begin{array}{ll}
=0 & \textup{ if } m\not|p \\ 
\in \C\cdot b_5 & \textup{ if }m|p
\end{array}\right. 
\end{eqnarray}
for $\sigma\in H_0''(f_t)$ with 
$\alpha(\sigma)>\alpha_2$.

\medskip
(c) In the five series with $r_I=2$ (see table \eqref{t5.1})
for $b\in \Z_{\geq 0}$
\begin{eqnarray}\label{9.49}
\alpha(s[y^{b+1}\omega_0]_0(t))
&=& \alpha_2+\frac{2b}{m_2} = \frac{2b+1}{m_2},\\
s(y^{b+1}\omega_0,\frac{2b+1}{m_2})(t)
&=& t_1^{-(2b+1)/m_2}\cdot s(y^{b+1}\omega_0,
\frac{2b+1}{m_2})(1,0). \hspace*{1cm}\label{9.50}
\end{eqnarray}
Especially, if $p=mr$ then 
$\frac{2r+1}{m_2}=\frac{1}{m}=\beta_2$, 
$b_5$ can be chosen as $b_5=s(y^{r+1}\omega_0,\beta_2)(1,0)$,
and 
\begin{eqnarray}\label{9.51}
s(y^{r+1}\omega_0,\beta_2)(t)
&=& t_1^{-1/m}\cdot b_5.
\end{eqnarray}

\medskip
(d) In the three subseries $W_{1,12r}^\sharp,
S_{1,10r}^\sharp, U_{1,9r}$ (i.e. the subseries with
$r_I=1$ and $m|p$), $b_5$ can be chosen such that
$b_5$ and $\omega$ in the following table \eqref{9.54}
satisfy
\begin{eqnarray}\label{9.52}
\alpha(s[\omega]_0(t)) 
&=& \beta_2 =\frac{1}{m},\\
s(\omega,\beta_2+1)(t)
&=& t_1^{-1/m}\cdot b_5.
\end{eqnarray}
\label{9.53}
\begin{eqnarray}\label{9.54}
\begin{array}{r|l}
 & \omega \\ \hline 
W_{1,12+24r}^\sharp,S_{1,10+20r}^\sharp & xy^r\omega_0\\
U_{1,9+18r} & y^rz\omega_0 \\
W_{1,24r}^\sharp,S_{1,20r}^\sharp,U_{1,18r} & 
y^{r+1}\omega_0
\end{array}
\end{eqnarray}
\end{theorem}

{\bf Proof:}
(a) Observe $\nu(\omega_0)-1=\alpha_1<s(f)$ and
$\nu(y\omega_0)-1=\alpha_2<s(f)$.
This, theorem \ref{t7.16} and corollary \ref{t7.17}
show \eqref{9.38}, \eqref{9.45}, $b_1\neq 0$ and 
$b_3\neq 0$. $b_2\neq 0$ will be shown below.
\eqref{9.40} (which will also be shown below) and
$K_f^{(-2)}(H_0''(f_t),H_0''(f_t))=0$
give especially
$$0=K_f^{(-2)}(s[\omega_0]_0(1,0),s[y\omega_0]_0(1,0))
=K_f^{(-2)}(b_1+b_2,b_3+b_4).$$
As $K_f^{(-2)}(b_2,b_3)\neq 0$, also $K_f^{(-2)}(b_1,b_4)\neq 0$
and $b_4\neq 0$ and in the case $m|p$
$C^{\beta_2}=\C\cdot b_4\oplus \C\cdot b_5$.

\medskip
(b)--(d) We restrict to the series $E_{3,p}$. 
The calculations for the series 
$Z_{1,p},Q_{2,p},W_{1,p}$ and $S_{1,p}$ are very similar.
The calculations for the series 
$W_{1,p}^\sharp, S_{1,p}^\sharp$ and $U_{1,p}$ are similar,
but require more case discussions.

The two compact faces $\sigma_1$ and $\sigma_2$ (remark \ref{t9.5})
of the Newton boundary give rise to the following two relations
\begin{eqnarray}\label{9.55}
\frac{1}{3}xf_x+\frac{1}{9}yf_y+\frac{1}{2}zf_z 
-\frac{p}{9}t_1y^{9+p}-\frac{p+1}{9}t_2y^{10+p} &=& f,\\
\frac{6+p}{2(9+p)}xf_x+\frac{2}{2(9+p)}yf_y+\frac{1}{2}zf_z &&\nonumber\\
-\frac{p}{2(9+p)}x^3-\frac{1}{9+p}t_2y^{10+p} &=& f.\label{9.56}
\end{eqnarray}
These relations and \eqref{7.24} give the following two values
for $\ppp_\tau\tau s[x^ay^b\omega_0]_0(t)$:
\begin{eqnarray}
&&\ppp_\tau\tau s[x^ay^b\omega_0]_0(t) \nonumber\\
&=& l_{\sigma_1}(a+1,b+1,1)\cdot s[x^ay^b\omega_0]_0(t)\label{9.57}\\
&&-\frac{p}{9}t_1\ppp_\tau s[x^ay^{b+9+p}\omega_0]_0(t)
-\frac{p+1}{9}t_2\ppp_\tau s[x^ay^{b+10+p}\omega_0]_0(t),\nonumber\\
&=& l_{\sigma_2}(a+1,b+1,1)\cdot s[x^ay^b\omega_0]_0(t)\label{9.58}\\
&&-\frac{p}{2(9+p)}\ppp_\tau s[x^{a+3}y^b\omega_0]_0(t)\
-\frac{1}{9+p}t_2\ppp_\tau s[x^ay^{b+10+p}\omega_0]_0(t).\nonumber
\end{eqnarray}
This gives for any $\beta$ with $\dim C^\beta\neq 0$ 
\begin{eqnarray}
&&(\beta+1-l_{\sigma_1}(a+1,b+1,1))s(x^ay^b\omega_0,\beta)(t)\nonumber \\
&=& -\frac{p}{9}t_1\ppp_\tau s(x^ay^{b+9+p}\omega_0,\beta+1)(t)\nonumber\\
&&- \frac{p+1}{9}t_2\ppp_\tau s(x^ay^{b+10+p}\omega_0,\beta+1)(t),
\label{9.59}\\
&&(\beta+1-l_{\sigma_2}(a+1,b+1,1))s(x^ay^b\omega_0,\beta)(t)\nonumber \\
&=& -\frac{p}{2(9+p)}\ppp_\tau s(x^{a+3}y^b\omega_0,\beta+1)(t
\nonumber\\
&&- \frac{1}{9+p}t_2\ppp_\tau s(x^ay^{b+10+p}\omega_0,\beta+1)(t).
\label{9.60}
\end{eqnarray}
Furthermore, \eqref{7.42} gives
\begin{eqnarray}\label{9.61}
\ppp_{t_1} s[x^ay^b\omega_0]_0(t) &=& 
(-\ppp_\tau) s[x^ay^{b+9+p}\omega_0]_0(t),\\
\ppp_{t_2} s[x^ay^b\omega_0]_0(t) &=& 
(-\ppp_\tau) s[x^ay^{b+10+p}\omega_0]_0(t)\nonumber\\
&=& \ppp_{t_1} s[x^ay^{b+1}\omega_0]_0(t).\label{9.62}
\end{eqnarray}
\eqref{9.59}--\eqref{9.62} give
\begin{eqnarray}
\Bigl(\frac{p}{9}t_1\ppp_{t_1} + \frac{p+1}{9}t_2\ppp_{t_2}
-(\beta+1)+l_{\sigma_1}(a+1,b+1,1)\Bigr)\nonumber \\
s(x^ay^b\omega_0,\beta)(t)=0, \label{9.63}\\
\Bigl(\frac{1}{9+p}t_2\ppp_{t_2}-(\beta+1)+l_{\sigma_2}(a+1,b+1,1)\Bigr)
s(x^ay^b\omega_0,\beta)(t) \nonumber \\
= \frac{p}{2(9+p)}\ppp_\tau s(x^{a+3}y^b\omega_0,\beta+1)(t).
\label{9.64}
\end{eqnarray}
\eqref{9.63} gives for $t_2=0$
\begin{eqnarray}\label{9.65}
s(x^ay^b\omega_0,\beta)(t_1,0) =
t_1^{\frac{9}{p}(\beta+1-l_{\sigma_1}(a+1,b+1,1))}\cdot
s(x^ay^b\omega_0,\beta)(1,0).
\end{eqnarray}
The following eight equations are special cases of \eqref{9.65}.
\begin{eqnarray}\label{9.66}
s(\omega_0,\alpha_1)(t_1,0) 
&=& b_1, \\
s(\omega_0,\beta_1)(t_1,0) 
&=& t_1^{1/m_2}\cdot b_2, \label{9.67} \\
s(\omega_0,\alpha_2)(t_1,0) 
&=& t_1^{-1/m_2+1/p}\cdot s(\omega_0,\alpha_2)(1,0), \label{9.68} \\
s(\omega_0,\beta_2)(t_1,0) 
&=& t_1^{1/p}\cdot s(\omega_0,\beta_2)(1,0), \label{9.69} 
\end{eqnarray}
\begin{eqnarray}
s(y^{b+1}\omega_0,\frac{2b+1}{m_2})(t_1,0) 
&=& t_1^{-(2b+1)/m_2}\cdot 
s(y^{b+1}\omega_0,\frac{2b+1}{m_2})(1,0), \label{9.70} \\
s(y\omega_0,\alpha_2)(t_1,0) 
&=& t_1^{-\alpha_2}\cdot b_3
=t_1^{-1/m_2}\cdot b_3 , \label{9.71} \\
s(y\omega_0,\beta_2)(t_1,0) 
&=& b_4, \label{9.72} \\
s(y^{r+1}\omega_0,\beta_2)(t_1,0) 
&=& t_1^{-1/m_2}\cdot s(y^{r+1}\omega_0,\beta_2)(1,0)
\quad\textup{if }p=18r. \hspace*{1cm}\label{9.73} 
\end{eqnarray}

\medskip
{\bf Claim:}
Fix some $b\in\Z_{\geq 0}$.
\begin{list}{}{}
\item[(i)] 
$\nu(y^{b+1}\omega_0)=\alpha_2+\frac{b}{9+p}=\frac{2b+1}{m_2}$.
\item[(ii)]
Any $(n+1)$-form $df\land d\eta$ which contains $y^{b+1}\omega_0$
as a summand, contains a summand $g\cdot \omega_0$ with 
$g$ a monomial (times a nonzero scalar) with
$\nu(g\cdot \omega_0)\leq \nu(y^{b+1}\omega_0)$.
\item[(iii)]
$\oooo{\nu}(y^{b+1}\omega_0) = \frac{2b+1}{m_2}$.
\end{list}

\medskip
{\bf Proof of the claim:}
(i) Trivial.  (iii) follows from (i) and (ii).

(iii) The only monomial differential $(n-1)$-forms $\eta$
such that $df\land d\eta$ contains $f_y\cdot y^c\cdot\omega_0$
are $\eta_1=-xy^c dz$ and $\eta_2=y^cz dx$, and
\begin{eqnarray*}
df\land d\eta_1 &=& f_y\cdot y^c\cdot \omega_0 
-f_x\cdot c\cdot xy^{c-1}\cdot \omega_0,\\
df\land d\eta_2 &=& f_y\cdot y^c\cdot \omega_0 
-f_z\cdot c\cdot y^{c-1}z\cdot \omega_0.
\end{eqnarray*}
These $(n+1)$-forms contain $(3-2c)x^2y^{c+2}\omega_0$ 
respectively $3x^2y^{c+2}\omega_0$, and
\begin{eqnarray*}
\nu(x^2y^{c+2}\omega_0) \leq \nu(y^{c+8+p}\omega_0).
\end{eqnarray*}
\hfill $(\Box)$

\medskip
The claim and theorem \ref{t7.16} imply
\begin{eqnarray}\label{9.74}
\alpha(s[y^{b+1}\omega_0]_0(t)) &=& 
\frac{2b+1}{m_2},\\
s(y^{b+1}\omega_0,\frac{2b+1}{m_2})(t) &\neq & 0.\label{9.75}
\end{eqnarray}
Especially, $b_3\neq 0$, and if $p=18r$ also
$s(y^{r+1}\omega_0,\beta_2)(t)\neq 0$.
In this case $p=18r$, the vanishing
$$K_f^{(-2)}(s[\omega_0]_0(1,0),s[y^{r+1}\omega_0]_0(1,0))=0$$
gives $K_f^{(-2)}(b_1,s(y^{r+1}\omega_0,\beta_2)(1,0))=0$.
Therefore in this case we can choose 
$b_5=s(y^{r+1}\omega_0,\beta_2)(1,0)$.

The elementary sections $s(y^{b+1}\omega_0,\frac{2b+1}{m_2})(t)$
are independent of $t_2$ because \eqref{9.62} gives
\begin{eqnarray*}
\ppp_{t_2}s(y^{b+1}\omega_0,\frac{2b+1}{m_2})(t)
= \ppp_{t_1}s(y^{b+2}\omega_0,\frac{2b+1}{m_2})(t)=0.
\end{eqnarray*}
Now part (c), i.e. \eqref{9.49}--\eqref{9.51},
and \eqref{9.46} are proved.

\eqref{9.62} gives also
\begin{eqnarray}\label{9.76}
\ppp_{t_2}s[\omega_0]_0(t)
&=& \ppp_{t_1}s[y\omega_0]_0(t),\\
\textup{so }s(\omega_0,\beta)(t)
&=& s(\omega_0,\beta)(t_1,0)\qquad\textup{for }
\alpha_1\leq \beta < \alpha_2.\nonumber
\end{eqnarray}
With \eqref{9.66} and \eqref{9.59} and \eqref{9.75} we obtain
\begin{eqnarray*}
s(\omega_0,\beta)(t_1,0) = \left\{\begin{array}{ll}
b_1 & \textup{ if }\beta=\alpha_1,\\
\frac{-p}{9(\beta-\alpha_1)}t_1\ppp_\tau 
s(y^{9+p}\omega_0,\beta+1)(t_1,0)=0 & 
\textup{ if }\alpha_1<\beta<\beta_1,\\
\frac{-p}{9(\beta_1-\alpha_1)}t_1\ppp_\tau 
s(y^{9+p}\omega_0,\beta_1+1)(t_1,0)\neq 0 & 
\textup{ if }\beta=\beta_1. \end{array}\right. 
\end{eqnarray*}
This gives $b_2\neq 0$ and 
(together with \eqref{9.66} and \eqref{9.67})
\eqref{9.39}--\eqref{9.41}. 

The argument in the proof of part (a) with 
$K_f^{(-2)}(H_0''(f_t),H_0''(f_t))=0$
gives $b_4\neq 0$ and \eqref{9.47} and \eqref{9.48}.

It rests to show \eqref{9.42}--\eqref{9.44}.
From \eqref{9.76}, \eqref{9.46} and \eqref{9.47} 
we obtain
\begin{eqnarray*}
\ppp_{t_2}s(\omega_0,\alpha_2)(t)
&=& \ppp_{t_1} s(y\omega_0,\alpha_2)(t) 
= \ppp_{t_1}(t_1^{-1/m_2}\cdot b_3),\\
\ppp_{t_2}s(\omega_0,\beta_2)(t)
&=& \ppp_{t_1} s(y\omega_0,\beta_2)(t) 
\left\{\begin{array}{ll} = 0 & \textup{ if }m\not|p,\\
\in\C\cdot b_5 & \textup{ if }m|p,\end{array}\right.
\end{eqnarray*}
which gives \eqref{9.42} and \eqref{9.43}.

For \eqref{9.44} observe the following.
The sections
\begin{eqnarray*}
s(y\omega_0,\alpha_2)(t_1,0)&=& t_1^{-1/m_2}\cdot b_3,\\
s(y\omega_0,\beta_2)(t_1,0)&=& b_4,\\
\textup{and in the case }m|p\quad 
s(y^{r+1}\omega_0,\beta_2)(t_1,0)&=& t_1^{-1/m}\cdot b_5
\end{eqnarray*}
are univalued nowhere vanishing sections in the bundles
$\bigcup_{t_1\in T}C^{\alpha_2}(t_1,0)$ and
$\bigcup_{t_1\in T}C^{\beta_2}(t_1,0)$,
and they generate these bundles.
Also $s(\omega_0,\alpha_2)(t_1,0)$ and 
$s(\omega_0,\beta_2)(t_1,0)$ are univalued sections in
these bundles. \eqref{9.68} and \eqref{9.69}
show for $p\geq 2$ that they are everywhere vanishing.
For $p=1$ they give the statement for
$E_{3,1}$ in the last line of table \eqref{9.44}.
This finishes the proof of the parts (b) and (c)
for the series $E_{3,p}$. \hfill$\Box$

\bigskip
{\bf Proof of $M_T^{m_2}=\id$:}

By theorem \ref{t9.6}, the following sections in the bundles
$\bigcup_{t_1\in T}C^{\beta}(t_1,0)$ for $\beta$ as 
in table \eqref{9.77} are univalued nowhere vanishing 
sections and generate these bundles
(in the case $\beta=\alpha_1$ only if $m\not|p$).
\begin{eqnarray}\label{9.77}
\begin{array}{l|l|l|l|l|l}
\textup{section} & b_1 & t_1^{1/m_2}\cdot b_2 & 
t_1^{-1/m_2}\cdot b_3 & b_4 & t_1^{-1/m}\cdot b_5
\textup{ if }m|p \\ \hline 
\beta & \alpha_1 & \beta_1 & \alpha_2 & \beta_2 & \beta_2 \\ \hline
\textup{eigenvalue of} & 1 & e^{-2\pi i/m_2} & 
e^{2\pi i /m_2} & 1 & e^{2\pi i/m}\\
M_T\textup{ on }\C\cdot b_j
\end{array}
\end{eqnarray}
Therefore $b_1$ and $b_4$ are univalued, and $b_2$ and $b_3$ 
(and $b_5$ if $m|p$) are multivalued flat sections
with eigenvalues of $M_T$ as in the table.
Thus $M_T^{m_2}$ is on $C^{\alpha_1},C^{\beta_1},C^{\alpha_2},
C^{\beta_2},Ml_\zeta$ and $Ml_{e^{2\pi i/m_2}}$
the identity. We will show that it is the identity on all of $Ml$.

Consider firstly the case $m\not|p $. 
Then by \eqref{9.24} $M_T^{m_2}$ is in 
\begin{eqnarray}\label{9.78}
&&\{\id\} \qquad \textup{in the cases }
W_{1,2q-1}^\sharp, S_{1,2q-1}^\sharp, U_{1,2q}, E_{3,p}, Z_{1,p},\\
&&\{\id,\id|_{B_1}\times (-M_h^{m_2/2})|_{B_2}\} 
\quad \textup{in the cases }
W_{1,2q}^\sharp, S_{1,2q}^\sharp, U_{1,2q-1},\nonumber \\
&&\{\id,(-M_h^{m/2})|_{B_1}\times \id|_{B_2}\} 
\qquad \textup{in the cases }
Q_{2,p}, W_{1,p}, S_{1,p}.\nonumber
\end{eqnarray}
On the other hand, in the cases in the second and third line
of \eqref{9.78}, $m_2=m+r_Ip$ is even, and $M_T$ itself is 
in $G_\Z$ which is given by \eqref{5.6} in theorem \ref{t5.1}.
Thus $M_T^{m_2}=\id$ also in the second and third line 
of \eqref{9.78}.

Consider secondly the case $m|p$, so $p=mr$. 
By \eqref{5.9} in theorem \ref{t5.1},
$M_T^{m_2}=\varepsilon\cdot M_h^k$ for some 
$\varepsilon\in\{\pm 1\}$ and some $k\in\Z$.
Then $\varepsilon\cdot \zeta^k=1$ and
$\varepsilon\cdot e^{2\pi i k/m_2}=1$.
If $\varepsilon=1$, then the two conditions boil down
to $m|k$ and $m_2|k$, so to $m_2|k$. 
Then $M_T^{m_2}=\id$. 
If $\varepsilon=-1$, we will come below to a 
contradiction. Then the two conditions require $m$ even
and $m_2$ even. 

For each eigenvalue $\lambda$ of $M_h$ on $Ml$ with 
$\dim Ml_\lambda =1$, an eigenvector in 
$Ml_{\lambda,\Z[\lambda]}$ exists.
Then $M_T$ has an eigenvalue in $\Eiw(\lambda)$
on this eigenvector, and $M_T^{m_2}$ has the 
eigenvalue 1 on this eigenvector.
Here $m_2$ even is used.
Therefore $M_T^{m_2}=\id$ on $Ml_\lambda$
for each 
$$\lambda\in\{\zeta,e^{2\pi i /m_2}\}
\cup\{\www\lambda\, |\, \dim Ml_{\www\lambda}=1\}.$$
Comparison with table \eqref{5.1} shows that no 
$k\in\Z$ with $-\lambda^k=1$ for all these $\lambda$ exists.
This gives a contradiction. The case 
$\varepsilon=-1$ is impossible. 
$M_T^{m_2}=\id$ is proved in all cases.\hfill $\Box$

\bigskip
{\bf Proof} that $M_T$ has the eigenvalues $1$ and 
$\oooo{\zeta}$ on $Ml_\zeta$ and on $C^{\alpha_1}$:

By table \eqref{9.77}, $M_T$ has on $C^{\beta_2}$
and on $H^\infty_{e^{-2\pi i \beta_2}}=H^\infty_{\oooo\zeta}$
the eigenvalues $1$ and $\zeta$. 
As $Ml_\zeta$ is dual to $H^\infty_{\oooo\zeta}$ and
$H^\infty_{\zeta}$ is complex conjugate to 
$H^\infty_{\oooo\zeta}$, $M_t$ has on $Ml_\zeta$,
$H^\infty_{\zeta}=H^\infty_{e^{-2\pi i \alpha_1}}$ 
and $C^{\alpha_1}$ the eigenvalues $1$ and $\oooo\zeta$.
\hfill$\Box$

\bigskip
{\bf Proof of theorem \ref{t9.2} (a)+(b)+(d):}

(a) This follows immediately from \eqref{9.7}, \eqref{9.9}
and lemma \eqref{9.4} (c).

(b) All of this follows by carefully putting together the results
in theorem \ref{t9.6}. 
Here $v_1^0=b_1$, $v_2^0=b_2$, $v_4^0\in\C^*\cdot b_4$ suitable,
and the section in the brackets on the right hand side 
of \eqref{9.10} is
\begin{eqnarray}\label{9.79}
s[\omega_0]_0(t)  &+& 
\left(\frac{1}{m}\frac{t_2}{t_1} + 
\left\{\begin{matrix} 0 \\ -c_1\cdot t_1^2 \\ -c_1\cdot t_1 \end{matrix}
\right\}\right)\cdot
s[y\omega_0]_0(t)\\
&&\mod \bigoplus_{\alpha_2<\beta<\beta_2}C^\beta 
\oplus\C\cdot b_5 \oplus V^{>\beta_2}.\nonumber
\end{eqnarray}
The three cases in $\{...\}$ correspond to the three lines in
\eqref{9.44}. 
The linear combination is chosen such that it has no part in 
$C^{\alpha_2}$. This section and the fact
$K_f^{(-2)}(H_0''(f_t),H_0''(f_t))=0$
determine $H_0''(f_t)$. 
By table \eqref{9.77}, $M_T$ has on $v_1^0=b_4$ the eigenvalue 1.

\medskip
(c) Consider the coordinate change
\begin{eqnarray}\label{9.80}
\varphi:(\C^3,0)\to (\C^3,0),\quad (x,y,z)\mapsto (x,y,-z).
\end{eqnarray}
We treat the cases $U_{1,9+18r}$ and $U_{1,18r}$ separately.

{\bf The case $U_{1,9+18r}$:} Then $\varphi\in G^{smar,gen}_\RR
\subset G^{smar}$, and 
\begin{eqnarray}\label{9.81}
\varphi^*(\omega_0) = -\omega_0,\quad
\varphi^*(y^rz\omega_0)=y^rz\omega_0.
\end{eqnarray}
Now compare \eqref{9.39} and \eqref{9.54}.
$\varphi$ induces an automorphism 
$(\varphi)_{coh}$ on $C^{\alpha_1}$ and 
$C^{\beta_2}$ with 
\begin{eqnarray}\label{9.82}
(\varphi)_{coh}(b_1)=-b_1,\quad 
(\varphi)_{coh}(b_4)=-b_4,\quad
(\varphi)_{coh}(b_5)=b_5.
\end{eqnarray}
One can choose $g_3 = -M_T\circ (\varphi)_{hom}\in G^{mar}$.

{\bf The case $U_{1,18r}$:} Because of \eqref{9.54} and \eqref{9.77},
instead of \eqref{9.81} the identities
\begin{eqnarray}\label{9.83}
\varphi^*(\omega_0) = -\omega_0,\quad
\varphi^*(y^{r+1}\omega_0)=-y^{r+1}\omega_0
\end{eqnarray}
are relevant. Now $(\varphi)_{coh}$ is because of \eqref{9.15}
an isomorphism
\begin{eqnarray*}
H_0''(f_{(t_1,0)})\to H_0''(f_{(-t_1,0)}),\quad
C^{\beta_2}(t_1,0)\to C^{\beta_2}(-t_1,0).
\end{eqnarray*}
The composition
\begin{eqnarray*}
(-\id)\circ (\textup{math. pos. flat shift from }
C^{\beta_2}(-t_1,0)\textup{ to }C^{\beta_2}(t_1,0))
\circ (\varphi)_{coh}
\end{eqnarray*}
acts on $C^{\beta_2}(t_1,0)$ and has because of \eqref{9.76}
the eigenvectors $b_4$ and $b_5$ with the eigenvalues 1
and $e^{\pi i/9}$:
\begin{eqnarray*}
\begin{array}{cccc}
b_4 & t_1^{-1/9}b_5 & & C^{\beta_2}(t_1,0) \\
\downarrow & \downarrow & (\varphi)_{coh} & \downarrow \\
-b_4 & -(e^{-\pi i}t_1)^{-1/9}b_5 & & C^{\beta_2}(-t_1,0) \\
\downarrow & \downarrow & \textup{shift} & \downarrow \\
-b_4 & -e^{\pi i /9}t_1^{-1/9}b_5 & & C^{\beta_2}(t_1,0)  \\
\downarrow & \downarrow & -\id & \downarrow \\
b_4 & e^{\pi i/9}t_1^{-1/9}b_5 & & C^{\beta_2}(t_1,0)
\end{array}
\end{eqnarray*}
The corresponding composition 
\begin{eqnarray*}
(-\id)\circ (\textup{math. pos. flat shift from }
Ml(f_{(-t_1,0)})\textup{ to }Ml(f_{(t_1,0)}))
\circ (\varphi)_{hom}
\end{eqnarray*}
is in $G^{mar}$ and can be chosen as $g_3$. 
\hfill $\Box$

\section{Period maps and Torelli results for the 
quadrangle singularities}\label{c10}
\setcounter{equation}{0}

\noindent
In this section we will prove for the quadrangle singularities
the strong global Torelli conjecture \ref{t8.11} (a), 
the conjectures \ref{t8.6} (b) $-\id\notin G^{smar}$ and (a)
$G_\Z=G^{mar}$. 
The Torelli conjecture for the unmarked singularities
had been proved in \cite{He93} (and the proof had been
sketched in \cite{He95}).
The main new ingredient for the Torelli result for marked
singularities is a much stronger control of the group
$G_\Z$, in theorem \ref{t6.1}.
But we will also recall the old ingredients from \cite{He93},
the space $D_{BL}$ and a period map for which we
need calculations of the Gauss-Manin connection.

The six bimodal families of quadrangle singularities 
have as surface singularities the normal forms $f_{(t_1,t_2)}$
in table \eqref{10.1}. These are not the normal forms
in \cite[15.1]{AGV85}. We will justify the normal forms
and explain their properties after theorem \ref{t10.1}.
The parameters $(t_1,t_2)$ are in 
$T^{(5)}:=(\C-\{0,1\})\times \C$.
Table \eqref{10.1} lists additionally weights $(w_x,w_y,w_z)$
such that $f_{(t_1,0)}$ is quasihomogeneous of weighted
degree 1 and two numbers $m_0$ and $m_\infty$
We set $m_1:=m_0$.
Observe $w_y=\frac{2}{m}<w_x\leq w_z$.

\begin{eqnarray}\label{10.1}
\begin{array}{lllll}
 & & (w_x,w_y,w_z) & m_0 & m_\infty \\ \hline 
W_{1,0}
& x^4+(4t_1-2)x^2y^3+y^6+t_2x^2y^4+z^2
& (\frac{1}{4},\frac{1}{6},\frac{1}{2}) & 12 & 6 \\
S_{1,0}
& x^2z+y^3z+yz^2+t_1x^2y^2+t_2x^2y^3
& (\frac{3}{10},\frac{2}{10},\frac{4}{10}) & 10 & 5 \\
U_{1,0}
& xz(x-z)+y^3(x-t_1z)+t_2y^4z
& (\frac{1}{3},\frac{2}{9},\frac{1}{3}) & 9 & 9 \\
E_{3,0}
& x(x-y^3)(x-t_1y^3)+t_2x^2y^4+z^2
& (\frac{1}{3},\frac{1}{9},\frac{1}{2}) & 9 & 9 \\
Z_{1,0}
& xy(x-y^2)(x-t_1y^2)+t_2x^2y^4+z^2
& (\frac{2}{7},\frac{1}{7},\frac{1}{2}) & 7 & 7 \\
Q_{2,0}
& x(x-y^2)(x-t_1y^2)+yz^2+t_2xz^2 
& (\frac{1}{3},\frac{1}{6},\frac{5}{12}) & 6 & 6 
\end{array}
\end{eqnarray}

Recall that table \eqref{6.1} lists for these singularities
the Milnor number $\mu$, the characteristic polynomials
$b_j$, $j\geq 1$, of the monodromy on the Orlik blocks
$B_j$ in theorem \ref{t5.1}, the order $m$ of the monodromy
and the index $r_I$.

For each 2-parameter family in table \eqref{10.1},
we choose $f_0:=f_{(i,0)}$ as reference singularity.
And as in section \ref{c9},
$M_\mu^{mar}$, $(M_\mu^{mar})^0$, $G_\Z$, $G^{mar}$, $Ml$, 
$H^\infty$ and $C^\alpha$ mean the objects for $f_0$.
As always, $\zeta:=e^{2\pi i /m}$.

We will construct branched coverings $c^{(2)}$ and $c^{(6)}$
and unbranched coverings $c^{(1)}$ and $c^{(5)}$ as follows.

\begin{eqnarray}\label{10.2}
\begin{xy}
  \xymatrix{
      T^{(3)}  \ar[d]_{c^{(1)}}  &\subset &  T^{(4)}:=\H \ar[d]^{c^{(2)}} \\
      T^{(1)}:=\C-\{0,1\}        &\subset  &   T^{(2)}:=\P^1\C   
  }
\end{xy}
\end{eqnarray}

\begin{eqnarray*}
\begin{xy}
  \xymatrix{
     T^{(7)}:=T^{(3)}\times\C  \ar[d]_{c^{(5)}:=c^{(1)}\times\id} &\subset &  T^{(8)}:=T^{(4)}\times \C \ar[d]^{c^{(6)}:=c^{(2)}\times\id}\\
     T^{(5)}=T^{(1)}\times\C         &\subset  &  T^{(6)}:=T^{(2)}\times\C  
  }
\end{xy}
\end{eqnarray*}

Let $\Gamma\subset PGL(2,\R)$ be a triangle group of type 
$(\frac{1}{m_0},\frac{1}{m_1},\frac{1}{m_\infty})$.
The quotient $\H/\Gamma$ is an orbifold with three orbifold
points of orders $m_0,m_1$ and $m_\infty$. They are the 
images of the elliptic fixed points of $\Gamma$ on 
$T^{(4)}=\H$ of orders $m_0,m_1$ and $m_\infty$.
As a manifold $\H/\Gamma\cong \P^1\C$.
Choose coordinates on $\H/\Gamma$ such that $0$ and $1$
are orbifold points of order $m_0=m_1$ and $\infty$
is an orbifold point of order $m_\infty$.
Denote by
\begin{eqnarray}\label{10.3}
c^{(2)}:T^{(4)}=\H\to T^{(2)}=\P^1\C
\end{eqnarray}
the quotient map. It is a branched covering.
Denote
\begin{eqnarray}
T^{(3)}&:=&T^{(4)}-(c^{(2)})^{-1}(\{0,1,\infty\}),\nonumber\\
\quad\textup{and}\quad
c^{(1)}&:=&c^{(2)}|_{T^{(4)}}:T^{(3)}\to T^{(1)}.\label{10.4}
\end{eqnarray}
It is a covering.

\begin{theorem}\label{t10.1}
Consider a bimodal family of quadrangle surface 
singularities in table \eqref{10.1}.

(a) There are canonical isomorphisms
\begin{eqnarray}\label{10.5}
T^{(7)}\to (M_\mu^{smar})^0\to (M_\mu^{mar})^0.
\end{eqnarray}

\medskip
(b) $-\id\notin G^{smar}$, where $G^{smar}$ is the
group for the singularities of multiplicity $\geq 3$,
namely the curve singularities
$W_{1,0}, E_{3,0}, Z_{1,0}$
and the surface singularities $S_{1,0}, U_{1,0},
Q_{2,0}$. So, conjecture \ref{t8.6} (b) is true.

\medskip
(c) $G_\Z=G^{mar}$. So, $M_\mu^{mar}=(M_\mu^{mar})^0$, and
conjecture \ref{t8.6} (a) is true.

\medskip
(d) The period map $BL:M_\mu^{mar}\to D_{BL}$ is an embedding.
So, the strong global Torelli conjecture \ref{t8.11} (a) is true.
\end{theorem}

The Torelli result for unmarked singularities
(the period map $M_\mu^{mar}/G_\Z\to D_{BL}/G_\Z$ 
is an embedding)
was proved already in \cite{He93}, and also that
there is a well defined period map $T^{(7)}\to D_{BL}$
and that it is an embedding. But we prefer to give an
independent account and recover these results.
The hardest part is in any case new. 
It is the precise control of $G_\Z$ in theorem \ref{t6.1}.

First we discuss the normal forms in table \eqref{10.1}
and the right equivalence classes in them.

Each bimodal family of quadrangle surface singularities
contains a 1-parameter subfamily of quasihomogeneous
singularities. The exceptional set of the minimal
good resolution of such a singularity 
consists of 5 smooth rational curves.
One, the central curve, intersects each of the other 4,
the branches, in one point. The right equivalence class
of one quasihomogeneous surface singularity is determined
by the central curve with the 4 intersection points and
the self intersection numbers of the 4 branches.
Table \eqref{10.6} lists these self intersection numbers.
\begin{eqnarray}\label{10.6}
\begin{array}{llllll}
W_{1,0} & S_{1,0} & U_{1,0} & E_{3,0} & Z_{1,0} & Q_{2,0} \\
(2,2,3,3) & (2,2,3,4) & (2,3,3,3) & (2,2,2,3) & (2,2,2,4) &
(2,2,2,5) 
\end{array}
\end{eqnarray}

In table \eqref{10.1}, the singularities with $t_2=0$ are
quasihomogeneous. Their normal forms are not taken from
\cite[15.1]{AGV85}, but from \cite[Anhang A2, p. 191]{Bi92}.
They are chosen such that the cross ratio
of the 4 intersection points on the central curve has
$j$-invariant 
$j=\frac{4}{27}\frac{(t_1^2-t_1+1)^3}{t_1^2(1-t_1)^2}$.
This fact implies that the families in table \eqref{10.1}
contain representatives of all right equivalence classes
in one $\mu$-homotopy class.

From the weights (or the spectral numbers, see below
theorem \ref{t10.6})
one deduces easily that any monomial basis of the Jacobi
algebra of one quasihomogeneous surface singularity
$f_{t_1,0}$ 
contains precisely one monomial $p_{>1}$ of weighted
degree $>1$ and that $\deg_w p_{>1}=1+\frac{2}{m}=1+w_y$.
\cite[12.6 Theorem]{AGV85} says here that any 
semiquasihomogeneous
singularity with quasihomogeneous part $f_{t_1,0}$
is right equivalent to $f_{t_1,0}+t_2\cdot p_{>1}$
for some $t_2\in\C$. In table \eqref{10.1} we have
chosen the monomial $p_{>1}$ such that it is part of a
monomial basis of the Jacobi algebra of $f_{t_1,0}$
for any $t_1\in T^{(1)}$.

\begin{remarks}\label{t10.2}
It is nontrivial (and slightly surprising) that such a 
monomial $p_{>1}$ exists simultaneously for all
$t_1\in T^{(1)}$.
In \cite{He93}\cite{He95} the second author had overlooked
this problem and had chosen in the four cases 
$S_{1,0},E_{3,0},Z_{1,0},Q_{2,0}$ a monomial which
does not work for special parameters $t_1\in T^{(1)}$.
The following table \eqref{10.8} lists for all 6 families
all monomials $\www p$ of weighted degree $1+\frac{2}{m}$
and for each of them the function $q(t_1)$ with
\begin{eqnarray}\label{10.7}
\www p \equiv q(t_1)\cdot p_{>1}\mod
(\textup{Jacobi ideal of }f_{t_1,0}),
\end{eqnarray}
where $p_{>1}=\frac{\ppp f_{(t_1,t_2)}}{\ppp t_2}$ is the
monomial chosen in table \eqref{10.1}.
\begin{eqnarray}\label{10.8}
\begin{array}{lllll}
 & p_{>1} & \www p\ :\ q(t_1) & \www p\ :\ q(t_1) 
 & \www p\ :\ q(t_1)  \\
W_{1,0} & x^2y^4 & x^4y:1-2t_1 & y^7:1-2t_1 & x^2yz:0\\
 & & y^4z:0 & yz^2:0 \\
S_{1,0} & x^2y^3 & x^2yz:-t_1  & y^4z:-t_1  & y^2z^2:t_1 \\
 & & y^6:2t_1-1 & x^4:2t_1-1 & z^3:t_1(2t_1-3) \\
U_{1,0} & y^4z & x^2yz:-t_1 & xyz^2:-t_1 & xy^4:t_1 \\
 & & x^3y:t_1(t_1-2) & yz^3:1-2t_1 & \\
E_{3,0} & x^2y^4 & x^3y:\frac{t_1+1}{2} & 
xy^7:\frac{t_1+1}{2t_1} & y^{10}:\frac{t_1^2-t_1+1}{t_1^2} \\
 & & yz^2:0 & & \\
Z_{1,0} & x^2y^4 & x^3y^2:\frac{t_1+1}{2} 
& xy^6:\frac{t_1+1}{2t_1} & y^8:\frac{t_1^2-t_1+1}{t_1^2} \\
 & & x^4:\frac{3}{2}t_1^2-2t_1+\frac{3}{2} & yz^2:0 &  \\
Q_{2,0} & xz^2 & x^2y^3:\frac{1}{(1-t_1)^2} & 
x^3y:\frac{t_1+1}{2(1-t_1)^2} & xy^5:\frac{t_1+1}{2t_1(1-t_1)^2}\\
 & & y^7:\frac{t_1^2-t_1+1}{t_1^2(1-t_1)^2} & y^2z^2:0
\end{array}
\end{eqnarray}
Thus $p_{>1}$ could be replaced in the normal form in table
\eqref{10.1} by any of the following monomials:
\begin{eqnarray}\label{10.9}
\begin{array}{llllll}
W_{1,0} & S_{1,0} & U_{1,0} & E_{3,0} & Z_{1,0} & Q_{2,0} \\
- & x^2yz, y^4z, y^2z^2 & x^2yz,xyz^2,xy^4 & - & - & x^2y^3 
\end{array}
\end{eqnarray}
\end{remarks}

We denote by $G_3$ and $G_2\subset G_3$ the groups of
automorphisms of $T^{(2)}=\P^1\C$
\begin{eqnarray}\label{10.10}
\begin{array}{rll}
G_3:=& \{t_1\mapsto t_1,1-t_1,\frac{1}{t_1},
\frac{t_1}{t_1-1},\frac{1}{1-t_1},\frac{t_1-1}{t_1}\}&
\cong S_3\textup{ as a group,}\\
G_2:=& \{t_1\mapsto t_1,1-t_1\}\subset G_3&
\cong S_2\textup{ as a group.}
\end{array}
\end{eqnarray}
They act also on $T^{(1)}=\C-\{0,1\}$.

\begin{theorem}\label{t10.3}
Consider a bimodal family of quadrangle surface singularities in
table \eqref{10.1}. A function
\begin{eqnarray}\label{10.11}
\kappa: G_2\times T^{(1)}\to\C^* & \textup{ for }
W_{1,0},S_{1,0},\\
\kappa: G_3\times T^{(1)}\to\C^* & \textup{ for }
U_{1,0},E_{3,0},Z_{1,0},Q_{2,0},\nonumber
\end{eqnarray}
with the following properties exists.
\begin{eqnarray}
f_{(t_1,t_2)}\sim_\RR f_{(\www t_1,\www t_2)}&\iff&
\exists\ g\in\left\{ \begin{array}{ll}
G_2& \textup{ for } W_{1,0},S_{1,0},\\
G_3& \textup{ for } U_{1,0},E_{3,0},Z_{1,0},Q_{2,0},
\end{array}\right. \nonumber\\
&& \textup{with }\www t_1=g(t_1),
\www t_2^{m_\infty} = \kappa(g,t_1)\cdot t_2^{m_\infty},\hspace*{1cm}
\label{10.12}\\
\kappa(\id,t_1)&=&1,\label{10.13}\\
\kappa(g_2g_1,t_1)&=& \kappa(g_1,t_1)\cdot \kappa(g_2,g_1(t_1)).
\label{10.14}
\end{eqnarray}
Table \eqref{10.15} lists $\kappa(g,t_1)$ for generators
$g$ of the group.
\begin{eqnarray}\label{10.15}
\begin{array}{lllllll}
 & W_{1,0} & S_{1,0} & U_{1,0} & E_{3,0} & Z_{1,0} & Q_{2,0} \\
t_1\mapsto 1-t_1 & 1 & -1 & 1 & 
\left(\frac{1-t_1}{t_1}\right)^{18} &
\left(\frac{1-t_1}{t_1}\right)^{14} & -1 \\
t_1\mapsto t_1^{-1} & - & - & -t_1^{-3} & t_1^{-12} & 
t_1^{-10} & t_1^3
\end{array}
\end{eqnarray}
\end{theorem}

{\bf Proof:}
\eqref{10.13}--\eqref{10.15} are consistent
(to check this is nontrivial only for $E_{3,0}$
and $Z_{1,0}$) and define a unique function $\kappa$ as 
in \eqref{10.11}.
We will show now that it satisfies $\Leftarrow$ in \eqref{10.12}.
We postpone the proof of $\Rightarrow$ in \eqref{10.12}
to the end of this section.

The equality
\begin{eqnarray}\label{10.16}
f_{(t_1,t_2)}(x\cdot e^{2\pi iw_x},y\cdot e^{2\pi i w_y},
z\cdot e^{2\pi i w_z})=f_{(t_1,t_2\cdot e^{2\pi i 2/m})}
\end{eqnarray}
gives $\Leftarrow$ in \eqref{10.12} for $g=\id$
and $\kappa(\id,t_1)=1$
(for $U_{1,0}$ $m=m_\infty=9$, in the other cases
$m_\infty=\frac{m}{2}$).
We list now coordinate changes 
$(x,y,z)\mapsto \varphi^{(1)}(x,y,z)$ and
$(x,y,z)\mapsto \varphi^{(2)}(x,y,z)$ with
\begin{eqnarray}\label{10.17}
f_{(t_1,t_2)}(\varphi^{(1)}(x,y,z)) &=&
f_{(1-t_1,0)}+t_2\cdot p^{(1)}(t_1,x,y,z) \nonumber\\
&&\textup{for all 6 cases},\\
f_{(t_1,t_2)}(\varphi^{(2)}(x,y,z)) &=& 
f_{(t_1^{-1},0)}(x,y,z)+t_2\cdot p^{(2)}(t_1,x,y,z)\nonumber\\
&&\textup{for }U_{1,0},E_{3,0},Z_{1,0},Q_{2,0}\label{10.18}
\end{eqnarray}
for certain quasihomogeneous polynomials $p^{(1)}$
and $p^{(2)}$ in $x,y,z$ with 
$\deg_wp^{(1)}=\deg_wp^{(2)}=1+\frac{2}{m}$.
\begin{eqnarray}\label{10.19}
\begin{array}{lllllll}
 & \varphi^{(1)}(x,y,z) & \varphi^{(2)}(x,y,z) \\
W_{1,0} & (x,-y,z) & - \\
S_{1,0} & (ix,y,-z-y^2) & - \\
U_{1,0} & (-x+z,-y,z) & (-z,t_1^{-1/3}y,-x) \\
E_{3,0} & (x-y^3,-y,z) & (x,t_1^{-1/3}y,z) \\
Z_{1,0} & (e^{-2\pi i/14}(x-y^2),i\cdot e^{-2\pi i/28}y,z) 
& (t_1^{1/7}x,t_1^{-3/7}y,z) \\
Q_{2,0} & (x-y^2,iy,e^{-2\pi i/8}z) & (x,t_1^{-1/2}y,t_1^{1/4}z) 
\end{array}
\end{eqnarray}
One can calculate $p^{(1)}$ and $p^{(2)}$ easily.
The proof of \cite[12.6 Lemma]{AGV85} implies here
\begin{eqnarray}
&&f_{(\www t_1,0)}+t_2\cdot \www p \ \sim_\RR\ 
f_{(\www t_1,\www t_2)}\nonumber \\ 
&&\textup{ where }
t_2\cdot\www p\equiv \www t_2\cdot 
p_{>1}\mod (\textup{Jacobi ideal of }f_{(\www t_1,0)}).
\hspace*{1cm}
\label{10.20}
\end{eqnarray}
With table \eqref{10.8} one finds $\www t_2$
with \eqref{10.20} for $\www p=p^{(1)}$
and for $\www p=p^{(2)}$. Then one verifies table
\eqref{10.15}.
\hfill $\Box$

\begin{remarks}\label{t10.4}
(i) For the quasihomogeneous singularities,
\eqref{10.12} becomes
\begin{eqnarray*}
f_{(t_1,0)}\sim_\RR f_{(\www t_1,0)}&\iff & 
\exists\ g\in G_2\textup{ resp. }G_3\textup{ with }
\www t_1=g(t_1).
\end{eqnarray*}
This is proved in \cite[Satz 1.5.2]{Bi92}
using the minimal good resolution.
Our proof of $\Rightarrow$ in \eqref{10.12} for all singularities
at the end of this section will be different.

\medskip
(ii) The right equivalence classes in $T^{(5)}$ are the orbits
of a group action on $T^{(5)}$ in the cases
$W_{1,0}$ and $S_{1,0}$. There the group is 
a central extension of $G_2$
by a cyclic group of order $m_\infty=\frac{m}{2}$,
\begin{eqnarray*}
1\to \left(\begin{array}{ll}\textup{cyclic group}\\ 
\textup{of order }m\end{array}\right)\to 
(\textup{group acting on }T^{(5)})\to G_2
\to 1.
\end{eqnarray*}
In the other cases $U_{1,0},E_{3,0},Z_{1,0}$ and $Q_{2,0}$,
an $m$-th root of $\kappa(t_1\to t_1^{-1},.):T^{(1)}\to \C^*$
is not uni-valued, but multi-valued.
There one has only a groupoid acting on $T^{(5)}$,
whose orbits are the right equivalence classes in $T^{(5)}$.

\medskip
(iii) In any case, the space 
$M_\mu^{mar}=(M_\mu^{mar})^0\cong T^{(7)}$ 
(by theorem \ref{t10.1})
will be more canonical than $T^{(5)}$, and there the 
right equivalence classes are the orbits of the action of
the group $G_\Z=G^{mar}$.
\end{remarks}

Now we come to the spectral numbers and the classifying space
$D_{BL}$.

\begin{lemma}\label{t10.5}
Consider a bimodal family of quadrangle surface singularities
in table \eqref{10.1}. Denote $\omega_0:=dxdydz$.

(a) The spectral numbers $\alpha_1,\ldots,\alpha_\mu$ with
$\alpha_1\leq \ldots\leq \alpha_\mu$ satisfy
\begin{eqnarray}\label{10.21}
\alpha_1=\frac{-1}{m}<\alpha_2=\frac{1}{m}<\alpha_3\leq \ldots
\leq \alpha_{\mu-2}\\
<\alpha_{\mu-1}=1-\frac{1}{m}
<\alpha_\mu=1+\frac{1}{m},\nonumber \\
\dim C^{\alpha_1}=\dim C^{\alpha_2}=2.\label{10.22}
\end{eqnarray}
The following picture illustrates this.

%\bigskip
%\includegraphics[width=1.0\textwidth]{pic-kap10-1.jpg}
\begin{scriptsize}
\begin{tikzpicture}[line cap=round,line join=round,>=triangle 45,x=0.78cm,y=0.8cm]
\clip(-8.09,-0.5) rectangle (8.10,2.100);
\fill[line width=0.8pt,color=ttqqqq,fill=ttqqqq,pattern=north east lines,pattern color=ttqqqq] (-5.6,0.) -- (-5.6,1.) -- (-5.4,1.) -- (-5.4,0.) -- cycle;
\fill[line width=0.8pt,color=ttqqqq,fill opacity=0] (-6.,0.) -- (-6.,1.) -- (-6.2,1.) -- (-6.2,0.) -- cycle;
\fill[line width=0.8pt,color=ttqqqq,fill opacity=0] (-6.6,0.) -- (-6.6,1.) -- (-6.8,1.) -- (-6.8,0.) -- cycle;
\fill[line width=0.8pt,color=ttqqqq,fill opacity=0] (-5.6,1.) -- (-5.4,1.) -- (-5.4,2.) -- (-5.6,2.) -- cycle;
\fill[line width=0.8pt,color=ttqqqq,fill=ttqqqq,pattern=north east lines,pattern color=ttqqqq] (-4.6,0.) -- (-4.4,0.) -- (-4.4,1.) -- (-4.6,1.) -- cycle;
\fill[line width=0.8pt,color=ttqqqq,fill opacity=0] (-4.6,1.) -- (-4.6,2.) -- (-4.4,2.) -- (-4.4,1.) -- cycle;
\fill[line width=0.8pt,color=ttqqqq,fill=ttqqqq,pattern=north east lines,pattern color=ttqqqq] (-4.,0.) -- (-4.,1.) -- (-3.8,1.) -- (-3.8,0.) -- cycle;
\fill[line width=0.8pt,color=ttqqqq,fill=ttqqqq,pattern=north east lines,pattern color=ttqqqq] (-3.4,1.) -- (-3.2,1.) -- (-3.2,0.) -- (-3.4,0.) -- cycle;
\fill[line width=0.8pt,color=ttqqqq,fill=ttqqqq,pattern=north east lines,pattern color=ttqqqq] (3.2,0.) -- (3.2,1.) -- (3.4,1.) -- (3.4,0.) -- cycle;
\fill[line width=0.8pt,color=ttqqqq,fill=ttqqqq,pattern=north east lines,pattern color=ttqqqq] (3.8,0.) -- (3.8,1.) -- (4.,1.) -- (4.,0.) -- cycle;
\fill[line width=0.8pt,color=ttqqqq,fill=ttqqqq,fill opacity=1.0] (4.4,1.) -- (4.4,0.) -- (4.6,0.) -- (4.6,1.) -- cycle;
\fill[line width=0.8pt,color=ttqqqq,fill=ttqqqq,pattern=north east lines,pattern color=ttqqqq] (4.4,1.) -- (4.4,2.) -- (4.6,2.) -- (4.6,1.) -- cycle;
\fill[line width=0.8pt,color=ttqqqq,fill=ttqqqq,fill opacity=1.0] (5.4,0.) -- (5.4,1.) -- (5.6,1.) -- (5.6,0.) -- cycle;
\fill[line width=0.8pt,color=ttqqqq,fill=ttqqqq,pattern=north east lines,pattern color=ttqqqq] (5.4,1.) -- (5.4,2.) -- (5.6,2.) -- (5.6,1.) -- cycle;
\fill[line width=0.8pt,color=ttqqqq,fill=ttqqqq,fill opacity=1.0] (6.,0.) -- (6.,1.) -- (6.2,1.) -- (6.2,0.) -- cycle;
\fill[line width=0.8pt,color=ttqqqq,fill=ttqqqq,fill opacity=1.0] (6.6,1.) -- (6.8,1.) -- (6.8,0.) -- (6.6,0.) -- cycle;
\draw [line width=1.2pt] (-8.,0.)-- (8.,0.);
\draw [line width=0.8pt,color=ttqqqq] (-5.6,0.)-- (-5.6,1.);
\draw [line width=0.8pt,color=ttqqqq] (-5.6,1.)-- (-5.4,1.);
\draw [line width=0.8pt,color=ttqqqq] (-5.4,1.)-- (-5.4,0.);
\draw [line width=0.8pt,color=ttqqqq] (-5.4,0.)-- (-5.6,0.);
\draw [line width=0.8pt,color=ttqqqq] (-6.,0.)-- (-6.,1.);
\draw [line width=0.8pt,color=ttqqqq] (-6.,1.)-- (-6.2,1.);
\draw [line width=0.8pt,color=ttqqqq] (-6.2,1.)-- (-6.2,0.);
\draw [line width=0.8pt,color=ttqqqq] (-6.2,0.)-- (-6.,0.);
\draw [line width=0.8pt,color=ttqqqq] (-6.6,0.)-- (-6.6,1.);
\draw [line width=0.8pt,color=ttqqqq] (-6.6,1.)-- (-6.8,1.);
\draw [line width=0.8pt,color=ttqqqq] (-6.8,1.)-- (-6.8,0.);
\draw [line width=0.8pt,color=ttqqqq] (-6.8,0.)-- (-6.6,0.);
\draw [line width=0.8pt,color=ttqqqq] (-5.6,1.)-- (-5.4,1.);
\draw [line width=0.8pt,color=ttqqqq] (-5.4,1.)-- (-5.4,2.);
\draw [line width=0.8pt,color=ttqqqq] (-5.4,2.)-- (-5.6,2.);
\draw [line width=0.8pt,color=ttqqqq] (-5.6,2.)-- (-5.6,1.);
\draw [line width=0.8pt,color=ttqqqq] (-4.6,0.)-- (-4.4,0.);
\draw [line width=0.8pt,color=ttqqqq] (-4.4,0.)-- (-4.4,1.);
\draw [line width=0.8pt,color=ttqqqq] (-4.4,1.)-- (-4.6,1.);
\draw [line width=0.8pt,color=ttqqqq] (-4.6,1.)-- (-4.6,0.);
\draw [line width=0.8pt,color=ttqqqq] (-4.6,1.)-- (-4.6,2.);
\draw [line width=0.8pt,color=ttqqqq] (-4.6,2.)-- (-4.4,2.);
\draw [line width=0.8pt,color=ttqqqq] (-4.4,2.)-- (-4.4,1.);
\draw [line width=0.8pt,color=ttqqqq] (-4.4,1.)-- (-4.6,1.);
\draw [line width=0.8pt,color=ttqqqq] (-4.,0.)-- (-4.,1.);
\draw [line width=0.8pt,color=ttqqqq] (-4.,1.)-- (-3.8,1.);
\draw [line width=0.8pt,color=ttqqqq] (-3.8,1.)-- (-3.8,0.);
\draw [line width=0.8pt,color=ttqqqq] (-3.8,0.)-- (-4.,0.);
\draw [line width=0.8pt,color=ttqqqq] (-3.4,1.)-- (-3.2,1.);
\draw [line width=0.8pt,color=ttqqqq] (-3.2,1.)-- (-3.2,0.);
\draw [line width=0.8pt,color=ttqqqq] (-3.2,0.)-- (-3.4,0.);
\draw [line width=0.8pt,color=ttqqqq] (-3.4,0.)-- (-3.4,1.);
\draw [line width=0.8pt,color=ttqqqq] (3.2,0.)-- (3.2,1.);
\draw [line width=0.8pt,color=ttqqqq] (3.2,1.)-- (3.4,1.);
\draw [line width=0.8pt,color=ttqqqq] (3.4,1.)-- (3.4,0.);
\draw [line width=0.8pt,color=ttqqqq] (3.4,0.)-- (3.2,0.);
\draw [line width=0.8pt,color=ttqqqq] (3.8,0.)-- (3.8,1.);
\draw [line width=0.8pt,color=ttqqqq] (3.8,1.)-- (4.,1.);
\draw [line width=0.8pt,color=ttqqqq] (4.,1.)-- (4.,0.);
\draw [line width=0.8pt,color=ttqqqq] (4.,0.)-- (3.8,0.);
\draw [line width=0.8pt,color=ttqqqq] (4.4,1.)-- (4.4,0.);
\draw [line width=0.8pt,color=ttqqqq] (4.4,0.)-- (4.6,0.);
\draw [line width=0.8pt,color=ttqqqq] (4.6,0.)-- (4.6,1.);
\draw [line width=0.8pt,color=ttqqqq] (4.6,1.)-- (4.4,1.);
\draw [line width=0.8pt,color=ttqqqq] (4.4,1.)-- (4.4,2.);
\draw [line width=0.8pt,color=ttqqqq] (4.4,2.)-- (4.6,2.);
\draw [line width=0.8pt,color=ttqqqq] (4.6,2.)-- (4.6,1.);
\draw [line width=0.8pt,color=ttqqqq] (4.6,1.)-- (4.4,1.);
\draw [line width=0.8pt,color=ttqqqq] (5.4,0.)-- (5.4,1.);
\draw [line width=0.8pt,color=ttqqqq] (5.4,1.)-- (5.6,1.);
\draw [line width=0.8pt,color=ttqqqq] (5.6,1.)-- (5.6,0.);
\draw [line width=0.8pt,color=ttqqqq] (5.6,0.)-- (5.4,0.);
\draw [line width=0.8pt,color=ttqqqq] (5.4,1.)-- (5.4,2.);
\draw [line width=0.8pt,color=ttqqqq] (5.4,2.)-- (5.6,2.);
\draw [line width=0.8pt,color=ttqqqq] (5.6,2.)-- (5.6,1.);
\draw [line width=0.8pt,color=ttqqqq] (5.6,1.)-- (5.4,1.);
\draw [line width=0.8pt,color=ttqqqq] (6.,0.)-- (6.,1.);
\draw [line width=0.8pt,color=ttqqqq] (6.,1.)-- (6.2,1.);
\draw [line width=0.8pt,color=ttqqqq] (6.2,1.)-- (6.2,0.);
\draw [line width=0.8pt,color=ttqqqq] (6.2,0.)-- (6.,0.);
\draw [line width=0.8pt,color=ttqqqq] (6.6,1.)-- (6.8,1.);
\draw [line width=0.8pt,color=ttqqqq] (6.8,1.)-- (6.8,0.);
\draw [line width=0.8pt,color=ttqqqq] (6.8,0.)-- (6.6,0.);
\draw [line width=0.8pt,color=ttqqqq] (6.6,0.)-- (6.6,1.);
\draw (-5.85,0.04) node[anchor=north west] {$\alpha_{1}$};
\draw (-4.8,0.04) node[anchor=north west] {$\alpha_{2}$};
\draw (5.25,0.04) node[anchor=north west] {$\alpha_{\mu}$};
\draw (4.1,0.04) node[anchor=north west] {$\alpha_{\mu-1}$};
\begin{scriptsize}
\draw [color=ttqqqq] (0.,0.)-- ++(-2.5pt,0 pt) -- ++(5.0pt,0 pt) ++(-2.5pt,-2.5pt) -- ++(0 pt,5.0pt);
\draw[color=ttqqqq] (0.012,0.32) node {1/2};
\draw [color=ttqqqq] (-5.,0.)-- ++(-2.5pt,0 pt) -- ++(5.0pt,0 pt) ++(-2.5pt,-2.5pt) -- ++(0 pt,5.0pt);
\draw[color=ttqqqq] (-5.01,0.32) node {0};
\draw [color=ttqqqq] (5.,0.)-- ++(-2.5pt,0 pt) -- ++(5.0pt,0 pt) ++(-2.5pt,-2.5pt) -- ++(0 pt,5.0pt);
\draw[color=ttqqqq] (5.00,0.32) node {1};
\draw [fill=ttqqqq] (-2.4,0.5) circle (1.2pt);
\draw [fill=ttqqqq] (-2.,0.5) circle (1.2pt);
\draw [fill=ttqqqq] (-1.6,0.5) circle (1.2pt);
\draw [fill=ttqqqq] (-7.25,0.5) circle (1.2pt);
\draw [fill=ttqqqq] (-7.65,0.5) circle (1.2pt);
\draw [fill=ttqqqq] (1.92035,0.5) circle (1.2pt);
\draw [fill=ttqqqq] (2.32035,0.5) circle (1.2pt);
\draw [fill=ttqqqq] (2.72035,0.5) circle (1.2pt);
\draw [fill=ttqqqq] (7.67831,0.5) circle (1.2pt);
\draw [fill=ttqqqq] (7.27831,0.5) circle (1.2pt);
\end{scriptsize}
\end{tikzpicture}
\end{scriptsize}

%\bigskip
We also have
\begin{eqnarray}\label{10.23}
V^{\alpha_1}(f_{(t_1,t_2)})&\supset& H_0''(f_{(t_1,t_2)}
\supset V^{>\alpha_2}(f_{(t_1,t_2)}),\\
H_0''(f_{(t_1,t_2)})&=&
\C\cdot\left(s(\omega_0,\alpha_1)(t_1,t_2)
+ s(\omega_0,\alpha_2)(t_1,t_2)\right) \nonumber\\
&+& \C\cdot s(y\omega_0,\alpha_2)(t_1,t_2)+
V^{>\alpha_2}(f_{(t_1,t_2)}).\label{10.24}
\end{eqnarray}

\medskip
(b) The polarizing form $S$ defines an indefinite form
$((a,b)\mapsto S(a,\oooo b))$ on $H^\infty_\zeta$.
We get a half-plane
\begin{eqnarray}\label{10.25}
\HH(C^{\alpha_1})&:=&
\{\C\cdot v\, |\, v\in C^{\alpha_1}\textup{ with }
S(\psi_{\alpha_1}^{-1}(v),\oooo{\psi_{\alpha_1}^{-1}(v)})<0\}\\
&\subset & \P^1(C^{\alpha_1}).\nonumber
\end{eqnarray}

(c) 
\begin{eqnarray}
D_{BL}&=& \{\C\cdot (v_1+v_2)\, |\, v_1\in C^{\alpha_1}-\{0\}
\textup{ with }[\C\cdot v_1]\in \HH(C^{\alpha_1}),\nonumber\\
&& \hspace*{2cm}
v_2\in \C\cdot \psi_{\alpha_2}(\oooo{\psi_{\alpha_1}^{-1}(v_1)})
\subset C^{\alpha_2}\}\label{10.25b}\\
&\cong & \HH(C^{\alpha_1})\times\C.\nonumber
\end{eqnarray}
\end{lemma}

{\bf Proof:}
(a) The spectral numbers are well known 
\cite[13.3.4, p. 389]{AGV88} and can be calculated in the
semiquasihomogeneous cases for example with the 
generating series (here $m=2$, $(w_0,w_1,w_2)=(w_x,w_y,w_z)$)
\begin{eqnarray}\label{10.26}
\prod_{j=0}^m \frac{t-t^{w_j}}{t^{w_j}-1} =
\sum_{i=1}^\mu t^{\alpha_i+1}.
\end{eqnarray}
\eqref{10.22} and \eqref{10.23} are obvious.
\eqref{10.24} follows from lemma \ref{7.20} and
$\deg_w(\omega_0)=\alpha_1+1,\deg_w(y\omega_0)=\alpha_2+1$
and $\deg_w(x^iy^jz^k\omega_0)>\alpha_2+1$ for any
other monomial $x^iy^jz^k$, because $w_y<w_x\leq w_z$.

\medskip
(b) This follows as in section \ref{c9} before theorem
\ref{t9.2}. It follows also from the fact that
$\Gr_V^\bullet H_0''(f_{(t_1,t_2)})$ and $S$ induce
as in \eqref{7.27} a polarized Hodge structure of weight 2
on $H^\infty(f_{(t_1,t_2)})$. Especially,
\begin{eqnarray}\label{10.27}
a_1(t_1,t_2):=& \psi_{\alpha_1}^{-1} 
s(\omega_0,\alpha_1)(t_1,t_2)\in H^\infty(f_{(t_1,t_2)})_\zeta,\\
a_2(t_1,t_2):=& \psi_{\alpha_2}^{-1} 
s(y\omega_0,\alpha_2)(t_1,t_2)\in 
H^\infty(f_{(t_1,t_2)})_{\oooo \zeta}\nonumber
\end{eqnarray}
satisfy
\begin{eqnarray}\label{10.28}
\textup{on }H^\infty(f_{(t_1,t_2)})_\zeta:&&
\C\cdot a_1 = H^{2,0}= F^2 \subset 
H^\infty(f_{(t_1,t_2)})_\zeta \\ 
&&=F^1=H^{2,0}\oplus H^{1,1}
=\C\cdot a_1\oplus \C\cdot\oooo{a_2},\nonumber\\
\textup{on }H^\infty(f_{(t_1,t_2)})_{\oooo\zeta}:&&
\C\cdot a_2 = H^{1,1}= F^1 \subset 
H^\infty(f_{(t_1,t_2)})_{\oooo\zeta} \label{10.29}\\
&&=F^0=H^{1,1}\oplus H^{0,2}
=\C\cdot a_2\oplus \C\cdot\oooo{a_1},\nonumber\\
0<i^{2-0}S(a_1,\oooo{a_1}),&&
0<i^{1-1}S(\oooo{a_2},a_2),\ 
0=S(a_1,a_2).\label{10.30}
\end{eqnarray}

(c) %The space $\C\cdot a_2$ is $S$-orthogonal to $\C\cdot a_1$.
%Therefore $\C\cdot s(y\omega_0,\alpha_2)(t_1,t_2)$
%is determined by $\C\cdot s(\omega_0,\alpha_1)(t_1,t_2)$.
%Thus $H_0''(f_{(t_1,t_2)})$ is determined by 
%$\C\cdot (v_1+v_2)$ where
%\begin{eqnarray*}
%v_1&=&s(\omega_0,\alpha_1)(t_1,t_2),\\
%v_2&\in&\C\cdot \psi_{\alpha_2}^{-1}
%(\oooo{\psi_{\alpha_1}(v_1)}),\\
%\textup{and }v_2&\equiv& s(\omega_0,\alpha_2)(t_1,t_2)
%\mod \C\cdot s(y\omega_0,\alpha_2)(t_1,t_2)
%\end{eqnarray*}
This follows as in lemma \ref{t9.4} (c) in the case $m|p$.
\hfill$\Box$

\bigskip

The multi-valued period map $BL_{T^{(5)}}:T^{(5)}\to D_{BL}$
had been calculated in \cite{He93}.
We recall the result and sketch the proof.
In part (e) of theorem \ref{t10.6} we add a formula for the
case $S_{1,0}$ which will be useful for the determination
of a transversal monodromy in theorem \ref{t10.7}.

\begin{theorem}\label{t10.6}
Consider a bimodal family of quadrangle surface singularities
in table \eqref{10.1}. 

\medskip
(a) $s(\omega_0,\alpha_1)(t_1,t_2)=s(\omega_0,\alpha_1)(t_1,0)
=s[\omega_0](t_1,0)$ is independent of $t_2$ and satisfies
the hypergeometric differential equation
\begin{eqnarray}\label{10.31}
0=\Bigl(t_1(1-t_1)\ppp_{t_1}^2 + (c-(a+b+1)t_1)\ppp_{t_1}-ab\Bigr) 
s[\omega_0](t_1,0)
\end{eqnarray}
with $(1-c,c-a-b,a-b)=(\frac{1}{m_0},\frac{1}{m_1},
\frac{1}{m_\infty})$.

\medskip
(b) The multi-valued period map
\begin{eqnarray}\label{10.32}
BL_{T^{(1)}}:T^{(1)}\to \HH(C^{\alpha_1}),\quad
t_1\mapsto \C\cdot s[\omega_0](t_1,0),
\end{eqnarray}
lifts to a uni-valued period map
\begin{eqnarray}\label{10.33}
BL_{T^{(3)}}:T^{(3)}\to \HH(C^{\alpha_1})
\end{eqnarray}
which is an open embedding and extends to an isomorphism
\begin{eqnarray}\label{10.34}
BL_{T^{(4)}}:T^{(4)}\to \HH(C^{\alpha_1}).
\end{eqnarray}

(c) 
\begin{eqnarray}\label{10.35}
s(\omega_0,\alpha_2)(t_1,t_2)=t_2\cdot(-\ppp_\tau)
s[p_{>1}\omega_0](t_1,0),\\
C^{\alpha_2}= \C\cdot s[y\omega_0](t_1,0)\oplus
\C\cdot \ppp_\tau s[p_{>1}\omega_0](t_1,0).
\label{10.36}
\end{eqnarray}

(d) The multi-valued period map
\begin{eqnarray}\label{10.37}
BL_{T^{(5)}}:T^{(5)}\to D_{BL}
\end{eqnarray}
is locally in $T^{(1)}$ and $\HH(C^{\alpha_1})$ an isomorphism
of line bundles and lifts to an open embedding of line bundles
\begin{eqnarray}\label{10.38}
BL_{T^{(7)}}:T^{(7)}\to D_{BL}.
\end{eqnarray}
(We do not know whether this extends to an isomorphism
of line bundles $T^{(8)}\to D_{BL}$, but we do not
expect it.)

\medskip
(e) In the case of $S_{1,0}$
\begin{eqnarray}\label{10.39}
\ppp_{t_1} s[x\omega_0](t_1,0) = 
\frac{2t_1-1}{5t_1(1-t_1)}\cdot s[x\omega_0](t_1,0).
\end{eqnarray}
\end{theorem}

{\bf Proof:} 
(a) We just sketch the ansatz for the calculations which
prove \eqref{10.31}.
$f_{(t_1,0)}$ and $\ppp_{t_1}f_{(t_1,0)}$ are quasihomogeneous
of weighted degree 1. 
List all monomials $d_1,\ldots,d_l$ in $x,y,z$ which turn up
in $f_{(t_1,0)}^2$, $f_{(t_1,0)}\cdot\ppp_{t_1}f_{(t_1,0)}$
and $(\ppp_{t_1}f_{(t_1,0)})^2$, 
find $l-2$ independent linear combinations of 
$d_1\omega_0,\ldots,d_l\omega_0$ in
$df_{(t_1,0)}\land d\Omega^1_{\C^3}$, and determine an equation
\begin{eqnarray}
&&p_1\cdot(\ppp_{t_1}f_{(t_1,0)})^2 \cdot\omega_0
+p_2\cdot f_{(t_1,0)}\cdot \ppp_{t_1}f_{(t_1,0)}\cdot\omega_0
+p_3\cdot f_{(t_1,0)}^2 \cdot\omega_0\nonumber \\
&&\equiv 0 
\mod  df_{(t_1,0)}\land d\Omega^1_{\C^3} \label{10.40}
\end{eqnarray}
with $p_1,p_2,p_3\in\Q[t_1]$. Then
\begin{eqnarray}\label{10.41}
\left(p_1\ppp_{t_1}^2 - (\alpha_1+2)p_2\ppp_{t_1}
+(\alpha_1+2)(\alpha_1+1)p_3\right) s[\omega_0](t_1,0).
\end{eqnarray}
Because of corollary \ref{t8.14} one can work in the cases
$W_{1,0},E_{3,0},Z_{1,0}$ with the curve singularities.
There the number $l$ of monomials is $l=5$.
In the other cases, the surfaces singularities
$S_{1,0},U_{1,0},Q_{2,0}$, it is $l=9$.

\medskip
(b) The period map $BL_{T^{(1)}}$ is not constant because
$s[\omega_0](t_1,0)$ and 
$\ppp_{t_1}s[\omega_0](t_1,0)
=(-\ppp_\tau)s[\ppp_{t_1}f_{(t_1,0)}\cdot\omega_0](t_1,0)$
are linearly independent because $\ppp_{t_1}f_{(t_1,0)}$ 
is not in the Jacobi ideal.
Therefore the multi-valued coefficient functions
$f_1(t_1)$ and $f_2(t_1)$ with
\begin{eqnarray}\label{10.42}
%s[\omega_0](t_1,0) = f_1(t_1)\cdot \psi_{\alpha_1}
%(\oooo{\psi_{\alpha_2}^{-1}({s[y\omega_0](i,0)})})
%+f_2(t_1)\cdot s[\omega_0](i,0)
s[\omega_0](t_1,0) = f_1(t_1)\cdot v_1^0
+f_2(t_1)\cdot v_2^0
\end{eqnarray}
for an arbitrary basis $v_1^0,v_2^0$ of $C^{\alpha_1}$ 
are linearly independent scalar solutions of the same
hypergeometric differential equation.
Their quotient
$(t_1\mapsto \frac{f_1(t_1)}{f_2(t_1)})$ is a Schwarzian
function \cite[sec. 113+114]{Fo51}, which maps
the closure of the upper half-plane to a hyperbolic
triangle with angles 
$\frac{\pi}{m_0},\frac{\pi}{m_1},\frac{\pi}{m_\infty}$.
The vertices are the images of $0,1,\infty$.
Therefore the multi-valued map 
$BL_{T^{(1)}}:T^{(1)}\to \HH(C^{\alpha_1})$
is an inverse of the quotient map
$c^{(1)}:T^{(3)}\to T^{(1)}$. This shows \eqref{10.33}
and \eqref{10.34}.

\medskip
(c) $s(\omega_0,\alpha_2)(t_1,0)=0$ because of 
formula \eqref{7.55} in lemma \ref{t7.20} (a).
\begin{eqnarray}
\ppp_{t_2}s(\omega_0,\alpha_2)(t_1,t_2) &=&
(-\ppp_\tau)s(p_{>1}\omega_0,\alpha_2+1)(t_1,t_2)\nonumber\\
&=& (-\ppp_\tau)s[p_{>1}\omega_0](t_1,0)\nonumber\\
\textup{thus }
s(\omega_0,\alpha_2)(t_1,t_2)&=& t_2\cdot (-\ppp_\tau)
s[p_{>1}\omega_0](t_1,0)\nonumber\\
&\equiv& t_2\cdot v_2 \mod \C\cdot s[y\omega_0](t_1,0)
\label{10.43}\\
\textup{with a suitable }v_2&\in& 
\psi_{\alpha_2}^{-1}(\oooo{\psi_{\alpha_1}(s[\omega_0](t_1,0))})
-\{0\}.\nonumber
\end{eqnarray}
Here $v_2\neq 0$ follows from \eqref{10.36} which is a 
consequence of the fact that $p_{>1}$ is not in the
Jacobi ideal of $f_{(t_1,0)}$.

\medskip
(d) This follows from \eqref{10.33} and part (c).

\medskip
(e) The proof is similar to the calculations which prove
part (a), but simpler.
\begin{eqnarray*}
&&\ppp_{t_1} s[x\omega_0](t_1,0)\\
&=& (-\ppp_\tau)s[\ppp_{t_1}f_{(t_1,0)}\cdot x\omega_0](t_1,0)
=(-\ppp_\tau)s[x^3y^2\omega_0](t_1,0)\\
&\stackrel{(*)}{=}& 
\frac{2t_1-1}{6t_1(t_1-1)}(-\ppp_\tau)
s[f_{(t_1,0)}\cdot x\omega_0](t_1,0)\\
&=& \frac{2t_1-1}{6t_1(t_1-1)}(-\ppp_\tau\tau)
s[x\omega_0](t_1,0)
=\frac{2t_1-1}{6t_1(t_1-1)}(-\frac{6}{5})s[x\omega_0](t_1,0)\\
&=& \frac{2t_1-1}{5t_1(1-t_1)}s[x\omega_0](t_1,0).
\end{eqnarray*}
For $\stackrel{(*)}{=}$ one has to find 3 relations in
$df_{(t_1,0)}\land d\Omega^1_{\C^3}$ between the 
monomial differential forms 
$x^3y^2\omega_0,xy^3z\omega_0,xyz^2\omega_0$ and $x^3z\omega_0$
in $f_{(t_1,0)}\cdot x\omega_0$ and $x^3y^2\omega_0$.
\hfill$\Box$

\bigskip
The last step before the proof of theorem \ref{t10.1}
is the following result on a transversal monodromy group.
Its proof uses formula \eqref{6.8} in theorem \ref{t6.1}.

\begin{theorem}\label{t10.7}
Consider a bimodal family of quadrangle surface singularities
in table \eqref{10.1}.
The pull back to $T^{(3)}$ with $c^{(1)}$ of the homology group
$\bigcup_{t_1\in T^{(1)}} Ml(f_{(t_1,0)})\to T^{(1)}$
comes equipped with a monodromy representation
$\pi^{(3)}:\pi_1(T^{(3)},\tau^{(3)})\to G_\Z$
(with $c^{(1)}(\tau^{(3)})=i$)
which is called {\it transversal monodromy group}.

\medskip
(a) The following table lists the local monodromies around 
elliptic fixed points in $(c^{(2)})^{-1}(0)$,
$(c^{(2)})^{-1}(1)$ and $(c^{(2)})^{-1}(\infty)$.
\begin{eqnarray}\label{10.44}
\begin{array}{lllllll}
 & W_{1,0} & S_{1,0} & U_{1,0} & E_{3,0} & Z_{1,0} & Q_{2,0}\\
(c^{(2)})^{-1}(\{0,1\}) & \id & \id & \id & \id & \id & \id \\
(c^{(2)})^{-1}(\infty) & \id & M_h^5 & \id & \id & \id & 
M_h^6
\end{array}
\end{eqnarray}
Therefore $\Imm(\pi^{(3)})=\{\id\}$ for 
$W_{1,0},U_{1,0},E_{3,0},Z_{1,0}$,
and $\Imm(\pi^{(3)})=\{\id,M_h^{m_\infty}\}$ for 
$S_{1,0}$ and $Q_{2,0}$. 

\medskip
(b) 
\begin{eqnarray}
&&\{g\in G_\Z\, |\, g\textup{ acts trivially on }D_{BL}\}
\nonumber \\
&=& \{g\in G_\Z\, |\, g=\pm\id \textup{ on }Ml_\zeta\}\nonumber\\
&=&\{\pm\id,\pm M_h^{m_\infty}\}\nonumber \\
&=& \left\{\begin{array}{ll}
\{\pm\id\} & \textup{for }U_{1,0},E_{3,0},Z_{1,0}\\
\{\pm\id,\pm M_h^{m_\infty}\} & \textup{for }
W_{1,0},S_{1,0},Q_{2,0}.\end{array}\right. \label{10.45}
\end{eqnarray}

(c) $G^{smar,gen}_\RR$ is here the group in \eqref{8.13}
for the singularities of multiplicity $\geq 3$, namely
the curve singularities $W_{1,0},E_{3,0},Z_{1,0}$ and the
surface singularities $S_{1,0},U_{1,0},Q_{2,0}$.
\begin{eqnarray}\label{10.46}
G^{smar,gen}_\RR &=& \left\{\begin{array}{ll}
\{\id\} & \textup{for }U_{1,0},E_{3,0},Z_{1,0},\\
\{\id,M_h^{m_\infty}\} &\textup{for } W_{1,0},S_{1,0},Q_{2,0}.
\end{array}\right. 
\end{eqnarray}
\end{theorem}

{\bf Proof:}
We start with part (b).
Suppose that $g\in G_\Z$ acts trivially on $D_{BL}$.
Then it acts trivially on $\HH(C^{\alpha_1})$, so 
$g=\lambda\cdot\id$ on $Ml_\zeta$ for some
$\lambda\in \C^*$. And $\C\cdot (v_1+v_2)=\C\cdot (\lambda v_1
+\oooo\lambda v_2)$, so $\lambda=\oooo\lambda\in\{\pm 1\}$.
This together with formula \eqref{6.8} and the set of
eigenvalues of $M_h$ gives \eqref{10.45}.

(a) The Papperitz-Riemann symbol
\begin{eqnarray}\label{10.47}
\left\{\begin{matrix} 
0 & 1 & \infty & \\ 0 & 0 & a & z \\ 1-c & c-a-b & b & 
\end{matrix}
\right\}
\end{eqnarray}
encodes the local behaviour near $0$, $1$ and $\infty$ of 
scalar solutions of the hypergeometric equation.
Locally suitable solutions have the following form
(h.o.t. = higher order terms):
\begin{eqnarray}\label{10.48}
\begin{array}{llll}
\textup{near }0:&
t_1^0+\textup{h.o.t.}&\textup{ and }&t_1^{1-c}+\textup{h.o.t.},\\
\textup{near }1:&
(t_1-1)^0+\textup{h.o.t.}&\textup{ and }&
(t_1-1)^{c-a-b}+\textup{h.o.t.},\\
\textup{near }\infty:&
t_1^{-a}+\textup{h.o.t.}&\textup{ and }&t_1^{-b}+\textup{h.o.t.} 
\end{array}
\end{eqnarray}
Especially, the local monodromy of the space of solutions
has the eigenvalues 
\begin{eqnarray}\label{10.49}
\begin{array}{llll}
\textup{around }0:&
1 &\textup{ and }& e^{2\pi i (1-c)},\\
\textup{around }1:&
1 &\textup{ and }& e^{2\pi i(c-a-b)},\\
\textup{around }\infty:&
e^{-2\pi i a} &\textup{ and } & e^{-2\pi i b}.
\end{array}
\end{eqnarray}
In our situation $(1-c,c-a-b,a-b)=(\frac{1}{m_0},\frac{1}{m_1},
\frac{1}{m_\infty})$,
\begin{eqnarray}\label{10.50}
\begin{array}{lllllll}
 & W_{1,0} & S_{1,0} & U_{1,0} & E_{3,0} & Z_{1,0} & Q_{2,0}\\
a & \frac{1}{2} & \frac{1}{2} & \frac{4}{9} & 
\frac{4}{9} & \frac{3}{7} & \frac{5}{12} \\
b & \frac{1}{3} & \frac{3}{10} & \frac{1}{3} & 
\frac{1}{3} & \frac{2}{7} & \frac{1}{4} \\
c & \frac{11}{12} & \frac{9}{10} & \frac{8}{9} & 
\frac{8}{9} & \frac{6}{7} & \frac{5}{6} 
\end{array}
\end{eqnarray}
The branched covering $c^{(2)}:T^{(4)}\to T^{(2)}$
has at elliptic fixed points the orders $m_0,m_1,m_\infty$.
Therefore the local monodromies of the pull back to $T^{(3)}$
of the solutions on $T^{(1)}=\C-\{0,1\}\subset T^{(2)}=\P^1\C$
become $+\id$ except around the elliptic fixed points
in $(c^{(2)})^{-1}(\infty)$ in the cases $S_{1,0}$
and $Q_{2,0}$ where they become $-\id$.

The same holds for the restrictions to $Ml_\zeta$ of the
local monodromies in $\pi^{(3)}$.

With \eqref{6.8} we obtain \eqref{10.44} for 
$U_{1,0},E_{3,0},Z_{1,0}$ and the following approximation
of \eqref{10.44} for $W_{1,0},S_{1,0},Q_{2,0}$.
\begin{eqnarray}\label{10.51}
\begin{array}{llll}
 & W_{1,0} & S_{1,0} & Q_{2,0}\\
(c^{(2)})^{-1}(\{0,1\}) & \id\textup{ or }-M_h^6 & 
\id\textup{ or }-M_h^5 & 
\id\textup{ or }-M_h^6 \\
(c^{(2)})^{-1}(\infty) & \id\textup{ or }-M_h^6 & 
-\id\textup{ or }M_h^5 & 
-\id\textup{ or }M_h^6
\end{array}
\end{eqnarray}

The case $W_{1,0}$: The sublattice $Ml_{-1,\Z}$ has rank 1.
Therefore the local transversal monodromies of the
homology bundle on $T^{(1)}$ around $0$, $1$ and $\infty$
have on $Ml_{-1,\Z}$ eigenvalues in $\{\pm 1\}$.
The branched covering $T^{(4)}\to T^{(2)}$ is at
the elliptic fixed points of even order. Thus $\pi^{(3)}$
restricts to the trivial monodromy on $Ml_{-1,\Z}$.
This excludes $-M_h^6$ in \eqref{10.51}.

The case $S_{1,0}$: The local transversal monodromies of the
homology bundle on $T^{(1)}$ around $0$, $1$ and $\infty$ 
have on $Ml_{e^{-2\pi i/5}}$ eigenvalues in $\Eiw(\zeta)$.
The branched covering is at the elliptic fixed points
in $(c^{(2)})^{-1}(\{0,1\})$ of order $10$.
Thus the local monodromies of $\pi^{(3)}$ around points
in $(c^{(2)})^{-1}(\{0,1\})$ are trivial on
$Ml_{e^{-2\pi i/5}}$. This excludes $-M_h^5$ in the
first line of \eqref{10.51}.
The branched covering is at the elliptic fixed points
in $(c^{(2)})^{-1}(\infty)$ of order $5$. 
Formula \eqref{10.39} in theorem \ref{t10.6} gives near $\infty$
\begin{eqnarray}\label{10.52}
s[x\omega_0](t_1,0) = 
(t_1^{-2/5}+\textup{h.o.t.})\cdot 
(\textup{a flat multi-valued section}).
\end{eqnarray}
Therefore also the local monodromy of $\pi^{(3)}$ around
points in $(c^{(2)})^{-1}(\infty)$ is trivial.
This excludes $-\id$ in the second line of \eqref{10.51}.

The case $Q_{2,0}$: The local transversal monodromies of the
homology bundle on $T^{(1)}$ around $0$, $1$ and $\infty$
have on $Ml_{e^{-2\pi i/3}}$ eigenvalues in 
$\Eiw(e^{2\pi i/6})$. 
The branched covering $T^{(4)}\to T^{(2)}$ is at
the elliptic fixed points of order $6$. Thus $\pi^{(3)}$
restricts to the trivial monodromy on $Ml_{e^{-2\pi i/3}}$.
This excludes $-M_h^6$ in the first line and $-\id$ in
the second line of \eqref{10.51}.

\medskip
(c) $-\id\notin G^{smar,gen}_\RR$ by theorem \ref{t8.8} (d).
$G^{smar,gen}_\RR$ fixes $BL(f,\pm\rho)$ for any
$(f,\pm\rho)\in M_\mu^{mar}$.
Because $T^{(7)}\to D_{BL}$ is an open embedding, 
$G^{smar,gen}_\RR$ fixes $D_{BL}$.
By part (b) $G^{smar,gen}_\RR=\{\id\}$ for
$U_{1,0},E_{3,0},Z_{1,0}$, and 
$G^{smar,gen}_\RR = \{\id\}$ or  $\{\id, M_h^{m_\infty}\}$ 
or $\{\id,-M_h^{m_\infty}\}$ for 
$W_{1,0},S_{1,0},Q_{2,0}$. 
The coordinate changes $\varphi$ of the curve singularities
$W_{1,0}$ and the surface singularities $S_{1,0}$ and 
$Q_{2,0}$ in the following table give a nontrivial element
of $G^{smar,gen}_\RR$.
\begin{eqnarray}\label{10.53}
\begin{array}{lll}
W_{1,0} & S_{1,0} & Q_{2,0} \\
(x,y)\mapsto (-x,y) & (x,y,z)\mapsto (-x,y,z) &
(x,y,z)\mapsto (x,y,-z)\end{array}
\end{eqnarray}
The coordinate change $\varphi$ maps $\omega_0$ to $-\omega_0$
and $s[\omega_0](t_1,0)$ to $-s[\omega_0](t_1,0)$.
Therefore $(\varphi)_{hom}|_{Ml_\zeta}=-\id$ and
$(\varphi)_{hom}=M_h^{m_\infty}$ 
(and not $-M_h^{m_\infty}$). This shows \eqref{10.45}
for $W_{1,0},S_{1,0},Q_{2,0}$.
\hfill$\Box$ 

\bigskip
Finally we come to the proof of theorem \ref{t10.1}.
Within this proof, we will also finish the proof
of theorem \ref{t6.1}.
After it, we will finish the proof of theorem \ref{t10.3}.

\bigskip
{\bf Proof of theorem \ref{t10.1}:}
By theorem \ref{t10.7} (a)+(c), the transversal monodromy
representation $\pi^{(7)}$ of the pull back to $T^{(7)}$
with $c^{(5)}$ of the homology bundle
$\bigcup_{(t_1,t_2)\in T^{(5)}} Ml(f_{(t_1,t_2)})\to T^{(5)}$
is trivial in the cases $W_{1,0},U_{1,0},E_{3,0},Z_{1,0}$
and has image in $G^{smar,gen}_\RR=\{\id,M_h^{m_\infty}\}$
in the cases $S_{1,0}$ and $Q_{2,0}$.
Thus the strong marking $+\id$ on $f_{(i,0)}$ induces for
each $f_{(t_1,t_2)}$ two strong markings in the same
right equivalence class in the cases $S_{1,0}$ and $Q_{2,0}$
and one strong marking in the other cases.
In any case, this gives a map $T^{(7)}\to (M_\mu^{smar})^0$.

The composition $T^{(7)}\to (M_\mu^{smar})^0\to D_{BL}$
is an open embedding by theorem \ref{t10.6}.
Also recall that $(M_\mu^{smar})^0\to D_{BL}$ is an immersion
and that all three spaces are 2-dimensional manifolds.
Therefore $T^{(7)}\to (M_\mu^{smar})^0$ and
$(M_\mu^{smar})^0\to D_{BL}$ are open embeddings.
We postpone the proof that the map
$T^{(7)}\to (M_\mu^{smar})^0$ is an isomorphism.

Part (b) follows now easily: Consider the case of singularities
of multiplicity $\geq 3$.
$-\id\in G_\Z$ acts trivially on $D_{BL}$. 
It acts nontrivially on $M_\mu^{smar}$ by theorem \ref{t8.5} (c).
The map $(M_\mu^{smar})^0\to D_{BL}$ is an embedding.
Therefore $-\id\in G_\Z$ does not act on $(M_\mu^{smar})^0$.
Therefore $-\id\notin G^{smar}$. This shows part (b).
In this case $(M_\mu^{smar})^0\cong (M_\mu^{mar})^0$
by theorem \ref{t8.5} (c).

In the case of singularities of multiplicity 2,
$M_\mu^{smar}=M_\mu^{mar}$ and 
$(M_\mu^{smar})^0=(M_\mu^{mar})^0$ hold anyway.

$c^{(2)}:T^{(4)}=\H\to T^{(2)}=\P^1\C$ is the branched covering
from an action of a triangle group $\Gamma$ of type
$(\frac{1}{m_0},\frac{1}{m_1},\frac{1}{m_\infty})$ on $\H$.
The group $\Gamma$ is a normal subgroup of index 2 respectively
6 of a triangle group $\Gamma^{qh}$ of type
$(2,2m,2m)$ for $W_{1,0}$ and $S_{1,0}$ and of type
$(2,3,2m)$ for $U_{1,0},E_{3,0},Z_{1,0}$ and $Q_{2,0}$
such that
$\Gamma^{qh}/\Gamma =(G_2\textup{ respectively }G_3)$.
The following pictures show hyperbolic triangles associated
to $\Gamma$ and $\Gamma^{qh}$. The symbols
$[0],[1],[\infty],[\frac{1}{2}],[2],[-1],[e^{2\pi i /6}]$
at special points indicate the images of these points
under $c^{(2)}$.

\bigskip
\includegraphics[width=1.0\textwidth]{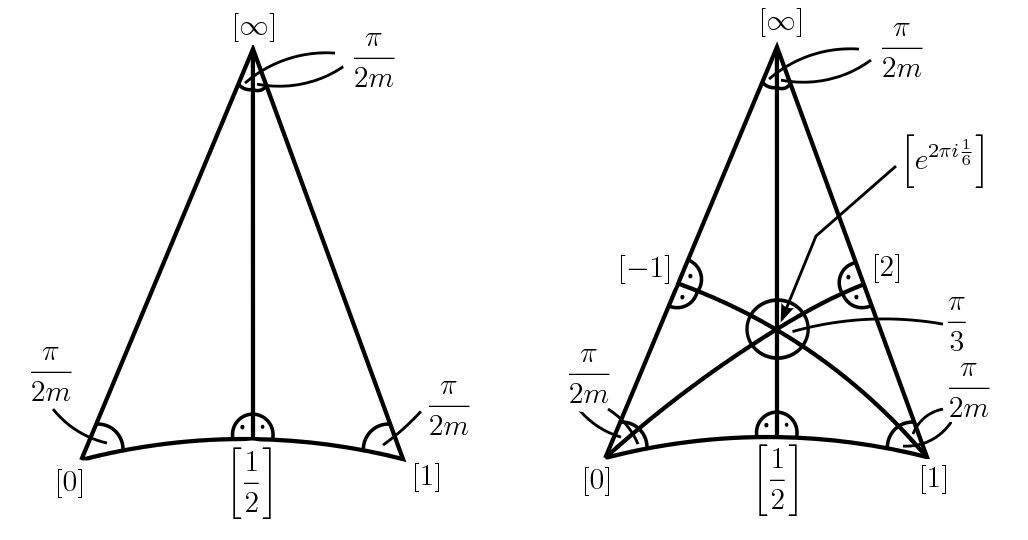}
The group $\Gamma^{qh}$ maps the set of elliptic fixed points
$(c^{(2)})^{-1}(\{0,1,\infty\})=T^{(4)}-T^{(3)}$ of $\Gamma$ 
to itself, so it acts on $T^{(3)}$.

By the proved implication $\Leftarrow$ in \eqref{10.12}
in theorem \ref{t10.3}, the orbits of $\Gamma^{qh}$
in $T^{(3)}$ are contained in the right equivalence classes
of quasihomogeneous singularities.
By the embedding $T^{(3)}\to \HH(C^{\alpha_1})$ in theorem 
\ref{t10.6}, $\Gamma^{qh}$ acts also on $\HH(C^{\alpha_1})$,
and the orbits are contained in the orbits of 
$\Psi(G^{mar})$, because the orbits of $G^{mar}$ on
$(M_\mu^{mar})^0$ are the right equivalence classes
in $(M_\mu^{mar})^0$.

Now compare the actions of $\Gamma^{qh}$ and $\Psi(G^{mar})$
on $\HH(C^{\alpha_1})$. 
$\Gamma^{qh}$ acts as a triangle group of type 
$(2,2m,2m)$ respectively $(2,3,2m)$,
and $\Psi(G^{mar})$ acts by theorem \ref{t6.1} (b)
as a subgroup of a triangle group of the same type.
And the orbits of $\Gamma^{qh}$ are contained  in the
orbits of $\Psi(G^{mar})$. Therefore the actions coincide,
and $\Psi(G^{mar})=\Psi(G_\Z)$ is a triangle group of 
the claimed type in \eqref{6.7}.
This gives the surjectivity in theorem \ref{t6.1}
and finishes the proof of theorem \ref{t6.1}.

It also shows that $G^{mar}$ acts on $T^{(3)}$.
Because $T^{(3)}$ contains representatives of the right
equivalence classes of all quasihomogeneous singularities
in the given $\mu$-homotopy family,
the marked quasihomogeneous singularities in 
$(M_\mu^{mar})^0$ must all be in $T^{(3)}$.
This proves that the open embedding
$T^{(7)}\to (M_\mu^{mar})^0$ is an isomorphism.

Next we will prove $G_\Z=G^{mar}$.
Consider an element $g_1\in G_\Z$. Because of
$\Psi(G^{mar})=\Psi(G_\Z)$, we can multiply it with
an element $g_2\in G^{mar}$ such that $g_3=g_1g_2$
satisfies $\Psi(g_3)=\id$. By formula \eqref{6.8}
in theorem \ref{t6.1} $g_3\in\{\pm M_h^k\, |\, k\in\Z\}
\subset G^{mar}$. This proves $G_\Z=G^{mar}$.

Now $M_\mu^{mar}=(M_\mu^{mar})^0$ holds.
Because $BL: (M_\mu^{mar})^0\to D_{BL}$ is an embedding,
$BL:M_\mu^{mar}\to D_{BL}$ is an embedding.
This finishes the proof of theorem \ref{t10.1}.
\hfill$\Box$

\bigskip
{\bf Proof of $\Rightarrow$ in \eqref{10.12}
in theorem \ref{t10.3}:}
$G_\Z$ acts as $\Gamma^{qh}$ on $\HH(C^{\alpha_1})$
and thus as $G_2$ respectively $G_3$ on $T^{(1)}$.
This shows $\Rightarrow$ in \eqref{10.12} for the 
quasihomogeneous singularities.

An element $g\in G_\Z$ which acts trivially on $T^{(3)}$
is in $\{\pm M_h^k\, |\, k\in\Z\}$ and restricts to 
$\lambda\cdot\id$ on $Ml_\zeta$ for some $\lambda\in\Eiw(\zeta)$.
Because of
\begin{eqnarray*}
g:\ \C\cdot (v_1+v_2)\mapsto \C(\lambda\cdot v_1+
\oooo\lambda\cdot v_2)
=\C\cdot(v_1+\oooo{\lambda}^2\cdot v_2)
\end{eqnarray*}
it acts on the fibers of the projection 
$D_{BL}\to\HH(C^{\alpha_1})$
by multiplication with $\oooo{\lambda}^2$,
and it acts in the same way on the fibers of the 
projection $T^{(7)}\to T^{(3)}$. 
But $(\oooo{\lambda}^2)^{m_\infty}=1$.
This shows $\Rightarrow$ in \eqref{10.12}
for all singularities.
\hfill$\Box$

%\nocite{*}
%\bibliographystyle{cdraifplain}
%\bibliography{xampl}
\end{document}